\documentclass[twoside,12pt]{amsart}

\usepackage{graphicx}
\usepackage{framed}
\usepackage{times}

\usepackage{color}

\usepackage[utf8]{inputenc}

\usepackage[all]{xy}
\usepackage{pb-diagram, pb-xy}

\usepackage{amssymb}
\usepackage{amsthm}

\usepackage{tikz}
\usepackage{pgffor}

\usepackage{etoolbox}

\usepackage{hyperref}

\usepackage{calc}

\usepackage[framemethod=tikz]{mdframed}

\usepackage{soul}

\usepackage{color}

\def\part{\@startsection{part}{0}%
\z@{\linespacing\@plus\linespacing}{.5\linespacing}%
{\large\bfseries\centering}}

\newtheoremstyle{saetze} 
    {5pt}                    
    {5pt}                    
    {\itshape}                   
    {12pt}                           
    {\bfseries}                   
    {.}                          
    {.5em}                       
    {}  

\theoremstyle{saetze}
\newtheorem{thm}{Theorem}[section]
\newtheorem{lem}[thm]{Lemma}
\newtheorem{cor}[thm]{Corollary}
\newtheorem{prop}[thm]{Proposition}

\newcommand{\bbZ}{{\mathbb Z}}

\newcommand{\g}{{\mathfrak g}}
\newcommand{\p}{{\mathfrak p}}
\newcommand{\gl}{{\mathfrak{gl}}}

\newcommand{\calC}{{\mathcal C}}

\newcommand{\calF}{{\mathcal F}}

\newcommand{\calN}{{\mathcal N}}

\newcommand{\calR}{{\mathcal R}}

\newcommand{\calT}{{\mathcal T}}

\newcommand{\End}{{\mathit{End}}}

\newcommand{\one}{{\mathbf{1}}}

\newcommand{\Z}{{\mathbf{Z}}}

\def\N{\mathbb{N}}
\newcommand{\A}{{\mathbb A}}

\DeclareMathSizes{12}{11}{8}{5} 
\usepackage{xypic,dsfont}

\hypersetup{
    colorlinks,
    citecolor=red,
    filecolor=black,
    linkcolor=blue,
    urlcolor=black
}

\DeclareRobustCommand{\gobblefive}[5]{}
\newcommand*{\SkipTocEntry}{\addtocontents{toc}{\gobblefive}}

\makeatletter
\def\part{\@startsection{part}{1}%
\z@{.7\linespacing\@plus\linespacing}{.5\linespacing}%
{\large\scshape\centering}}
\makeatother

\begin{document}

\title[Cohomological tensor functors]{Cohomological tensor functors on representations of the General Linear Supergroup}
\author{{\rm Th. Heidersdorf, R. Weissauer}}
\date{}


\begin{abstract} We define and study cohomological tensor functors from the category $T_n$ of finite-dimensional representations of the supergroup $Gl(n|n)$ into $T_{n-r}$ for $0 <r \leq n$. In the case $DS: T_n \to T_{n-1}$ we prove a formula $DS(L) = \bigoplus \Pi^{n_i} L_i$ for the image of an arbitrary irreducible representation. In particular $DS(L)$ is semisimple and multiplicity free. We derive a few applications of this theorem such as the degeneration of certain spectral sequences and a formula for the modified superdimension of an irreducible representation.
\end{abstract}


\maketitle

\setcounter{tocdepth}{1}

\tableofcontents 


\SkipTocEntry\section*{Introduction}


Little is known about the decomposition of tensor products between finite-dimensional representations of the general linear supergroup $Gl(m|n)$ over an algebraically closed field of characteristic 0. 
In this article we define and study \textit{cohomological tensor functors} from the category $T_n= Rep(Gl(n|n))$ of finite-dimensional representations of $Gl(n|n)$ to $T_{n-r}$ for $0 < r \leq n$. One of our aims is to reduce questions about tensor products between irreducible representations by means of these functors to lower rank cases so that these can hopefully be inductively understood. This is indeed the case for small $n$ as the $Gl(1|1)$-case has been completely worked out in \cite{Goetz-Quella-Schomerus} and the $Gl(2|2)$-case is partially controlled by the theory of mixed tensors \cite{Heidersdorf-mixed-tensors} \cite{Heidersdorf-Weissauer-gl-2-2}. Along the way we obtain formulas for the (modified) superdimensions of irreducible representations.

\bigskip

The tensor functors that we study are variants and generalizations of a construction due to Duflo-Serganova \cite{Duflo-Serganova} and Serganova \cite{Serganova-kw}. For any $x \in X =
	 \{ x \in \mathfrak{g}_1 \ | \ [x,x] = 0\},$ where $\mathfrak{g}_1$ denotes the odd part of the underlying Lie superalgebra $\mathfrak{\gl}(m|n)$, the cohomology of the complex associated to $(V,\rho) \in T_n$ \[ \xymatrix{ \ldots \ar[r]^-{\rho(x)} & V \ar[r]^-{\rho(x) } &  V  \ar[r]^-{\rho(x) } & V \ar[r]^-{\rho(x) } & \ldots } \] defines a functor $V \mapsto V_x: T_n \to T_{n-r}$ (where $r$ is the so-called rank of $x$) which preserves tensor products. The category $T_n$ splits in two abelian subcategories $T_n = \calR_n \oplus \Pi \calR_n$  where $\Pi$ denotes the parity shift (lemma \ref{thm:decomposition}). We therefore focus on the $\calR_n$-case and fix a special $x$ of rank 1 in section \ref{DF} and denote the corresponding tensor functor $DS: \calR_n \to T_{n-1}$. Later in section \ref{sec:cohomology-functors} we refine this construction  to define for any $V \in \calR_n$ a complex \[ \xymatrix{ \ldots \ar[r]^-{\partial} & \Pi (V_{2\ell-1}) \ar[r]^-{\partial} &  V_{2\ell}  \ar[r]^-{\partial} &  \Pi (V_{2\ell+1})  \ar[r]^-{\partial} & \ldots }\ \] whose cohomology in degree $\ell$ is denoted by $H^\ell(V)$. The representation $DS(V)$ is naturally $\Z$-graded and we have a direct sum decomposition \[DS(V)
 = \bigoplus_{\ell \in \Z} \ \Pi^\ell (H^\ell(V)\] for $Gl(n\! -\! 1 | n\! -\! 1)$-modules $H^\ell(V))$ in $\calR_{n-1}$. The definition of $DS$ can be easily generalized to the case $x \in X$ of higher rank $r>1$, and we denote the corresponding tensor functors by $DS_{n,n-r}: T_n \to T_{n-r}$. Like the BGG category the category $T_n$ has two different duality functors, the ordinary dual $()^{\vee}$ and the contragredient dual $()^*$. The tensor functors $DS$ and $DS_{n,n-r}$ are not $*$-invariant in the sense that $DS(V^*) \not\simeq DS(V)^*$. We therefore define an analog $D$ of the Dirac operator and we denote the corresponding Dirac cohomology groups by \[ H_D (V) = ker(D:M \to M)/ Im(D: M \to M)\] for a certain module $M\in T_{n-1}$ attached to $V$ in section \ref{DDirac}. 
 This defines a $*$-invariant tensor functor. It agrees with $DS$ on irreducible modules, but in general gives rise to an analog of Hodge decomposition (proposition \ref{H1}). The definition of Dirac cohomology generalizes easily to define $\Z$-graded tensor functors  $\omega_{n,n-r}= \bigoplus_{\ell\in\mathbb Z} \omega_{n,n-r}^\ell$ whose graded pieces are functors $\omega_{n,n-r}^\ell: T_n \to T_{n-r}$ that are described in section \ref{m}.

\bigskip

The second part is devoted to the main theorem \ref{mainthm}. In the main theorem we give an explicit formula for the image of an irreducible representation $L = L(\lambda)$ in $\calR_n$ of atypicality $j$ (for $0 < j \leq n$) under the functor $DS$. Surprisingly, with its natural $\mathbb Z$-gradation, 
the representation \[ DS(L) = \bigoplus_i L_i [-\delta_i] \] decomposes completely into a finite direct sum of irreducible representations. 
Here, for certain integers $\delta_i\in \mathbb Z$, the summands are attached to irreducible $Gl(n\! -\! 1 | n\! -\! 1)$-modules $L_i\in {\calR}_{n-1}$, where $L_i [-\delta_i]$ denotes the module $\Pi^{\delta_i}(L_i) \in T_{n-1}$ 
concentrated in degree $\delta_i$ with respect to the $\mathbb Z$-graduation of $DS(L)$. If we ignore the $\mathbb Z$-graduation, the module $DS(L) \in T_{n-1}$ always is semisimple and multiplicity free for irreducible $L$. This makes the main theorem into an effective tool to reduce questions about tensor products or superdimensions to lower rank cases in the absence of any known branching laws.

\bigskip

To analyse the $\Z$-graded object $DS(L)$ in more detail, we can assume that $L$ is a representation in the principal block containing the trivial representation (the \textit{maximally atypical} case). In fact, one can inductively reduce the general case to this special case. The irreducible maximal atypical representations $L\in {\calR_n}$ can be described in different ways. For the moment it may be sufficient that up to isomorphism they uniquely correspond  to spaced forests of rank $n$ in a natural way. By definition, such spaced forests ${\calF}$ are defined by data $$(d_0,{\calT}_1,d_1,{\calT}_2, \cdots, d_{k-1},{\calT}_k)$$ where the $\calT_i$ for $i=1,...,k$  are rooted planar trees positioned on points of the numberline  from left to right. The integer $d_0$ specifies the absolute position of the leftmost tree ${\calT}_1$
and the natural numbers $d_{i}$ for $i=1,...,k-1$ indicate the distances between the position of the trees ${\calT}_{i}$ and ${\calT}_{i+1}$. Here we allow $d_i=0$, i.e. some trees may be positioned at the same point of the numberline. The absolute positions $\delta_i = \sum_{j<i} d_i \in\mathbb Z$ of the planar trees $\calT_i$ therefore satisfy $$\delta_1 \leq \delta_2 \leq \cdots \leq \delta_k\ .$$ In particular,
$\delta_1$ describes the absolute position of the leftmost tree ${\calT}_1$ of the forest and $\delta_k$ describes the absolute position of the rightmost tree $\calT_k$ of this forest. 
Each tree ${\calT}_i$ is a planar tree with say $r_i$ nodes, among which is the distinguished node defined by the root of the tree. By definition, the rank of the forest ${\calF}$ is the sum
$\sum_{i=1}^k r_i$ of the nodes of all trees. Since in the equivalence above the rank $n$ is fixed, only
forests with at most $k\leq n$ trees occur.

\bigskip
This being said, we are now able to describe the summands  of the decomposition 
of $DS(L)$ mentioned above. For simplicity, we still assume $L$ to be maximal atypical. If $L$ corresponds to the spaced forest with trees ${\calT_1}, {\calT_2},...,{\calT}_k$ in the sense above with the positions at $\delta_1,...,\delta_k$, then $DS(L)$ has precisely $k$ irreducible constituents $L_i[-\delta_i]$ for $i=1,...,k$, so that $L_i$ corresponds to the spaced forest
$\calT_1, ... ,\calT_{i-1}, \partial{\calT}_i, ... ,{\calT}_k$ of rank $n-1$ where $\partial{\calT}_i$ denotes the forest of planar trees obtained from $\calT_i$
by removing its root. The trees are now at the new positions $\delta_1-1,...,\delta_{i-1}-1,\delta_i,...,\delta_i,\delta_{i+1}+1,...,\delta_k +1$ where we use the convention that $\delta_i$ denotes the common position of all the trees in $\partial \calT_i$. In the special case where ${\calT}_i$ has only one node,  $\partial {\calT}_i$ is not defined and will be discarded (together with $\delta_i$). In other words, in this case the new spaced forest has only $k-1$ trees.

\bigskip

This description of the $\Z$-graded object $DS(L)$ follows from the results in sections \ref{duals} - \ref{koh3}. We introduce spaced forests in section \ref{duals} where we describe the dual of an irreducible representation. The $\Z$-grading of $DS(L)$ for maximally atypical $L$ is then obtained in proposition \ref{hproof} and in the general case in proposition \ref{hproof-2}. These results follow from the main theorem and its proof by a careful bookkeeping, but they are considerably stronger and in particular incorporate theorem \ref{mainthm} as a special case. 

\bigskip

We show $DS_{n,0}(L) \cong
\bigoplus \omega_{n,0}^\ell(L)[-\ell]$ for irreducible maximal atypical representations  $L$. From this,
as an application of the main theorem, we obtain  in theorem \ref{thm:forest-formula} a nice explicit formula for the Laurent polynomial \[ \sum_{\ell \in \Z} sdim(\omega_{n,0}^{\ell}(L)) \cdot t^{\ell}. \] (Hilbert polynomial) 
attached to the Dirac cohomology tensor functors $$\omega_{n,0}^\ell: T_n \to T_0 \ $$
in the case of an irreducible maximal atypical representation  $L$.  

\bigskip

As already mentioned, the main theorem does not require $L$ to be in the principal block.
Applying $DS$  repeatedly $k$-times to an irreducible representation $L = L(\lambda)$ of atypicality $i$ we obtain an isotypical typical representation $m(\lambda) L^{core}$ in $T_{n-i}$, and $L^{core}$ only depends on the block of $L$ (section \ref{sec:main}).  We derive a closed formula for the multiplicity $m(\lambda)$ in section \ref{sec:main}. The multiplicity $m(\lambda)$ can be expressed as \[ m(\lambda) = \frac{|\mathcal{F}(\lambda)|!}{\mathcal{F}(\lambda)!} \] where $\mathcal{F}(\lambda)$ is the spaced forest associated to $L(\lambda)$, $|\mathcal{F}(\lambda)|$ is the number of its nodes and $\mathcal{F}(\lambda)!$ is the forest factorial \ref{sec:main}. This not only implies that the so called modified superdimension of $L$ does not vanish (i.e. the \textit{generalized Kac-Wakimoto conjecture}), but moreover gives a closed formula for it. The main theorem has a number of other useful applications and we refer the reader to the list given after theorem \ref{mainthm}.

\bigskip


The proof of the main theorem occupies the entire second part. Build on an involved induction using translation functors, carried out in sections \ref{sec:loewy-length} - \ref{sec:moves}, the proof is reduced to the case of ground-states; these are rather specific irreducible modules in a block. For instance, ground states of the principal block are powers of the Berezin determinant. Then, for a block of atypicality $k < n$, we prove in section \ref{stable0} that every ground state is a Berezin twist of a mixed tensor, an irreducible direct summand in an iterated tensor $X^{\otimes r} \otimes (X^{\vee})^{\otimes s}$ where $X$ denotes the standard representation of $Gl(n|n)$. It is easy to verify directly that the main theorem holds for Berezin powers and mixed tensors. In section \ref{sec:loewy-length} we study the Loewy structure of translation functors applied to irreducible representations and their behavior under $DS$. We also explain why we can restrict to the maximally atypical case for the proof of the main theorem. In section \ref{sec:inductive} we prove both parts of the main theorem (semisimplicity and determination of the constituents) under certain assumption on translation functors which are verified in the later section \ref{sec:moves}.

\bigskip

In section \ref{sec:chevalley-eilenberg} we discuss the cohomology ring $H_{DS_n}^{\bullet}(V(\one))$ for the tensor functor $DS_{n,0}$. Although the description of the composition factors of an arbitrary Kac module $V(\lambda)$ is much more complicated than that of $V(\one)$, we show in lemma \ref{lem:kac=kac} that there is an isomorphism \[ H^{\bullet}_{DS_n}(V(\one)) \cong H_{DS_n}^{\bullet + deg(\rho)} (V(\lambda)) \] of $I$-modules. In fact the cohomology ring of $V(\one)$ can be identified with the Lie algebra homology ring $H_{\bullet}(\mathfrak{gl}(n))$ and defines an exterior algebra $I$ on primitive elements $f_1,f_3, \ldots, f_{2n-1}$ so that $I$ acts on the graded cohomology $H^\bullet_{DS_n}$ of finite dimensional $\frak g$-modules.
 We also discuss the relationship between the cohomology of a Kac module and its irreducible quotient. We show in theorem \ref{cohom-proj} that the induced homomorphism \[ H^{\nu}_{DS_n}(pr): H_{DS_n}^{\nu}(V(\lambda)) \to H_{DS_n}^{\nu}(L(\lambda)) \] is an isomorphism in the top degree and trivial in all lower degrees. In section \ref{sec:primitive} we describe the elements of $I \cong H^{\bullet}_{DS_n}V(\one))$ in terms of the representation
theory of the superlinear group $Gl(n\vert n)$.

\bigskip

Since the image of an irreducible representation under $DS$ is therefore understood, it is natural to look at the image $DS(I)$ of an indecomposable representations $I$. The kernel of $DS$ is the tensor ideal of representations with a filtration by anti-Kac modules by results in section \ref{sec:support}. If $R(\lambda)$ is a mixed tensor we can easily compute $DS(R(\lambda))$. In other cases it is rather complicated to determine $DS$. As an example for the importance of this problem consider the computation of the tensor product between two irreducible representations $L_1 \otimes L_2 = \bigoplus_i I_i$ in indecomposable summands in $\calR_n$. The decomposition of $DS(L_1)$ and $DS(L_2)$ gives estimates on the number of possible direct summands, but these are rather weak unless something is known about $DS(I_i)$. For an easy example of the use of the cohomological tensor functors in this setting see \cite{Heidersdorf-Weissauer-gl-2-2}. In the last sections \ref{kac-module-of-one} - \ref{hooks} we give a cohomological criterion \ref{splitting1} for an indecomposable representation to be equal to the trivial representation. We call an epimorphism $q: V \to W$ strict, if the induced morphism $\omega(q): \omega(V) \to \omega(W)$ for the tensor functor $\omega = \omega_{n,0}:T_n \to svec_k$ is surjective. We prove in corollary \ref{splitting1} that if $Z$ is an indecomposable module with cosocle $\one$ such that the quotient map $q: Z \to \one$ is strict, then $Z \simeq \one$. Any such representation $Z$ contains extensions of the trivial representations $\one$ with the other irreducible constituents in the second upper Loewy layer. This leads us to study the cohomology $H^i$ of extensions of the trivial representation \[ \xymatrix{ 0 \ar[r] & S_{\nu} \ar[r] & V \ar[r]^{q_V}  & \one \ar[r] & 0} \] for irreducible representations $S_{\nu}$. We show in the key lemma \ref{trivialextension} that in this case the map $\omega^0(q_V)$ vanishes. This is a contradiction to our analysis in section \ref{strictmorphisms} if we suppose that $Z$ is not irreducible. 

\bigskip
Most of the results in this article can be rephrased for representations of the supergroup $Gl(m|n)$ where $m \neq n$. This will be discussed elsewhere.

\bigskip

\subsection*{Acknowledgements} The authors are grateful to the referee for providing useful suggestions.

\bigskip



\part{Cohomological tensor functors}



\section{ The superlinear groups}\label{2}

\medskip

Let $k$ be an algebraically closed field of characteristic zero.
A super vectorspace $V$ over $k$ is a $\mathbb Z/2\mathbb Z$-graded $k$-vectorspace $V=V_{\overline 0} \oplus V_{\overline 1}$. Its superdimension is  $sdim(V)= \dim(V_{\overline 0})\! -\!
\dim(V_{\overline 1})$.
The parity shift functor $\Pi$ on the category of super vectorspaces over $k$ 
is defined by $\Pi(V)_{\overline 0} = V_{\overline 1}$ and $\Pi(V)_{\overline 1} =V_{\overline 0}$ and
the parity endomorphism of $V$ is $p_V=id_{V_{\overline 0}}\oplus -id_{V_{\overline 1}}$ in $End_k(V)$.  

\medskip

{\it Conventions on gradings}. For $\mathbb Z$-graded object $M=\bigoplus_i M_i$ with objects $M_i$ in an additive category $\calC$  one has the shifted  $\mathbb Z$-graded objects $M\langle j\rangle$
defined by $(M\langle j\rangle)_i = M_{i+j}$. If $\calC$ carries a super structure defined by a functor $\Pi: \calC \to \calC$ such that $\Pi \circ \Pi$ is the identity functor, we mainly use the $\mathbb Z$-graded objects 
 $M[j]$ defined by $(M[j])_i := \Pi^j(M_{i+j})$. Considering objects $L$ in $\calC$ as graded objects concentrated in one degree, we often consider the $\mathbb Z$-graded objects $L[-\ell]$
 concentrated in degree $\ell$. In this context, forgetting the $\mathbb Z$-grading of $L[\ell]$ for $L\in \calC$ and $\ell\in \mathbb Z$ gives the object $\Pi^\ell(L)$ in $\calC$.

\medskip
{\it The categories $F$ and $T$}. Let  $\g=\gl(m\vert  n) = \g_{\overline{0}} \oplus \g_{\overline{1}}$
be the general Lie superalgebra. The even part $\g_{\overline{0}} = \gl(m) \oplus \gl(n)$ of $\gl(m\vert n)$ can be considered as the Lie algebra of the classical subgroup
$G_{\overline 0}=Gl(m)\times Gl(n)$ in $G=Gl(m\vert n)$.
By definition a finite-dimensional representation $\rho$ of $\mathfrak{gl}(m\vert n)$ 
defines a representation $\rho$ of $Gl(m \vert n)$, if its restriction to $\g_{\overline{0}}$ 
comes from an algebraic representation of $G_{\overline 0}$, also denoted $\rho$.
For the linear supergroup $G=Gl(m\vert n)$ over $k$
let $F$ be the category of the 
super representations $\rho$ of $Gl(m\vert n)$ on finite dimensional super vectorspaces over $k$.
If $(V,\rho)$ is in $F$, so is $\Pi(V,\rho)$.
The morphisms in the category $F$ are the $G$-linear maps $f:V \to W$ between super representations, 
where we allow even and odd morphisms  with respect to the gradings on $V$ and $W$, i.e
morphisms with $f \circ p_V =\pm p_W \circ f$.
For $M,N \in F$ we have $Hom_F(M,N) = Hom_F(M,N)_{\overline 0} \oplus Hom_F(M,N)_{\overline 1}$, where $Hom_F(M,N)_{\overline 0}$ are the even morphisms. Let $T=sRep_\Lambda(G)$ be the subcategory of $F$ with the same objects as $F$ and $Hom_T(M,N)=Hom_F(M,N)_{\overline 0}$ . 
Then $T$ is an abelian category, whereas $F$ is not.

\medskip
{\it The category ${\calR}$}. Fix the morphism $\varepsilon: \mathbb Z/2\mathbb Z \to G_{\overline 0}=Gl(m)\times Gl(n)$ which maps $-1$ to the element 
$diag(E_m,-E_n)\in Gl(m)\times Gl(n)$ denoted $\epsilon_{mn}$. We write $\epsilon_n = \epsilon_{nn}$. Notice that
$Ad(\epsilon_{mn})$ induces the parity morphism on the Lie superalgebra $\gl(m|n)$ of $G$. 
We define the abelian subcategory
$\calR = sRep(G,\varepsilon)$ of $T$ as the full subcategory of all objects $(V,\rho)$ in $T$  
with the property $  p_V = \rho(\epsilon_{mn})$; here $\rho$ denotes the underlying homomorphism $\rho: Gl(m)\times Gl(n) 
\to Gl(V)$ of algebraic groups over $k$.
The subcategory ${\calR}$ is stable under the dualities ${}^\vee$ and $^*$. 
For $G=Gl(n\vert n)$ we usually write $T_n$ instead of $T$, and ${\calR}_n$ instead of $\calR$,
to indicate the dependency on $n$. 

\medskip{\it The duality $*$}.
The Lie superalgebra $\g=\mathfrak{gl}(m\vert  n)$ has a consistent \cite{Kac-Rep} $\mathbb Z$-grading $\g = \g_{(-1)} \oplus \g_{(0)} \oplus \g_{(1)}$, where $\g_{\overline{0}} = \g_{(0)}$ and where $\g_{\overline{1}} = \g_{(-1)} \oplus \g_{(1)}$ is defined by the upper triangular block matrices $\g_{(1)}$ and $\g_{(-1)}$ by the lower triangular block matrices.  The supertranspose $x^T$ (see \cite{Scheunert}, (3.35) and (4.14)]) of a graded endomorphism $x\in \End(k^{m\vert n})$ is 
defined by \[  x=\begin{pmatrix} m_1 & m_2 \\ m_3 & m_4 \end{pmatrix} \ \mapsto \ x^T = \begin{pmatrix} m_1^t & -m_3^t \\ m_2^t & m_4^t \end{pmatrix}, \ \ x \in \g \ \] where $m_i^t$ denotes the ordinary transpose of the matrices $m_i$. If we identify $\g$ and $End(k^{m\vert n})$, then $\tau(x)= - x^T$ defines an automorphism of the Lie superalgebra $\g$ such that $\tau(\g_{(i)}) = \g_{(-i)}$ holds for $i=-1,0,1$. For a representation $M=(V,\rho)$ in $T_n$ and homogenous $x $ in $\g$ the Tannaka dual representation $M^\vee=(V^\vee,\rho^\vee)$ is the representation
$x\mapsto - \rho(x)^T$ on $V$,  using the supertranpose $\rho(x)^T$ of $\rho(x)$ in $End(V)$.
Finally we define the representation $M^* =(V^{\vee},\rho^\vee \circ \tau) $, where $\tau(x)=-x^T$ is the automorphism of $\g$ defined by the supertranspose on $\g$. 
See also \cite{BKN-complexity}, 3.4 using a different convention. $V \in {\calR}_n$ (see below) implies $V^*\in {\calR}_n$ by \cite{Brundan-Kazhdan}, lemma 4.43. For simple and for projective objects $V$ of $T_n$ furthermore $V^*\cong V$. Also $V^*\vert_{G_{\overline 0}} \cong
V\vert_{G_{\overline 0}}$ for all $V$ in $T_n$. Notice that both $\vee$ and $*$
define contravariant functors on $T_n$.

\medskip
{\it Weights}.
Consider the standard Borel subalgebra $\mathfrak{b}$ of upper triangular matrices in $\g$
and its unipotent radical $\mathfrak{u}$.  
The basis $\Delta$ of positive roots associated to $\mathfrak{b}$ is given by the basis of the positive roots
associated to $\mathfrak{b} \cap \g_{\overline 0}$ for the Lie algebra $\g_{\overline 0}$
and a single odd root $x$ whose weight will be called $\mu$. If we denote by $e_{i,i}$, $i=1,\ldots,2n$, the linear form which sends a diagonal element $(t_1,\ldots,t_{2n})$ to $t_i$, then the simple roots in this basis are given by the set $\{e_{1,1} - e_{2,2}, \ldots, e_{2n-1,2n-1} - e_{2n,2n} \}$ with $\mu = e_{n,n} - e_{n+1,n+1}$. The diagonal elements $t=diag(t_1,...,t_n,t_{n+1},...,t_{2n})$ in $G_{\overline 0}$ act by semisimple matrices on $V$ for any representation $(V,\rho)$ in $T_n$. Hence $V$ decomposes into a direct sum of eigenspaces $V = \bigoplus_\lambda V_\lambda$
for certain characters $t^\lambda= t_1^{\lambda_1} \cdots t_n^{\lambda_n}
(t_{n+1})^{\lambda_{n+1}} \cdots (t_{2n})^{\lambda_{2n}}$. Then write $\lambda=(\lambda_1,...,\lambda_n ;
\lambda_{n+1}, \cdots, \lambda_{2n})$. A {\it primitive} weight vector $v$ (of weight $\lambda$) in a representation $(V,\rho)$
of $g$ is a nonzero vector in $V$ with the property $\rho(X)v=0$ for $X\in u$ and
$\rho(t)v = t^\lambda$. An irreducible representation $L$ has a unique primitive weight vector (up to a scalar), the highest weight vector. Its weight $\lambda$ 
uniquely determines the irreducible module $L$ up to isomorphism in $\calR_n$. Therefore we write $L=L(\lambda)$.

\medskip
{\it Kac modules}. We put $\mathfrak{p}_{\pm} = \g_{(0)} \oplus \g_{(\pm1)}$. We consider a simple $\g_{(0)}$-module as a $\mathfrak{p}_{\pm}$-module in which $\g_{(1)}$ respectively $\g_{(-1)}$ acts trivially. We then define the Kac module $V(\lambda)$ and the AntiKac module $V'(\lambda)$ via \[ V(\lambda)  = Ind_{\mathfrak{p}_+}^{\g} L_0(\lambda) \ , \ V'(\lambda)  = Ind_{\mathfrak{p}_-}^{\g} L_0(\lambda) \] where $L_0(\lambda)$ is the simple $\g_{(0)}$-module with highest weight $\lambda$. The Kac-modules are universal highest weight modules. $V(\lambda)$ has a unique maximal submodule $I(\lambda)$ and $L(\lambda) = V(\lambda)/I(\lambda)$ \cite{Kac-Rep}, prop.2.4.

\medskip
{\it The Berezin}. The Berezin determinant of the supergroup $G=G_n$
defines a one dimensional representation $Ber=Ber_n$. Its weight is
is given by $\lambda_i=1$ and $\lambda_{n+i} =-1$ for $i=1,..,n$.
The representation space of $Ber_n$ has the superparity $(-1)^n$.
We denote the trival representation $Ber^0$ by $\one$.

\medskip
{\it Ground states}. Each $i$-atypical block of $\calR_n$ contains irreducible representations $L(\lambda)$ of the form
$$  \lambda = (\lambda_1,...,\lambda_{n-i},\lambda_n,...,\lambda_n\ ;\
-\lambda_n,...,-\lambda_n,\lambda_{n+1+i}, ..., \lambda_{2n}) \ .$$
with $\lambda_n \leq \min(\lambda_{n-i}, - \lambda_{n+1+i})$. We call these the ground states of the block. They will play a major role in our computation of $DS(L)$ in theorem \ref{mainthm}.

\medskip
{\it Equivalence}. Two irreducible representations $M,N$ on $T$ are said to be
equivalent $M \sim N$, if either $M \cong Ber^r \otimes N$ or $M^\vee \cong Ber^r \otimes N$ holds for some $r\in \mathbb Z$. This obviously defines an equivalence relation on the set of isomorphism classes of irreducible representations of $T$. A self-equivalence of $M$ is given by an isomorphism $f: M \cong Ber^r \otimes M$ (which implies $r=0$ and $f$ to be a scalar multiple of the identity) respectively an isomorphism $f: M^\vee \cong Ber^r \otimes M$.
If it exists, such an isomorphism uniquely determines $r$ and is unique up to a scalar
and we say $M$ is of type (SD). Otherwise we say $M$ is of type (NSD).

\medskip
{\it Negligible objects}. An object $M\in T_n$ is called negligible if it is the direct
sum of indecomposable objects $M_i$ in $T_n$ with superdimensions $sdim(M_i)=0$. The tensor ideal of negligible objects is denotes $\calN$ or $\calN_n$.

\medskip



\section{ The Duflo-Serganova functor $DS$}\label{DF}

\medskip {\it An embedding}. Fix some $1\leq m\leq n$. We view $G_{n-m}= Gl(n-m\vert  n-m)$ as an \lq{outer  block matrix}\rq\ in $G_n=Gl(n\vert n)$ and $G_1$ as the \lq{inner  block matrix}\rq\  as below. Here $G_0$ is the empty group. We fix some invertible $m\times m$-matrix $J$ 
with the property $J=J^t = J^{-1}$. 
For example take $J$ to be the identity matrix $E$ , or the matrix with nonzero entries equal to 1 only in the antidiagonal. 
We furthermore fix the 
embedding
$$  \varphi_{n,m}: G_{n-m} \times G_1 \hookrightarrow G_n \  $$
defined by $$\begin{pmatrix} A & B \\ C & D \end{pmatrix} \times \begin{pmatrix} a & b \\ c  & d\end{pmatrix} \mapsto \begin{pmatrix} A & 0 & 0 & B \\ 0 & a  E  & b  J  & 0\\ 0 & c J  &  d  E   & 0\\ C & 0 & 0 & D \end{pmatrix}$$ 
We use this embedding to identify elements in $G_{n-m}$ and $G_1$ with elements
in $G_n$.
In this sense $\epsilon_n = \epsilon_{n-m} \epsilon_1$ holds in $G_n$, for the corresponding elements $\epsilon_{n-m}$ and 
$\epsilon_1$ in $G_{n-m}$ resp. $G_1$, defined in section \ref{2}.

\medskip {\it Two functors}. 
One has a functor $(V,\rho) \mapsto V^+ =\{ v \in V\ \vert \ \rho(\epsilon_1)(v)=v \}$
$$ {}^+: {\calR}_n \to {\calR}_{n-m}$$ 
where $V^+$ is considered as a $G_{n-m}$-module using $\rho(\epsilon_1) \rho(g) = \rho(g) \rho(\epsilon_1)$
for $g\in G_{n-m}$. Indeed $Ad(\epsilon_1)(g)=g$ holds for all $g\in G_{n-m}$.
The grading on $V$ induces a grading on $V^+$
by $(V^+)_{\overline 0}= V_{\overline 0} \cap V^+$ and $(V^+)_{\overline 1}= V_{\overline 1} \cap V^+$. For this grading
the decomposition $V^+ = (V^+)_{\overline 0} \oplus (V^ +)_{\overline 1}$ is induced by the parity morphism
$\rho(\epsilon_n)$ or equivalently $\rho(\epsilon_{n-1})$. With this grading on $V^ +$ the restriction of $\rho$
to $G_{n-m}$ preserves $V^+$ and defines a representation $(V^+,\rho)$ of $G_{n-m}$ in ${\calR}_{n-m}$.

\medskip
Similarly define $V^- =\{ v \in V\ \vert \ \rho(\epsilon_1)(v)=-v \}$. With the grading induced from $V=V_{\overline 0}\oplus V_{\overline 1}$
this defines a representation $V^-$  of $G_{n-m}$  in $\Pi {\calR}_{n-m}$. Obviously
$$   (V,\rho)\vert_{G_{n-m}} \ =\ V^+  \ \oplus \ V^-  \ .$$  


\medskip {\it The exact hexagon}.  Fix the following element $x\in \g_n$  $$x = \begin{pmatrix} 0 & y \\ 0 & 0 \end{pmatrix} \in \g_{n} \ \text{ for } \ y = \begin{pmatrix} 0 & 0 & \ldots & 0 \\ 0 & 0 & \ldots & 0 \\ \ldots & & \ldots &  \\ J & 0  & 0 & 0 \\ \end{pmatrix} $$
for the fixed invertible $m\times m$-matrix $J$. Since $x$ is an odd element with $[x,x]=0$, we get $$2 \cdot \rho(x)^2 =[\rho(x),\rho(x)] =\rho([x,x]) =0 $$ for any representation $(V,\rho)$ of $G_n$ in ${\calR}_n$. Notice $d= \rho(x)$ supercommutes with $\rho(G_{n-m})$. 
Furthermore $\rho(x): V^{\pm} \to V^{\mp}$ holds as a $k$-linear map, an immediate consequence of $d\rho(\varepsilon_{1}) = - \rho(\varepsilon_{1})d$, i.e. of $Ad(\varepsilon_1)(x)=-x$. 
Since $\rho(x) \in Hom_F(V,V)_1$ is an {\it odd} morphism,
$\rho(x)$ induces the following {\it even} morphisms (morphisms in ${\calR}_{n-m}$)
$$ \rho(x): V^+ \to \Pi(V^-) \quad \text{ and } \quad \rho(x): \Pi(V^-) \to V^+ \ .$$
The $k$-linear map $\partial=\rho(x): V\to V$ is a differential and commutes with the action of $G_{n-m}$ on $(V,\rho)$. Therefore $\partial$ defines a complex
in ${\calR}_{n-m}$ 
$$ \xymatrix{ \ar[r]^-{\partial} & V^+ \ar[r]^-{\partial} &  \Pi(V^-) \ar[r]^-{\partial} & V^+ \ar[r]^-{\partial} & \cdots } $$
Since this complex is periodic,  it has essentially only two cohomology groups denoted $H^+(V,\rho)$ and $H^-(V,\rho)$ in the following. This defines two functors $(V,\rho) \mapsto D_{n,n-m}^\pm(V,\rho)=H^{\pm}(V,\rho)$ 
$$  \fbox{$ D_{n,n-m}^\pm: {\calR}_n \to {\calR}_{n-m} $} \ .$$ 
It is obvious that an exact sequence \[ \xymatrix{ 0 \ar[r] & A \ar[r]^\alpha & B \ar[r]^\beta & C \ar[r] & 0}\] in ${\calR}_n$ gives rise to an exact sequences of complexes in ${\calR}_{n-m}$.
Hence 

\begin{lem} \label{hex} The long exact cohomology sequence defines
an exact hexagon in ${\calR}_{n-m}$ \[ \xymatrix{ & H^+ (A) \ar[r]^{H^+(\alpha)} & H^+(B) \ar[dr]^{H^+(\beta)} & \\ H^-(C) \ar[ur]^\delta & & & H^+(C) \ar[dl]^{\delta} \\ & H^-(B) \ar[ul]^{H^-(\beta)} & H^-(A) \ar[l]^{H^-(\alpha)} & }\]
\end{lem}


\medskip\noindent  
{\it Alternative point of view}. For the categories $T=T_n$ resp. $T_{n-m}$ (for the groups
$G_n$ resp. $G_{n-m}$) consider the  
tensor functor of Duflo and Serganova in \cite{Duflo-Serganova} $$ DS_{n,n-m}: T_n \to T_{n-m} $$
defined by $DS_{n,n-m}(V,\rho)= V_x:=Kern(\rho(x))/Im(\rho(x))$. For $(V,\rho)\in \calR_n$ we obtain 
$$   H^+(V,\rho) \oplus \Pi (H^-(V,\rho)) = DS_{n,n-m}(V) \ .$$
Indeed, the left side is $DS_{n,n-m}(V)=V_x$ for the $k$-linear 
map $\partial=\rho(x)$ on $V=V^+ \oplus V^-$. Hence $H^ +$  is the functor obtained by composing the tensor functor
$$ DS_{n,n-m}: {\calR}_n \to T_{n-m} $$
with the functor 
$$ T_{n-1} \to {\calR}_{n-m} $$
that projects the abelian category $T_{n-m}$ onto ${\calR}_{n-m}$ using

\begin{lem} \label{thm:decomposition} Every object $M \in T_n$ decomposes uniquely as $M = M_0 \oplus M_1$ with $M_0 \in {\calR}_n$ and $M_1 \in \Pi({\calR}_n)$. This defines a block decomposition of the abelian category $$\fbox{$ T = {\calR}_n \oplus \Pi ({\calR}_n) $} \ .$$
\end{lem}

\medskip{\it Proof}. For any $M, N \in {\calR}_n$ the $\mathbb Z_2$-graded space $Ext_T^i (M,N)$ is concentrated in degree zero \cite{Brundan-Kazhdan}, Cor. 4.44. \qed

\medskip
{\it Tensor property}. As a graded module over $R=k[x]/x^2$ any representation $V$ decomposes into a direct sum of a trivial representation $T$ and copies of $R$ (ignoring shifts by $\Pi$). To show that
$DS_{n,m}(V)= R_x \oplus T_x = T$ is a tensor functor, it suffices that $(R\otimes R)_x=0$, see also \cite{Serganova-kw}. For this we use that the underlying tensor product is the supertensor product.
Indeed for $R=V_{\overline 0}\oplus V_{\overline 1}$ and $V_{\overline 0}=k\cdot 1$ and $V_{\overline 1}=k\cdot x$ 
we have $x(e_1)=e_2$ and $x(e_2)=0$. The induced {\it superderivation} $d$ on $R\otimes R$ satisfies $d(1\otimes 1)=x\otimes 1+ 1\otimes x$, $d(x\otimes 1)= -x\otimes x$,
$d(1\otimes x)=x\otimes x$ and $d(x\otimes x)=0$. Hence $Im(d)= Ker(d)= k\cdot (1\otimes x + x \otimes 1) \oplus k \cdot x \otimes x$
and therefore $(R\otimes R)_x=0$.

\medskip



\section{ Cohomology Functors}\label{sec:cohomology-functors}

\medskip
In this section we assume $V \in T_n$ and $m=1$. In the following
let $DS$ be the functor $DS_{n,n-1}$ (for $J=1$). 

\medskip
{\it  Enriched weight structure}. The maximal torus of diagonal matrices in $G_n$
naturally acts on $DS(V)$ so that $DS(V)$ decomposes into weight
spaces $DS(V) = \bigoplus_\lambda DS(V)_\lambda$ for $\lambda$ in the
weight lattice $X(n)$ of $\g_n$. Indeed for the weight decomposition $V= \bigoplus_\lambda V_\lambda$
every $v\in V$ has the form $v = \sum_\lambda v_\lambda$
for $v_\lambda\in V_\lambda$.  Now $\partial v=0$ if and only if 
$\partial v_\lambda=0$ holds for all $\lambda$, since $\partial(V_\lambda) \subseteq
V_{\lambda + \mu}$ for the odd simple weight $\mu$
(ignoring parities on $V$). Similarly $v = \partial w$ if and only if $v_\lambda = \partial w_\lambda$
for all $\lambda$, since we can always project on the weight eigenspaces. 
This trivial remark shows that $DS(V)$ naturally carries a weight decomposition 
with respect to the weight lattice $X(n)$ of $\g_n$. The weight structure for $\g_{n-1}$
is obtained by restriction.
The kernel of the restriction $X(n) \to X(n-1)$ of weights, denoted by $$\lambda \mapsto \overline\lambda\ ,$$ are the multiples $\mathbb Z \cdot \mu$ of the odd simple root $\mu=e_{n,n} - e_{n+1,n+1}$. We may therefore view $DS(V)$ as endowed with the richer weight structure coming from the $G_n$-module $V$. This decomposition induces a natural decomposition of $DS(V)$
into eigenspaces $DS(V) = \bigoplus_\ell \ DS(V)_\ell$.
To make this more convenient consider the torus of elements
$diag(1,...,1,1;t^{-1},1,...,1)$ for $t\in k^*$, called the {\it small torus}. These elements commute with $G_{n-1}$ and their eigenvalue decomposition gives a decomposition $$ V = \bigoplus_{\ell\in\mathbb Z} \ V_{\ell} $$ into $G_{n-1}$-modules $V_{\ell}$. Here $V_\ell \subseteq V$ denotes the subspace defined by all vectors in $V$ on which the above elements of the small torus acts by multiplication with $t^\ell$.
Obviously $V_\ell=0$ for $\ell\notin [\ell_0, \ell_1]$ and suitable $\ell_0,\ell_1$. For the odd morphism $\partial = \rho(x)$
the properties $\mu(diag(1,...,1,1;t^{-1},1,...,1))=
t$ and $\partial(V_\lambda) \subseteq V_{\lambda+\mu}$ show that
$$  \xymatrix{ \ar[r]^-\partial & \Pi(V_{2\ell - 1}) \ar[r]^-\partial & V_{2\ell} 
\ar[r]^-\partial & \Pi(V_{2\ell + 1}) \ar[r]^-\partial & V_{2\ell + 2} \ar[r]^-\partial &  }$$
defines a complex. Its cohomology is denoted 
$   H^\ell(V)$. Obviously 
\[   \Pi^\ell(H^\ell(V)) = DS(V)_\ell \] and hence we obtain a decomposition 
of $DS(V,\rho)$ into a direct sum of $G_{n-1}$-modules
$$   DS(V,\rho) \ = \ \bigoplus_{\ell \in\mathbb Z}  \ \Pi^\ell(H^\ell(V)) \ ,$$
If we want to emphasize the $\mathbb Z$-grading, we also
write this in the form
$$  \fbox{$ DS(V,\rho) \ = \ \bigoplus_{\ell \in\mathbb Z}  \ H^\ell(V)[-\ell] $} \ .$$
We will calculate $DS(L) \in T_{n-1}$ for irreducible $L$ in theorem \ref{mainthm} and we will compute its $\Z$-grading in proposition \ref{hproof} and proposition \ref{hproof-2}.

An exact sequence \[ \xymatrix{ 0 \ar[r] & A \ar[r]^\alpha & B \ar[r]^\beta & C \ar[r] & 0}\] in ${\calR}_n$ then gives rise to a long exact sequence in ${\calR}_{n-1}$
$$  \xymatrix@C=1.4em{\ar[r] & H^{\ell-1}(C)  \ar[r] & H^{\ell}(A) \ar[r] & H^{\ell}(B) 
\ar[r] & H^{\ell}(C) \ar[r] & H^{\ell+1}(A) \ar[r] & . }$$

\begin{lem}\label{-ell} 
For $V$ in $T_n$ we have 
$H^\ell(Ber_n \otimes V) = Ber_{n-1} \otimes H^{\ell -1}(V)$. 
For the Tannaka dual $V^\vee$ of $V$ 
$H^\ell(V)^\vee \cong H^{-\ell}(V^\vee)$ holds for all $\ell\in\mathbb Z$
(isomorphisms of $G_{n-1}$-modules).
\end{lem}

\medskip
{\bf Proof}. The first property follows from $DS(Ber_n)=Ber_{n-1}[-1]$ and the fact that 
$DS$ is a tensor functor. 
Furthermore $DS(V)^\vee \cong DS(V^\vee)$, since $DS$ is a tensor functor.
Hence the second claim follows from $(V^\vee)_{-\ell} = (V_{\ell})^\vee$,
since $\Pi^2$ is the identity 
and duality \lq{commutes}\rq\ with the parity shift $\Pi$. \qed

\medskip
Note that for $V_\ell \in T_{n-1}$ the module $(V_\ell)^* \in T_{n-1}$
is isomorphic to $(V^*)_\ell$. 

\medskip
Finally, for $(V,\rho)\in \calR_n$ we get  
$ V^+ = \bigoplus_{\ell \in 2\mathbb Z} V_\ell$ and $\Pi(V^-) = \bigoplus_{\ell \in1+2\mathbb Z} V_\ell $. Hence we obtain the next lemma.

\begin{lem} \label{+ell} For $V$ in $\calR_n$ the following holds
$$  \fbox{$ H^{+}(V) \ = \ \bigoplus_{\ell \in 2\mathbb Z} H^\ell(V) \quad , \quad  H^-(V) = \bigoplus_{\ell \in1+2\mathbb Z} H^\ell(V). $} $$
\end{lem}

\medskip



\section{ Support varieties and the kernel of $DS$}\label{sec:support}

\medskip

We show that the kernel of $DS$ consists of the modules which have a filtration by AntiKac modules.

\medskip
{\it Support varieties}. We review results from \cite{BKN-1}, \cite{BKN-2} and \cite{BKN-complexity} on support varieties. Recall the decomposition $\g = \g_{(-1)}\oplus \g_{(0)} \oplus \g_{(-1)}$. The support varieties are defined by 
\[ V_{\g_{(\pm 1)}}(M) = \{ \xi \in \g_{(\pm 1)} \ 
| \ M \text{ not projective as a } U(\langle \xi \rangle)- \text{module} \} \cup \{0\} \ . \]  
Notice that $\xi \in \g_{(\pm 1)}$ generates an odd abelian Lie superalgebra 
$\langle\xi \rangle$ with $[\xi,\xi] = 0$, which up to isomorphisms has only two indecomposable modules: The trivial module and its projective cover $U(\langle\xi \rangle)$. 
 By \cite{BKN-1}, prop 6.3.1
  \[ V_{\g_{(\pm 1)}}(M \otimes N) = V_{\g_{(\pm 1)}}(M) \cap V_{\g_{(\pm 1)}}(N). \] 
The associated variety of Duflo and Serganova is 
defined as \[ X_M  = \{ \xi \in X \ | \ M_\xi \neq 0 \} \]
where $X$ is the cone $X = \{ \xi \in \g_{\bar{1}} \ | [\xi,\xi] = 0 \}$. For $\xi \in X$ the condition $M_\xi \neq 0$ is equivalent by \cite{BKN-complexity}, 3.6.1, to the condition that $M$ is not projective as a $U(\langle \xi \rangle)$-module.  Hence
$ X_M$ is the set of all $\xi \in X$ such that $M$  is not projective as a  $U(\langle \xi \rangle)$-module
together with $\xi=0$. Thus  
$$ V_{\g_{(-1)}}(M) \cup V_{\g_{(1)}}(M) \subseteq X_M \quad , \quad  V_{\g_{(\pm 1)}}(M)  = X_M \cap \g_{(\pm 1)} \ . $$

\bigskip
{\it Kac and anti-Kac objects}. We denote by $\calC^+$ the tensor ideal of modules with a filtration by Kac modules in $\calR_n$ and by $\calC^-$ the tensor ideal of modules with a filtration by anti-Kac modules in $\calR_n$ and quote from \cite{BKN-complexity}, thm 3.3.1, thm 3.3.2 
$$ M \in \calC^+  \Leftrightarrow V_{\g_{(1)}}(M) = 0 \quad , \quad   M \in \calC^-  \Leftrightarrow V_{\g_{(-1)}}(M) = 0 \ .$$ Hence $M$ is projective if and only if $V_{\g_{(1)}}(M)\! =\! V_{\g_{(-1)}}(M)\! =\! 0$ holds.

\bigskip

{\it Vanishing criterion}. For any $\xi \in X$ there exists $g \in Gl(n) \times Gl(n)$ and isotropic mutually orthogonal linearly independent roots $\alpha_1, \ldots, \alpha_k$ such that $Ad_g(\xi) = \xi_1 + \ldots + \xi_m$ with $\xi_i \in g_{\alpha_i}$. The number $m=r(\xi)$ is called the rank of $\xi$ \cite{Serganova-kw}. 
The orbits for the action of $Gl(n)\times Gl(n)$ on $\g_{(1)}$ are \cite{BKN-complexity}, 3.8.1 \[ (\g_{(1)})_m = \{ \xi \in \g_{(1)} \ | \ r(\xi) = m\} \ \ \text{ for } \ \  0 \leq m \leq n \ .\]
By a minimal orbit for the adjoint action of $Gl(n) \times Gl(n)$ on $g_{(\pm 1)}$ we mean a minimal non-zero orbit with respect to the partial order given by containment in closures. 
The unique  minimal orbit  
$(\g_{(1)})_1$ is the orbit of the element $x$ defined earlier. The situation is analogous for $\g_{(-1)}$, where $\overline x=\tau(x)$ generates the corresponding minimal orbit. A slight modification of \cite{BKN-complexity}, thm 3.7.1 and its proof gives  

\begin{thm} \label{kernel}  For $\xi \in \g_{(1)}$ and  $M\in \calC^-$ we have $M_\xi=0$.  For $\xi \in \g_{(-1)}$ and  $M\in \calC^-$ we have $M_\xi=0$.  For $\xi=x$ we have $DS(M) = M_x = 0$ if and only if $M\in \calC^-$
and $M_{\overline x}=0$ if and only if $M\in \calC^+$.  
\end{thm} 

\medskip{\it Proof}. Let $M \in \calC^-$. Then the definition of Kac objects implies $V_{\g_{(1)}}(M) = 0$. Hence $ \{\xi \in \g_{(1)} \ | \ M_\xi \neq 0\} = 0$. Conversely assume $M_x = 0$. 
Since
$V_{\g_{(1)}}(M)$ is a closed $Gl(n) \times Gl(n)$-stable variety, it contains a closed orbit. Since the orbits $(\g_{(1)})_m$ are closed only for $m=1$, unless $V_{\g_{(1)}}(M)$ is empty, it must contain $(\g_{(1)})_1$. But this would imply $M_x \neq 0$, a contradiction. Hence $V_{\g_{(1)}}(M) = \emptyset$. \qed 

\begin{cor} \label{van} For our fixed $x \in (\g_{(1)})_1$ 
\begin{enumerate}
	\item $M$ is projective if and only if $M_x = 0$ and $M_{\tau x} = 0$. 
	\item $M$ is projective if and only if $M_x = 0$ and $M^*_x = 0$
	\item If $M = M^*$, then $M$ is projective if and only if $M_x = 0$.
\end{enumerate}
\end{cor} 

\medskip{\it Proof}.  $M_x = 0$ implies $V_{\g_{(1)}}(M) = 0$ and $M_{\tau(x)} = 0$ implies $V_{\g_{(-1)}}(M) = 0$, hence (1). Now (2) and (3) follow from \cite{BKN-complexity}, 3.4.1 
using \[ V_{\g_{(\pm 1)}}(M^*) = \tau( V_{\g_{(\mp 1)}} (M)).\] \qed

\medskip



\section{ The tensor functor $D$} \label{DDirac}

\medskip
In this section we construct another tensor functor $H_D:T_n \to T_{n-1}$ which is defined as the cohomology of a complex given by a Dirac operator $D$. This tensor functor has the advantage that it is compatible with the twisted duality $^*$.

\medskip
In this section we assume $V \in T_n$. 
For $t\in k^*$ the diagonal matrices $$diag(E_{n-1},t,t,E_{n-1}) \in G_{\overline 0}$$ define a one dimensional torus, the center of $G_1$; for this recall
the embeddings $G_1=id \times G_1 \hookrightarrow G_{n-1}\times G_1 \hookrightarrow G_n$.  The center of $G_1$ commutes with $G_{n-1} \times id \subset G_n$. 
Hence the center of $G_1$
naturally acts  on $DS(V)$ in a semisimple way for any representation $(V,\rho) \in T$.
Hence the underlying vectorspace $V$ decomposes into $H$-eigenspaces for $H=diag(0_{n-1},1,1,0_{n-1})$ 
in $\g_n=Lie(G_n)$ which generates the Lie algebra of the torus.   

\medskip
Let  $x\in \g_n$ be the fixed nilpotent element specified in section \ref{DF}.
Let $\overline x = x^T$ denote the supertranspose of $x$.
Now $Ad(\epsilon_1)(H)=H$ and
$[H,x]=[H,\overline x]=0$ imply that the operators $\partial=\rho(x)$ and $\overline\partial= c \cdot \rho(\overline x)$ (for any $c \in k^*$)
commute with $H$. 
Furthermore $[x,\overline x] =  H$
for the odd elements $x$ and $\overline x$ implies 
$$ \partial \overline\partial + \overline\partial \partial = c \cdot \rho(H) \ .$$   
Since $H$ commutes with $x$, the operator $\rho(H)$ acts on $V_x$.
Since $H$ commutes with $\varepsilon_1$, the
grading $V^\pm$ is compatible with taking invariants \[ V^H = \{ v \in V \ | \ \rho(H) = 0 \}.\]
Similarly we denote the space of coinvariants by $V_H$. On $V$ the odd operator $\overline\partial$ defines a  homotopy
 of the complex
 $$  \xymatrix{ \ar[r]^-\partial & \Pi(V_{2\ell - 1}) \ar[r]^-\partial & V_{2\ell} 
 \ar[dl]_-{\overline\partial}\ar[r]^-\partial & \Pi(V_{2\ell + 1})  \ar[dl]_-{\overline\partial}\ar[r]^-\partial & V_{2\ell + 2}  \ar[dl]_-{\overline\partial}\ar[r]^-\partial &  \cr
\ar[r]^-\partial & \Pi(V_{2\ell - 1}) \ar[r]^-\partial & V_{2\ell} 
\ar[r]^-\partial & \Pi(V_{2\ell + 1}) \ar[r]^-\partial & V_{2\ell + 2} \ar[r]^-\partial & }$$
Hence $c \cdot \rho(H)$
is homotopic to zero. In particular, the natural action of
$\rho(H)$ on the cohomology modules $H^\ell(V)$ is trivial. Therefore 

\begin{lem}\label{homotopy}
$\rho(H)$ acts trivially
on the cohomology $DS(V)=V_x$.
\end{lem} 

\medskip
Since $H$ acts in a semisimple way, taking $H$-invariants $V \mapsto V^H$ is an exact functor and commutes
with the cohomology functor $V\mapsto V_x$. Thus
$$  DS(V) = M_x  \quad \text{ for } \quad M = V^H \ $$
and similarly $H^\pm(V) = H^\pm(V^H)$ etc. Notice $(V^H)^\pm = (V^\pm)^H$.
Since  the operators $\partial$ and $\overline\partial$ commute with $H$, they preserve $M=V^H$
and anti-commute on $M$. In this way we obtain a double complex for $M=V^H$ defined by
 $$ \xymatrix{ 
 \ar[r]^-{\overline\partial} & M^+ \ar[r]^-{\overline\partial} & \Pi(M^-) \ar[r]^-{\overline\partial} & \cr
  \ar[r]^-{\overline\partial} & \Pi(M^-) \ar[r]^-{\overline\partial} \ar[u]^-{\partial} & M^+ \ar[r]^-{\overline\partial} \ar[u]^-{\partial} &  \cr
 \ar[r]^-{\overline\partial}  & M^+ \ar[r]^-{\overline\partial} \ar[u]^-{\partial}& \Pi(M^-) \ar[r]^-{\overline\partial} \ar[u]^{\partial}&.  } $$

\medskip
{\it The Dirac operator}. This double complex is related to the complex 
$$ \xymatrix{ \cdots \ar[r]^-D &  M^+ \ar[r]^-D &  \Pi(M^-) \ar[r]^-D &  M^+  \ar[r]^-D &  \Pi(M^-) \ar[r]^-D & \cdots} $$
for $M = V^H$ attached to the Dirac operator
$$   D = \partial + \overline\partial \ .$$ 
Since $M = M^+ \oplus \Pi(M^-)$,
the two cohomology modules $H_D^+(V)$ and  $H_D^-(V)$ 
of this periodic complex compute
$$  H_D(V) = Kern(D: M \to M)/Im(D: M \to M) \ $$
in the sense that
$$H_D(V)= H_D^+(V) \oplus \Pi(H_D^-(V))\ $$
gives the decomposition of  $H_D(V)$ into its
$\calR_n$ and $\Pi(\calR_n)$-part. 

\medskip
{\bf Remark}. Note that $D$ commutes with $\rho(H)$. 
Hence the operator $D$ respects the eigenspaces of $H$ on $V$.   
Since $D^2 = \partial^2 + (\partial\overline \partial + \partial \overline\partial)
+ \overline\partial^2 = (\partial\overline \partial + \partial \overline\partial) = c \cdot \rho(H)$,
we have $Ker(D:V\to V)=Kern(D:V^H \to V^H)$. However $D(V)$ is in general different from
$D(V^H)$, although both spaces have the same intersection with $V^H$.  

\begin{lem} \label{*1} For $c = i$ there exist natural isomorphisms $H_D(V^*,\rho^*) \cong H_D(V,\rho)^*$, $H_D(V)^{\vee} \cong H_D(V^{\vee})$, $H_D^{\pm}(V^*) \cong H_D^{\pm}(V)^*$ and $H_D^{\pm}(V^{\vee}) \cong H_D^{\pm}(V)^{\vee}$ of $G_{n-1}$-modules. For short exact sequences in $\calR_n$ one obtains an exact hexagon in $\calR_{n-1}$ for
the functors $H^\pm_D$. 
\end{lem}

{\it Proof}. The assertion $H_D(V)^{\vee} = H_D(V^{\vee})$ follows since $H_D$ is a tensor functor by lemma \ref{H_D-tensor}. We calculate $\tau(x + i\overline x)= -(\overline x - i x)= i (x+ i\overline x)$,
since $\tau^2(x)=-x$. Now recall that $\rho^*(D)= \rho^\vee(\tau(D))= i\rho^\vee(D)$
is defined as endomorphism on $V^*=V^\vee$. Hence $H_D(V^*,\rho^*)$, by definition the cohomology of $\rho^*(D)$ on $(V^*)^H$, can be identified with
the space $$ Ker(i\rho^\vee(D): (V^\vee)^H \to (V^\vee)^H)/Im(i\rho^\vee(D):(V^\vee)^H \to (V^\vee)^H)\ .$$  Of course we can ignore the factor $i$, and identify this representation
with the representation on 
$$\bigl(Ker(\rho(D): V_H \to V_H)/Im(\rho(D):(V_H \to V_H)\bigr)^\vee$$ or hence
with $$H_D(V,\rho)^\vee= \bigl(Ker(\rho(D): V^H \to V^H)/Im(\rho(D):(V^H \to V^H)\bigr)^\vee \ ,$$ using the dual $(V_H)^\vee \to (V^H)^\vee$ of the natural morphism $V^H \to V_H$, which is  an isomorphism by the semisimplicity of $H$. Finally recall $H_D(V,\rho)^\vee = H_D(V,\rho)^*$ for the underlying representation spaces. This is an isomorphism of $G_{n-1}$-modules since $\tau$ restricts to the corresponding $\tau$ on $G_{n-1}$. \qed

\medskip
So from now on assume $c=i$. 
Then, in contrast to lemma \ref{-ell}, we obtain

\begin{lem} \label{*} There exist natural isomorphisms of functors
$\calR_n \to \calR_{n-1}$ 
$$ \fbox{$ \mu_V: H_D^\pm(V^*) \cong H_D^\pm(V)^* $} \ .$$   
\end{lem}

\medskip
{\it Proof}. It remains to show that the isomorphism $ \mu_V: H_D^\pm(V^*) \cong H_D^\pm(V)^*$
given above defines a natural transformation. For a $G_n$-linear map $f: V \to W$
the induced map $f^*: W^* \to V^*$ is nothing but the morphism $f^\vee: W^\vee \to V^\vee$,
using $V^*=V^\vee$ and $W^* = W^\vee$. This now easily shows that the above identifications
$\mu_V, \mu_W$ induce a commutative  diagram                                                                                                                                                                            
$$  \xymatrix{  H_D(W)^*  \ar[r]^{H_D(f)^*} &   H_D(V)^*  \cr          
  H_D(W^*) \ar[u]_{\mu_W} \ar[r]^{H_D(f^*)} &   H_D(V^*) \ar[u]_{\mu_V}  \cr}   $$     
\qed

{\bf Example}. Let $V$ be the Kac module $V({\bf 1})$ in $\calR_1$.
Then $DS(V^*)=0$ and $DS(V) = {\bf 1} \oplus \Pi({\bf 1})$. On the other hand $H_D(V)=0$
and $H_D(V^*)=0$.

\medskip
{\bf Remark}. It is not a priori clear how to define a Dirac analog of 
the modules $H^\ell(V)$.
Indeed
$\overline\partial$ and $\partial$ (in the sense of odd morphisms) 
satisfy $\overline\partial : V_\lambda \to V_{\lambda - \mu}$
and $\partial : V_\lambda \to V_{\lambda + \mu}$
for the odd simple weight $\mu$. Hence $\overline\partial : V_\ell \to V_{\ell-1}$
and $\partial : V_\ell \to V_{\ell+1}$ and therefore $D= \partial + \overline\partial$
does not simply shift the grading.  We adress this question in section \ref{Hodge}.

\medskip
{\it $H_D$ as a tensor functor}. Although taking $H$-invariants $V\mapsto M=V^H$ is not a tensor functor, 
$H_D$ is nevertheless a tensor functor.
To show this it is enough to restrict the representations $(V,\rho)$ to $G_1 \hookrightarrow G_n$.
Hence it suffices to show that the functor
$$ H_D: T_1 \to T_0 = svec_k \ $$
is a tensor functor.
$H$ generates the center of $\mathfrak{gl}(1|1)$ and $D^2= \rho(H)$. 
Hence $Kern(D) \subset V^H$. Since $H$ is semisimple, the Jordan blocks
of $D$ on $V$ (ignoring the grading!) are Jordan blocks $B_\lambda$ of length 1 except for the eigenvalue $\lambda=0$, where
they are either Jordan blocks $B_0$ of length 1 or Jordan blocks $R$ of lenght 2. 
Indeed the square of an indecomposable Jordan block of length $a$ and eigenvalue $\lambda$
is again an indecomposable Jordan block of length $a$ for $\lambda\neq 0$.
Since $D^2=\rho(H)$ is semisimple, this implies $\lambda=0$ and $a\leq 2$ for $a>1$.
By definition, for $V= \bigoplus_\lambda k_\lambda(V) \cdot B_\lambda \oplus k(V)\cdot R $ 
we have $H_D(V) = k_0(V) \cdot B_0$, if we ignore the grading. Now
$B_\lambda \otimes B_{\lambda'} = B_{\lambda
\pm\lambda'}$, where the sign depends on the parity of $B_\lambda$. 
Furthermore the characteristic polynomial of $D$ on  $R\otimes B_\lambda$ 
is $X^2 - \lambda^2$, hence $D$ has eigenvalue $0$ on  $R\otimes B_\lambda$
only for $\lambda=0$, in which case $R\otimes B_\lambda$ is isomorphic to $R$.   
Finally $R\otimes R \cong R^2$. Hence the only possible deviation from the tensor functor property for $H_D$ might come from tensor products $B_\lambda \otimes B_{\lambda'}$ where $\lambda\pm \lambda'=0$. In this case $H=\lambda^2 \cdot id$ 
on $B_\lambda$ and $B_{\lambda'}$, hence $H= 2 \lambda^2 \cdot id$ on $B_\lambda \otimes
B_{\lambda'}$. But the even operator $D^2$ then acts by $  2 \lambda^2 \cdot id$ on $B_\lambda \otimes
B_{\lambda'}$. Hence $D$ does not have the eigenvalue zero on $B_\lambda \otimes
B_{\lambda'}$ unless $\lambda=\lambda'=0$. Therefore $B_0\otimes B_0 \cong B_0$
is the only relevant case. Hence $H_D(V \otimes W) = k_0(V)k_0(W)\cdot B_0 =
k_0(V)\cdot B_0 \otimes k_0(W)\cdot B_0 = H_D(V) \otimes H_D(W)$. This remains true
if we also take into account gradings.

\begin{lem}\label{H_D-tensor}
$H_D: T_n \to T_{n-1}$ is a tensor functor.
\end{lem}

\section{ The relation between $DS(V)$ and $D(V)$}\label{DS-vs-D}

\medskip\noindent
For $(V,\rho)\in T_n$ the eigenvalue decomposition with respect to the small torus gives
a decomposition $$ V = \bigoplus_{\ell\in\mathbb Z} V_{\ell} $$ into $G_{n-1}$-modules $V_{\ell}$. 
Furthermore $\overline\partial$ and $\partial$ (in the sense of odd morphisms) 
satisfy $\overline\partial : V_\ell \to V_{\ell-1}$
and $\partial : V_\ell \to V_{\ell+1}$. In other words, they give rise
to morphisms $\overline\partial : \Pi^\ell(V_\ell) \to \Pi^{\ell -1}(V_{\ell-1})$
and $\partial : \Pi^{\ell }(V_\ell) \to \Pi^{\ell +1}(V_{\ell+1})$, hence
induce morphisms on $\bigoplus_{\ell\in\mathbb Z} H^\ell(V)$ which
shift the grading by $-1$ resp. $+1$.

\medskip
Since the generator $H$ of the center of $Lie(G_1)$ commutes with the small torus,
we obtain an induced decomposition for the invariant subspace $M= V^H \subseteq V$
$$  M \ = \ \bigoplus_\ell \ \Pi^{\ell}(M_\ell)  $$
for $\Pi^{\ell}(M_\ell) = M \cap V_\ell = (V_\ell)^H$.
Notice $M = M^+ \oplus \Pi(M^-)$ for $(V,\rho)\in \calR_n$, with $M^+$ and $M^-$ defined in $\calR_n$ by
$$ M^+ = \bigoplus_{\ell \in 2\mathbb Z} M_\ell \quad , \quad M^- = \bigoplus_{\ell \in1+2\mathbb Z} M_\ell \ .$$
The spaces $M_\ell$ are $G_{n-1}$-modules. 

\medskip
On $M$ the operators $\partial$ and $\overline\partial$ define even morphisms and 
they anticommute in the diagram below. Hence
we get a double complex $K=K^{\bullet,\bullet}$ in $T_{n-1}$ attached to $(V,\rho)$
$$ \xymatrix{ M_{\ell+2} \ar[r]^{\overline\partial} & M_{\ell+1} \ar[r]^{\overline\partial} & M_\ell \ar[r]^{\overline\partial} & M_{\ell-1} \cr
M_{\ell+1} \ar[u]^{\partial}\ar[r]^{\overline\partial}  & M_\ell \ar[u]^{\partial}\ar[r]^{\overline\partial} & M_{\ell-1} \ar[u]^{\partial}\ar[r]^{\overline\partial} & M_{\ell-2} \ar[u]^{\partial} \cr
M_{\ell}\ar[u]^{\partial}\ar[r]^{\overline\partial} & M_{\ell-1}\ar[u]^{\partial}\ar[r]^{\overline\partial} &  M_{\ell-2} \ar[u]^{\partial}\ar[r]^{\overline\partial} & M_{\ell-3} \ar[u]^{\partial}  } $$
with $K^{i,j} = M_{j-i}$. This double complex is periodic with respect
to $(i,j) \mapsto (i+1,j+1)$. The modules $K^{i,j}$ vanish for $j-i \notin [\ell_0,\ell_1]$
and certain $\ell_0,\ell_1 \in \mathbb Z$. 

\medskip
The associated single complex $(Tot(K),D)$ has the objects $Tot(K)^n =
\bigoplus_{i\in \mathbb Z} M_{n + 2i}$ and the differential $D=\partial + \overline\partial$.
The total complex therefore is periodic with $Tot^0(K)= M^+$ and $Tot^1(K)=\Pi(M^-)$
and computes the cohomology $H^n(Tot(K),D)= H_D^+(V)$ for $n\in 2\mathbb Z$
and $H^n(Tot(K),D)= H_D^-(V)$ for $n\in 1+2\mathbb Z$.

\medskip
On the total complex $(Tot(K),D)$ we have a decreasing filtration
defined by $F^p Tot^n(K) = \bigoplus_{r+s=n, r\geq p} K^{r,s}$.
This filtration induces decreasing filtrations on the cohomology
of the total complex 
$$  ... \supseteq F^p(H_D^\pm(V)) \supseteq F^{p+1}(H_D^\pm(V)) \supseteq ... $$
and a spectral sequence $(E_r^{p,q},d_r)$ converging to
$$   E_\infty^{p,q} \  = \  gr^p H^{p+q}(Tot(K),D) \ .$$
Indeed the convergence of the sequence follows from the fact that
the higher differentials $d_r: E_r^{pq} \to E_r^{p+r,q-r+1}$ vanish
for $2r - (q-p+1) > \ell_1 - \ell_0$. The $E_1$-complex of the spectral sequence
is the direct sum over all $q$ of the horizontal
complexes $E_1^{p,q} = (H^q_\partial(K^{p,\bullet}),\overline\partial)$.
For the various $q$ these complexes are the same up to a shift of the complex.
So, if we ignore this shift, these complexes are given by the natural
action of $\overline\partial$ on $\bigoplus_{\ell\in\mathbb Z} H^\ell(V)$
defining the complex 
$$        \xymatrix{ ... \ar[r]^-{\overline\partial} & H^{q+1}(V)
\ar[r]^-{\overline\partial} & H^q(V)
\ar[r]^-{\overline\partial} & H^{q-1}(V) \ar[r]^-{\overline\partial} & ... } \ .$$ 
The decreasing filtration $F^p$ induced on 
$$E_1(K)^n \ =\ \bigoplus_{i\in \mathbb Z} H^{n+2i}(V) $$
has graded terms $gr^p(E_1(K)^n) = H_\partial(K^{p,n-p})= H_\partial(M_{n-2p})
= H^{n-2p}(V)$. We now define the subquotient 
$  H_D^{n-2p}(V) \ := \ gr^p(E_\infty(K)^n) $
of $H^{n-2p}(V)$, hence
$$  H_D^{\ell}(V) \ := \ gr^p(E_\infty(K)^{\ell + 2p}). $$
 Note that this definition does not 
depend on the choice of $p$.
We thus obtain

\begin{lem} For $T\in T_n$
the cohomology modules $H^\pm_D(V)$ admit canonical decreasing
filtrations $F^p$ whose graded pieces are the $G_{n-1}$-modules $H_D^{-2p}(V)$
for $H^+_D(V)$ and $H_D^{-2p-1}(V)$ for $H^-_D(V)$.
\end{lem}

\medskip
{\bf Condition {\tt T}}. {\it We say that condition {\tt T} holds for $(V,\rho)$ 
in $T_n$ if the  natural operation of the operator $\overline\partial=\rho(\tau(x))$ on $DS(V,\rho)$ is trivial}. 

\medskip{\bf Example}. The standard representation  $X=X_{st}$ of $G_n$ on $k^{n\vert n}$ satisfies condition {\tt T}.

\medskip
{\bf Remark}. If $\tau(x)$ act trivially both on $DS(V)$ and $DS(W)$ for some 
$V,W \in T_n$, then $\tau(x)$ acts trivially on $DS(V\otimes W)= DS(V) \otimes DS(W)$.
If $\tau(x)$ acts trivially on $DS(U)$ for 
$U \in T_n$, then $\tau(x)$ act trivially on every retract of $DS(U)$.
Hence condition {\tt T} for $(V,\rho)=L(\lambda)$ implies condition {\tt T} for every retract
$U$ of $DS(V \otimes W)$. Thus the subcategory of objects in $\calR_n$ satisfying condition {\tt T} is closed under tensor products and retracts.

\medskip
Now consider the following conditions for $(V,\rho)$:
\begin{enumerate}
\item $(V,\rho)$ is irreducible.
\item $H^+(V) \oplus H^-(V)$ is multiplicity free. 
\item $H^+(V)$ and $H^-(V)$ do not have common constituents.
\item Condition $\tt T$ holds.
\item $\overline\partial$ acts trivially on $DS(V)$.
\item The $E_1^{p,q}$ and the $E_2^{p,q}$ terms 
of the spectral sequence coincide $$H^{}_{\overline\partial}\bigl(H^{\ell}(V)\bigr) =
H^{\ell}(V)\ $$
where $\ell:=n-2p=q-p$.  
\end{enumerate}
Later in theorem \ref{mainthm} we prove that (1) implies (2). Furthermore it is trivial that $(2) \Longrightarrow (3) 
\Longrightarrow (4) \Longrightarrow (5) \Longrightarrow (6)$.

\begin{prop}
If condition (3) holds, then the spectral sequence degenerates at the $E_1$-level
and $H_D^{\pm}(V)$ is naturally isomorphic
to $ H^\pm(V)$.
\end{prop}

\medskip
{\it Proof}. The differentials of the spectral sequence $d_r: E_r^{pq} \to E_r^{p+r,q-r+1}$ 
define maps from the subquotient $E_r^{pq}$ of 
$H^{n-2p}(V)$ (for $n=p+q$) to the subquotient $E_r^{p+r,q-r+1}$
of $H^{n-2p- 2r+1}(V)$. If $H^{n-2p}(V)$
contributes to $H^\pm(V)$, then $H^{n-2p- 2r+1}(V)$ 
contributes to $H^\mp(V)$. Since all the higher differentials are $G_{n-1}$-linear,
condition (3) forces all differentials $d_r$ to be zero for $r\geq 1$. 
Hence the spectral sequence degenerates at the $E_1$-level.
\qed

\begin{prop}\label{abutment}
The spectral sequence always degenerates at the $E_2$-level, i.e. for all objects $(V,\rho)$ in $T_n$ we have
$$ \fbox{$ H_{\overline\partial}(H^\ell(V)) \ \cong \ H^{\ell}_{D}(V)  $} \ . $$
\end{prop}

\begin{cor} The kernel of $H_D:T_n \to T_{n-1}$ contains $\mathcal{C}^+ \cup \mathcal{C}^-$.
\end{cor}

\medskip {\bf Remark}. It seems plausible that the kernel equals $\mathcal{C}^+ \cup \mathcal{C}^-$.

\medskip{\it Proof}. This is a general assertion on spectral
sequences arising from a double complex   
$K$ such that
$K^{i,j}=M_{j-i}$ for maps $\partial: M_\ell \to M_{\ell +1}$ and $\overline\partial: M_{\ell}
\to M_{\ell -1}$ between finite dimensional $k$-vectorspaces $M_\ell, \ell\in\mathbb Z$
so that $M_\ell=0$ for almost all $\ell$. Indeed, any such double complex $K$ can be viewed
as an object in the category $T_1$ via the embedding $\varphi_{n,m}$ of section \ref{DF}.. Using $T_1= \calR_1 \oplus \Pi(\calR_1)$ we can decompose and assume without restriction of generality that it is an object in $\calR_1$.
However, then it defines a  maximal atypical object in the category $\calR_1^1 \subset \calR_1$. 
For this notice that $\calR_1^1$ can be identified
with the category of objects in $\calR_1$ with trivial central character.
Note that this condition on the central character for a representation $(V,\rho)$ of $G_1$
simply means $V = V^H=M$, since $H$ generates the center of $Lie(G_1)$.
This reduces  our claim to the special case $n=1$
for $(V,\rho)$ in $\calR_1^1$. Obviously we can assume
that $(V,\rho)$ is indecomposable. 

\medskip
The indecomposable objects $V$ in $\calR_1^1$ were classified by Germoni \cite{Germoni-sl}.
Either $V\in \calC^{+}$ (Kac object), or 
$V\in \calC^-$ or there exists an object $U \subset V, U\in\calC^-$ with irreducible
quotient $L$ or there exists a quotient $Q$ of $V$ in $\calC^-$ with irreducible kernel $L'$.
Since $DS(N)=0$ for all objects in $\calC^-$ (theorem \ref{kernel}), we conclude from the long exact sequence of $H^\ell$-cohomology that  
we can either assume $V\in \calC^+$ or that $V$ is irreducible, since in the remaining cases $DS(V)=0$ or $DS(V)\cong DS(L)$ or $DS(L') \cong DS(V)$. 
As already mentioned, by the later theorem \ref{mainthm} for irreducible $V$,
the spectral sequence already abuts. For $r=1$ however this is obvious anyway, since
any atypical irreducible $L$ is isomorphic to a Berezin power $L \cong Ber^m$. Hence $H_D^\nu(L) = H^\nu(L)= k$ for $\nu=m$ and  $H_D^\nu(L) = H^\nu(L)= 0$ otherwise.  

\medskip
So it remains to consider the case of indecomposable
Kac objects $V\in\calC^+ $ in $\calR_1^1$. Unless $V\in \calC^+ \cap \calC^-$, by Germoni's
results $V \cong V(i;m)$ for $i\in \mathbb Z$ and $m\in\mathbb N$ is a successive extension 
$$  0 \to V(i-2;m-1) \to V(i;m) \to V(Ber^i) \to 0 $$
of the Kac objects with $V(i;1)=V(Ber^i)$. Furthermore the Kac module $V(Ber^i)$  is an extension of
Berezin modules
$$        0\to   Ber^{i-1} \to             V(Ber^i) \to Ber^i \to 0 \ ,$$
hence $H^\ell(V(Ber^i))\cong k$ for $\ell= i,i-1$
and is zero otherwise. From the long exact cohomology sequence and induction we obtain $\dim(H^\nu(V(i;m))=1$ for $\nu\in\{i,i-1,...,i-2m+1\}$, and $H^\nu(V(i;m))=0$
otherwise. So (for fixed $q$) the complexes in the $E_1$-term of the spectral
sequence for $V=V(i;m)$ have the form 
$$ 0 \to H^i(V) \to H^{i-1}(V) \to \cdots \to H^{i-2m+2}(V) \to H^{i-2m+1}(V) \to 0 \ $$ 
with differentials $\overline\partial$ and $H^\nu(V)$ of dimension one for $\nu=i,i-1,...,i-2m+1$. 
We have to show that these complexes are acyclic for all $V=V(i;m)$. 
For this it suffices that the first differential $\overline\partial: H^i(V)\to H^{i-1}(V)$, 
the third differential $\overline\partial: H^{i-2}(V)\to H^{i-3}(V)$ and so on, are injective.
By dimension reasons the differentials  $\bar{\partial}: H^{i-1}(V)  \to H^{i-2}(V)$, $\bar{\partial}: H^{i-3}(V)  \to H^{i-4}(V)$ etc. are then isomorphisms and the differentials $\bar{\partial}: H^{i}(V)  \to H^{i-1}(V)$, $\bar{\partial}: H^{i-2}(V)  \to H^{i-3}(V)$ etc. are zero. Hence the cohomology of this complex vanishes and the $E_2$-term of the spectral sequence is zero.. Hence the spectral sequence abuts at $r=2$,
which proves our claim. 

\medskip
To prove the injectivity  for the first, third and so on differential $\overline\partial$ 
we use induction on $m$. For $m=1$
and $V=V(Ber^i)\in\calC^+$ we know $H_D(V)=0$ by theorem \ref{kernel}
and lemma \ref{*} .
Since $H^\nu(V)=0$ for $\nu\neq i, i-1$, all higher differentials $d_r$  for $r\geq 2$ are
zero by degree reasons. Hence $\overline\partial: H^i(V) \to H^{i-1}(V)$
must be an isomorphism.

\medskip
For the induction step put $V_i:=V(i,1)$ and 
$N=V(i-2,m-1)$; then $V/N \cong V_i$. Hence
we get a commutative diagram
with horizontal exact sequences  $$
\xymatrix@C-0.5cm{ ... \ar[r] & H^{\nu -1}(N) \ar[r]\ar[d]_{\overline\partial} & H^{\nu -1}(V) \ar[d]_{\overline\partial}\ar[r] & H^{\nu -1}(V_i)\ar[d]_{\overline\partial}  \ar[r] &  H^\nu(N) \ar[d]_{\overline\partial} \ar[r] & ... \cr
... \ar[r] &  H^{\nu-2}(N) \ar[r] & H^{\nu-2}(V) \ar[r] & H^{\nu-2}(V_i) \ar[r] &  H^{\nu-1}(N) \ar[r] & ...  } $$
Since $H^\nu(N)=0$ for $\nu >i-2$ and $H^\nu(V_i)=0$ for $\nu\neq i,i-1$
 $$
\xymatrix@C-0.5cm{ 0\ar[r] & H^{i-1}(N) \ar[r]\ar[d]_{\overline\partial} & H^{i-1}(V) \ar[d]_{\overline\partial}\ar[r] & H^{i-1}(V_i)\ar[d]_{\overline\partial}  \ar[r] &  0 \ar[d]_{\overline\partial}\ar[r] & H^i(V) \ar@{^{(}->}[d]_{\overline\partial}\ar[r]^\sim & H^i(V_i)\ar@{^{(}->}[d]_{\overline\partial} \ar[r] & 0 \cr
0\ar[r] &  H^{i-2}(N) \ar[r] & H^{i-2}(V) \ar[r] & H^{i-2}(V_i) \ar[r] &  0 \ar[r] & H^{i-1}(V) \ar[r]^\sim & H^{i-1}(V_i) \ar[r] & 0 } $$
Thus $\overline\partial: H^i(V) \to H^{i-1}(V)$ is injective by a comparison with $V_i$.
The assertion for the third, fifth and so on differential $\overline\partial$
follows from the induction assumption on $N$, since $H^\nu(V) \cong H^\nu(N)$
for $\nu\leq i-2$.  \qed

\medskip\noindent



\section{ Hodge decomposition}\label{Hodge}

\medskip
We show in proposition \ref{H1} that the groups $H_D^{\pm}(V)$ satisfy a Hodge decomposition.
\medskip
Put $F_p = F^{-p}$.
This defines a decreasing filtration of $G_{n-1}$-modules $F_p(H_D^\pm(V))$  on $H_D^\pm(V)$ as in the last section for $V\in T_n$. Here 
$$ F_p(H_D^\pm(V)) \ = \ Im\bigl( (\bigoplus_{\ell \leq 2p} M^\pm_\ell ) \cap Ker(D) \to H_D^\pm(V) \bigr) \ .$$
One has also  a decreasing filtration of $G_{n-1}$-modules $\overline F_q(H_D^\pm(V))$  on $H_D^\pm(V)$ defined by  the second filtration of the cohomology of  $(Tot(K),D)$ for
the double complex $K^{\bullet,\bullet}$ defined in the last section.
It is defined by the subcomplexes $\overline F^q(Tot(K)^n) = \bigoplus_{r+s=n, s\geq q} K^{r,s}$
of $(Tot(K),D)$. Notice that 
$$ \overline F^q(H_D^\pm(V)) $$
is the image of the $D$-cohomology of this subcomplex in $H_D(K)$.
This filtration has analogous properties. In particular
$$  \overline H^{n-2q}_D(V) := \overline F^q(Tot^n(K))/\overline F^{q+1}(Tot^n(K)) $$
by an analog of proposition \ref{abutment} is isomorphic to
$$  \overline H^\ell_D(V) \cong H_\partial(\overline H^\ell(V)) $$
where $\overline H^\ell(V)$ is defined as $H^\ell(V)$, only by using $\overline\partial$
instead of $\partial$.  

\medskip
We remark that both
filtrations are functorial with respect to morphisms $f: V \to W$ in $T_n$.
Hence also $V \mapsto \overline F^q(H^n_D(V)) \cap F_p(H^n_D(V))$ defines a functor
from $T_n$ to $T_{n-1}$.

\begin{prop}\label{H1}
For all objects $V$ in $T_n$ we have a canonical decomposition of $H_D^\pm(V)$ into
$G_{n-1}$-modules
$$   H_D^\pm(V) \ = \ \bigoplus_{\nu\in \mathbb Z} H_D^\nu(V) $$
where for $\varepsilon = (-1)^\nu$
$$ H_D^\nu(V) :=  F_\nu(H_D^\varepsilon(V)) \cap  \overline F_\nu(H_D^\varepsilon(V)) \ .$$
Furthermore for $\mu > \nu$ we have  
$$ F_\nu(H_D^\pm(V)) \cap  \overline F_\mu(H_D^\pm(V)) = 0\ .$$
\end{prop}
\smallskip

\begin{cor}\label{H2}
For a short exact sequence $0\to A \to B\to C\to 0$
in $T_n$ the sequences
$$   H_D^\nu(A) \to H_D^\nu(B) \to H_D^\nu(C) $$
are exact for all $\nu$. 
\end{cor}

\medskip
{\bf Remark}. As shown after lemma \ref{Ia} these halfexact sequences can not be
extended to long exact sequences! 

\medskip
{\it Proof}. If  $x \in H_D^\nu(B)$ maps to zero in $H_D^\nu(C) \subset H_D(C)$
there exist $y\in H_D(A)$ such that $x$ is the image of $y$ by the exact hexagon
for $H_D^\pm$. But then, for the decomposition $y = \sum_\nu y_\nu$
and $y_\nu \in H_D^\nu(C)$ given in proposition \ref{Hodge}, the components $y_\nu$ also maps to $x$ by the functoriality of $H_D^\nu(.)$. \qed

\medskip
{\it Proof of proposition \ref{H1}}. 
As in the proof of proposition \ref{abutment} we can reduce to the case
of an indecomposable object $V$ in $\calR_1^1$. For such  $V$ either 
$H_D(V)=0$, in which case the assertion is trivial, or $V$ is of the form 
$$        0 \to    L \to        V \to Q  \to 0 $$
with irreducible $L$ and $Q\in \calC^-$
or of the form
$$ 0 \to U \to V \to L \to 0 $$
with irreducible $L$ and $U\in \calC^-$.
These two situation are duals of each other. So we 
restrict ourselves to the first case. The irreducible module
$L$ is isomorphic to $Ber^m$ for some $m\in \mathbb Z$.
Then according to \cite{Germoni-sl} the quotient module $Q$
has socle and cosocle
$$  socle(Q) = \bigoplus_{i=1}^s  Ber^{m+2i} $$
$$ cosocle(Q)  = \bigoplus_{i=1}^s  Ber^{m+2i-1} \ .$$ 
Recall $H^m(V) = H^\bullet(V)$, and hence $H_D(V)=H^m(V)$. 
Hence by the abutment of the spectral sequences
$$    H_{\overline\partial}(H^\nu(V)) = H^\nu(V) \cong k $$  
for $\nu = m$ and is zero otherwise. 
Hence $H^m(V) \cong H_D(V)$, since the filtration $F^q$
only jumps for $p=m$. Similarly
$$    H_\partial(\overline H^\nu(V)) = \overline H^\nu(V) \cong  k $$  
for $\nu = m$ and is zero otherwise. 
Hence $\overline H^m(V) \cong H_D(V) $, since the filtration $\overline F^q$
only jumps for $q=m$. This simultaneous jump shows
$$ H_D^m(V) = \overline F^m(H_D(V)) \cap F_m(H_D(V)) $$
and also for $q > p$.
$$ \overline F^q(H_D(V)) \cap F_p(H_D(V)) \ = \ 0 \ .$$
 \qed



\section{ The case $m>1$}\label{m}

\medskip
As the diligent reader may have observed, the results 
obtained in the last sections on the functor $DS$ carry over to the case of the
more general functors $DS_{n,n-m}$. For this fix $m \geq 1$. 
The enriched weight structure of $DS_{n,n-m}$ (which depends on $m$)
is obtained from the decomposition of $(V,\rho) \in T_n$ into eigenspaces with respect to
the eigenvalues $t^\ell$ under the elements $\varphi_{n,m}(E \times diag(1,t^{-1}))$ of the small torus. This allows to
give a decomposition
$$  DS_{n,n-m}(V) = \bigoplus_\ell \ DS_{n,n-m}^\ell(V)[-\ell] \ $$
into eigenspaces $\Pi^\ell(DS_{n,n-m}^\ell(V))$ and gives long exact sequences in $T_{n-m}$ attached to  short exact sequences
in $T_n$ as in section \ref{DF} . Furthermore
lemma \ref{-ell} and lemma \ref{+ell} carry over verbatim. Notice,
$$DS_{n,n-m}^\ell(Ber_n)= Ber_{n-m} \quad , \quad \mbox{ for } \ell=m $$ and it is zero for $\ell\neq m$.
Indeed,  $\varphi_{n,m}(E \times diag(1,t^{-1}))$ acts on $Ber_n$ by $t^m$.
Since $\epsilon_n =  \varphi_{n,m}(E \times diag(1,-1)) \epsilon_{n-m}$, the restriction
of $Ber_n$ to $G_{n-m}$ via $\varphi_{n,m}$ defines the module $\Pi^m(Ber_{m-n})$. 

\medskip

{\bf Remark}. Note that $DS_{n,n-m}^\ell(V)[-\ell]  =  (DS_{n,n-m}(V))_\ell$. Here upper indices denote graduations without twist, lower indices graduations with twist. This is consistent with $DS_{n,n-1}^\ell(V)=H^\ell(V)$ in section \ref{sec:cohomology-functors}. Note that it is essential to have non twisted objects such as $H^\bullet(V)$ or more generally $DS^\bullet(V)$ due to the comparison with $H_D$ respectively $\omega_{n,n-m}$ (see below) which don't have any twists. 

\medskip
For $n-m_1=n_1$ and $n_1-m_2=n_2$ the functors
$DS_{n,n_1}: T_n \to T_{n_1}$ and $DS_{n_1,n_2}: T_{n_1} \to T_{n_2}$
are related to the functor $DS_{n,n_2}: T_n \to T_{n_2}$ by a Leray type
spectral sequence with the $E_2$-terms
$$ \fbox{$ \bigoplus_{p+q=k} DS^p_{n_1,n_2}(DS_{n,n_1}^q(V)) \ \Longrightarrow
\ DS^k_{n,n_2}(V) $} \ .$$
To be more precise, choose matrices
$$  J \ = \ \begin{pmatrix} 0 & J_2 \cr J_1 & 0 \end{pmatrix} $$ and 
$m_i\times m_i$-matrices $J_i, i=1,2$ with zero enties except for the entries 1 in the antidiagonal. 
Then $J$ and $J_1$ define functors $DS_{n,n_2}(V,\rho)=(V,\rho)_x$ resp. $DS_{n,n_1}(V,\rho)=(V,\rho)_{x_1}$
and $J_2$ defines a functor $DS_{n_1,n_2}(W,\psi)=(W,\psi)_{x_2}$. Obviously
we have $x=x_1 + x_2 \in \g_n$ and $x_2\in \g_{n_1} \subset \g_n$ such that $[x_1,x_2]=0$.  
Then indeed $\partial = \partial_1 + \partial_2$ and $\partial_1 \partial_2 +
\partial_2 \partial_1 =0$ for $\partial = \rho(x), \partial_1=\rho(x_1)$ and $\partial_2=\psi(x_2)$.
Consider the weight (eigenvalue) decomposition $$  V \ = \ \bigoplus_{p,q\in\mathbb Z} V^{p,q} $$
of $(V,\rho)$ with respect to the matrices 
$$g(t_1,t_2)=diag(1,..,1;t_1^{-1},...,t_1^{-1},t_2^{-1},..,t_2^{-1},1,..,1)$$ 
in $G_n$ ($m_1$ entries $t_1^{-1}$ and $m_2$ entries $t_2^{-1}$) so that
$v\in V^{p,q}$ if and only if $ g(t_1,t_2)v = t_1^q t_2^p \cdot v$. (We now write indices on top
to avoid confusion with the lower indices $n$ and $n-m$).
 Then  $\partial_2: V^{p,q} \to V^{p+1,q}$ and
$\partial_1: V^{p,q} \to V^{p,q+1}$. Hence the Leray type spectral sequence
is obtained by the spectral sequence of this double complex. For this note that the the functors $D_{n,n_2}^k$ are defined by the eigenvalues $t^k$ of the elements $g(t,t)$.  

\begin{prop} \label{Leray}
For irreducible maximal atypical objects $L$ in $T_n$ the Leray type
spectral sequence degenerates:
$$  \fbox{$ DS_{n,n_2}(L) \ \cong \ DS_{n_1,n_2}(DS_{n,n_1}(L)) $} \ .$$
\end{prop}

\medskip
{\it Proof}. Up to a parity shift,  we can replace $L=L(\lambda)$ by $X_\lambda$ in $T_n$, so that $sdim(X_\lambda) >0$ using that $sdim (X) \neq 0$ \cite{Serganova-kw}, \cite{Weissauer-gl}. Then it suffices to prove inductively (for $DS$ applied $m$ times)
$$   (DS\circ DS .... \circ DS)(X_\lambda) \ \cong \ DS_{n,n-m}(X_\lambda) \ .$$
The case $m=1$ is obvious by definition, since $DS_{n,n-1}=DS$.  
Suppose this assertion holds for $m$. Let us show that it then also holds for $m$ replaced by 
$m+1$. Indeed, the $E_2$-term of the spectral sequence 
$$  DS\circ (DS\circ DS .... \circ DS)(X_{\lambda}) \Longrightarrow D_{n,n-m-1}(X_\lambda) \ $$
are of the form
$$  DS\circ (DS\circ DS .... \circ DS)(X_{\lambda}) \ \cong \bigoplus_\mu X_\mu $$
for irreducible representations $X_\mu$ in $T_{n-m-1}$ of superdimension $sdim(X_\mu)>0$.
Indeed this follows by repeatedly applying the later theorem \ref{mainthm}, which implies $DS(X_\lambda) \cong \bigoplus_{i=1}^k X_{\lambda_i}$ for irreducible maximal atypical objects $X_{\lambda_i}$ in $T_{n-1}$ with $sdim(X_{\lambda_i}) >0$. 
Now $DS$ is a tensor functor, and hence preserves superdimensions.
Hence $sdim(X_\lambda) = \sum_\mu sdim(X_\mu)$. If the spectral sequence would not degenerate at $E_2$-level, then the $E_\infty$-term is a proper subquotient of the semisimple
$E_2$-term. Hence $$ sdim(DS_{n-m-1})(X_\lambda) < \sum_\mu sdim(X_\mu)\ ,$$ since $sdim(X_\mu) >0$. This would imply $$sdim(DS_{n-m-1})(X_\lambda) < sdim(X_\lambda) \ .$$
However this is a contradiction, since $D_{n-m-1}$ is a tensor functor
and hence $sdim(DS_{n-m-1}(X_\lambda)) = sdim(X_\lambda)$. Hence the spectral sequence degenerates and $DS_{n,n_2}(X_{\lambda})$ has a filtration with graded pieces that are computed by appropriate $DS_{n_1,n_2}^p DS_{n,n_1}^q (X_{\lambda})$. In order to prove that $DS_{n,n_2}(X_{\lambda})  = DS_{n_1,n_2}(DS_{n,n_1}(X_{\lambda}))$ we show that this filtration splits and $DS_{n,n_2}(X_{\lambda})$ is semisimple.
This follows from the sign rules of the main theorem \ref{mainthm}. Indeed for $X_{\lambda}$ with $\varepsilon(\lambda) = 1$, the constituents of $DS(X_{\lambda})$ in $\calR_n$ have the same sign $\varepsilon = 1$  and the constituents of $DS(X_{\lambda})$  in $\Pi \calR_n$ have sign $\varepsilon = -1$. Now use that $Ext^1(\calR_n,\Pi \calR_n) = 0$ by lemma \label{thm:decomposition} and $Ext^1(L(\lambda),L(\mu)) = 0$ if $\varepsilon(\lambda) = \varepsilon(\mu)$ by corollary \ref{semisimple-sign}. Hence there are no extensions between the constituents of $DS(X_{\lambda})$. Repeated application of $DS$ gives again constituents which are either in $\calR_n$ with sign $\varepsilon = 1$ or constituents in $\Pi \calR_n$ with sign $\varepsilon = -1$. Since the constituents of $DS_{n,n_2}(L)$ are given by the constituents of the graded pieces, the semisimplicity follows. \qed   
 
\medskip
We have seen in the last proposition that the Leray type spectral sequence
degenerates at the $E_2$-level for irreducible maximal atypical objects. 
Let $F^p$ be the decending first (or second) filtration of the total complex.
Due to the degeneration we can make use of the following lemma.

\begin{lem}
Suppose given a finite double complex $(K^{\bullet,\bullet}, d_{hor},d_{vert})$ with associated total
complex $K^\bullet = Tot(K^{\bullet,\bullet})$ and total differential $d$. Suppose
the associated spectral sequence for the first (second) filtration degenerates 
at the $E_2$-level and suppose
$x\in F^p(K^\bullet)$ is a boundary in $K^\bullet$. Then there exists $y\in F^{p-1}(K^\bullet)$
such that $x=dy$. 
\end{lem}

\medskip
{\it Proof}. We can assume that $x=\sum_{i=p}^{\infty} x_{p,n-p}$ has fixed degree
$n$.  The spectral sequence degenerates at $E_2$ and $[x]=0$ in $F^pH^n(K^\bullet)$.
Hence the class of $x$ in $Gr^p(H^n(K^\bullet)) = H_{hor}^p(H^{n-p}_{vert}(K^{\bullet,\bullet}))$
vanishes. In other words
there exists $v\in K^{p,n-p-1}$ and $u\in K^{p-1,n-p}$ such that $d_{vert}(u)=0$ and such that
$d_{hor}(u) + d_{vert}(v) = x_{p,q}$. Hence $x- d(u+v) \in F^{p+1}(K^\bullet)$ with $u+v \in F^{p-1}(K^\bullet)$ again is closed.
Iterating this argument we conclude that for any $r$ large enough we find $y\in F^{p-1}(K^\bullet)$ such that $x - dy \in F^r(K^\bullet)$.  If $r$ is large enough, then $F^r(K^\bullet)=0$ and hence the claim follows. \qed
 
\medskip
{\it Dirac cohomology}. Similarly the results of section \ref{DDirac} hold verbatim for $\partial = \rho(x)$ and $\overline\partial = c \cdot \rho(\overline x)$
and $D=\partial + \overline\partial$. In particular,  for a generator $z$ of the Lie algebra of the center of $G_1$, let $H$ denote its image $\varphi_{n,m}(z)\in \g_n$.
The $D$-cohomology of the fixed space $V^H$ then 
gives objects $\omega_{n,n-m}(V,\rho)$ so that
so that 
$$  \omega_{n,n-m}: T_n \to T_{n-m} $$
defines a tensor functor generalizing $H_D=\omega_{n,n-1}$. Note that $\omega_{n,n-m}$ restricts to a tensor functor $\omega_{n,n-m}:\mathcal{R}_n \to \mathcal{R}_{n-m}$ unlike $DS_{n,n-m}$.

\medskip
As in section \ref{DS-vs-D} there is a spectral sequence that allows to
define a filtration on $\omega_{n,n-m}(V)$ whose graded pieces
are $$ \fbox{$ \omega_{n,n-m}^\ell(V) \ \cong \ H_{\overline\partial}(DS_{n,n-m}^\ell(V)) $}  \ .$$ 
This generalizes proposition \ref{abutment}. Furthermore
the results of proposition \ref{H1} and corollary \ref{H2} of section \ref{Hodge}
carry over and define a Hodge decomposition for $\omega_{n,n-m}$ in terms 
of the functors $\omega_{n,n-m}^\ell$. Finally the same
argument used in the proof of proposition \ref{Leray} also shows

\begin{prop} \label{Leray2}
For irreducible maximal atypical objects $L$ in $T_n$ the
spectral sequence  above degenerates, i.e. for all $\ell$
$$  \fbox{$  \omega_{n,n-m}^\ell(L)\ \cong \  DS_{n,n-m}^\ell(L) $} \ .$$
\end{prop}
 
Now consider the $\mathbb Z$-graded object $DS_{n,n-m}^\bullet(L)
=\bigoplus_{\ell\in \Z}\ DS_{n,n-m}^\ell(L)$ (which is different from $DS_{n,n-m}(L)$ if we forget the graduation)  to compare with  
$\bigoplus_{\ell\in \Z}\ \omega_{n,n-m}^\ell(L)$ (that is $\omega_{n,n-m}(L)$
after forgetting the graduation).

\begin{lem}\label{insteps}
Suppose for irreducible $V\in T_n$ that $DS_{n,n_1}^\bullet(V) \cong \omega_{n,n_1}(V)$. Then 
$$  \omega_{n,n_2}^\bullet(V) \cong \omega_{n_1,n_2}^\bullet(DS_{n,n_1}^\bullet(V)) $$
holds. 
\end{lem} 
 
\medskip
{\it Proof}. Use that $\overline \partial = \overline\partial_1 + \overline\partial_2$.
By the assumption $DS_{n,n_1}^\bullet(V) \cong \omega_{n,n_1}^\bullet(V)$ the differential
$\overline\partial_1$ is trivial on $DS_{n,n_1}^\bullet(V)$, hence
trivial on $DS^\bullet_{n_1,n_2}(DS^\bullet_{n,n_1}(V)) \cong DS^\bullet_{n,n_2}(V)$. Therefore
the $\overline\partial$-homology of $DS_{n,n_2}^\bullet(V)$ is the same as the
$\overline\partial_2$-homology attached to $
DS_{n_1,n_2}^\bullet(DS_{n,n_1}(V))^\bullet$. \qed

\medskip
This implies $\omega_{n,n-m}^\bullet(L)\ \cong \  DS_{n,n-m}^\bullet(L)$ 
for any irreducible $L$ in $T_n$. We prove this
by induction on $m$. For $m=1$ this follows from the fact that
irreducible representations satisfy property $\tt T$. Now we use
$ \omega_{n,n-m-1}^\bullet(L) \cong \omega_{n-m,n-m-1}^\bullet(DS_{n,n-m}^\bullet(L))$
from lemma \ref{insteps}. Since $DS_{n,n-m}^\bullet(L)$ is semisimple
by proposition \ref{Leray} (as iteration of $m$ times $DS^\bullet$), we 
have $$ \omega_{n-m,n-m-1}^\bullet(DS_{n,n-m}^\bullet(L)) = DS^\bullet(DS_{n,n-m}^\bullet(L))
= DS_{n,n-m-1}^\bullet(L)\ .$$ This implies

\begin{prop} \label{Leray3}
For all irreducible objects $L$ in $T_n$ and all $\ell$ we have 
$$  \fbox{$  \omega_{n,n-m}^\ell(L)\ \cong \  DS_{n,n-m}^\ell(L) $} \ .$$
\end{prop}
 
\bigskip\noindent
The case $m=n$ is of particular interest. Notice that $T_0$ is the
category $svec_k$ of finite dimensional super $k$-vectorspaces.
Hence
$$  \omega=\omega_{n,0}: T_n \ \longrightarrow \ svec_k \ .$$

\medskip
{\it The tori $A_i$}.
Let $A_i\subseteq G_n$ denote the diagonal torus of all elements
of the form $diag(1,\ldots,1,t_{n-i+1},\ldots,t_n \ | \ t_n,\ldots, t_{n-i+1},1,\ldots,1)$. In particular $A_1$ is the torus of section 5 and $H = H_{n,n-1}$. It commutes with all
operators $\partial_{n,n-i}, \overline\partial_{n,n-i}$ and $D_{n-i}$
and hence acts on $DS_{n,n-i}(V)$ respectively $\omega_{n,n-i}(V)$.
We claim

\begin{lem} The action of $A_i$ on $DS_{n,n-i}(V)$ and $\omega_{n,n-i}(V)$
is trivial. 
\end{lem}

\medskip
{\it Proof}. 
For this we can assume without loss of generality
that $i=n$. The $H_{n,n-i}$ for $i=1,...,n$ generate the Lie algebra
of the torus $A$. Hence it suffices that all $H_{n,n-i}$ act trivially. This follows
from the Leray type spectral sequence $DS_{n-i,0}\circ DS_{n,n-i} \Longrightarrow
DS_{n,0}$. As in the proof of lemma \ref{homotopy} one shows that $H_{n,n-i}$ acts trivially on $DS_{n,n-i}(V)$. Hence by the spectral sequence
$H_{n,n-i}$ acts by a nilpotent matrix on $D_{n,0}(V)$. On the other hand $A$, and hence
$H_{n,n-i} \in Lie(A)$, acts in a semisimple way. This proves the claim. \qed 

\medskip\noindent



\section{ Boundary maps}

\medskip
Suppose given a module $S$ in $\calR_n$.
Consider $D_{tot} = D+D'$ for $D=D_{n,n-i}$ and $D'=D_{n-i,0}$.
Notice that $DD' = -D'D$ and $D^2= c \rho(H)$, $(D')^2 = c \rho(H')$
and $D_{tot}^2 = c \rho(H_{tot})$. 

\medskip
For fixed $i$ we write $A = A_i$. We have $H_D(S)= Kern(D:S\to S)/(S^H \cap Im(D:S\to S))$.
We have also shown that this is equal to 
$$ H_D(S)= Kern(D:S^A\to S^A)/Im(D:S^A\to S^A)$$
for the torus $A$ whose Lie algebra is generated by all $H_{n,n-j}$
for $j=1,...,n$. In a similar way
$$ H_{D_{tot}}(S)= Kern(D_{tot}:S^A\to S^A)/Im(D_{tot}:S^A\to S^A)\ .$$
Recall that $A$ commutes with $D_{tot},D,D'$ and acts 
in a semisimple way.

\medskip
Let $U\subseteq S$ denote the image of $D':S^A\to S^A$.
Then $U$ and $\overline S^A= S^A/U$ are stable under $D$ and $D'$.
If $s\in S^A$ is in $Kern(D_{tot})$, then $Ds = - D's \in U$. 
Hence $s\mapsto s + U$ defines a map
from $Kern(D_{tot}:S^A\to S^A)/Im(D_{tot}:S^A\to S^A)$ to $Kern(D:\overline S^A\to \overline S^A)/Im(D_{tot}:\overline S^A\to \overline S^A)$, hence a map 
or
$$  \sigma_S: H_{D_{tot}}(S) \longrightarrow H_D(\overline S) \ .$$

\bigskip\noindent
Suppose given modules $S, V, L$ in $\calR_n$ defining an extension
$$ 0 \to S \to V \to L \to 0 \ .$$
We get a boundary map
$$ \delta_{tot}: H^\pm_{D_{tot}}(L) \longrightarrow   
H^\mp_{D_{tot}}(S) $$
defined as usually by
$ Kern(D_{tot}:L\to L) \ni \overline v  \ \mapsto \  [s] \ , \ s = D_{tot} v \in L $.
Here $v\in V^A$ is any lift of $\overline v \in L^A$ (it exists by the semisimple action of $A$ on $V$). Obviously $D_{tot}(s)=0$, since $D_{tot}^2=0$ on the space of $A$-invariant vectors. Therefore the class $[s]$ of $s$ in  $H^\mp_{D_{tot}}(S)$ is well defined. 

\medskip
In a completely similar way one defines the boundary map
$$ \delta: H^\pm_{D}(L) \longrightarrow   
H^\mp_{D}(S) \ .$$
We claim that there exists a commutative diagram
$$ \xymatrix{  H_D^\pm(L) \ar[r]^\delta & H_D^\mp(S) \cr
H^\pm_{D_{tot}}(L)\ar[u]_{\sigma_L}  \ar[r]^{\delta_{tot}} & H_{D_{tot}}^\mp(S) \ar[u]_{\sigma_S} }  $$
In fact on the level of representatives $v\in V^A$ it amounts to the assertion
$$ \xymatrix{  \overline v = v \ mod \ U  \ar[r] & \overline s = D\overline v  \cr
v \ar[u]_{\sigma_L}  \ar[r] &  s= D_{tot}v \ar[u]_{\sigma_S} }  $$
using $D_{tot}v \equiv Dv\ mod \ U_L$.

\medskip
We now consider two extension $(S,V,L)$ and $(S,\tilde V, \tilde L)$.
Then the commutative diagram 
$$ \xymatrix{  H^+_D(\tilde L) \ar[r]^{\tilde \delta} & H^-_D(S) & H^+_D(L) \ar[l]_{\delta} \cr  H^+_{D_{tot}}(\tilde L) \ar[u]_{\sigma_{\tilde L}}\ar[r]^{\tilde \delta_{tot}} & H^-_{D_{tot}}(S) \ar[u]_{\sigma_{S}} & H^+_{D_{tot}}(L) \ar[u]_{\sigma_{L}}\ar[l]_{\delta_{tot}} } $$
implies

\begin{lem} Suppose $L\cong{\bf 1}$. Then $Im(\delta_{tot})$ is not contained in $Im(\tilde\delta_{tot})$,
if there exists an integer $i$ for $1\leq i\leq n$ such that
\begin{itemize}
\item $H^+_D(\tilde L)$ does not contain $\bf 1$ as a $G_{n-i}$-module.
\item $\delta({\bf 1}) \neq 0 $ in $H_D^-(S)$.
\end{itemize}
for $D=D_{n,n-i}$.
\end{lem}

\medskip
{\it Remark}. If $S=L_n(j)$ for some $1\leq j\leq n$ is one of the hook representation discussed in section \ref{hooks},
then the smallest integer $i$ for which 
$\delta({\bf 1}) \neq 0 $ in $H_D^-(S)$ holds is given by $n-i=j-1$. Indeed we later
show that although  $0 \to D_{n,j}^0(S) \to D_{n,j}^0(V) \to {\bf 1}\to 0$ is exact,
the map $D_{n,j-1}^0(V) \to {\bf 1}$ is not surjective any longer.
So in the later applications, to apply the last lemma, we have to check
whether $H^+_D(\tilde L)$ contains the trivial $G_{j-1}$-module in this case
or not.

\medskip



\section{Highest Weight Modules}\label{HighW}

\medskip{\it \! Irreducible representations}.
The irreducible $Gl(n\vert n)$-modules $L$ in ${\calR}_n$ are uniquely determined up to isomorphism by their highest weights $\lambda$. These highest weights $\lambda$ are in the set $X^+(n)$ of dominant weights, where $\lambda$ is in $X^+(n)$ if and only if $\lambda$ is of the form 
$$\lambda = (\lambda_1,\lambda_2,\ldots,\lambda_n \ ; \ \lambda_{n+1},\ldots, \lambda_{2n})$$ with integers $\lambda_1 \geq \lambda_2 \geq \ldots \geq \lambda_n$ and $\lambda_{n+1} \geq \lambda_{n+2} \geq \ldots \geq \lambda_{2n}$.

\medskip
We remark that
the condition $$  \lambda_n = -\lambda_{n+1} \ $$ for $\lambda$  is equivalent to the condition $\lambda(H)=0$.
In the language of Brundan and Stroppel in section \ref{BS} the condition $\lambda(H)=0$ is tantamount to the condition that the irreducible representation $L(\lambda)$ is not projective and the smallest $\vee$-hook is to the left of all $\times$'s and $\circ$'s. Any at least 1-atypical
block contains such $L(\lambda)$. If these equivalent conditions hold we
write $$\overline\lambda =(\lambda_1,...,\lambda_{n-1} \ ; \ \lambda_{n+2},\ldots, \lambda_{2n})$$
defining an irreducible representation $L(\overline\lambda)$ in $\calR_{n-1}$.

\medskip
Using the notation of \cite{Drouot}
the irreducible {\it maximally atypical} $Gl(n\vert n)$-modules $L$ in ${\calR}_n$ are given by highest weights $\lambda$ of the form 
$$\lambda = (\lambda_1,\lambda_2,\ldots,\lambda_n \ ; \ -\lambda_n,\ldots, -\lambda_1) 
$$ with integers $\lambda_1 \geq \lambda_2 \geq \ldots \geq \lambda_n$. We abbreviate this by writing $[\lambda_1,\ldots,\lambda_n]$ for the corresponding irreducible representation.
The full subcategory of $\calR_n$ generated by these will be denoted $\calR_n^n$.

\medskip
{\it Highest weight modules}. Recall that
a vector $v\neq 0$ in a module $(V,\rho)$ in $\calR_n$ is called {\it primitive}, if $\rho(X)v=0$
holds for all $X$ in the standard Borel subalgebra $\mathfrak{b}$ of $\g=\g_n$.   
A highest weight vector of a module $V$ (of weight $\lambda$) in $\calR_n$ is a vector $v\in V$ that is a primitive
eigenvector of $\mathfrak{b}$ (of the weight $\lambda$) 
generating the module $V$. In this case $V$ is called a highest weight module (of weight $\lambda$).  Every irreducible representation $L(\lambda)$ in $\calR_n$ is a highest weight module of weight $\lambda$. Every highest weight module $V$ of weight $\lambda$ has 
cosocle isomorphic to $L(\lambda)$. 

\begin{lem}\label{GEN}
For $(V,\rho)=L(\lambda)$, or more generally a cyclic representation generated by a highest weight vector of weight $\lambda$, the weight space in $V$ of weight $\lambda-\mu$
is generated by $\rho(\overline x)v$, where $v$ is a highest weight vector of $(V,\rho)$.  
\end{lem} 

\bigskip\noindent
{\it Proof}. For the simple positive roots $\Delta=\{\alpha_1,...,\alpha_r\}$, i.e. the union of the odd simple root $\{\mu\}$ and the even simple roots in $\g_{\overline 0}$ with respect to the standard Borel subalgebra of upper triangular matrices, choose generators $X_\alpha\in \mathfrak{u}$. Put $\tau(X_\alpha)= Y_{-\alpha}$
and $V_0=F\cdot v$. Recursively define $V_i = V_{i-1}+\sum_{\alpha\in \Delta}
\rho(Y_{-\alpha})(V_{i-1}).$ We claim that $V_\infty=\bigcup_{i=0}^\infty V_i$ is
a $\g$-submodule of $V$, hence equal to $V$. This claim also implies that
the weight space $V_{\lambda-\mu}$ is generated by $\rho(\tau(x))v$.

\medskip
$V_\infty$ is invariant under all $\rho(Y_{-\alpha}), \alpha\in \Delta=\{\alpha_1,...,\alpha_r\}$. Each $V_i$ obviously is invariant under $\rho(X)$ for diagonal $X\in \mathfrak{g}$. Indeed each $V_i/V_{i-1}$ decomposes in weight spaces for weights
$$  \lambda - \sum_{j=1}^r n_j \alpha_j \quad , \quad \sum_{j=1}^r  
n_j = i  \quad  \quad (n_j \in \mathbb N_{\geq 0}) \ .$$ Note $\rho(X_\alpha)\rho(Y_{-\beta}) \pm \rho(Y_{-\beta})\rho(X_\alpha) = \rho(H_\alpha)$ for $\alpha=\beta\in \Delta$ and $\rho(X_\alpha)\rho(Y_{-\beta}) \pm \rho(Y_{-\beta})\rho(X_\alpha) = \rho([X_\alpha, Y_{-\beta}]) =0$ for $\alpha,\beta\in \Delta$ and $\alpha\neq \beta$
[since $\alpha - \beta \notin \Phi^+ \cup \Phi^-$ for $\alpha,\beta\in\Delta$].  
Hence $V_\infty$ is invariant under $g$, since $Y_{-\beta}, \beta\in\Delta$ and diagonal $ X$
and $X_\alpha, \alpha\in \Delta$, generate $\mathfrak{g}$ as a Lie superalgebra. \qed

\begin{lem} \label{stable}\label{0}
 Suppose $\lambda=(\lambda_1,...,\lambda_{n-1},\lambda_n\ ; \ \lambda_{n+1}, \lambda_{n+2},\ldots, \lambda_{2n})$ satisfies $\lambda_n = -\lambda_{n+1} $. If $V$ is a highest weight representation generated by a highest weight vector $v$ of weight $\lambda$, the module $H^{\lambda_n}(V)$ contains a highest weight submodule of weight $\overline\lambda$ generated by the image of $v$ with parity $(-1)^{\lambda_n}$. In particular the  representation $\Pi^{\lambda_n} L(\overline\lambda)$ in $\calR_{n-1}$ is a 
Jordan-H\"older constituent of $H^{\lambda_n}(V)$. 
\end{lem}

\medskip
{\it Proof of the lemma}. 
The highest weight vector $v$ of $V$ is a highest weight vector of the restriction
of $V$ to the subgroup $G_{n-1}$ of $G_n$ and is annihilated by $\rho(x)$.
By our assumption on the weight $\lambda$ furthermore
$   v \in V^H $.
To prove our claim it suffices to show that $v$ is not contained in $Im(\rho(x))$.
Suppose $v = \rho(x)(w)$. Since the weight of $x$ is $\mu$,  we can assume that the weight of $w$ is $\lambda - \mu$. Since $V$ is a highest weight representation, by lemma \ref{GEN} then  $w$ is proportional to $\rho(\overline x)v$. So that to show $\rho(x)w=0$ and to finish our proof, it suffices that by $ [x, \overline x] =  H$
$$  \rho(x) \rho(\overline x) v  = - \rho(\overline x) \rho(x) v + \rho(H) v = 0$$
vanishes, since $\rho(x)v=0$ and $v\in V^H$. \qed

\medskip
Note that (for the notation used see section \ref{sec:Casimir}) $$ z_i= [x_i,\overline x] = x_i \overline x + \overline x x_i \quad , \quad z'_i= [x'_i,\overline x] = x_i' \overline x + \overline x x_i' $$
are in the unipotent Lie algebra $\mathfrak{u}_{\overline 0} \subset \mathfrak{u}$ of the standard Borel $\mathfrak{b}_{\overline 0}$ of $\g_{\overline 0}$ for all $i=1,..,n-1$.
Suppose $(V,\rho)$
is a representation of $G_n$. If $(V,\rho)$ has a highest weight vector $v$, then $\rho(X)v=0$
holds for all $X$ in the unipotent radical $\mathfrak{u}$ of the standard Borel of $\g$.
In particular $$\rho(x)v=0, \rho(x_i)v=0, \rho(x'_i)v=0, \rho(z_i)v=0, \rho(z'_i)v=0 $$ 
 and hence by the commutation relations above this implies for $i=1,...,n-1$ also $$\rho(x_i)\rho(\overline x)v =0
 \quad , \quad \rho(x'_i)\rho(\overline x)v =0 \ .$$
Now also suppose $v\in V^H$ and put $w=\rho(\overline x)v$. 
Then $\rho(x)w=0$, as shown in the proof of lemma \ref{stable}. 
Similarly one can show $\rho(x_i)w=0$ (since $\rho(x_i)v=\rho(z_i)v=0$)
and $\rho(x'_i)w=0$. All elements $\mathfrak{u} \cap \g_{n-1}$ commute with $\rho(\overline x)$ and annihilate $v$, hence annihilate $w$.
Finally, since $\rho(\overline x)$ and $\rho(x_i),\rho(x'_i)$ annihilate $w$, also $\rho(z_i)$ and $\rho(z'_i)$ annihilate $w$. It follows that $\rho(X)w=0$ for all $X\in \mathfrak{u}$,
since $\mathfrak{u}$ is spanned by $\mathfrak{u}\cap \g_{n-1}$ and the $x,x_i,x'_i,z_i,z'_i$,. This implies
that $w$ is a highest weight vector in $(V,\rho)$ of weight $\lambda -\mu$, 
if $w\neq 0$. Hence

\begin{cor}\label{companion}
If $(V,\rho)$ is a highest weight representation with  highest weight vector $v$ and highest weight $\lambda$ so that $\lambda(H)=0$,
then $w=\rho(\overline x)v$ defines a highest weight vector  of weight $\lambda-\mu$ in $V$ if $w\neq 0$. 
\end{cor}

\medskip
In the situation of the last corollary, the following
conditions are equivalent
\begin{enumerate}
\item $w=0$
\item $D(v)=0$
\item $D(v)=0$ and $v$ defines a nonvanishing cohomology class in $H_D(V)$.
\end{enumerate}
Indeed $D(v)= i\rho(\overline x) v + \rho(x)v= i w$.
Furthermore, if $v=D(\tilde w)$, then
$v = i\rho(\overline x) \tilde w_1 + \rho(x) \tilde w_2$
for $w_1\in V_{\lambda+\mu}$ and $w_2\in V_{\lambda-\mu}$.
Since $\lambda$ is highest weight, therefore $V_{\lambda+\mu}=0$.
Furthermore $V_{\lambda-\mu}$ is generated by $w$, and $\rho(x)w=0$.
Hence $v \notin D(V)$. 

\medskip
A highest weight representation $V$ of weight $\lambda$ 
canonically admits the irreducible representation $L=L(\lambda)$ as a quotient. Let $q : V \to L$
denote the quotient map.

\begin{cor}\label{companion2}
In the highest weight situation of corollary \ref{companion} 
the following holds for the representation $V$:
\begin{enumerate}
\item If $V$ contains
a highest weight subrepresentation $W\neq 0$ of weight $\lambda-\mu$, then 
$H_D^{\lambda_n}(V)$ has trivial weight space $H_D^{\lambda_n}(V)_{\overline \lambda} \subseteq H_D^{\lambda_n}(V)$.
\item If the natural map $H_D(q): H_D(V) \to H_D(L)$ is surjective, then $V$ does
not contain a highest weight subrepresentation $W\neq 0$ of weight $\lambda-\mu$.
\end{enumerate}
\end{cor} 

\medskip
{\it Proof}. For the first assertion, notice that $D(v)=0$ implies $w=0$
and $w$ generates $V_{\lambda-\mu}$. For the second assertion  
notice that the highest weight vector $v\in V$ maps to the highest weight
vector $q(v)$ of $L$. By the first assertion and lemma \ref{stable},
applied for $L$, the vector $q(v)$ is $D$-closed
and defines a nonzero class in $H_D^{\lambda_n}(L)_{\overline\lambda}$. Since now 
$H_D(q)$ is surjective by assumption, corollary \ref{H2} implies that 
this class is the image of a nonzero cohomology class $\eta$ in $H_D^{\lambda_n}(V)$.
This class is representated by a nonzero
$\overline\partial$ closed class in $H^{\lambda_n}(V)=DS_{\lambda_n}$
in the weight space $\overline\lambda$. Hence this class has a $D$-closed representative
$v'$ in $V_{\lambda}$, since the enriched weight structure on $DS(V)$ allows to recover the
weight structure of $V$. Since $V$ is a highest weight representation, the space $V_\lambda$
has dimension one and therefore $v'$ is proportional to $v$. Thus $D(v)=0$. But, as explained above, this implies $w=0$ and hence $V_{\lambda-\mu}=0$. \qed

\medskip
Since Kac modules $V(\lambda)$ are highest weight modules of weight $\lambda$
with $H_D(V(\lambda))=0$,
lemma \ref{stable} and its corollaries above imply

\begin{lem}\label{KAC} 
For $\lambda$ in $X^+$ with $\lambda_n = \lambda_{n+1}=0$ 
the cohomology $H^{0}(V(\lambda))$ of the Kac module $V(\lambda)$ contains a highest
weight module of weight $\overline\lambda$. Furthermore $V(\lambda)$ contains
a nontrivial highest weight representation of weight $\lambda -\mu$.
\end{lem}

\smallskip{\bf Example}. Let $(V,\rho)=V({\bf 1})$ in $\calR_2$ be the Kac module of the trivial representation. Then $DS(V^*)=0$ and $DS(V)\neq 0$, since $V$ is not projective.
The module $V$ is a cyclic module generated by it highest weight vector of weight $\lambda=0$ (this is not true for the anti-Kac module $V^*$).
Furthermore $V$ has Loewy length 3 with Loewy series $(Ber_2^{-2},Ber_2^{-1}S^1,{\bf 1})$ where $S^1 = [1,0]$.
We claim $$DS(V) = (Ber_1^{-2}\oplus {\bf 1}) \otimes \bigl({\bf 1} \oplus {\Pi(\bf 1)}\bigr)\ .$$
This follows from the later results, e.g. lemma \ref{hex} and theorem \ref{mainthm}:  $d(Ber_2) = - Ber_1$ and $d(S^1)= Ber_1^{-1} + Ber_1$ imply $d(V)=0$, hence $DS(V)$ has at most 4 Jordan-H\"older
constituents $Ber_1^{-2}, \Pi(Ber_1^{-2}), {\bf 1}, \Pi({\bf 1})$. By lemma \ref{stable} the constituent
${\bf 1}$ occurs. By duality then also the constituent $Ber_1^{-2}$ must occur.
Since $d(V)=0$ the constituent $\Pi (Ber_1^{-2}), {\Pi (\bf 1)}$ must occur. Finally apply proposition \ref{ext-0}. This example shows that $DS$ in general does not preserve negligible objects.

\medskip
{\it Highest weights}. 
Suppose $(V,\rho)$ is a highest weight module of weight $\lambda$ such that $\lambda(H)=0$. Let $\nu$ be a weight of $V$.
Then 
$$   \nu = \lambda - \sum_{\alpha\in\Delta_n} \N_{\geq 0} \cdot \alpha $$
for the set $\Delta_n$ of simple positive roots $\alpha$ of $G_n$.

\medskip
Now suppose $\nu$ contributes to $DS(V)$. Then $\nu$ is a weight of $V^H$
and hence $\nu(H)=0$. 
Notice $\Delta_n$ is the union of $\Delta^+_n=\{e_1-e_2,...,e_{n-1}-e_n, e_{n+1}-e_{n+2},...,
e_{2n-1}-e_{2n}\}$ and $\Delta^-_n =\{ e_n - e_{n+1} \}$.
The restriction of the simple roots $\alpha \in \Delta_n$ are in $\Delta_{n-1}$ (i.e. simple root of $G_{n-1}$) except for the even simple roots $\alpha = e_{n-1}-e_n$, $\alpha= e_{n+1}-e_{n+2}$
and the odd simple root $\alpha = e_n-e_{n-2}$. A linear combination $\sum_{\alpha\in \Delta_n}
n_{\alpha}\alpha $ annihilates $H$ if and only if the coefficient, say $m$, of $e_{n-1}-e_n$ and $e_{n+1}-e_{n+2}$
coincides; hence this holds iff $\nu$ is of the form $\sum_{\alpha\in \Delta_{n-1}^+} n_\alpha \alpha +
(n_\mu - m) \cdot (e_n - e_{n+1}) +  m\cdot (e_{n-1} - e_{n+2})$. Notice that $\mu=(e_n - e_{n+1})$ is trivial
on the maximal torus of $G_{n-1}$ and that $(e_{n-1} - e_{n+2})$ defines the new odd simple
root in $\Delta_{n-1}^-$. Hence the restriction of $\nu \in V^H$ is of the form 
$$   \nu\vert_{\mathfrak{b}\cap \g_{n-1}} \ \in \ \lambda\vert_{\mathfrak{b}\cap \g_{n-1}} - \sum_{\alpha\in\Delta_{n-1}}
\N_{\geq 0} \cdot \alpha \ $$ under our assumptions above. 
Notice for $V_\lambda \subset V^H \subset V$ we have 
$$  \ell = \lambda(diag(1,..,1;t^{-1},1,..,1)) = n_\mu - m = \lambda'_n \ .$$
The discussion above implies

\begin{lem} \label{Hi}
For a highest weight module $(V,\rho)$ in $T_n$ of weight $\lambda$ with $\lambda(H)=0$
the module $DS(V,\rho)$ has its weights $\nu$ in $\lambda - \sum_{\alpha\in\Delta_{n-1}} \N_{\geq 0} \cdot \alpha$.   
\end{lem} 

\begin{cor} Given $(V,\rho)\in T_n$,
suppose $L(\lambda)$ is a Jordan-H\"older constituent of $(V,\rho)$ 
such that for all Jordan-H\"older constituents $L(\nu)$ of $(V,\rho)$
we have $\nu \in \lambda - \sum_{\alpha\in\Delta_n} \N_{\geq 0} \cdot \alpha$
and $\nu(H)=0$. Then $L(\overline\lambda)$ appears in $DS(V,\rho)$ and  
all other irreducible constituents $L(\nu')$ or $\Pi L(\nu')$ of $DS(V,\rho)$ satisfy
$\nu' \in \lambda - \sum_{\alpha\in \Delta_{n-1}} \N_{\geq 0} \cdot \alpha$.
\end{cor}

\medskip
{\it Proof}. This follows from the last lemma and the 
weak exactness of the functor $DS$. \qed



\section{The Casimir}\label{sec:Casimir}
\medskip
We study the operation of the Casimir $C_n$ on $DS(V)$. This will be used in section \ref{sec:loewy-length} when we study the effect of $DS$ on translation functors $F_i(L_{\times \circ})$.

\medskip 
Consider the fixed element $x\in \g_n$  $$x = \begin{pmatrix} 0 & y \\ 0 & 0 \end{pmatrix} \in \g_{n} \ \text{ for } \ y = \begin{pmatrix} 0 & 0 & \ldots & 0 \\ 0 & 0 & \ldots & 0 \\ \ldots & & \ldots &  \\ 1 & 0  & 0 & 0 \\ \end{pmatrix} $$
Similarly we define $$x_i\ ,\ x_i' \ \text{ for } i=1,..,n-1$$ for matrices $y=y_i$ resp. $y_i'$ with
a unique entry 1 in the first column resp. last row at positions different from the entry 1 in the above $y$
$$  \begin{pmatrix} * & 0 & \ldots & 0 \\ * & 0 & \ldots & 0 \\ \ldots & & \ldots &  \\ 0 & *  & * & * \\ \end{pmatrix} $$
Then $x, x_i, x'_i$ are in $u$ for $i=1,...,n-1$.
The elements $x_i,x'_i$ satisfy $[x_i,x]=0=[x'_i,x]$.

\medskip
Using Brundan-Stroppel's notations \cite{Brundan-Stroppel-4}, (2.14), let $e_{r,s}\in \g_n$ be the $rs$-matrix unit.
Then the  Casimir operator $C_n= \sum_{r,s=1}^n (-1)^{\overline s} e_{r,s} 
e_{s,r}$ of the super Lie algebra $\g_n=Lie(G_n)$
is recursively given by 
$$ C_n = C_{n-1} + C_1 +   2(\overline z_1 z_1 + \cdots + \overline z_{n-1} z_{n-1}) + (e_{1,1} + \cdots + e_{n-1,n-1} - (n-1)e_{n,n})  $$
$$ -2(\overline z'_1 z'_1 + \cdots + \overline z'_{n-1} z'_{n-1}) - ( - e_{n+2,n+2} - \cdots - e_{2n,2n} + (n-1) e_{n+1,n+1})  $$
$$ + 2(\overline x_1 x_1 + \cdots + \overline x_{n-1} x_{n-1}) - (e_{1,1} + \cdots + e_{n-1,n-1} + (n-1)e_{n+1,n+1}) $$
$$ + 2(\overline x'_1 x'_1 + \cdots + \overline x'_{n-1} x'_{n-1}) - (e_{n+2,n+2} + \cdots + e_{2n,2n} + (n-1)e_{n,n})  $$
with the notations $x_i = e_{i,n+1}$, 
$x'_i = e_{n,2n+1-i}$, $z_i=e_{i,n}$ and $z'_i= e_{n+1,2n+1-i}$. Furthermore $\overline x_i, \overline x'_i,
\overline z_i$ and $\overline z'_i$ denote the supertransposed of $x_i, x'_i,
z_i$ and  $z'_i$. 
Hence
$$ C_n = C_{n-1} + C_1 + 2(\overline z_1 z_1 + \cdots + \overline z_{n-1} z_{n-1}
-\overline z'_1 z'_1 - \cdots - \overline z'_{n-1} z'_{n-1}) $$ $$ + 2(\overline x_1 x_1 + \cdots + \overline x_{n-1} x_{n-1} + \overline x'_1 x'_1 + \cdots + \overline x'_{n-1} x'_{n-1})
- 2(n-1)H \ $$
using $[\tau(x),x_i] = z_i$ and
$[\tau(x),x'_i]= z'_i$ and $$[z_i,\overline z_i] = e_{i,i} - e_{n,n}\ \quad , \quad [z'_i,\overline z'_i] = e_{n+1,n+1} - e_{2n+1-i,2n+1-i}$$
and $[\overline x_i,x_i] = e_{i,i} + e_{n+1,n+1}$ and $[\overline x'_i,x'_i] = e_{2n+1-i,2n+1-i} + e_{n,n}$. Notice $\overline x_i x_i - x_i \overline x_i = 2 \overline x_i x_i - e_{i,i} - e_{n+1,n+1}$
and $\overline x'_i x'_i - x'_i \overline x'_i = 2 \overline x'_i x'_i - e_{2n+1-i,2n+1-i} - e_{n,n}$.
Finally $C_1 = e_{n,n}^2 - e_{n+1,n+1}^2 - x \overline x + \overline x x
= e_{n,n}^2 - e_{n+1,n+1}^2 + 2 \overline x x - H$.

\medskip
{\it Representations}. Suppose $(V,\rho)$ is a representation of $\g_n$. On $DS(V,\rho)$
we have $\rho(H)=0$ and $\rho(x)=0$. Since 
$$[x,x_i] = [x,x'_i] = [x,z_i] = [x,z'_i]= 0 \ ,$$
the elements $x_i,x'_i,z_i,z'_i$  naturally act on the cohomology $DS(V)=V_x$. 
Since $x$ commutes with $H$, the spaces $Kern(x)$ and its subspace
$Im(x)$ decompose into $H$-eigenspaces $Kern(x)(j)$ and $Im(x)(j)$ for $j\in \mathbb Z$. 
By lemma \ref{homotopy} however $Kern(x)(j)= Im(x)(j)$, expect for the zero-eigenspace of
$H$. Now, although $x,\overline x$ commute with $H$, the operators $y\in\{x_i,x'_i,z_i,z'_i\}$ satisfy  $[H,y]=\pm y$ and hence map the zero eigenspace $M=V^H$ into the $\pm 1$-eigenspace of $H$ on $V$. Since the $j=\pm 1$-eigenspaces do not give a nonzero 
contribution to the cohomology $DS(V) = V_x$, this implies 

\begin{lem}\label{T}
The natural action of $\rho(x_i), \rho(x'_i), \rho(z_i), \rho(z'_i)$ and 
 $\rho(x), \rho(H)$
on $DS(V,\rho)$ is trivial.
\end{lem}

\medskip
Notice that $C_n$ commutes with all elements in $\g_n$,
hence induces a linear map
on $DS(V,\rho)$ that commutes with the action of $G_{n-1}$ on $DS(V,\rho)$.

\begin{lem} \label{Cas}
The restriction of the Casimir $C_n$ acts on $DS(V,\rho)$ like the Casimir $C_{n-1}$ of $T_{n-1}$ acts on $DS(V,\rho)\in T_{n-1}$.  
\end{lem}

\medskip{\it Proof}. By lemma \ref{T} the restriction of $C_n$ to 
$DS(V,\rho)$ is the sum of $C_{n-1}$ and the 
operator $C_1 = e_{n,n}^2 - e_{n+1,n+1}^2 $. 
Now consider a weight space 
of $DS(V,\rho)$ with eigenvalue $\overline\lambda$. Then $\overline\lambda$
is the restriction of an eigenvalues $\lambda$ of the weight decomposition of $(V,\rho)$.
Since $DS(V,\rho)$ is represented by elements in $M=V^H$, the condition $\lambda(H)=0$ implies $\lambda_n = - \lambda_{n+1}$ and hence $\lambda_n^2 - \lambda_{n+1}^2=0$. Therefore $C_1$ acts trivially on $DS(V,\rho)$. \qed

\medskip\noindent
{\it Remark}. As the referee pointed out, there is a more conceptual proof of lemma \ref{Cas}. For a module $M \in T_n$, the Casimir map is the composition \[ \xymatrix{ C_M: M \ar[r] & \mathfrak{gl}(n|n) \otimes \mathfrak{gl}(n|n)^* \otimes M \ar[r] & \mathfrak{gl}(n|n) \otimes M \ar[r] & M } \] where the first map is the coevaluation map for the adjoint representation of $\g$ and the last two are the action maps $\mathfrak{gl}(n|n) \otimes M \to M$. Since $DS$ maps the standard representation to the standard representation, it preserves the adjoint representation as well, and hence preserves the Casimir map in the sense that $DS(C_M)$ is the Casimir on $DS(M)$.

\bigskip



\part{The Main Theorem and its proof}


\medskip

In this part we prove the main theorem, stating that $DS(L) = \bigoplus_i \Pi^{n_i} L(\lambda_i)$ in $T_{n-1}$ for any irreducible representation $L$. We have seen that $DS(L)$ is actually a $\Z$-graded object in $T_{n-1}$; and we calculate the $\Z$-grading for any $L$ in the propositions \ref{hproof} \ref{hproof-2}. These statements contain the main theorem as a special case. Their proofs however depend on the main theorem and its proof. We will prove the main theorem first for special irreducible $L$, called ground states, and reduce the general question to these by means of translation functors.


\section{The language of Brundan and Stroppel}\label{BS}

\medskip
By the work of Brundan-Stroppel \cite{Brundan-Stroppel-4} the block combinatoric and notably the $Ext^1$ between irreducible representations can be described in terms of weight and cup diagrams associated to any irreducible $L()\lambda)$.

\medskip

{\it Weight diagrams}. Consider a weight $\lambda=(\lambda_1,...,\lambda_n ; \lambda_{n+1}, \cdots, \lambda_{2n})$. Then $\lambda_1 \geq ... \geq \lambda_n$ and $\lambda_{n+1} \geq ... \geq \lambda_{2n}$ are integers, and every $\lambda\in {\mathbb Z}^{m+n}$ satisfying these inequalities occurs as the highest weight of an irreducible representation $L(\lambda)$. The set of highest weights will be denoted by $X^+=X^+(n)$. Following \cite{Brundan-Stroppel-4} to each highest weight $\lambda\in X^+(n)$  we associate  two subsets of cardinality $n$ of the numberline $\mathbb Z$
\begin{align*} I_\times(\lambda)\ & =\ \{ \lambda_1  , \lambda_2 - 1, ... , \lambda_n - n +1 \} \\
 I_\circ(\lambda)\ & = \ \{ 1 - n - \lambda_{n+1}  , 2 - n - \lambda_{2n-1} , ... ,  - \lambda_{2n}  \}. \end{align*}

We now define a labeling of the numberline $\mathbb Z$.
The integers in $ I_\times(\lambda) \cap I_\circ(\lambda) $ are labeled by $\vee$, the remaining ones in $I_\times(\lambda)$ resp. $I_\circ(\lambda)$ are labeled by $\times$ respectively $\circ$. All other integers are labeled by $\wedge$. 
This labeling of the numberline uniquely characterizes the weight vector $\lambda$. If the label $\vee$ occurs $r$ times in the labeling, then $r=atyp(\lambda)$ is called the {\it degree of atypicality} of $\lambda$. Notice $0 \leq r \leq n$, and for $r=n$ the weight $\lambda$ is called
{\it maximal atypical}. 

\medskip
{\it Blocks}. A block $\Gamma$ of $X^+(n)$ is a connected component of the Ext-quiver of ${\calR}_n$. Let ${\calR}_{\Gamma}$ (or by abuse of notation $\Gamma$) be the full subcategory of objects of ${\calR}_n$ such that all composition factors are in $\Gamma$. This gives a decomposition ${\calR}_n = \bigoplus_{\Gamma} {\calR}_{\Gamma}$ of the abelian category. Two irreducible representations $L(\lambda)$ and $L(\mu)$ are in the same block if and only if the weights $\lambda$ and $\mu$ define labelings with the same position of the labels $\times$ and $\circ$. The degree of atypicality is a block invariant, and the blocks $\Lambda$ of atypicality $r$ are in 1-1 correspondence with pairs of disjoint subsets of $\bbZ$ of cardinality $n-r$ resp. $n-r$. Let $\calR_n^i$ be the full subcategory of $\calR_n$ defined by the blocks of atypicity $n-i$. In particular $\calR_n$ has a unique maximally atypical block, and any block of atypicality $i$ in $\calR_n$ is equivalent to the maximally atypical  block in $\calR_i$.

\medskip
{\it Cups}. To each weight diagram we associate a cup diagram as in \cite{Brundan-Stroppel-1}. Here a cup is a lower semi-circle joining two points in $\mathbb Z$. To construct the cup diagram go from left to right throught the weight diagram until one finds a pair of vertices $\vee \ \ \ \wedge$ such that there only $x$'s, $\circ$'s or vertices which are already joined by cups between them. Then join $\vee \ \ \wedge$ by a cup. This procedure will result in a weight diagram with $r$ cups.  
For example

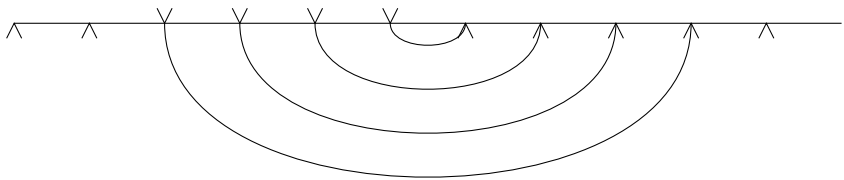
\begin{center}
\medskip

\begin{tikzpicture}
 \draw (-5,0) -- (6,0);
\foreach \x in {0,-1,-2,-3} 
     \draw (\x-.1, .2) -- (\x,0) -- (\x +.1, .2);
\foreach \x in {-5,-4,1,2,3,4,5} 
     \draw (\x-.1, -.2) -- (\x,0) -- (\x +.1, -.2);
\foreach \x in {} 
     \draw (\x-.1, .1) -- (\x +.1, -.1) (\x-.1, -.1) -- (\x +.1, .1);

\draw [-,black,out=270, in=270](0,0) to (1,0);
\draw [-,black,out=270, in=270](-1,0) to (2,0);
\draw [-,black,out=270, in=270](-2,0) to (3,0);
\draw [-,black,out=270, in=270](-3,0) to (4,0);

\end{tikzpicture}

\end{center}
is the labelled cup diagram ($n=4$) of the trivial representation
attached to the weight $\lambda=(0,\ldots,0 \vert 0, \ldots,0)$.

{\it Sectors and segments}. For the purpose of this paragraph we assume $\lambda\in X^+$ to be in a maximal atypical block, so that the weight diagram does not have labels $\times$ or $\circ$. Some of the $r$ cups of a cup diagram may be nested. If we remove all inner parts of the nested cups we obtain a cup diagram defined by the (remaining)  outer cups. We enumerate these  cups from left to right. The starting point of the $j$-th lower cup is denoted $a_j$ and its endpoint is denoted $b_j$. 
Then there is a label $\vee$ at the position $a_j$ and a label $\wedge$ at position $b_j$.
The interval $[a_j,b_j]$ of $\bbZ$ will be called the $j$-th sector of the cup diagram.
Adjacent sectors, i.e with $b_j=a_{j+1} -1$ will be grouped together into segments. The segments 
again define intervals in the numberline. Let $s_j$ be the starting point of the $j$-th segment
and $t_j$ the endpoint of the $j$-th segment. Between any two segments there is a distance at least $\geq 1$.  In the following case the weight diagram has 2 segments and 3 sectors  

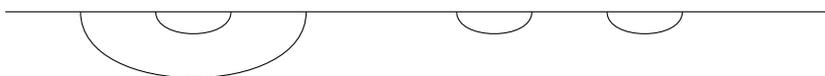
\begin{center}
\bigskip

\begin{tikzpicture}
 \draw (-5,0) -- (6,0);
\foreach \x in {} 
     \draw (\x-.1, .2) -- (\x,0) -- (\x +.1, .2);
\foreach \x in {} 
     \draw (\x-.1, -.2) -- (\x,0) -- (\x +.1, -.2);
\foreach \x in {} 
     \draw (\x-.1, .1) -- (\x +.1, -.1) (\x-.1, -.1) -- (\x +.1, .1);

\draw [-,black,out=270, in=270](-4,0) to (-1,0);
\draw [-,black,out=270, in=270](-3,0) to (-2,0);
\draw [-,black,out=270, in=270](1,0) to (2,0);
\draw [-,black,out=270, in=270](3,0) to (4,0);
%

\end{tikzpicture}
\end{center}

whereas the following weight diagram has 1 segment and 1 sector. 

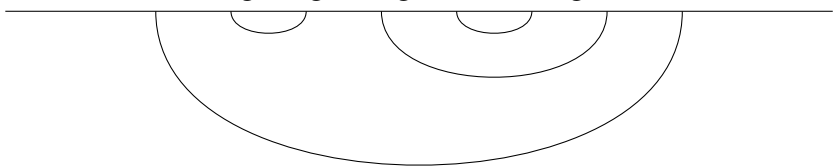
\begin{center}
\begin{tikzpicture}
 \draw (-5,0) -- (6,0);
\foreach \x in {} 
     \draw (\x-.1, .2) -- (\x,0) -- (\x +.1, .2);
\foreach \x in {} 
     \draw (\x-.1, -.2) -- (\x,0) -- (\x +.1, -.2);
\foreach \x in {} 
     \draw (\x-.1, .1) -- (\x +.1, -.1) (\x-.1, -.1) -- (\x +.1, .1);

\draw [-,black,out=270, in=270](-3,0) to (4,0);
\draw [-,black,out=270, in=270](-2,0) to (-1,0);
\draw [-,black,out=270, in=270](0,0) to (3,0);
\draw [-,black,out=270, in=270](1,0) to (2,0);
%

\end{tikzpicture}
\end{center}

Removing the outer circle would result in a cup diagram with two sectors and one segment. We can also define the notion of a sector or segment for blocks which are not maximally atypical. In this case we say that two sectors are adjacent (and belong to the same segment) if they are only separated by $\times$ or $\circ$'s. For our purpose the $\times$ and $\circ$'s will not play a role and we will often implicitly assume that we are in the maximally atypical block.

\medskip
{\it Important invariants}.
Note that the segment and sector structure of a weight diagram is completely encoded by the positions of the $\vee$'s. Hence any finite subset of $\bbZ$ defines a unique weight diagram in a given block. This will lead to the notion of a {\it plot} in the next section
where we associate to a maximal atypical highest weight the following  invariants: 

\begin{itemize}
\item the type (SD) resp. (NSD),
\item the number $k=k(\lambda)$ of sectors of $\lambda$, 
\item the sectors $S_\nu=(I_\nu,K_\nu)$ from left to right (for $\nu=1,...,k$),
\item the ranks $r_\nu = r(S_\nu)$, so that $\# I_\nu = 2r_\nu$, 
\item the distances $d_\nu$ between the sectors (for $\nu=1,...,k-1$), 
\item the total
shift factor $d_0=\lambda_n$ 
\item and the added distances $\delta_i = \sum_{\nu=0}^{i-1} d_{\nu}$. \end{itemize}
If convenient, $k$ sometimes may also denote the number of segments, 
but hopefully no confusion will arise from this. 

\medskip
A maximally atypical weight is called basic if $\lambda_\nu = -\lambda_{n+\nu}$ holds for $\nu=1,...,n$ such that $[\lambda]:=(\lambda_1,...,\lambda_n)$
defines a decreasing sequence $\lambda_1 \geq \cdots \geq \lambda_{n-1} \geq \lambda_n=0$
with the property $n-\nu \geq \lambda_\nu$ for all $\nu=1,...,n$. The total number
of such {\it basic weights} in $X^+(n)$ is the Catalan number $C_n$. Reflecting the graph of such a sequence
$[\lambda]$ at the diagonal, one obtains another basic weight $[\lambda]^*$. We will show that a basic weight $\lambda$ is of  type (SD) if and only if $[\lambda]^* = [\lambda]$ holds. 
To every maximal atypical highest 
weight $\lambda$ is attached a unique maximal atypical highest 
weight $\lambda_{basic}$ 
$$  \lambda \mapsto \lambda_{basic} \ $$
having the same invariants as $\lambda$, except that
$d_0 = d_1=\cdots = d_{k-1}=0$ holds for $\lambda_{basic}$. For example, the basic weight attached to the irreducible representation $[5,4,-1]$ in $\calR_3$ with cup diagram

\begin{center}
\medskip

\begin{tikzpicture}
 \draw (-5,0) -- (6,0);
\foreach \x in {} 
     \draw (\x-.1, .2) -- (\x,0) -- (\x +.1, .2);
\foreach \x in {} 
     \draw (\x-.1, -.2) -- (\x,0) -- (\x +.1, -.2);
\foreach \x in {} 
     \draw (\x-.1, .1) -- (\x +.1, -.1) (\x-.1, -.1) -- (\x +.1, .1);

\draw [-,black,out=270, in=270](-4,0) to (-3,0);
\draw [-,black,out=270, in=270](2,0) to (3,0);
\draw [-,black,out=270, in=270](4,0) to (5,0);
%

\end{tikzpicture}

\end{center}

is the basic representation $[2,1,0]$ with weight diagram

\medskip
\begin{center}

\begin{tikzpicture}
 \draw (-5,0) -- (6,0);
\foreach \x in {} 
     \draw (\x-.1, .2) -- (\x,0) -- (\x +.1, .2);
\foreach \x in {} 
     \draw (\x-.1, -.2) -- (\x,0) -- (\x +.1, -.2);
\foreach \x in {} 
     \draw (\x-.1, .1) -- (\x +.1, -.1) (\x-.1, -.1) -- (\x +.1, .1);

\draw [-,black,out=270, in=270](-2,0) to (-1,0);
\draw [-,black,out=270, in=270](0,0) to (1,0);
\draw [-,black,out=270, in=270](2,0) to (3,0);
%

\end{tikzpicture}

\end{center}

\medskip



\section{On segments, sectors and plots}\label{derivat}

\medskip
If $\lambda$ is a maximally atypical weight in $\calR_n$, it is completely encoded by the $n$ $\vee$'s in its weight diagram. We change the point of view and regard it as a map (a plot) $\lambda: \Z \to \{ \boxplus, \boxminus\}$ where the $\boxminus$ correspond to the $\vee$'s. If $\lambda$ is not maximal atypical, its weight diagram has crosses and circles. These do not play any role in the combinatorial arguments, and we can still describe $\lambda$ by a plot if we just ignore and remove the crosses and circles from the weight diagram. 
 
\medskip
A plot $\lambda$ is a map
$$ \lambda: \mathbb Z \to \{\boxplus,\boxminus\}\ $$ such that  
the cardinality $r$ of the fiber $\lambda^{-1}(\boxplus)$ is finite.
Then by definition $r=r(\lambda)$ is the degree and $\lambda^{-1}(\boxplus)$ is the support
of $\lambda$. 
As usual an interval $[a,b] \subset \mathbb Z$  is the set  $\{x\in\mathbb Z\ \vert \
a \leq x \leq b\}$.
Replacing $\boxplus$ by $1$ and $\boxminus$ by $-1$ we may view $\lambda(x)$
as a real valued function extended by $\lambda(x):= \lambda([x])$ to a function 
on $\mathbb R$  for $[x] = \max_{n\in \mathbb Z} \{ n\leq x\}$.

\medskip
{\bf Segments and sectors}. 
An interval $I=[a,b]$ of even cardinality $2r$ and a subset $K$ of cardinality
of rank $r$ defines a plot $\lambda$ of rank $r$ with support $K$. We call
$(I,K)$ a {\it segment}, if $f(x) = \int_a^x \lambda(x) dx$ is 
nonnegative on $I$. Notice, then $a\in K$ but $b\notin K$.

\medskip{\it Factorization}.
For a given plot $\lambda$ put $a=\min(supp(\lambda))$ and for the first zero $x_0>a$ of the function $f(x)=\int_{a}^x \lambda(x)dx$
put $b=x_0-1$. This defines an interval $I=[a,b]$ of even length, such that $\lambda\vert_I$
(now again viewed as a function on $I\cap \mathbb Z$)
admits the values $1$ and $-1$ equally often. 
If  $supp(\lambda)\subset I$, then $\lambda$ is called a {\it prime}
plot. If $\lambda$ is not a prime plot, the plot $\lambda_1$  
with support $I\cap supp(\lambda)$ 
defines a prime plot. It is called the first sector of the plot $\lambda$.
Now replace the plot $\lambda$ by the plot, where the support $K_1$ of the first sector $I=I_1$
is removed from the support $K$ of $\lambda$. Repeating the process above, we obtain
a prime plot $\lambda_2$ with support $K_2$ defining a segment $(I_2,K_2)$. This segment is called the second sector of $\lambda$. Obviously $I_1$ is an interval in $\mathbb Z$ on the right
of $I_1$, hence in particular they are disjoint. Continuing with this process, one defines finitely
many prime plots $\lambda_1,...,\lambda_k$ attached to a given plot defining disjoint
segments $S_1=(I_1,K_1)$,..., $S_k(I_k,K_k)$. These segments $S_\nu$
are called the sectors of the plot
$\lambda$. Let
$$  d_\nu = dist(I_\nu, I_{\nu+1}) \quad , \quad \nu=1,..., k-1$$
denote the distances between these sectors $S_\nu$, i.e. $d_\nu = \min(S_{\nu+1}) - \max(S_\nu)$.

For disjoint segments $(I_1,K_1)$ and $(I_1,K_2)$  the union $(I,K)=(I_1\cup I_2,K_1\cup K_2)$
again is a segment, provided $I=I_1\cup I_2$ is an interval in $\mathbb Z$.

Grouping together adjacent sectors of $\lambda$ with distances $d_\nu=0$ 
defines the segments of $\lambda$. In other words, the union
of the intervals $I_\nu$ of the sectors $S_\nu$ of the $\lambda_\nu$ can be written a disjoint union of intervals $I$ of maximal length.
These intervals $I$ define the {\it segments of $\lambda$} as $(I, T \cap supp(\lambda))$.

\medskip
We consider formal finite linear combinations $\sum_i n_i\cdot \lambda_i$
of plots with integer coefficients. This defines an abelian group $R = \bigoplus_{r=0}^\infty R_r$
(graduation by rank $r$).
We define a commutative ring structure on $R$ so that the product of two plots $\lambda_1$ and $\lambda_2$ is zero unless the segments of $\lambda_1$ and $\lambda_2$ are disjoint,
in which case the support of the product becomes the union of the supports.
A plot $\lambda$ that can not be written in the form $\lambda_1 \cdot \lambda_2$ for
plots $\lambda_i$ of rank $r_i>0$ is called a prime plot.

\begin{lem} \label{primef}  Every plot can be written
as a product of prime plots uniquely up to permutation of the factors.
\end{lem}

\medskip
Of course this prime factorization of a given plot $\lambda$ is given by the
prime factors $\lambda_\nu$ attached to the sectors $S_\nu, \nu=1,..,k$
of $\lambda$. Hence 
for $\lambda= \prod_i \lambda_i$ with prime plots $\lambda_i$, the 
{\it sectors of $\lambda$} are the segments attached to the prime factors
$\lambda_i$.  
The interval $I=[a,...,b]$ attached to a prime plot $\lambda$ is the unique sector or the unique segment of the prime plot $\lambda$. It has cardinality $2r(\lambda)$, and the support $K$ of $\lambda$ defines a subset of the sector
$I$ of cardinality $r$. Recall $a\in K$ but $b\notin K$.

\medskip
{\bf Differentiation}.
We define a derivation on $R$ called derivative.
Indeed the derivative induces an additive map
$$  \partial: R_{n} \to R_{n-1} \ .$$
To {\it differentiate} a plot of rank $n>0$, or a segment,
 we  use the formula
$$  \partial(\prod_i \lambda_i) = \sum_i \partial \lambda_i \cdot \prod_{j\neq i} \lambda_j $$ 
in the ring $R$ to reduce the definition to the case of a prime plot $\lambda$.
For prime $\lambda$ let $(I,K)$ be its associated sector. Then $I=[a,b]$. Using $a\in K$, $b\notin K$,  for a sector $(I,K)$ of a prime plot $\lambda$  it is easy to verify by the integral criterion 
that
$$\partial(I,K) = (I,K)' = (I',K')$$ for $I'=[a+1,b-1]$ and $K'=I\cap K$
again defines the sector of a prime plot $\partial \lambda$ of rank $r(\lambda')=r(\lambda)-1$.
Then for prime plots $\lambda$ of rank $n$ with sector $(I,K)$ we define $\partial \lambda$ in $R$
by   
$$ \fbox{$ \partial \lambda = \partial(I,K) $} \quad , \quad  I=[a,b] \ .$$

\medskip
{\bf Integration}. For a segment $(I,K)$ with $I=[a,b]$ put
$$ \int(I,K) =  ([a-1,b+1],K\cup\{a-1\}) \ $$
increasing the rank by 1.
Observe, that the integral criterion implies that 
$([a-1,b+1],K\cup\{a-1\})$ always defines a prime segment.
Obviously 
$$ \partial \int (I,K) = (I,K) \ . $$
Similarly $\int \partial (I,K) =(I,K)$ for a prime segment $(I,K)$ of rank $>0$. 

\medskip
{\bf Lowering sectors}. For a sector $S=(I,K)$ with $I=[a,b]$
define 
$$  S^{low} \ =\ ([a-1,a],\{a-1\}) \cup  \partial(S) \ .$$ 
Notice that $S^{low}$ is a segment with interval $[a-1,b-1]$. 

\medskip
{\bf Melting sectors}. Suppose $\lambda_1$ and $\lambda_2$ are prime
plots. Let $(I_1,K_1)$ and $(I_2,K_2)$ be their defining sectors. 
Assume that $(I,K)=(I_1\cup I_2,K_1\cup K_2)$ defines a segment
with plot $\lambda$. Hence $I_1=[a,i]$ and $I_2=[i+1,b]$ for some $i\in \mathbb Z$ and
$i\notin K_1$ and $i+1\in K_2$. Then by the
integral criterion 
$$   (I,K)^{melt} \ = \ (I_1\cup I_2, K_1 \cup \{i \} \cup (K_2 - \{i+1\} ) $$
defines a prime plot with $I=[a,b]$. 

\medskip
{\it Example.} We can represent plots with labelled cup diagrams. A plot of rank $r$ has $r$ cups. For instance the irreducible representation $[3,3,1,1] \in \calR_4$ has the cup diagram

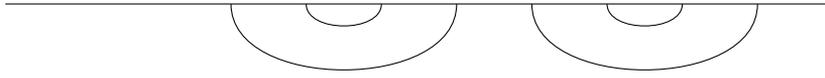
\begin{center}
\bigskip

\begin{tikzpicture}
 \draw (-5,0) -- (6,0);
\foreach \x in {} 
     \draw (\x-.1, .2) -- (\x,0) -- (\x +.1, .2);
\foreach \x in {} 
     \draw (\x-.1, -.2) -- (\x,0) -- (\x +.1, -.2);
\foreach \x in {} 
     \draw (\x-.1, .1) -- (\x +.1, -.1) (\x-.1, -.1) -- (\x +.1, .1);

\draw [-,black,out=270, in=270](-2,0) to (1,0);
\draw [-,black,out=270, in=270](-1,0) to (-0,0);
\draw [-,black,out=270, in=270](2,0) to (5,0);
\draw [-,black,out=270, in=270](3,0) to (4,0);
%

\end{tikzpicture}
\end{center} 

The corresponding plot is defined by its support $\{-2,-1,2,3\}$. Its derivative is the sum of two plots of rank 3 corresponding to the two cup diagrams

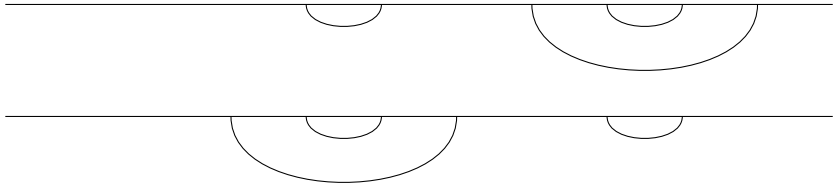
\begin{center}
\bigskip

\begin{tikzpicture}
 \draw (-5,0) -- (6,0);
\foreach \x in {} 
     \draw (\x-.1, .2) -- (\x,0) -- (\x +.1, .2);
\foreach \x in {} 
     \draw (\x-.1, -.2) -- (\x,0) -- (\x +.1, -.2);
\foreach \x in {} 
     \draw (\x-.1, .1) -- (\x +.1, -.1) (\x-.1, -.1) -- (\x +.1, .1);

\draw [-,black,out=270, in=270](-1,0) to (0,0);
\draw [-,black,out=270, in=270](2,0) to (5,0);
\draw [-,black,out=270, in=270](3,0) to (4,0);
%

\end{tikzpicture}

\end{center} 

\begin{center}
\medskip

\begin{tikzpicture}
 \draw (-5,0) -- (6,0);
\foreach \x in {} 
     \draw (\x-.1, .2) -- (\x,0) -- (\x +.1, .2);
\foreach \x in {} 
     \draw (\x-.1, -.2) -- (\x,0) -- (\x +.1, -.2);
\foreach \x in {} 
     \draw (\x-.1, .1) -- (\x +.1, -.1) (\x-.1, -.1) -- (\x +.1, .1);

\draw [-,black,out=270, in=270](-2,0) to (1,0);
\draw [-,black,out=270, in=270](-1,0) to (0,0);
\draw [-,black,out=270, in=270](3,0) to (4,0);
%

\end{tikzpicture}

\end{center} 

If we integrate the first segment of the plot we get the plot of rank 5 with support $\{-3,-2,-1,2,3\}$ with corresponding cup diagram

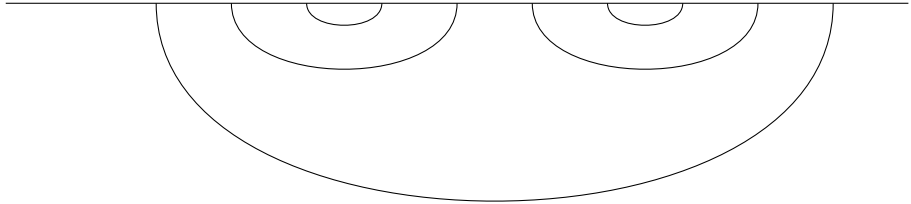
\begin{center}
\bigskip

\begin{tikzpicture}
 \draw (-5,0) -- (7,0);
\foreach \x in {} 
     \draw (\x-.1, .2) -- (\x,0) -- (\x +.1, .2);
\foreach \x in {} 
     \draw (\x-.1, -.2) -- (\x,0) -- (\x +.1, -.2);
\foreach \x in {} 
     \draw (\x-.1, .1) -- (\x +.1, -.1) (\x-.1, -.1) -- (\x +.1, .1);

\draw [-,black,out=270, in=270](-2,0) to (1,0);
\draw [-,black,out=270, in=270](-1,0) to (-0,0);
\draw [-,black,out=270, in=270](2,0) to (5,0);
\draw [-,black,out=270, in=270](3,0) to (4,0);
\draw [-,black,out=270, in=270](-3,0) to (6,0);
%

\end{tikzpicture}
\end{center} 

The plot of $[3,3,1,1]$ has two adjacent sectors. Melting these two gives the plot with support \{ -2,-1,1,3\} with cup diagram 

\bigskip
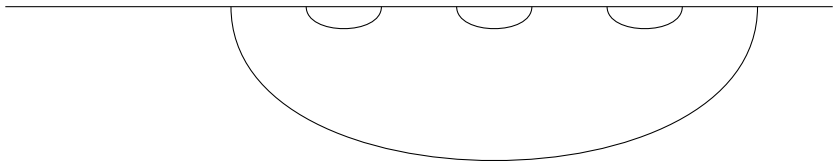
\begin{center}

\begin{tikzpicture}
 \draw (-5,0) -- (6,0);
\foreach \x in {} 
     \draw (\x-.1, .2) -- (\x,0) -- (\x +.1, .2);
\foreach \x in {} 
     \draw (\x-.1, -.2) -- (\x,0) -- (\x +.1, -.2);
\foreach \x in {} 
     \draw (\x-.1, .1) -- (\x +.1, -.1) (\x-.1, -.1) -- (\x +.1, .1);

\draw [-,black,out=270, in=270](-2,0) to (5,0);
\draw [-,black,out=270, in=270](-1,0) to (0,0);
\draw [-,black,out=270, in=270](1,0) to (2,0);
\draw [-,black,out=270, in=270](3,0) to (4,0);
%

\end{tikzpicture}
\end{center} 

{\it The not maximally atypical case}. As with sectors and segments we can define the notion of a plot for representations which are not maximally atypical. We fix the block of the irreducible representations, ie. the positions of the $\times$'s (say at the vertices $x_1, \ldots,x_r$) and the positions of the $\circ$'s (say at the vertices $\circ_1,\ldots, \circ_r$). Once these are fixed we define $\Z_{\times \circ} := \Z \setminus ( \{x_1,\ldots,x_r\}  \cup \{\circ_1, \ldots, \circ_r\})$. Then a plot is a map $\lambda: \Z_{\times \circ} \to \{\boxplus,\boxminus\}$. The reader can convince himself that all the previous definitions and operations on plots (factorization, derivatives etc) can be adapted easily to this more general setting. However this amounts in practice only to fixing the positions of the $\times$ and $\circ$'s and then ignoring them. We will associate in section \ref{sec:loewy-length} to every weight $\lambda$ of atypicality $i$ a plot $\phi(\lambda)$ of rank $i$ (without $\times$'s and $\circ$'s) and work with these instead.

\bigskip




\section{Mixed tensors and ground states} \label{stable0}

\medskip\noindent

We compute $DS(L) \in T_{n-1}$ for special irreducible representations $L$ in a block $\Gamma$, the so-called ground states. The general case for arbitrary $L$ will then be reduced to this case in later sections.

\medskip
Let $MT$ denote the full subcategory of mixed tensors in ${\calR_n}$ whose objects are direct sums of the indecomposable objects in ${\calR}_n$ that appear in a decomposition $X_{st}^{\otimes r} \otimes (X_{st}^{\vee})^{\otimes s}$ for some natural numbers $r,s \geq 0$, where $X_{st} \in {\calR_n}$ denotes the standard representation. By \cite{Brundan-Stroppel-5} and \cite{Comes-Wilson} the indecomposable objects in $MT$ are parametrized by $(n|n)$-cross bipartitions. Let $R_n(\lambda^L,\lambda^R)$ denote the indecomposable representation in ${\calR}_n$ corresponding to the bipartition $\lambda = (\lambda^L ,\lambda^R)$ under this parametrization. 

\medskip

To any bipartition $\lambda$ we attach a weight diagram in the sense of \cite{Brundan-Stroppel-1}, ie. a labelling of the numberline $\mathbb Z$ according to the following dictionary. Put \[ I_{\wedge}(\lambda)  := \{ \lambda_1^L, \lambda_2^L - 1, \lambda_3^L - 2, \ldots \} \quad \text{and}\quad  I_{\vee}(\lambda)  := \{1  -\lambda_1^R, 2 - \lambda_2^R, \ldots \}\ . \] Now label the integer vertices $i$ on the numberline by the symbols $\wedge, \vee, \circ, \times$ according to the rule \[ \begin{cases} \circ \quad \text{ if } \ i \  \notin I_{\wedge} \cup I_{\vee}, \\ \wedge \quad \text{ if } \ i \in I_{\wedge}, \ i \notin I_{\vee}, \\ \vee \quad \text{ if } \ i \in I_{\vee}, \ i \notin I_{\wedge}, \\ \times \quad \text{ if } \ i \in I_{\wedge} \cap I_{\vee}. \end{cases} \]  To any such data one 
attaches a cup-diagram as in section \ref{BS} or \cite{Brundan-Stroppel-1} and we define the following three invariants 

\begin{align*} a(\lambda) & = \text{ number of crosses } \\ d(\lambda) & = \text{ number of cups } \\ k(\lambda) & = a(\lambda) + d(\lambda) \end{align*} A bipartition is {\it $(n|n)$-cross} if and only if $k(\lambda) \leq n$. By \cite{Brundan-Stroppel-5} the modules $R( \lambda^L, \lambda^R)$ have irreducible socle and cosocle equal to $L(\lambda^{\dagger})$ where the highest weight $\lambda^{\dagger}$ can be obtained by a combinatorial algorithm from $\lambda$. Let  $$\theta: \Lambda \to X^+(n) \ $$ denote the resulting map $\lambda \mapsto \lambda^\dagger$ between the set of $(n|n)$-cross bipartitions $\Lambda$ and the set $X^+(n)$ of highest weights of $\calR_n$. 

\begin{thm} \cite{Heidersdorf-mixed-tensors} $R = R(\lambda)$ is an indecomposable module of Loewy length $2 d(\lambda) + 1$. It is projective if and only if $k(\lambda) = n$, in which case we have $R = P(\lambda^{\dagger})$. 
\end{thm}

\medskip
Hence $R$ is irreducible if and only if $d(\lambda) = 0$, and then $R = L(\lambda^{\dagger})$.

\medskip
{\it Deligne's interpolating category}. For every $t \in k$ there exists the category $Rep(Gl_t)$ defined in  \cite{Deligne-interpolation}. This is a $k$-linear pseudoabelian rigid tensor category. By construction it contains an object $st$ of dimension $t$, called the standard representation. Given any $k$-linear pseudoabelian tensor category $C$ with unit object and a tensor functor $$F: Rep(Gl_t) \to C$$ the functor $F \to F(st)$ is an equivalence between the category of $\otimes$-functors of $Rep(Gl_t)$ to $C$ with the category of $t$-dimensional dualisable objects $X \in C$ and their isomorphisms. 

\medskip
In particular, given a dualizable object $X$ of dimension $t$ in a $k$-linear pseudoabelian tensor category, a unique tensor functor $F_X: Rep(Gl_t) \to C$ exists mapping $st$ to $X$. Hence, for our categories ${\calR}_n$ and $t=0$, we get a tensor functor $F_n: Rep(Gl_0) \to {\calR}_n$ by mapping the standard representation of $Rep(Gl_0)$ to the standard representation of $Gl(n\vert n)$ in ${\calR}_n$.  Every mixed tensor is in the image of this tensor functor ( \cite{Comes-Wilson}). The indecomposable elements in Deligne's category are parametrized by the set of all bipartitions. The kernel of $F_n$ contains those indecomposables labelled by bipartitions that are not $(n|n)$-cross. Any $(n|n)$-cross bipartion $\lambda$ defines  an indecomposable object  in the image of $Rep(Gl_0)$. We write $R_n(\lambda)$ for $F_n(R(\lambda))$. By the universal property of Deligne's category any tensor functor from $Rep(Gl_0)$ to a tensor category $C$ is fixed up to isomorphism by the choice of an image of the standard representation of $Rep(Gl_0)$. 

\begin{lem}\label{stable2}\cite{Heidersdorf-mixed-tensors} $DS(R_n(\lambda)) = R_{n-1}(\lambda) $ holds unless $R_n(\lambda)$ is projective, in which case $DS(R_n(\lambda)) =0$.  
\end{lem} 

Note that the vanishing of $R_n(\lambda)_x$ in the projective case is just a special case of lemma \ref{van} (i) and (ii).

\medskip{\it Proof}. An easy computation shows that under the Duflo-Serganova functor the standard representation of $\g_n$ is mapped to the standard representation of $\g_{n-1}$. Since any indecomposable mixed tensor module is in the image of a tensor functor from Deligne's category $Rep(Gl_0)$ \cite{Comes-Wilson} the result follows from the commutative diagram \[ \xymatrix{ Rep (Gl_{0}) \ar[d]_{F_{n}} \ar[dr]^{F_{n-1}} & \\ {\calR}_{n} \ar[r]_-{DS} & {\calR}_{n-1} }. \] The kernel of $F_n$ consists of the $R(\lambda)$ with $k(\lambda) > n$, the kernel of $F_{n-1}$consists of the $R(\lambda)$ with $k(\lambda) \geq n$ Hence $R(\lambda) \in ker(DS)$ if and only if $k(\lambda) = n$ which is equivalent to $R(\lambda)$ projective. \qed

\medskip{\bf Example}.  As in section \ref{sec:strategy} put $\A_{S^i} :=R((i),(1^i)) \in {\calR}_n$. By lemma \ref{stable2} we have $(\A_{S^i})_x = \A_{S^i}$ for all $i\geq 1$.

\begin{cor}\label{projective-image} Every indecomposable projective module of ${\calR}_{n-1}$ is in the image of $DS$.
\end{cor}

\medskip{\it Proof}. The indecomposable projective modules are precisely the modules $DS(R(\lambda^L,\lambda^R))$ with $k(\lambda) = n-1$. Note that every indecomposable projective module is a mixed tensor \cite{Heidersdorf-mixed-tensors}.\qed


\medskip
{\it Irreducible mixed tensors}. By the results above the map $\theta: \Lambda \to X^+(n)$ is injective if restricted to bipartitions with $d(\lambda) = 0$. We denote by $\theta^{-1}$ its partial inverse. A closer inspection \cite{Heidersdorf-mixed-tensors} of the assignment $\theta: \lambda \mapsto \lambda^{\dagger}$  shows that $\theta$ and $\theta^{-1}$ are given by the following simplified rule: Define \begin{align*} m  &= \text{ maximal coordinate of a } \times \text{ or } \circ  \\ t & = max( k(\lambda) + 1, m + 1) \\ s & = \begin{cases} 0 \quad \quad \quad \quad \ m + 1 \leq k(\lambda) + 1 \\ m - k(\lambda)  \quad m + 1 > k(\lambda) + 1 \end{cases} \end{align*} The weight diagram of $\lambda^{\dagger}$ is obtained from the weight diagram of $\lambda$ by switching all $\vee$'s to $\wedge$'s and vice versa at positions $\geq t$ and switching the first $s + n - k(\lambda)$ $\vee$'s at positions $< t$ to $\vee$'s and vice versa. The numbers labelled by a $\wedge$ or $\vee$ will be called free positions. 
Conversely if 
$L(\lambda^{\dagger})$ is some irreducible representation in $MT$, the corresponding bipartition with $\theta^{-1}(\lambda^{\dagger}) = \lambda$ is obtained in the same way: Define $t, m, s$ as above and apply the same switching rules to the weight diagram of $\lambda^{\dagger}$.

\begin{prop}\label{mixed-tensor-derivative} Let $$L = L(\lambda^{\dagger}) = L(\lambda_1, \ldots, \lambda_{n-i}, 0,\ldots, 0\ ;\ 0, \ldots,0, \lambda_{n+i+1},\ldots, \lambda_{2n})$$ be an irreducible $i$-fold atypical representation. Then $L$ is a mixed tensor $L = R(\lambda)$ for a unique bipartition of defect 0 and $rk = n-i$. Then \[ DS(L) = R_{n-1}(\lambda) = L(\bar{\lambda}^{\dagger})\ ,\] where $\bar{\lambda}^{\dagger}$ is obtained from $\lambda^{\dagger}$ by removing the two innermost zeros corresponding to $\lambda^{\dagger}_n$ and $\lambda^{\dagger}_{n+1}$.                                                                                                                                                                                            \end{prop}

\medskip\noindent
{\it Proof}. We apply $\theta^{-1}$ to $\lambda^{\dagger}$. It transforms the weight diagram of $\lambda^{\dagger}$ into some other weight diagram which might not be the weight diagram of a bipartition. However if the resulting weight diagram is the weight diagram of an $(n|n)$-cross bipartition of defect 0, then $\theta(\lambda) = \lambda^{\dagger}$ and $R(\lambda) = L(\lambda^{\dagger})$. For $\lambda^{\dagger}$ \begin{align*} I_{\times} & = \{ \lambda_1,\lambda_2 -1,\ldots,\lambda_{n-i} - (n-i) + 1, - n + i, \ldots, -n + 1 \} \\ I_{\circ} & = \{ 1-n,2-n,\ldots,i-n,i+1-n - \lambda_{n+1+i},\ldots, - \lambda_{2n} \}.\end{align*} Then $I_{\times} \cap I_{\circ} = \{-n+1, \ldots,-n+i\}$ (since the atypicality is $i$) and the $n-i$ crosses are at the positions $\lambda_1,\lambda_2 -1,\ldots,\lambda_{n-i} - (n-i) + 1$ and the $n-i$ circles at the positions $i+1-n - \lambda_{n+1+i},\ldots, - \lambda_{2n}$. Define $m,t,s$ as above. Note that $k(\lambda) = n-i$. We distinguish two cases, either $t = n-i+1$ or $t=m+1$. Assume first $m + 1 
\leq n-i + 1$. Switch all free labels at positions $ \geq t$ and the first $n-(n-i) = i$ free labels at positions $<t$. By assumption the $2n- 2i$ crosses and circles lie at positions $> i- n$ and $< n-i+1$. However there are exactly $2n-2i$ such positions. Hence the switches at positions $<t$ turn exactly the $i$ $\vee$'s at positions $i-n,\ldots,1-n$ into $\wedge$'s. In the second case $t = m + 1 > n-i+1$ switch the first $m + n - 2 (n-i)$ free labels at positions $<t$. There are exactly $m + n -i$ positions between $m$ and $i-n$, $m - n +2i$ switches and $2n - 2i$ crosses and circles between $i-n$ and $t$. This results in $m - n + i$ free positions between $i-n$ and $t$. The remaining $i$ switches transform the $i$ $\vee$'s to $\wedge$'s. Hence in both cases $\theta^{-1}$ transforms the weight diagram of $\lambda^{\dagger}$ into a weight diagram where the rightmost $\wedge$ is at position $i-n$ and the leftmost $\vee$ is at the first free position $> i-n$ and all labels at positions $\geq t$ are given by $\vee$'s. This is the weigth diagram of a bipartition of defect 0 and rank $n-i$. Indeed the labelling defines the two sets $I_{\wedge}$ and $I_{\vee}$ and this defines two tuples $\lambda^L = (\lambda^L_1,\lambda^L_2,\ldots)$ and $\lambda^R = (\lambda^R_1,\lambda^R_2,\ldots)$. The positioning of the $\wedge$'s implies that $\lambda^L_{n-i+1} = 0$ and the positioning of the $\vee$'s implies $\lambda^R_t = 0$. Clearly $\lambda_1^L = \lambda_1 > 0$ and $\lambda^R_1 \geq 0$. Hence the pair $\lambda:= (\lambda^L,\lambda^R)$ is a bipartition (of defect 0 and rank $n-i$) and $\theta(\lambda) = \lambda^{\dagger}$. It remains to compute the highest weight of $R_{n-1}(\lambda)$. The two sets $I_{\vee}$ and $I_{\wedge}$ and accordingly the weight diagram of $\lambda$ do not depend on $n$. Neither do $t, m, s$ and the switches at positions $\geq t$. To get $\lambda^{\dagger}$ in ${\calR}_n$ from $\lambda$ we switch the first $s + n - (n-i)$ free labels $<t$. To get $\lambda^{\dagger}$ in ${\calR}_{n-1}$ from $\lambda$ we switch the first $s + (n-1) - (n-i)$ free labels $<t$. This results in removing the leftmost $\vee$ at position $1-n$. \qed


\bigskip
{\it Ground states}. 
Let ${\calR}_n^{i} \subset {\calR}_n$ denotes the full subcategory of $i$-atypical objects.
Every block in $\calR_n^i$ contains irreducible objects with the property
that all $i$ labels $\vee$ are adjacents and to the left of all $n-i$ labels $\times$ and all
$n-i$ labels $\circ$. We call such an irreducible object a groundstate of the corresponding 
block in $\calR_n^i$. Each block in $\calR_n^i$ uniquely defines its groundstate up to a simultaneous shift of the $i$ adjacent labels $\vee$. The weight $\lambda$ of such a groundstate $L(\lambda)$ is of the form
$$  \lambda = (\lambda_1,...,\lambda_{n-i},\lambda_n,...,\lambda_n\ ;\
-\lambda_n,...,-\lambda_n,\lambda_{n+1+i}, ..., \lambda_{2n}) \ .$$
with $\lambda_n \leq \min(\lambda_{n-i}, - \lambda_{n+1+i})$ (here $\lambda_n \mapsto 
\lambda_n - 1$ corresponds to the shift of the $i$ adjacent labels $\vee$).
The coefficients $\lambda_1,...,\lambda_{n-i}, \lambda_{n+1+i}, ..., \lambda_{2n}$
determine and are determined by the position of the labels $\times$ and $\circ$ defining
the given block in $\calR_n^i$. We define $$\overline \lambda
= (\lambda_1,...,\lambda_{n-i},\lambda_n,...\ ;\
...,-\lambda_n,\lambda_{n+1+i}, ..., \lambda_{2n}) \ $$ by omitting the innermost  $\lambda_n ; -\lambda_n$ pair.
Then $ L(\overline\lambda) \in \calR_n^{i-1} \subset T_{n-1} $.

\medskip  
{\it Berezin twists}. Twisting with $Ber=Ber_n$ induces an endofunctor of $\calR_n^i$
and permutes blocks. By a suitable twist one can replace a given block in $\calR_n^i$ such that
it contains the groundstate 
$$   \lambda' = (\lambda_1-\lambda_n,...,\lambda_{n-i}-\lambda_n,0,...,0\ ;\
0,...,0,\lambda_{n+1+i}+\lambda_n, ..., \lambda_{2n}+\lambda_n) \ .$$

\begin{prop}\label{13}
For a groundstate $L=L(\lambda)$ of a block in $\calR_n^i \subset \calR_n$ the image 
$DS(L)$ in $T_{n-1}$ of $L$ under the Duflo-Serganova functor is
$$   DS(L(\lambda)) = \Pi^{-\lambda_n}L(\overline\lambda) \ $$
for $i>0$ or $DS(L)=0$ for $i=0$.
\end{prop}

\medskip\noindent
In particular therefore theorem \ref{mainthm} holds for the groundstates $L=L(\lambda)$ of blocks 
in $\calR_n^i \subset \calR_n$.

\medskip\noindent
{\it Proof}. We can assume $i>0$. Then we can assume 
$\lambda_{n}=\lambda_{n+1}=0$ by a suitable Berezin twist. Hence
$$    L = R_n(\lambda^L,\lambda^R)  $$
for an $(n|n)$-cross bipartition $(\lambda^L,\lambda^R)$  and therefore
$$  DS(L) = R_{n-1}(\lambda^L,\lambda^R) $$
is irreducible of weight $\overline\lambda$, i.e. $DS(L(\lambda))=L(\overline\lambda)$.
This proves the claim, since by assumption now $\lambda_{n+1}=0$. \qed


\medskip\noindent



\section{Sign normalizations}\label{signs}

\medskip\noindent

The main theorem \ref{mainthm} asserts in particular that $DS(L)$ is semisimple. In order to show this we define a sign $\varepsilon(L) \in \{ \pm 1 \}$ for every irreducible representation $L$ with the property that $Ext^1_{\calR_n}(L(\lambda),L(\mu)) = 0$ if $\varepsilon(\lambda) = \varepsilon(\mu)$. More precisely, this sign should satisfy the following conditions: 

\begin{enumerate}
\item Let $\calR_n(\varepsilon)$ denote the full subcategory of all objects whose irreducible constituents $X$ have sign $\varepsilon(X)$. Similarly define the full subcategories $\Gamma(\varepsilon)$ for a block $\Gamma$. Then we require that the categories $\calR(+)$ and $\calR(-)$ are semisimple categories.
\end{enumerate}

\medskip

Clearly it is enough to require that the categories $\Gamma(\varepsilon)$ are semisimple. Note that any such sign function is unique on a block $\Gamma$ up to a global sign $\pm 1$. Hence if we normalize the sign by $\varepsilon(L) = +1$ where $L$ is an irreducible representation in a block $\Gamma$, the sign is uniquely determined on $\Gamma$. Our second condition is: 

\begin{enumerate}
\setcounter{enumi}{1}
\item $\varepsilon(L(\lambda)) = 1$ if $L(\lambda)$ is an irreducible mixed tensor (see section \ref{stable0}) and $\varepsilon(\one) = 1$.
\end{enumerate}

If $L(\lambda)$ is maximal atypical, we put \[ \varepsilon(L(\lambda)) = (-1)^{p(\lambda)}\] for the parity $p(\lambda) = \sum_{i=1}^n \lambda_{n+i}$. In the maximal atypical case we have $Ext^1_{\calR_n}(L(\lambda), L(\mu)) = 0$ if $p(\lambda)  \equiv p(\mu) \ mod \ 2$ by \cite{Weissauer-gl}. Hence the categories $\Gamma_n(\pm)$ are semisimple. This determines the sign $\varepsilon$ up to a global $\pm 1$ on each block $\Gamma$ of atypicality $i$. Indeed by \cite{Serganova-blocks} any block of atypicality $i$ is equivalent to the maximal atypical block $\Gamma_i$ of $\calR_i$. We fix once and for all a particular equivalence denote it by $\tilde{\phi}_n^i$. We describe the effect of $\tilde{\phi}_n^i$ on an irreducible module $L(\lambda)$ of atypicality $i$ \cite{Serganova-blocks} \cite{Gruson-Serganova}. For  an $i$-fold atypical weight in $ X^+(n)$  its weight diagram  has $n-i$ vertices labelled with $\times$ and $n-i$ vertices labelled with $\circ$. Let $j$ be the leftmost vertex labelled  either by $\times$ or $\circ$. By removing this vertex and shifting all vertices at  the positions $> j$ one position to the left, recursively we remove all vertices labelled by $\times$ or $\circ$ from the given weight diagram. The remaining finite subset $K$ of labels $\vee$ has cardinality $i$ and the weight diagram so obtained defines a unique irreducible maximally atypical module in $\calR_{i}$. Under $\tilde{\phi}_n^i$ the irreducible representation $L(\lambda)$ maps to the irreducible $Gl(i|i)$-representation described by the removal of crosses and circles above. We denote the weight of this irreducible representation by $\phi_n^i(\lambda)$. We make the preliminary definition $\varepsilon(L(\lambda)) = (-1)^{p(\phi_n^i(L(\lambda)))}$ for a weight $\lambda$ of atypicality $i$. We claim that this sign satisfies condition 1. Indeed in the maximal atypical case we have $Ext^1_{\calR_i}(L(\lambda), L(\mu)) = 0$ if $p(\lambda)  \equiv p(\mu) \ \text{ mod } \ 2$ by \cite{Weissauer-gl}. Hence the categories $\Gamma_i(\pm)$ are semisimple. The equivalence $(\tilde{\phi}_n^i)^{-1}$ between the two abelian categories $\Gamma_i$ and $\Gamma$ is exact and sends the semisimple category $\Gamma_i(\pm) \subset \calR_i$ to the category $\Gamma(\pm)$ by the definition of the sign $\varepsilon$; and $\Gamma(\pm)$ is semisimple.

\medskip

This preliminary definition however doesn't satisfy condition 2. If $L  = L(\lambda_1, \ldots, \lambda_{n-i}, 0,\ldots, 0\ ;\ 0, \ldots,0, \lambda_{n+i+1},\ldots, \lambda_{2n})$ is the $i$-atypical mixed tensor of section \ref{stable0}, its weight diagram has $n-i$ circles at the vertices $\lambda_1,\lambda_2 -1,\ldots,\lambda_{n-i} - (n-i) + 1$ and $n-i$ circles at the vertices $i+1-n - \lambda_{n+1+i},\ldots, - \lambda_{2n}$. The $i$ $\vee$'s are to the left of the crosses and circles. Applying $\tilde{\phi}_n^i$ removes the crosses and circles but leaves the $\vee$'s unchanged at the vertices $-n + 1, \ldots, -n+i$. The irreducible representation of $Gl(i|i)$ so obtained is $Ber^{-n+i}$ with $p([-n+i,\ldots,-n+i]) = i(-n+i)$. Hence the preliminary sign $\varepsilon(L(\lambda))$ of a mixed tensor of atypicality $i$ is $\varepsilon(L(\lambda)) = (-1)^{i(-n+i)}$. In order to satisfy condition 2) we have to normalize the sign by the additional factor $(-1)^{i(-n+i)}$ for an $i$-atypical weight, and we define \[ \varepsilon(\lambda) = (-1)^{i(-n+i)}  (-1)^{p(\tilde{\phi}_n^i(L(\lambda)))}\] where $p$ is the parity in the maximal atypical block of $Gl(i|i)$. This sign satisfies condition 1) and 2) by construction, and it is the unique sign with these properties. Note that our definition implies that the sign of a typical weight in $\calR_n$ is always positive.

\medskip

The additional sign factor can be understood as follows. The unique irreducible mixed tensor should play an analogous rule to the trivial representation $\one$ in $\calR_i$. W can modify the block equivalences $\tilde{\phi}_n^i$ as follows: Since the mixed tensor $L(\lambda)$ maps to $Ber^{-n+i}$ we twist with the inverse and define the normalized equivalence \[ \phi_n^i(L) = Ber^{n-i} \otimes \tilde{\phi}_n^i(L).\] Then we obtain $(-1)^{p(\phi_n^i(L(\lambda)))} = (-1)^{-i(-n+i)} (-1)^{p(\tilde{\phi}_n^i(L(\lambda)))}$. Hence \[ \varepsilon(\lambda) = (-1)^{p(\phi_n^i(L(\lambda)))}.\] 

\begin{cor} \label{semisimple-sign} The categories $\calR_n(\varepsilon)$ are semisimple categories.
\end{cor}

The sign $\varepsilon$ will automatically have the following important property: The translation functors of section \ref{sec:loewy-length} $\A = F_i L(\lambda_{\times \circ})$ have Loewy structure $(L,A,L)$ with $L \in \calR_n(\pm \varepsilon)$ and $A \in \calR_n(\mp \varepsilon)$. This is required in our axioms in section \ref{sec:inductive}. This property follows  immediately from the maximal atypical case \cite{Weissauer-gl} due to the description of the composition factors of $F^i L(\lambda_{\times \circ})$ given in section \ref{sec:loewy-length}. 

\medskip

The significance of the sign function is the following: $DS(L)$ is in general a representation in $T_{n-1}$, not $\calR_{n-1}$. The sign factor regulates whether an irreducible summand of $DS(L)$ is in $\calR_{n-1}$ or $\Pi \calR_{n-1}$, see theorem \ref{mainthm}.

\medskip

In the proof of the main theorem we use the language of plots of section \ref{derivat} for uniform bookkeeping. Each maximally atypical irreducible representation $L$ defines a plot: The segments of the cup diagram of $L$ define the segments of the plot, and the $\vee$ in the weight diagram give the support of the plot. If $L$ is not maximally atypical, we associate to it a plot via the map $\phi$ of section \ref{sec:loewy-length}. For each plot we defined its derivative $\partial \lambda$ in section \ref{derivat}. If $\lambda$ is a prime plot given by the sector $(I,K)$, $I = [a,b]$, then $\partial(I,K) = (I,K)' = (I',K')$ for $I'=[a+1,b-1]$ and $K'=I'\cap K = [a,b]$ defines the derivative. We normalize the derivative of section \ref{derivat} and put \[ \lambda' = (-1)^{a+n-1} \partial (I,K) \] where $I = [a,b]$ (i.e. the leftmost $\vee$ is at $a$). The reason for this is as follows. The sign has to be normalized in such a way that for objects $X=L(\lambda)$ in the stable range of the given block we get $d(X) = X'$ for the map $d$ of section \ref{sec:inductive}. Assume first that we are in the maximally atypical $\calR_n$-case and consider a weight with associated prime plot $\lambda$. The parity of the weight $\lambda$ is $p(\lambda) = \sum_{i=1}^n \lambda_{n+i}$. Applying $DS$ removes the $\vee$ in the outer cup. The parity of the resulting weight in $T_{n-1}$ is given by $p(\lambda') = \sum_{i=1}^{n-1} \lambda_{n+i}$, hence $p(\lambda) - p(\lambda') = \lambda_n$ and we get a shift by $n_i \equiv (-1)^{\lambda_n}$ according to theorem \ref{mainthm}. The leftmost $\vee$ is at the vertex $a = \lambda_n - n +1$, hence $(-1)^{a+n-1} = (-1)^{\lambda_n}$ and the two shifts agree.

\medskip

Let us now assume $at(L(\lambda)) = k <n$ and that the weight defines a prime plot of rank $k$. Here we have to use the normalized plot associated to the weight $\lambda$ by the map $\phi$ from section \ref{sec:loewy-length} in which case the two shifts agree again. We may pass to the maximally atypical case due to the lemmas \ref{shift-1}, \ref{shift-2}, \ref{shift-3} which allow us to shift all the circles and crosses sufficiently far to the right.

\medskip

 
\section{The Main theorem}\label{sec:main}

\medskip
In the main theorem we calculate $DS(L) \in T_{n-1}$ for any irreducible $L$. We refine this in section \ref{kohl-2}, \ref{koh3} and compute the $\Z$-grading of $DS(L)$.

\begin{thm} \label{mainthm}  Suppose $L(\lambda)\in \calR_n$ is an irreducible
atypical representation, so that $\lambda$ corresponds to
a cup diagram $$ \bigcup_{j=1}^r \ \ [a_j,b_j] $$ with $r$ sectors
$[a_j,b_j]$ for $j=1,...,r$. Then $$DS(L(\lambda)) \ \cong\ \bigoplus_{i=1}^r  \ \Pi^{n_i} L(\lambda_i)$$ is the direct sum of irreducible atypical
representations $L(\lambda_i)$ in $\calR_{n-1}$ with shift $n_i \equiv \varepsilon(\lambda)
- \varepsilon(\lambda_i)$ modulo 2. The representation $L(\lambda_i)$ is uniquely defined by
the property that its cup diagram is $$ [a_i +1, b_i-1] \ \ \ \cup \ \ \bigcup_{j=1, j\neq i}^r \ \ [a_j,b_j] \ ,$$ the union of the sectors $[a_j,b_j]$ for $1\leq j\neq i \leq r$ and (the sectors occuring in) the segment $[a_i+1,b_i-1]$.
\end{thm}
  
\medskip
{\it Consequence}. {\it In particular this implies that for irreducible representation $(V,\rho)$ the $G_{n-1}$-module $H^+(V)\oplus H^-(V)$ is semisimple in $\calR_{n-1}$ and multiplicity free.
Furthermore the sign of the constituents in $H^\pm(V)$ is $\pm sign(V)$.}

\medskip

If we use the language of plots, the main theorem says that the irreducible summands of $DS(L)$ are given by the derivatives of the sectors of the plot associated to $\lambda$.

\medskip
{\bf Example.} The maximally atypical weight $[3,0,0]$ has cup diagram 

\medskip
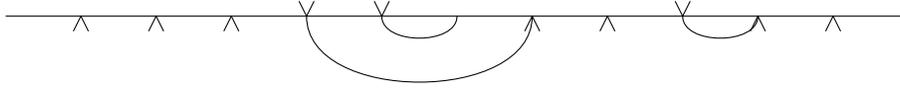
\begin{center}

  \begin{tikzpicture}
 \draw (-6,0) -- (6,0);
\foreach \x in {3,-1,-2} 
     \draw (\x-.1, .2) -- (\x,0) -- (\x +.1, .2);
\foreach \x in {-5,-4,-3,1,2,4,5} 
     \draw (\x-.1, -.2) -- (\x,0) -- (\x +.1, -.2);
\foreach \x in {} 
     \draw (\x-.1, .1) -- (\x +.1, -.1) (\x-.1, -.1) -- (\x +.1, .1);
\foreach \x in {} 
     \draw  node at (1,0) [fill=white,draw,circle,inner sep=0pt,minimum size=6pt]{};

\draw [-,black,out=270,in=270](-2,0) to (1,0);
\draw [-,black,out=270,in=270](-1,0) to (0,0);
\draw [-,black,out=270,in=270](3,0) to (4,0);

\end{tikzpicture}
\end{center} 

It splits into the two irreducible representations $[3,0]$ and $\Pi[-1,-1]$ in $\calR_2 \oplus \Pi\calR_2$. We will later compute its cohomology in proposition \ref{hproof} and obtain $H^{\bullet} ([3,0,0]) = S^2\langle0\rangle \oplus Ber^{-1}\langle-1\rangle$.

\medskip
{\bf Example}. Denote by $S^i$ the irreducible representation $[i,0,\ldots,0]$. Consider a nontrivial
extension 
$0 \to S^2 \to E \to Ber(S^2)^\vee \to 0$ in $\calR_3$  (such extensions exist).
Then $sdim(E)= 0$ and $E$ is indecomposable, hence negligible.
The derivative of $S^2=[2,0,0]$ (in the sense of plots) is
$  (S^2)' =   [2,0]   +   Ber^{-1}$
and the derivative
of $Ber(S^2)^\vee = [2,2,1]$ is   
$[2,2,1]' = (Ber[1,1,0])' = -Ber ([1,1,0]') = - [2,2]   -   [2,0]$.
From $ [2,2]=Ber^2$ then 
$H^+(E) = Ber^{-1}  \oplus    ? $
and $   H^-(E)  = Ber^2      \oplus    ?$
where ? is either $[2,0]$ or zero. 
Hence $E$ is negligible in $\calR_3$, 
but $D^2(E)= D(D(E)) \neq 0$.
In particular, $D(\calN_3)$ {\it is not contained} in $\calN_2$.

\medskip

{\it Block equivalences}. Applying $DS$ is compatible with taking the block equivalences $\phi_n^i$ and $\tilde{\phi}_n^i$ in the sense $DS(\phi_n^i(L)) = \phi_n^i(DS(L))$ for irreducible $L$ by the main theorem. This can be extended to arbitrary modules $M$ in a block. By \cite{Serganova-blocks} and \cite{Kujawa-generalized-kac-wakimoto} the block equivalence $\tilde{\phi}_n^i$ between an $i$-atypical block $\Gamma$ and the unique maximal atypical block of $Gl(i|i)$ is obtained as a series of translation functors, a restriction and the projection onto a weight space. Call a weight $\lambda$ in $\Gamma$ stable if all the $\vee$'s are to the left of all crosses and circles. By \cite{Serganova-blocks} we can apply a suitable sequence of translation functors to any indecomposable module in $\Gamma$ until all its composition factors are stable. We recall now the definition of $\tilde{\phi}_n^i$ as in \cite{Kujawa-generalized-kac-wakimoto} on an indecomposable module $M$. Embed $\mathfrak{gl}(k|k)$ as an inner block matrix in $\mathfrak{gl}(n|n)$. Let $\mathfrak{l} = \mathfrak{gl}(k|k) + \mathfrak{h}$ where $\mathfrak{h}$ are the diagonal matrices. Then choose $\mathfrak{h'} \subset \mathfrak{h}$ such that $\mathfrak{h'}$ is a central subalgebra of $\mathfrak{l}$ and $\mathfrak{l} = \mathfrak{gl}(k|k) \oplus \mathfrak{h'}$. We denote the restriction of a weight $\lambda$ to $\mathfrak{h'}$ by $\lambda'$.  Now move $M$ by a suitable sequence of translation functors until its composition factors are stable. The block $\Gamma$ is the full subcategory of modules admitting some central character $\chi_{\mu}$. Now define $Res_{\mu'}(M)  = \{  m \in M \ | \ h'm = \mu'(h') m \text{ for all } h' \in \mathfrak{h}' \}$.  Then, on a module $M$ with stable composition factors, the functor $\tilde{\phi}_n^i(M)$ is given by $Res_{\mu'}(M)$. Alternatively we could first restrict to $\mathfrak{l}$ and then project on the $\mu'$-eigenspace. By \cite{Serganova-kw}, cor 4.4, and the main theorem, $DS$ induces a bijection between the blocks in $T_n$ and $T_{n-1}$, and for any $M$ in $T_n$ $DS(F_i(M)) = F_i(DS(M))$. Since our fixed $x$, which we choose in the definition of $DS$, is contained in the embedded $\mathfrak{gl}(k|k)$, the operation of $\rho(x)$ on $Res(M)$ or on its $\lambda'$-eigenspace is the same as of $\rho(x)$ on $M$.
Hence $DS$ is clearly compatible with restriction, but it also doesn't matter whether we first apply $DS$ and project onto the $\lambda'$-eigenspace or first project to the $\lambda'$-eigenspace and then apply $DS$ since $\rho(x)$ commutes with $\mathfrak{h}'$. Hence $DS(\tilde{\phi_n^i}(M)) = \tilde{\phi}_n^i(DS(M))$ holds for any $M$, and the analogous statement for $\phi_n^i$ follows immediately. 
To summarize: If $\Gamma'$ denotes the unique block obtained from the $i$-atypical $\Gamma$ via $DS$, we obtain a commutative diagram \[ \xymatrix@+1,5cm{ \Gamma \ar^{\phi_n^i}[r] \ar_{DS}[d] & \Gamma_i \ar^{DS}[d] \\ \Gamma' \ar^{\phi_{n-1}^{i-1}}[r]& \Gamma_{i-1}}.\] \\ 

The main theorem has a number of useful consequences:

\medskip

\textit{Cohomology.} The main theorem permits us to compute the cohomology $H^i(L)$ of irreducible modules $L$ in section \ref{kohl-2} and \ref{koh3}. Although the calculation of the $\Z$-grading of $DS(L)$ is much stronger than the $\Z_2$ version of theorem \ref{mainthm}, it should be noted that the proof is based on the main theorem and a careful bookkeeping of the moves in section \ref{sec:moves}. 
\medskip

\textit{Spectral sequences.} The main theorem also shows the degeneration of the spectral sequences from section \ref{m} and shows \[ DS_{n,n_2}(L) \simeq DS_{n_1,n_2}(DS_{n,n_1}(L)).\] The degeneration can be extended in a similar way to the not maximally atypical case, see below.

\medskip

\textit{Tensor products.} The main theorem allows us to reduce some questions about tensor products of irreducible representations to lower rank. Since $DS$ is a tensor functor we have $DS(L(\lambda) \otimes L(\mu)) = DS(L(\lambda)) \otimes DS(L(\mu)) = \bigoplus_{i,j} (\Pi^{n_i} L(\lambda_i)) \otimes (\Pi^{n_j} L(\mu_j))$. If we inductively understood the tensor product in $T_{n-1}$, we would obtain estimates about the number of indecomposable summands and composition factors in this way. We use this method to calculate the tensor product of two maximal atypical representations of $Gl(2|2)$ in \cite{Heidersdorf-Weissauer-gl-2-2}, see also \cite{Heidersdorf-semisimple-quotient}.

\medskip

\textit{Negligible modules and branching laws.} The functor $DS$ does not preserve negligible modules as the example above shows. However when we restrict $DS$ to the full subcategory $\mathcal{RI}_n$ of modules which arise in iterated tensor products of irreducible representations, $DS$ induces a functor $DS: \mathcal{RI}_n/\mathcal{N} \to \mathcal{RI}_{n-1}/\mathcal{N}$. We show in \cite{Heidersdorf-Weissauer-tannaka} \cite{Heidersdorf-semisimple-quotient} that $\mathcal{RI}_n/\mathcal{N}$ is equivalent as a tensor category to the representation category of a proreductive group $H_n$. We also show that there is an embedding $H_{n-1} \to H_n$, and $DS$ can be identified with the restriction functor with respect to this embedding. In other words $DS$ gives us the branching laws for the restriction of the image of $L(\lambda)$ in $Rep(H_n)$ to the subgroup $H_{n-1}$. 

\medskip

\textit{Superdimensions and modified superdimensions.} The main theorem can be used to reprove parts of the generalized Kac-Wakimoto conjecture on modified superdimensions \cite{Serganova-kw}. In fact we derive a closed formula for the modified superdimension. We sketch this and prove the analog of proposition \ref{Leray2}. 

\medskip

\textit{A superdimension formula}. Assume $L$ maximally atypical. If $sdim(L) > 0$, $$DS(L(\lambda)) \ \cong\ \bigoplus_{i=1}^r  \ \Pi^{n_i}( L(\lambda_i) )$$ splits into a direct sum of irreducible modules of positive superdimension. Indeed the parity shift $\Pi^{n_i}$ occurs if and only if $p(\lambda) \not\equiv p(\lambda_i) \ mod \ 2$. Hence $DS^{n-1}(L)$ splits into a direct sum of irreducible representations of superdimension 1. Applying $DS$ $n$-times gives a functor $DS^n: \calR_n \to svec$, hence $DS^n(L) \simeq m \ k \oplus m' \Pi k$ for positive integers $m,m'$, hence $m = 0$ if and only if $sdim(L) < 0$ and $m' = 0$ if and only $sdim(L) > 0$. By \cite{Weissauer-gl} the superdimension of a maximal atypical irreducible representation in $\calR_n$ is given by \[ sdim(L(\lambda)) =( -1)^{p(\lambda)} m(\lambda)\] for a positive integer $m(\lambda)$ (see below for the definition). In particular \[ m(\lambda) = \begin{cases} m   & p(\lambda) \equiv 0 \text{ mod } 2 \\  m'  & p(\lambda) \equiv 1 \text{ mod } 2. \end{cases} \] By proposition \ref{Leray2} this also holds for $DS_{n,0}:\calR_n \to svec$: If $DS_{n,0} (L) \simeq m\ k \oplus m' \Pi k$, we get that either $m$ or $m'$ is zero.

\medskip

The positive integer $m(\lambda)$ for a maximally atypical weight can be computed as follows. We refer to \cite{Weissauer-gl}, but it would be an easy exercise to deduce this from the main theorem. We let $\underline{\lambda}$ be the associated oriented cup diagram to the weight $\lambda$ as defined in section \ref{BS}. To each such cup diagram we can associate a forest $\mathcal{F}(\lambda)$ with $n$ nodes, i.e. a disjoint union of rooted trees as in \cite{Weissauer-gl}. Each sector of the cup diagram corresponds to one rooted planar tree. We read the nesting structure of the sector from the bottom to the top such that the outer cup corresponds to the root of the tree. If the following is a sector of a cup diagram

\medskip
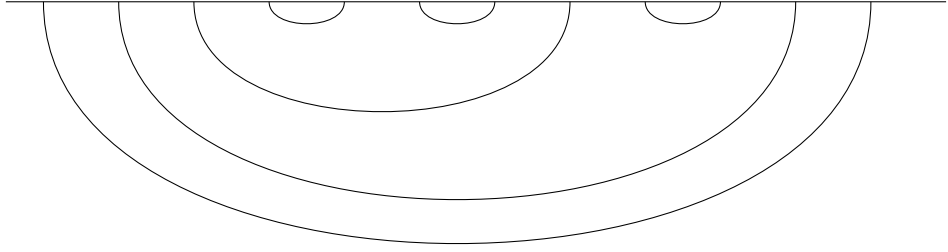
\begin{center}

  \begin{tikzpicture}
 \draw (-6.5,0) -- (6,0);
\foreach \x in {} 
     \draw (\x-.1, .2) -- (\x,0) -- (\x +.1, .2);
\foreach \x in {} 
     \draw (\x-.1, -.2) -- (\x,0) -- (\x +.1, -.2);
\foreach \x in {} 
     \draw (\x-.1, .1) -- (\x +.1, -.1) (\x-.1, -.1) -- (\x +.1, .1);

\draw [-,black,out=270,in=270](-5,0) to (4,0);
\draw [-,black,out=270,in=270](-4,0) to (1,0);
\draw [-,black,out=270,in=270](-3,0) to (-2,0);
\draw [-,black,out=270,in=270](-1,0) to (0,0);
\draw [-,black,out=270,in=270](2,0) to (3,0);
\draw [-,black,out=270,in=270](-6,0) to (5,0);

\end{tikzpicture}

\end{center}

then the associated planar rooted tree is 
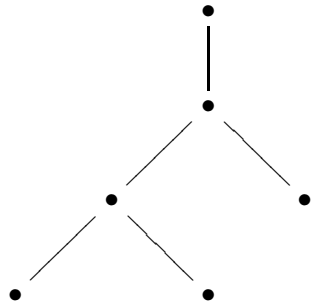
\[ \xymatrix{ & & \bullet \ar@{-}[d] & \\ & & \bullet \ar@{-}[dr] \ar@{-}[dl] & \\ &  \bullet \ar@{-}[dr] \ar@{-}[dl] & & \bullet \\ \bullet &  & \bullet &   } \]

If $\mathcal{F}$ is a forest let $|\mathcal{F}|$ the number of its nodes. We define the forest factorial $\mathcal{F}!$ as the the product $\prod_{x \in \mathcal{F}} |\mathcal{F}_x|$ where $\mathcal{F}_x$ for a node $x \in \mathcal{F}$ denotes the subtree of $\mathcal{F}$ rooted at the node $x$. Then the multiplicity is given by  \[ m(\lambda) = \frac{ |\mathcal{F}(\lambda)|!}{\mathcal{F}(\lambda)!}.\] For example $m(\lambda)$ for irreducible module in $\calR_4$ with cup diagram

\medskip
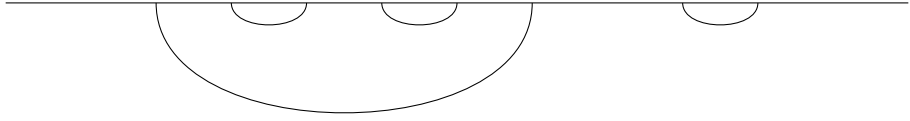
\begin{center}

  \begin{tikzpicture}
 \draw (-6,0) -- (6,0);
\foreach \x in {} 
     \draw (\x-.1, .2) -- (\x,0) -- (\x +.1, .2);
\foreach \x in {} 
     \draw (\x-.1, -.2) -- (\x,0) -- (\x +.1, -.2);
\foreach \x in {} 
     \draw (\x-.1, .1) -- (\x +.1, -.1) (\x-.1, -.1) -- (\x +.1, .1);

\draw [-,black,out=270,in=270](-4,0) to (1,0);
\draw [-,black,out=270,in=270](-3,0) to (-2,0);
\draw [-,black,out=270,in=270](-1,0) to (0,0);
\draw [-,black,out=270,in=270](3,0) to (4,0);

\end{tikzpicture}

\end{center}

is computed as follows: The associated planar forest is 
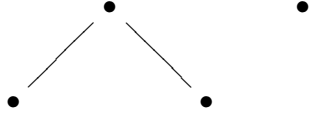
\[ \xymatrix{ & \bullet \ar@{-}[dr] \ar@{-}[dl] & & \bullet \\ \bullet &  & \bullet & &  } \] Hence \[sdim (L(\lambda)) = \frac{24}{3 \cdot 1 \cdot 1 \cdot 1} = 8.\]

\medskip

\textit{Modified Superdimensions} If $at(L(\lambda)) <n$, $sdim(L) = 0$. However one can define a modified superdimension for $L$ as follows. We recall some definitions and results from \cite{Kujawa-generalized-kac-wakimoto}, \cite{Geer-Kujawa-Patureau-Mirand} and \cite{Serganova-kw}. Denote by $c_{V,W}: V \otimes W \to W \otimes V$ the usual flip $v \otimes w \mapsto (-1)^{p(v)p(w)} w \otimes v$.  Put $ev'_V = ev_V \circ c_{V,V^{\vee}}$ and $coev'_V = c_{V,V^{\vee}} \circ coev_{V}$ for the usual evaluation and coevaluation map in the tensor categories $\calR_n$ and $T_n$. For any pair of objects $V,W$ and an endomorphism $f:V\otimes W \to V \otimes W$ we define \begin{align*} tr_L(f) & = (ev_V \otimes id_W) \circ (id_{V^{\vee}} \otimes f) \circ (coev'_{V} \circ id_w)  \in End_T(W) \\ tr_R(f) & = (id_V \otimes ev'_W ) \circ (f \otimes id_{W^{\vee}} ) \circ (id_V \otimes coev_{W})  \in End_T(V) \end{align*} 

For an object $J \in {\calR}_n$ let $I_J$ be the tensor ideal generated by $J$. A trace on $I_J$ is by definition a family of linear functions
\[t = \{t_V:\End_{{\calR}_n}(V)\rightarrow k \}\]
where $V$ runs over all objects of $I_J$ such that following two conditions hold.
\begin{enumerate}
\item  If $U\in I_{J}$ and $W$ is an object of ${\calR}_n$, then for any $f\in \End_{{\calR}_n}(U\otimes W)$ we have
\[t_{U\otimes W}\left(f \right)=t_U \left( t_R(f)\right). \]
\item  If $U,V\in I$ then for any morphisms $f:V\rightarrow U $ and $g:U\rightarrow V$  in ${\calR}_n$ we have 
\[t_V(g\circ f)=t_U(f \circ g).\]
\end{enumerate}

By Kujawa \cite{Kujawa-generalized-kac-wakimoto}, thm 2.3.1, the trace on the ideal $I_{L}$, $L$ irreducible, is unique up to multiplication by an element of $k$. Given a trace on $I_{J}$, $\{t_{V} \}_{V \in I_{J}}$, $J \in \calR_n$, define the modified dimension function on objects of $I_{J}$ as the modified trace of the identity morphism:
\begin{equation*}
d_{J}\left(V \right) =t_{V}(id_{V}).
\end{equation*}

We reprove the essential part of the generalized Kac-Wakimoto conjecture: We prove that there exists a nontrivial trace on the ideal of any $i$-atypical irreducible $L$, and we deduce a formula for the resulting modified superdimension.

{\it Tensor ideals}. By \cite{Serganova-kw} any two irreducible object of atypicality $k$ generate the same tensor ideal. Therefore write $I_i$ for the tensor ideal generated by any irreducible object of atypicality $i$. Clearly $I_0 = Proj$ and $I_n = T_n$ since it contains the identity. This gives the following filtration \[Proj = I_0 \subsetneq I_1 \subsetneq \ldots I_{n-1} \subsetneq I_n = T_n\] with strict inclusions by \cite{Serganova-kw} and \cite{Kujawa-generalized-kac-wakimoto}. We use this in the following. However it is not necessary for the results about the modified superdimension. We could simply consider consider the ideal $<L>$ generated by an $i$-atypical irreducible representation instead the ideal $I_i$. 

\medskip{\it The projective case}. Denote by $\Delta_0^+$ the positive even roots and by $\Delta_1^+$ the positive odd roots for our choice of Borel algebra. The half sums of the positive even roots is denoted $\rho_0$, the half-sum of the positive odd roots by $\rho_1$ and we put $\rho = \rho_0 - \rho_1$. We define a bilinear form $(,)$ on $\mathfrak{h}^*$ as follows: We put $(\epsilon_i, \epsilon_j) = \delta_{ij}$ for $i,j \leq m$, $(\epsilon_i,\epsilon_j) = - \delta_{ij}$ for $i,j \geq m+1$ and $(\epsilon_i,\epsilon_j) = 0$ for $i \leq m$ and $j >m$. Define for any typical module the following function \[ d(L(\lambda)) = \prod_{\alpha \in \Delta_0^+} \frac{(\lambda + \rho, \alpha)}{\rho,\alpha} / \prod_{\alpha \in \Delta_1^+} (\lambda + \rho,\alpha).\] Then $d(L(\lambda)) \neq 0$ for every typical $L(\lambda)$. By \cite{Geer-Kujawa-Patureau-Mirand}, 6.2.2 for typical $L$ \[ d_J(L) = \frac{d(L)}{d(J)}.\] Since the ideal $I_0$ is independent of the choice of a particular $J$ and any ambidextrous trace is unique up to a scalar, we normalize and define the modified normalized superdimension on $I_0$ to be \[ sdim_0 (L(\lambda)) := d(L(\lambda)).\]

\medskip{\it A formula for the modified superdimension}. Applying $DS$ iteratively $i$-times to a module of atypicality $i$ we obtain the functor \[ DS^i := DS \circ \ldots \circ DS: {\calR}_n \to T_{n-i}\] which sends $M$ with $atyp(M)  = i$ to a direct sum of typical modules.  

\medskip 

We show that there exists a nontrivial trace on $I_i$ similar to \cite{Kujawa-generalized-kac-wakimoto}, but without invoking Serganovas results. Denote by $t^P$ the normalized (such that we get $sdim_0$ from above) trace on $I_0 = Proj$. Now we define for $M \in I_i$ \[ t_M (f) : = t^p_{DS^i(M)} f_{DS^i(M)}: End_{{\calR}_n} (M) \to k \] where $f_{DS^i(M)}$ is the image of $f$ under the functor $DS^i$. We claim that this defines a nontrivial trace on $I_i$: Let $M = L$ be irreducible and put \[ t_L (id_L) := t_{DS^i(L)}^p (id_{DS^i(L)}). \] Now we compute $DS^i(L)$. By the main theorem the irreducible summands in $DS(L)$ are obtained by removing one of the outer cups of each sector. Applying $DS$ $i$-times gives then the typical module in $T_{n-i}$ given by the cup diagram of $L$ with all $\vee$'s removed. Applying $DS^i$ to any other irreducible module in the same block will result in the same typical weight. Following Serganova \cite{Serganova-kw} we call this unique irreducible module the core of the block $L^{core}$. Hence $DS^i(L) = m (L)\cdot L^{core} \oplus m'(L)\cdot \Pi L^{core}$. Since the positive integers $m$ and $m'$ only depend on the nesting structure of the cup diagram $\underline{\lambda}$, we may compute them in the maximally atypical case. By a comparison with the maximal atypical case $\calR_i$-case either $m$ or $m'$ is zero. As in the maximally atypical case a parity shift happens in $DS(L(\lambda))$ if and only if $\varepsilon(\lambda) \not\equiv \varepsilon(\lambda_i) \ mod \ 2$. Hence \begin{align*} m(\lambda) = \begin{cases}  m \ \ & \varepsilon(\lambda) \equiv 0 \ mod \ 2 \\  m' & \varepsilon(\lambda) \equiv 1 \ mod \ 2. \end{cases} \end{align*} This shows that the trace $t_L$ does not vanish: Indeed \[t_L (id_L) := t_{DS^i(L)}^p (id_{DS^i(L)}) = m(\lambda) t^P_{L^{core}}(id_{L^{core}}) \neq 0\] since $t^P$ is nontrivial. 

\medskip

Using our particular choice for $sdim_0$ on $I_0 = Proj$, we define the normalized modified superdimension as \begin{align*} sdim_i (L(\lambda)) & = sdim_0 (DS^i(L)) = sdim_0 (m L^{core} \oplus m'  \Pi L^{core}) \\ & = (-1)^{\varepsilon(\lambda)} m(\lambda) sdim_0(L^{core}) \end{align*}  In particular the modified super dimension does not vanish. Consider for example the irreducible 4-fold atypical representation in $\calR_6$ with cup diagram


\medskip
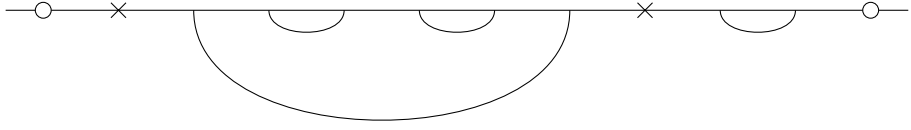
\begin{center}

  \begin{tikzpicture}
 \draw (-6.5,0) -- (5.5,0);
\foreach \x in {} 
     \draw (\x-.1, .2) -- (\x,0) -- (\x +.1, .2);
\foreach \x in {} 
     \draw (\x-.1, -.2) -- (\x,0) -- (\x +.1, -.2);
\foreach \x in {2,-5} 
     \draw (\x-.1, .1) -- (\x +.1, -.1) (\x-.1, -.1) -- (\x +.1, .1);
\foreach \x in {5} 
     \draw  node at (5,0) [fill=white,draw,circle,inner sep=0pt,minimum size=6pt]{};
\foreach \x in {-6} 
     \draw  node at (-6,0) [fill=white,draw,circle,inner sep=0pt,minimum size=6pt]{};

\draw [-,black,out=270,in=270](-4,0) to (1,0);
\draw [-,black,out=270,in=270](-3,0) to (-2,0);
\draw [-,black,out=270,in=270](-1,0) to (0,0);
\draw [-,black,out=270,in=270](3,0) to (4,0);

\end{tikzpicture}

\end{center}

We have already seen above that $m(\lambda) = 8$ in this case.  The core is given by the typical representation $L(3,-4 | 5,-5)$.

\medskip

As a consequence of our construction and the sign rule of the main theorem we get

\begin{cor} If $L$ is irreducible of atypicality $k$, then $sdim_k(L) = sdim_{k-1}(DS(L))$. If $sdim_k(L) > 0$, then all summands in $DS(L)$ have $sdim_{k-1}(L) > 0$.
\end{cor}

We can now copy the proof of proposition \ref{Leray} to get

\begin{cor}  \label{Leray-2}
For irreducible  atypical objects $L$ in $T_n$ the Leray type
spectral sequence degenerates:
$$  \fbox{$ DS_{n,n_2}(L) \ \cong \ DS_{n_1,n_2}(DS_{n,n_1}(L)) $} \ .$$
\end{cor}

\medskip\noindent



\section{Strategy of the proof}\label{sec:strategy}

\medskip
We have already proved the Main Theorem for the groundstates of each block. Recall that a groundstate is a weight with completely nested cup diagram such that all the vertices labelled $\times$ or $\circ$ are to the right of the cups. In the maximally atypical case the ground state are just the Berezin-powers. In the lower atypical cases every ground state is a Berezin-twist of a mixed tensor and we have already seen that these satisfy the main theorem in section \ref{stable0}. The proof of the general case will be a reduction to the case of groundstates.

\medskip

In the singly atypical case we just have to move the unique label $\vee$ to the left of all of the crosses and circles. We will see in section \ref{sec:loewy-length} that we can always move $\vee$'s to the left of $\circ$'s or $\times$. The proof of the general case will induct on the degree of atypicality, hence we will always assume that the theorem is proven for irreducible modules of lower atypicality. Hence for the purpose of explaining the strategy of the proof we will focus on the maximally atypical case.

\medskip{\it The modules $S^i$}.  Let us consider the following special maximally atypical case. Let $Ber\simeq [1,\ldots,1]\in \calR_n$ be the Berezin representation. Let $S^i$ denote the irreducible representation $[i,0,\ldots,0]$. Every $S^{i-1}$ occurs as the socle and cosocel of a mixed tensor denoted $\A_{S^{i+1}}$ \cite{Heidersdorf-mixed-tensors}. The Loewy structure of the modules $\A_{S^i} :=R((i),(1^i)) \in {\calR}_n$ is the following: $$ \A_{S^i} \ = \ 
(S^{i-1}, S^i \oplus S^{i-2}, S^{i-1}) \ $$

for  $i \neq n$ and $i\geq 1$ and $n\geq 2$ where we use $S^{-1} = 0$. Furthermore  
$$ \A_{S^n} \ = \ (S^{n-1}, S^n \oplus Ber^{-1} \oplus S^{n-2}, S^{n-1}) \ .$$

We saw in \ref{stable0} that for all mixed tensors $DS(R(\lambda^L,\lambda^R)) = R(\lambda^L,\lambda^R)$ holds, so we have $DS(\A_{S^i}) = \A_{S^i}$ for all $i\geq 1$.
Notice that by abuse of notation we view $S^i$ and also $\A_{S^i}$ as objects of $\calR_n$ for all $n$.

\medskip

The image $S^i \mapsto DS(S^i)$ can be computed recursively from the two exact sequences in $\calR_n$  \[ \xymatrix{ 0 \ar[r] & K_n^i \ar[r] & \A_{S^i} \ar[r]^p & S^{i-1} \ar[r] & 0 \\ 0 \ar[r] & S^{i-1} \ar[r]^j & K_n^i \ar[r] & S^i \oplus ? \oplus S^{i-2} \ar[r] & 0 } \] induced by projection $p$ onto the cosocle and the inclusion $j$ of the socle. According to the main theorem we should get for $n \geq 2$
$(Ber_{n})_x = \Pi Ber_{n-1}$ and 
\begin{enumerate}
\item $DS(S^i) = S^i\ $ for $i<  n-1$, 
\item $DS(S^{i}) = S^{i} \oplus \Pi^{n-1-i} Ber^{-1}\ $  for $i\geq n-1$.
\end{enumerate}

We proof this for $i \leq n-1$. First notice $H^-(\A_{S^i})=0$ and $H^+(\A_{S^i})=\A_{S^i}$.
Suppose $i\leq n-1$ and that $H^-(S^j)=0$, $H^+(S^j)=S^j$ already holds for $j<i$
by induction. This is justified since $S^0 = k$ equals the trival module. Then the exact hexagons give
\[ \xymatrix@-4mm{ H^+(K_n^i) \ar[r] & \A_{S^i} \ar[r] & S^{i-1} \ar[d] \\ 0 \ar[u] & 0 \ar[l] & H^-(K_n^i) \ar[l]  }\]
and
\[ \xymatrix@-4mm{ S^{i-1} \ar[r] & H^+(K_n^i) \ar[r] & H^+(S^i\oplus ?)\oplus S^{i-2} \ar[d] \\ H^-(S^i\oplus ?) \ar[u] & H^-(K_n^i) \ar[l] & 0 \ar[l] }\]
If $H^+(p)=0$, then $H^+(K_n^i)\cong \A_{S^i}$. Hence $H^+(K_n^i) \twoheadrightarrow H^+(S^i\oplus ?)\oplus S^{i-2}$ composed with the projection to $S^{i-2}$ is zero, 
since the cosocle of $H^+(K_n^i)\cong \A_{S^i}$ is $S^{i-1}$. This implies $S^{i-1}=0$, which is absurd. Hence $H^+(p)$ is surjective. Therefore $H^-(K_n^i)=0$ and $H^+(K_n^i)=K_{n-1}^i$, and in particular then $$H^+(K_n^i)=K_{n-1}^i$$ is indecomposable. Hence $K_{n-1}^i \to H^+(S^i \oplus ?)\oplus S^{i-2}$ is surjective, and $H^-(S^i)=0$.
Furthermore
$$H^+(S^i)=S^i   \quad , \quad i < n-1 $$ 
and
$$ H^+(S^i) = S^i \oplus Ber^{-1} \quad , \quad i=n-1 \ .$$ The proof for the cases $i \geq n$ is similar.

\medskip

The method described in the $S^i$-case doesn't work in general. In the general case we do not have exact analogs of the $\A_{S^i}$ - mixed tensors with the property $DS(\A) = \A$.  In section \ref{sec:loewy-length} we associate to every irreducible module three representations, the weight $L$, the \textit{auxiliary} representation $L^{aux}$ and the representation $L^{\times \circ}$ and an indecomposable rigid module $F_i (L^{\times \circ})$ of Loewy length $3$ with Loewy structure $(L, A, L)$ such that the irreducible module we started with and which we denote $L^{up}$ for reasons to be explained later is one of the composition factors of $A$. If we apply this construction to irreducible modules of the form $S^i = [i,0,\ldots,0]$ we recover the modules $\A_{S^i}$. Our aim is to use these indecomposable modules as a replacement for the modules $\A_{S^i}$.  

\medskip

In the $S^i$-case we reduced the computation of $DS(S^i)$ by means of the indecomposable modules $\A_{S^i}$ to the trivial case $DS(\one) = \one$. In the general case we will reduce the computation of $DS(L)$ by means of the indecomposable modules $F_i (L^{\times \circ})$ to the case of ground states. For that we define an order on the set of cup diagrams for a fixed block such that the completely nested cup diagrams (for which the Main Theorem holds) are the minimal elements. We prove the general case by induction on this order and will accordingly assume that the main theorem holds for all irreducible modules of lower order then a given module $L$. The key point is that for a given module $L^{up}$ we can always choose our weights $L^{aux}$ and $L^{\times \circ} = F_i(L(\lambda_{\times \circ}))$ such that all other composition factors of $F_i (L^{\times \circ})$ are of lower order then $L^{up}$. Hence the Main Theorem holds for all composition factors of $F_i (L^{\times \circ})$ except possibly $L^{up}$. This setup is similar to the $\A_{
S^i}$-case where we assumed by induction on $i$ that the Main Theorem held for all composition factors of $\A_{S^i} = (S^{i-1}, S^{i-2} + S^i, S^{i-1})$ except possibly $S^i$.

\medskip

Unlike the $\A_{S^i}$ the indecomposable modules $F_i ( L^{\times \circ})$ are not mixed tensors and hence we do not know a priori their behaviour under $DS$. However assuming that the Main Theorem holds for all composition factors except possibly $L^{up}$ we prove in section \ref{sec:loewy-length} a formula for $DS(F_i (L^{\times \circ}))$. In section \ref{sec:inductive} we show that under certain axioms on the modules $F_i (L^{\times \circ})$ and their image under $DS$ the module $DS(L^{up})$ is semisimple. These axioms are verified in section \ref{sec:moves}. Here it is very important that we can control the composition factors of  the $F_i (L^{\times \circ})$. The composition factors in the middle Loewy layer will be called $\textit{moves}$ since they can be obtained from the labelled cup diagram of $L$ by moving certain $\vee$'s in a natural way. The moves are described in detail in section \ref{sec:loewy-length}.

\medskip 

We still have to explain how the induction process works, i.e. how we relate a given irreducible module to irreducible modules with lesser number of segments respectively sectors. This is done by the so-called Algorithms I and II described first in \cite{Weissauer-gl}. As above for a given module $L^{up}$ all other composition factors of $F_i (L^{\times \circ})$ are of lower order then $L^{up}$. For $L^{up}$ with more then one segment we can choose $i$ and the representations $L^{aux}$ and $L^{\times \circ}$ in such a way that all composition factors have one segment less then $L^{up}$. We can now apply the same procedure to all the composition factors of $F_i(L^{\times \circ})$ with more then one segment - i.e. we choose for each of these (new) weights $L^{aux}$ and $L^{\times \circ}$ such that the composition factors of the (new) associated indecomposable modules have less segments then them. Iterating this we finally end up with a finite number of indecomposable modules where all composition factors have weight 
diagrams with only one segment. This procedure is called Algorithm I. In Algorithm II we decrease the number of sectors in the same way: If we have a weight with only one segment but more then one sector we can choose $i$ and the weights $L^{aux}$ and $L^{\times \circ}$ such that the composition factors of $F_i (L^{\times \circ})$ have less sectors then $L^{up}$. Applying this procedure to the composition factors of $F_i(L^{\times \circ})$ and iterating we finally relate the cup diagram of $L^{up}$ to a finite number of cup diagrams with only one sector.

\medskip

Hence after finitely many iterations we have reduced everything to irreducible modules with one segment and one sector. This sector might not be completely nested, e.g. we might end up with weights with labelled cup diagrams of the type

\medskip
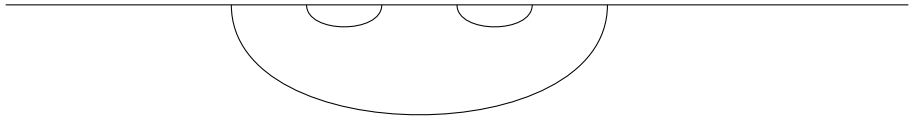
\begin{center}

  \begin{tikzpicture}
 \draw (-6,0) -- (6,0);
\foreach \x in {} 
     \draw (\x-.1, .2) -- (\x,0) -- (\x +.1, .2);
\foreach \x in {} 
     \draw (\x-.1, -.2) -- (\x,0) -- (\x +.1, -.2);
\foreach \x in {} 
     \draw (\x-.1, .1) -- (\x +.1, -.1) (\x-.1, -.1) -- (\x +.1, .1);

\draw [-,black,out=270,in=270](-3,0) to (2,0);
\draw [-,black,out=270,in=270](-2,0) to (-1,0);
\draw [-,black,out=270,in=270](0,0) to (1,0);

\end{tikzpicture}

\end{center}

In this case we can apply Algorithm II to the internal cup diagram having one segment enclosed by the outer cup. If we iterate this procedure we will finally end up in a collection of Kostant weights (i.e. weights with completely nested cup diagrams) of this block. 

\medskip

We still have to find the decomposition of the semisimple module $DS(L^{up})$ into its simple summands. Since we know the semisimplicity, we can compute $DS(L^{up})$ on the level of Grothendieck groups. Essentially we compute this in the following way: using the notation $\A = F_i(L(\lambda_{\times \circ}))$, we compute \[ d(\A) = H^+(\A) - H^-(\A) = 2d(L)  + d(A)  = 2d(L) + d(L^{up}) + d(A - L^{up}) \] in $K_0(\calR_{n-1})$ where we do not know $d(L^{up})$ and compare this to the known composition factors of $\tilde{\A} = DS(\A)$.  For this we need the so-called commutation rules for Algorithm I and Algorithm II. Using that the main theorem holds for all composition factors of $\A$ except possibly $L^{up}$ we can cancel most composition factors. The remaining factors have to be the simple factors of $DS(L^{up})$ and these factors are exactly those given by the derivative of $L^{up}$ (seen as a plot), finally proving the theorem. This is done in section \ref{sec:inductive}.

\medskip\noindent

{\it The case $[2,2,0]$}. We illustrate the above strategy with an example. In this part we ignore systematically all signs or parity shifts. The module $[2,2,0]$ has the labelled cup diagram 

\medskip
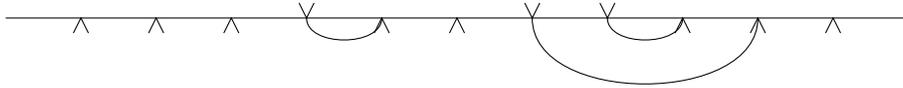
\begin{center}

  \begin{tikzpicture}
 \draw (-6,0) -- (6,0);
\foreach \x in {-2,1,2} 
     \draw (\x-.1, .2) -- (\x,0) -- (\x +.1, .2);
\foreach \x in {-5,-4,-3,-1,0,3,4,5} 
     \draw (\x-.1, -.2) -- (\x,0) -- (\x +.1, -.2);
\foreach \x in {} 
     \draw (\x-.1, .1) -- (\x +.1, -.1) (\x-.1, -.1) -- (\x +.1, .1);

\draw [-,black,out=270,in=270](-2,0) to (-1,0);
\draw [-,black,out=270,in=270](1,0) to (4,0);
\draw [-,black,out=270,in=270](2,0) to (3,0);

\end{tikzpicture}

\end{center}

hence it has two segments and two sectors. We will associate to $[2,2,0]$ an auxiliary weight $L$ and a twofold atypical weight $L^{\times \circ}$ in $T_3$ such that $[2,2,0]$ is of the form $L^{up}$ in the indecomposable module $F_i(L^{\times \circ})$. The auxiliary weight is in this case $[2,1,0]$ with labelled cup diagram 

\medskip
\begin{center}

  \begin{tikzpicture}
 \draw (-6,0) -- (6,0);
\foreach \x in {-2,0,2} 
     \draw (\x-.1, .2) -- (\x,0) -- (\x +.1, .2);
\foreach \x in {-5,-4,-3,-1,1,3,4,5} 
     \draw (\x-.1, -.2) -- (\x,0) -- (\x +.1, -.2);
\foreach \x in {} 
     \draw (\x-.1, .1) -- (\x +.1, -.1) (\x-.1, -.1) -- (\x +.1, .1);

\draw [-,black,out=270,in=270](-2,0) to (-1,0);
\draw [-,black,out=270,in=270](0,0) to (1,0);
\draw [-,black,out=270,in=270](2,0) to (3,0);

\end{tikzpicture}

\end{center}

with one segment and three sectors. The weight $\lambda_{\times \circ}$ is obtained from $[2,1,0]$ by replacing the $\vee \wedge$ at the vertices 0 and 1 by $\times \circ$

\medskip
\begin{center}

  \begin{tikzpicture}
 \draw (-6,0) -- (6,0);
\foreach \x in {-2,2} 
     \draw (\x-.1, .2) -- (\x,0) -- (\x +.1, .2);
\foreach \x in {-5,-4,-3,-1,3,4,5} 
     \draw (\x-.1, -.2) -- (\x,0) -- (\x +.1, -.2);
\foreach \x in {0} 
     \draw (\x-.1, .1) -- (\x +.1, -.1) (\x-.1, -.1) -- (\x +.1, .1);
\foreach \x in {1} 
     \draw  node at (1,0) [fill=white,draw,circle,inner sep=0pt,minimum size=6pt]{};

\draw [-,black,out=270,in=270](-2,0) to (-1,0);
\draw [-,black,out=270,in=270](2,0) to (3,0);

\end{tikzpicture}

\end{center}

The module $F_0 (L^{\times \circ})$ is $*$-selfdual of Loewy length $3$ and socle and cosocle $[2,1,0]$. It contains the module $[2,2,0]$ with multiplicity 1 in the middle Loewy layer. The rules of section \ref{sec:loewy-length} give the following composition factors (\textit{moves}) in the middle Loewy layer. In the labelled cup diagram of $[2,1,0]$ there is one internal upper sector $[2,3]$. The internal upper sector move gives the labelled cup diagram

\medskip
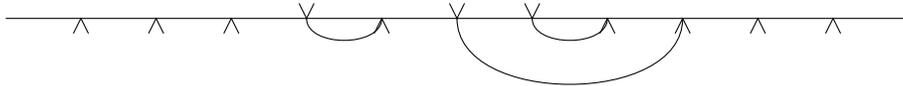
\begin{center}

  \begin{tikzpicture}
 \draw (-6,0) -- (6,0);
\foreach \x in {-2,0,1} 
     \draw (\x-.1, .2) -- (\x,0) -- (\x +.1, .2);
\foreach \x in {-5,-4,-3,-1,2,3,4,5} 
     \draw (\x-.1, -.2) -- (\x,0) -- (\x +.1, -.2);
\foreach \x in {} 
     \draw (\x-.1, .1) -- (\x +.1, -.1) (\x-.1, -.1) -- (\x +.1, .1);

\draw [-,black,out=270,in=270](-2,0) to (-1,0);
\draw [-,black,out=270,in=270](0,0) to (3,0);
\draw [-,black,out=270,in=270](1,0) to (2,0);

\end{tikzpicture}

\end{center}
hence the composition factor $[1,1,0]$. The labelled cup diagram of $[2,1,0]$ has one internal lower sector, namely the interval $[-2,-1]$. The associated internal lower sector move gives the labelled cup diagram

\medskip
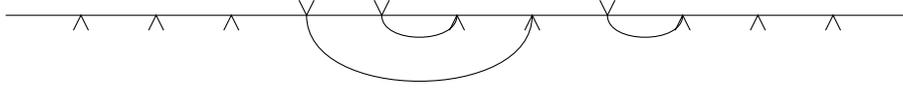
\begin{center}

  \begin{tikzpicture}
 \draw (-6,0) -- (6,0);
\foreach \x in {-2,-1,2} 
     \draw (\x-.1, .2) -- (\x,0) -- (\x +.1, .2);
\foreach \x in {-5,-4,-3,0,1,3,4,5} 
     \draw (\x-.1, -.2) -- (\x,0) -- (\x +.1, -.2);
\foreach \x in {} 
     \draw (\x-.1, .1) -- (\x +.1, -.1) (\x-.1, -.1) -- (\x +.1, .1);

\draw [-,black,out=270,in=270](-2,0) to (1,0);
\draw [-,black,out=270,in=270](-1,0) to (0,0);
\draw [-,black,out=270,in=270](2,0) to (3,0);

\end{tikzpicture}

\end{center}

The sector $[0,1]$ is unencapsulated, it is in the middle of the segment $[-2,3]$. Hence we also have the unencapsulated boundary move, i.e. we move the $\vee$ at the vertex $0$ to the vertex -3, resulting in the labelled cup diagram

\medskip
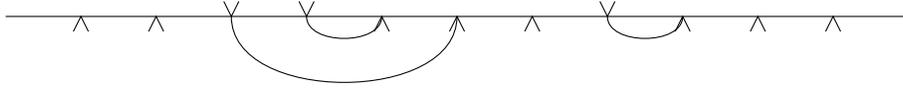
\begin{center}

  \begin{tikzpicture}
 \draw (-6,0) -- (6,0);
\foreach \x in {-2,-3,2} 
     \draw (\x-.1, .2) -- (\x,0) -- (\x +.1, .2);
\foreach \x in {-5,-4,-1,0,1,3,4,5} 
     \draw (\x-.1, -.2) -- (\x,0) -- (\x +.1, -.2);
\foreach \x in {} 
     \draw (\x-.1, .1) -- (\x +.1, -.1) (\x-.1, -.1) -- (\x +.1, .1);

\draw [-,black,out=270,in=270](-3,0) to (0,0);
\draw [-,black,out=270,in=270](-2,0) to (-1,0);
\draw [-,black,out=270,in=270](2,0) to (3,0);

\end{tikzpicture}

\end{center}

giving the composition factor $[2,-1,-1]$. The upward move of $[2,1,0]$ gives the composition factor $L^{up} = [2,2,0]$. Hence the Loewy structure of the indecomposable module $F_0(L^{aux})$ is \[ \begin{pmatrix} [2,1,0] \\ [2,-1,-1] + [1,1,0] + [2,0,0] + [2,2,0] \\ [2,1,0] \end{pmatrix}. \] We remark that all the composition factors have only one segment, hence we will not have to apply Algorithm I any more. Since the proof inducts on the degree of atypicality we know $DS(L^{\times \circ})$ and we can apply \ref{tildeA} to conclude $DS(F_i (L^{\times \circ})) = F_i (DS(L^{\times \circ})) = F_i (L_1 \oplus L_2)$ for two irreducible module obtained by applying $DS$ to $L^{\times \circ}$. By the main theorem $DS(L^{\times \circ})$ gives the modules 

\medskip
\begin{center}

  \begin{tikzpicture}
 \draw (-6,0) -- (6,0);
\foreach \x in {2} 
     \draw (\x-.1, .2) -- (\x,0) -- (\x +.1, .2);
\foreach \x in {-5,-4,-3,-2,-1,3,4,5} 
     \draw (\x-.1, -.2) -- (\x,0) -- (\x +.1, -.2);
\foreach \x in {0} 
     \draw (\x-.1, .1) -- (\x +.1, -.1) (\x-.1, -.1) -- (\x +.1, .1);
\foreach \x in {1} 
     \draw  node at (1,0) [fill=white,draw,circle,inner sep=0pt,minimum size=6pt]{};

\draw [-,black,out=270,in=270](2,0) to (3,0);

\end{tikzpicture}

\end{center}

and 

\medskip
\begin{center}

  \begin{tikzpicture}
 \draw (-6,0) -- (6,0);
\foreach \x in {-2} 
     \draw (\x-.1, .2) -- (\x,0) -- (\x +.1, .2);
\foreach \x in {-5,-4,-3,2,-1,3,4,5} 
     \draw (\x-.1, -.2) -- (\x,0) -- (\x +.1, -.2);
\foreach \x in {0} 
     \draw (\x-.1, .1) -- (\x +.1, -.1) (\x-.1, -.1) -- (\x +.1, .1);
\foreach \x in {1} 
     \draw  node at (1,0) [fill=white,draw,circle,inner sep=0pt,minimum size=6pt]{};

\draw [-,black,out=270,in=270](-2,0) to (-1,0);

\end{tikzpicture}

\end{center}

Applying $F_0$ to the first summand gives the module $\A_1$ with socle and cosocle $[2,1]$. The upward move gives the composition factor $[2,2]$. The unique internal upper sector move gives the composition factor $[1,1]$. We do not have any lower sector moves. The non-encapsulated boundary move gives the composition factor $[2,0]$. This results in the Loewy structures of $\A_1 = F_0(L_1)$ and $A_2 = F_0 (L_2)$  \[ \A_1 = \begin{pmatrix} [2,1] \\  [1,1] + [2,0] + [2,2] \\ [2,1] \end{pmatrix}, \ \  \A_2 =  \begin{pmatrix} [0,-1] \\  [1,-1] + [-1,-1] + [-2,-2] \\ [0,-1] \end{pmatrix}. \] 

The irreducible modules in the middle Loewy layers give the module $\tilde{A}$. We compare $\tilde{A}$ and $A'$ in $K_0$: Taking the derivative of $A = [2,-1,-1] + [1,1,0] + [2,0,0] + [2,2,0]$ gives \begin{align*} A' & = [2,-1] + [-2,-2] + [1,-1] + [1,1] + [2,0] \\ & + [-1,-1] + [2,-1] + [2,2] \end{align*} with the module $[2,-1] = L^{aux}$ appearing twice. The computation above of $\A_1$ and $\A_2$ gives \[ \tilde{A} =  [-2,-2] + [1,-1] + [1,1] + [2,0]  + [-1,-1]  + [2,2].\] This shows the following commutation rule in this example \[ A' = \tilde{A} + 2 (-1)^{i+n} L^{aux} \ \ \ in \ \ \ K_0(\calR_{n-1}).\] We remark that the composition factors $[2,0]$ in $\A_1$ and $[-1,-1]$ are detecting objects in the sense of section \ref{sec:moves}.

\medskip

We will prove in section \ref{sec:moves} that the properties of the modules $\A, \A_1$ and $\A_2$ imply that $DS(L^{up})$ is semisimple. Hence we can compute $DS(L^{up})$ by looking at $K_0$.  

\medskip

In Algorithm II we reduce everything to a single sector. Take one of the composition factors of $F_0 (L^{\times \circ})$ with more then one sector, eg. $[2,1,0]$ with one segment and three sectors. The associated auxiliary weight is in this case the weight $[2,0,0]$ with the twofold atypical weight $L^{\times \circ}$ given by the labelled cup diagram

\medskip
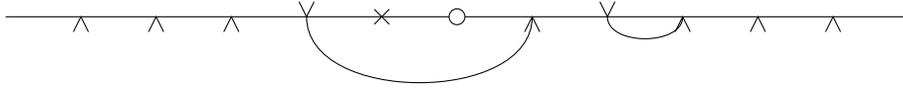
\begin{center}

  \begin{tikzpicture}
 \draw (-6,0) -- (6,0);
\foreach \x in {2,-2} 
     \draw (\x-.1, .2) -- (\x,0) -- (\x +.1, .2);
\foreach \x in {-5,-4,-3,1,3,4,5} 
     \draw (\x-.1, -.2) -- (\x,0) -- (\x +.1, -.2);
\foreach \x in {-1} 
     \draw (\x-.1, .1) -- (\x +.1, -.1) (\x-.1, -.1) -- (\x +.1, .1);
\foreach \x in {0} 
     \draw  node at (0,0) [fill=white,draw,circle,inner sep=0pt,minimum size=6pt]{};

\draw [-,black,out=270,in=270](-2,0) to (1,0);
\draw [-,black,out=270,in=270](2,0) to (3,0);

\end{tikzpicture}

\end{center} 

The module $F_{-1}(L^{\times \circ})$ has socle and cosocle $[2,0,0]$ and the followowing modules in the middle Loewy layer: The upward move gives $[2,1,0]$ and the upper sector move of the upper sector $[2,3]$ gives the weight $[0,0,0]$.
 
 There are no non-encapsulated boundary moves and no internal lower sector moves, hence we get the Loewy structure \[ F_{-1}(L^{\times \circ}) = \begin{pmatrix} [2,0,0] \\ [0,0,0] + [2,1,0] \\ [2,0,0] \end{pmatrix}.\] We compute $DS(F_{-1}(L^{\times \circ}))$ (using $DS(F_i (L^{\times \circ})) = F_i (DS(L^{\times \circ})) \ )$ (lemma \ref{tildeA}). By the main theorem $DS(L^{\times \circ})$ splits into two direct summands.
  
\medskip

Applying $F_{-1}$ to the first and second summand gives the indecomposable modules \[ \A_1 = \begin{pmatrix} [2,0] \\ [2,1] + [2,-1] + [0,0] \\ [2,0] \end{pmatrix}, \ \  \A_2 = \begin{pmatrix} [-1,-1] \\ [0,-1] \\ [-1,-1] \end{pmatrix} \] We remark that all the factors in the middle Loewy layers are detecting objects in the sense of section \ref{sec:inductive}. As shown in section \ref{sec:inductive} these properties already imply that $DS([2,1,0])$ is semisimple. To compute it we need the commutation rules for Algorithm II, i.e. we compare the derivative $A'$ of the middle Loewy layer of $F_{-1}(L^{\times \circ})$ with the modules $\tilde{A} = A_1 + A_2$ in the middle Loewy layers of $\A_1$ and $\A_2$. In both cases we get $[2,1] + [2,-1] + [0,0] + [0,-1]$, hence the commutation rule \[ \tilde{A} = A'.\] The general case is proven in lemma \ref{comII}.

\medskip\noindent





\section{Modules of Loewy length 3}\label{sec:loewy-length}

\medskip\noindent

As described in section \ref{sec:strategy} we reduce the main theorem to the case of ground states by means of translation functors $F_i(\ldots)$. In this section we describe the Loewy layers and composition factors of the objects $F_i(L_{\times \circ})$ and study their behaviour under $DS$.

\medskip

{\it Khovanov algebras}. We review some facts from the articles by Brundan and Stroppel  \cite{Brundan-Stroppel-1}, \cite{Brundan-Stroppel-2}, \cite{Brundan-Stroppel-4},  \cite{Brundan-Stroppel-5}. We denote the Khovanov-algebra of \cite{Brundan-Stroppel-4} associated to $Gl(m|n)$ by $K(m,n)$. These algebras are naturally graded. For $K(m,n)$ we have a set of weights or weight diagrams which parametrise the irreducible modules (up to a grading shift). This set of weights is again denoted $X^+$. For each weight $\lambda \in X^+$ we have the irreducible module $L(\lambda)$, the indecomposable projective module $P(\lambda)$ with top $L(\lambda)$ and the standard or cell module $V(\lambda)$. If we forget the grading structure on the $K(m,n)$-modules, the main result of \cite{Brundan-Stroppel-4} is:

\begin{thm} There is an equivalence of categories $E$ from $\calR_{m|n}$ to the category of finite-dimensional left-$K(m,n)$-modules such that $EL(\lambda) = L(\lambda)$, $EP(\lambda) = P(\lambda)$ and $EK(\lambda) = V(\lambda)$ for $\lambda \in X^+$.
\end{thm}

$E$ is a Morita equivalence, hence $E$ will preserve the Loewy structure of indecomposable modules. This will enable us to study questions regarding extensions or Loewy structures in the category of Khovanov modules. We will use freely the terminology of \cite{Brundan-Stroppel-1}, \cite{Brundan-Stroppel-2}, \cite{Brundan-Stroppel-4},  \cite{Brundan-Stroppel-5}. The notion of cups, caps, cup and cap diagrams are introduced in \cite{Brundan-Stroppel-1}. For the notion of \textit{matching} between a cup and a cap diagram see \cite{Brundan-Stroppel-2}, section 2. For the notion of $\Gamma$-\textit{admissible} see \cite{Brundan-Stroppel-4}, section 2. 

\medskip

Let $\lambda$ in ${\calR}_n$ be any atypical weight with a $\vee\wedge$-pair in its weight diagram, i.e. such that there exists an index $i$ labelled by $\vee$ and the index $i+1$ is labelled by $\wedge$. Fix such an index $i$ and replace $(\vee \wedge)$ by the labelling $(\times,\circ)$. This defines a new weight $\lambda_{\times \circ}$ of atypicality $atyp(\lambda)-1$. We denote by $F_i$, $i \in \Z$, the endofunctor from \cite{Brundan-Stroppel-4}, (2.13). The functor $F_i$ has an avatar $F_i$ on the side of Khovanov-modules. This projective functor $F_i$ is defined by $F_i := \bigoplus K_{(\Gamma - \alpha_i)\Gamma}^{t_i(\Gamma)} \otimes_K -$, see \cite{Brundan-Stroppel-4}, (2.3), for summation rules and also \cite{Brundan-Stroppel-2}, (4.1). Since by loc. cit. lemma 2.4, $F_i L(\lambda_{\times \circ})$ is indecomposable, $F_i L(\lambda_{\times \circ}) = K_{(\Gamma - \alpha_i)\Gamma}^{t_i(\Gamma)} \otimes_K -$ for one specific $i$-admissible $\Gamma$

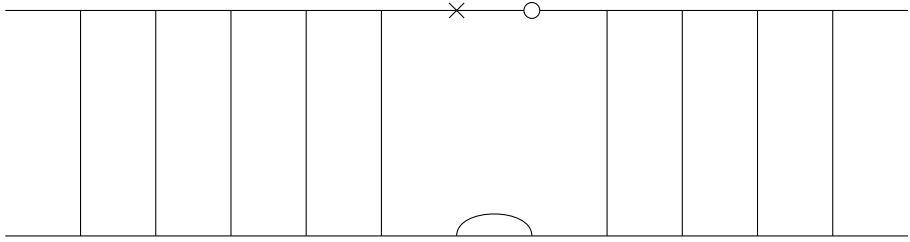
\medskip

  \begin{tikzpicture}
 \draw (-6,0) -- (6,0);
\foreach \x in {} 
     \draw (\x-.1, .2) -- (\x,0) -- (\x +.1, .2);
\foreach \x in {} 
     \draw (\x-.1, -.2) -- (\x,0) -- (\x +.1, -.2);
\foreach \x in {0} 
     \draw (\x-.1, .1) -- (\x +.1, -.1) (\x-.1, -.1) -- (\x +.1, .1);
     \draw  node at (1,0) [fill=white,draw,circle,inner sep=0pt,minimum size=6pt]{};

\begin{scope} [yshift = -3 cm]

 \draw (-6,0) -- (6,0);
\foreach \x in {} 
     \draw (\x-.1, .2) -- (\x,0) -- (\x +.1, .2);
\foreach \x in {} 
     \draw (\x-.1, -.2) -- (\x,0) -- (\x +.1, -.2);
\foreach \x in {} 
     \draw (\x-.1, .1) -- (\x +.1, -.1) (\x-.1, -.1) -- (\x +.1, .1);
\end{scope}

\draw [-,black,out=270,in=90](-5,0) to (-5,-3);
\draw [-,black,out=270,in=90](-4,0) to (-4,-3);
\draw [-,black,out=270,in=90](-3,0) to (-3,-3);
\draw [-,black,out=270,in=90](-2,0) to (-2,-3);
\draw [-,black,out=270,in=90](-1,0) to (-1,-3);
\draw [-,black,out=270,in=90](3,0) to (3,-3);
\draw [-,black,out=270,in=90](4,0) to (4,-3);
\draw [-,black,out=270,in=90](2,0) to (2,-3);
\draw [-,black,out=270,in=90](5,0) to (5,-3);

\draw [-,black,out=90, in=90](0,-3) to (1,-3);

\end{tikzpicture}
\medskip

Here the matching between $(\Gamma - \alpha_i)$ and $\Gamma$ is given by the diagram above and the rule that all other vertices, except those labelled by $\times$ or $\circ$, are connected by a vertical identity line segment. We want to determine its composition factors and Loewy layers. For that one considers the modules $F_i L(\lambda_{\times \circ})$ as modules in the \textit{graded} category of $K=K(n,n)$-modules where $K(n,n)$ is the Khovanov algebra from \cite{Brundan-Stroppel-4}. We recall some facts from \cite{Brundan-Stroppel-1} and \cite{Brundan-Stroppel-4}, see also \cite{Heidersdorf-mixed-tensors}. 

\medskip
Let $\Lambda$ be any block in the category of graded $K$-modules. For  a graded $K$-module $M = \bigoplus_{j \in \Z} M_j$, we write $M\langle j \rangle$ for the same module with the new grading $M\langle j \rangle_i := M_{i-j}$. Then the modules $\{ L(\lambda)\langle j \rangle \ | \ \lambda \in \Lambda, \ j \in \Z \}$ give a complete set of isomorphism classes of irreducible graded $K_{\Lambda}$-modules.  For the full subcategory $Rep(K_{\Lambda})$ of $Mod_{lf}(K_{\Lambda})$ consisting of finite-dimensional modules, the Grothendieck group is the free $\Z$-module with basis given by the $L(\lambda) \langle j \rangle$. Viewing it  as a $\Z[q,q^{-1}]$-module, so that by definition $q^j [M] : = [M\langle j \rangle]$ holds, $K_0(Rep(K_{\Lambda}))$ is the free $\Z[q,q^{-1}]$-module with basis $\{ L(\lambda) \ | \ \lambda \in \Lambda \}$. We refer to \cite{Brundan-Stroppel-2}, section 2, for the definition of the functors $G_{\Lambda \Gamma}^t$. For terminology used in the statement of the next theorem see loc.cit or section \ref{sec:n^n}.  We quote from \cite{Brundan-Stroppel-2}, thm 4.11

\begin{thm} \label{compos} Let $t$ be a proper $\Lambda \Gamma$-matching and $\gamma \in \Gamma$. Then in the graded Grothendieck group \[ [ G_{\Lambda \Gamma}^t L(\gamma) ] = \sum_{\mu} (q + q^{-1})^{n_{\mu}} [L(\mu)] \] where $n_{\mu}$ denotes the number of lower circles in $\underline{\mu}t$ and the sum is over all $\mu \in \Lambda$ such that a) $\underline{\gamma}$ is the lower reduction of $\underline{\mu}t$ and b) the rays of each lower line in $\underline{\mu}\mu t$ are oriented so that exactly one is $\vee$ and one is $\wedge$.
\end{thm}

Up to a grading shift  by $-caps(t)$ we have $F_i L(\lambda_{\times \circ}) = G_{(\Gamma-\alpha_i) \Gamma}^t L(\gamma)$ for some $\gamma$ and we may apply the theorem above to compute their Loewy structure. By \cite{Brundan-Stroppel-4}, lemma 2.4.v, $F_i L(\lambda_{\times \circ})$ is indecomposable with irreducible socle and head isomorphic to $L(\lambda)$. 

\begin{prop} $F_i L(\lambda_{\times \circ})$ has a three step Loewy filtration \[ F_i L(\lambda_{\times \circ}) = \begin{pmatrix} L(\lambda) \\ F \\ L(\lambda) \end{pmatrix} \] where all irreducible constituents in (the semisimple) module $F$ occur with multiplicity 1.
\end{prop}    

\medskip{\it Proof}. Let $F(j)$ be the submodule of $F_i L(\lambda_{\times \circ})$ spanned by all graded pieces of degree $\geq j$. Let $k$ be large enough so that all constituents of $F_i L(\lambda_{\circ \times})$  have degree $\geq -k$ and $\leq k$.  Then \[ F = F(- k) \supset F(- k +1) \supset \ldots \supset  F(k)\] with successive semisimple quotients $F(j)/F(j+1)$ in degree $j$. In our case we take $k=1$, since the irreducible socle and top $L(\lambda) = L(\lambda_{\vee\wedge})$ satisfies $n_{\lambda} = 1$. Then all other composition factors $L(\mu)$ necessarily satisfy $n_{\mu} = 0$ (we ignore the shift by $\langle -caps(t)\rangle $ here). The grading filtration thus gives our three step Loewy filtration. The statement about the multiplicity follows since the multiplicity of $L(\mu)$ in $F$ is given by $2^{n_{\mu}}$. The Loewy filtration of $F_i L(\lambda_{ \times \circ})$ is preserved by the Morita equivalence $E^{-1}$ of $K(n,n)\text{\it -mod}$ with ${\calR}_n$. \qed

\begin{lem} $F_i L(\lambda_{\times \circ})$ is $*$-invariant.
\end{lem}

\medskip{\it Proof}. Since $X \otimes L_{\times,\circ}$ is $^*$-invariant, $^*$ permutes its indecomposable summands. The indecomposable summands are either irreducible or are of the form $F_
j  L(\lambda_{\times \circ})$ for some $j$ with labeling $(\times,\circ)$ at position $(j,j+1)$. Since 
$^*$ preserves irreducible modules, the indecomposable summands corresponding to the $(\times,\circ)$-pairs in $\lambda_{x\circ}$ are permuted amongst themselves. Since $^*$ preserves irreducible modules $[M^*] = [M]$ in $K_0$. However all the non-irreducible $F_j (L(\lambda_{\times \circ}))$ lie in different blocks for $j \neq j'$ by the rules of \cite{Brundan-Stroppel-4}, lemma 2.4.\qed


\medskip
{\it Composition factors}. We describe the composition factors of $F_i (L^{\times \circ})$. We can restrict ourselves to the maximally atypical block (i.e. we can ignore $\times$'s and $\circ$'s).

\medskip

Let $\lambda$ be $i$-fold atypical. Since $F_i(L(\lambda_{\times \circ}))$ is indecomposable, any highest weight of a composition factor $\mu$ has the same positioning of the $n-i$ crosses and $n-i$ circles as $\lambda$. In particular it has the same positioning of the circles and crosses as $\lambda_{\times \circ}$ except at the position $(i,i+1)$. Let $F_i (L(\lambda_{\times \circ}))$ be given by a matching $t$ as follows

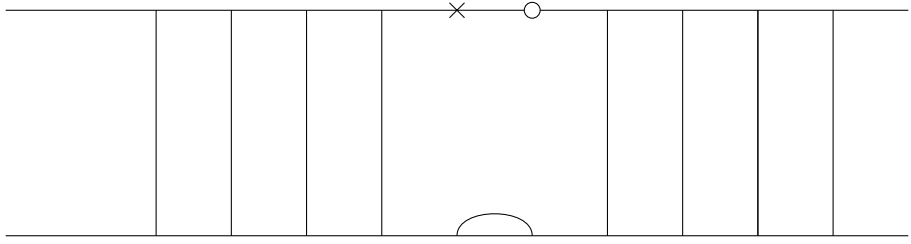
\begin{center}
\bigskip

  \begin{tikzpicture}
 \draw (-6,0) -- (6,0);
\foreach \x in {} 
     \draw (\x-.1, .2) -- (\x,0) -- (\x +.1, .2);
\foreach \x in {} 
     \draw (\x-.1, -.2) -- (\x,0) -- (\x +.1, -.2);
\foreach \x in {0} 
     \draw (\x-.1, .1) -- (\x +.1, -.1) (\x-.1, -.1) -- (\x +.1, .1);
     \draw  node at (1,0) [fill=white,draw,circle,inner sep=0pt,minimum size=6pt]{};

\begin{scope} [yshift = -3 cm]

 \draw (-6,0) -- (6,0);
\foreach \x in {} 
     \draw (\x-.1, .2) -- (\x,0) -- (\x +.1, .2);
\foreach \x in {} 
     \draw (\x-.1, -.2) -- (\x,0) -- (\x +.1, -.2);
\foreach \x in {} 
     \draw (\x-.1, .1) -- (\x +.1, -.1) (\x-.1, -.1) -- (\x +.1, .1);
\end{scope}

\draw [-,black,out=270,in=90](-4,0) to (-4,-3);
\draw [-,black,out=270,in=90](-3,0) to (-3,-3);
\draw [-,black,out=270,in=90](-2,0) to (-2,-3);
\draw [-,black,out=270,in=90](-1,0) to (-1,-3);
\draw [-,black,out=270,in=90](4,0) to (4,-3);
\draw [-,black,out=270,in=90](2,0) to (2,-3);
\draw [-,black,out=270,in=90](3,0) to (3,-3);
\draw [-,black,out=270,in=90](4,0) to (4,-3);
\draw [-,black,out=270,in=90](5,0) to (5,-3);

\draw [-,black,out=90, in=90](0,-3) to (1,-3);

\end{tikzpicture}
\end{center}
\bigskip

The crosses and the circles are now fixed. Since the composition factors depend only on the nesting structure and the matching $t$ as in theorem \ref{compos} we can fix them and assume that we are in the maximally atypical block of $Gl(i|i)$. In this case the composition factors can be determined from the segment and sector structure of $\lambda$ as in \cite{Weissauer-gl}. For symbols $x,y \in \{\circ,\wedge,\vee,\times\}$
we write $\lambda_{xy}$ for the diagram obtained from $\lambda$ with the $i$th and $(i+1)$th vertices relabeled by $x$ and $y$, respectively.

\begin{itemize}
\item {\it Socle and cosocle}. They are defined by $L(\mu)$ for $\mu=\lambda_{\vee\wedge}$.
\item {\it The upward move}. It corresponds to the weight $\mu = \lambda_{\wedge\vee}$ which
is obtained from $\lambda_{\vee\wedge}$ by switching $\vee$ and
$\wedge$ at the places $i$ and $i+1$. It is of type
$\lambda_{\wedge\vee}$.
\item {\it The nonencapsulated boundary move}. It only occurs in the nonencapsulated
case. It moves the $\vee$ in $\lambda_{\vee\wedge}$ from position
$i$ to the left boundary position $a$. The resulting weight $\mu$ is of type
$\lambda_{\wedge\wedge}$.
\item {\it The internal upper sector moves}. For every internal upper sector
$[a_j,b_j]$ (i.e. to the right of $[i,i+1]$) there is a summand
whose weight is obtained from $\lambda_{\vee\wedge}$ by moving the
label $\vee$ at $a_j$ to the position $i+1$. These moves define new weights $\mu$ of
type $\lambda_{\vee\vee}$.
\item {\it The internal lower sector moves}. For every internal lower sector
$[a_j,b_j]$ (i.e. to the left of $[i,i+1]$) there is a summand
whose weight is obtained from $\lambda_{\vee\wedge}$ by moving the
label $\vee$ from the position $i$ to the position $b_j$. These
moves define new weights $\mu$ of type $\lambda_{\wedge\wedge}$.
\end{itemize}

For examples see \cite{Weissauer-gl} or section \ref{sec:strategy}. It follows from the maximal atypical case and the definition of our sign $\varepsilon(L)$ that we have $F_i L(\lambda_{\circ \times}) = (L,F,L)$ with $L \in \calR_n(\pm \varepsilon)$ and $F \in \calR_n(\mp \varepsilon)$. For the following lemma see also \cite{Serganova-kw}, thm. 2.1 and cor. 4.4.

 \begin{lem} \label{tildeA}
Suppose theorem \ref{mainthm} holds for the irreducible representation
$L^{\times \circ}=L(\lambda_{\times \circ})$ in the block $\Gamma$ of $\calR_n$. Suppose 
$i\in \mathbb Z$ is $\Gamma$-admissible in the sense of \cite{Brundan-Stroppel-4}, p.6.
Then for the special projective functor $F_i$ the following holds:
$$  \fbox{$  DS( F_i L_{\times \circ}) = F_i DS(L_{\times \circ}) $} \ .$$
\end{lem}     

\medskip{\it Proof}. Given $(V,\rho)$ in $\calR_n$ the Casimir $C_n$ of $\calR_n$ 
restricts on $DS(V,\rho)$ to the Casimir $C_{n-1}$ of $\calR_{n-1}$ by lemma \ref{Cas}.
On irreducible representations $V$ the Casimir acts by a scalar $c(V)$.
Given representations $V_1,V_2$ in $\calR_n$, such that $C_n$ acts by $c(V_i) \cdot id_{V_i}$
on $V_i$,  and $v \in V_1\otimes V_2$, then
$C_n(v) = (c(V_1)+ c(V_2))\cdot v + 2\Omega_n(v)$ for 
$\Omega_n= \sum_{r,s=1}^n (-1)^{\overline s} e_{r,s} \otimes
e_{s,r} \!\in\! \g_n \otimes \g_n$.

Note $F_i(V) = pr_{\Gamma - \alpha_i} \circ ( V \otimes X_{st})
\circ pr_\Gamma$, so $F_i  L(\lambda_{\times \circ}) =  pr_{\Gamma - \alpha_i} ( L(\lambda_{\times \circ}) \otimes X_{st})$. By \cite{Brundan-Stroppel-4}, lemma 2.10, this is  also
the generalized $i$-eigenspace of $\Omega_n$ on $L(\lambda_{\times \circ}) \otimes X_{st}$.
Put $c=c(L(\lambda_{\times \circ})) + c(X_{st}) +2i$. Then $F_i  L(\lambda_{\times \circ}) $ is
the generalized $c$-eigenspace of $C_n$ on $L(\lambda_{\times \circ}) \otimes X_{st}$.
Hence $DS( F_i L(\lambda_{\times \circ}))$   is the generalized $c$-eigenspace of $C_{n-1}$
 on $DS( L(\lambda_{\times \circ}) \otimes X_{st}) = DS( L(\lambda_{\times \circ})) \otimes DS(X_{st}) =
 DS( L(\lambda_{\times \circ})) \otimes X_{st,n-1}$. Observe that $c(DS(V_1))+ c(DS(V_2)) = c(V_1)+ c(V_2)$, since $C_n$ induces $C_{n-1}$ on $DS(V_i)$.

 By the main theorem \ref{mainthm} (using induction over degree of atypicity)
 $DS(L(\lambda_{\times \circ}))$  is in a unique block $\overline\Gamma$. 
 So $F_i DS(L(\lambda_{\times \circ})) = pr_{\overline\Gamma - \alpha_i} \circ (? \otimes X_{st,n-1}) \circ pr_{\overline\Gamma} DS(L(\lambda_{\times \circ})) = pr_{\overline\Gamma - \alpha_i}  (DS(L(\lambda_{\times \circ})) \otimes X_{st,n-1})$, and again by \cite{Brundan-Stroppel-4}, lemma 2.10, this is the generalized $c$-eigenspace of the Casimir $C_{n-1}$ on $DS(L(\lambda_{\times \circ})) \otimes X_{st,n-1}$. Thus $DS( F_i L(\lambda_{\times \circ})) \cong F_i DS(L(\lambda_{\times \circ}))$.
 \qed


\medskip\noindent

{\it Weights, sectors, segments}. Let $L(\lambda)$ be $i$-atypical in a block $\Gamma$. Let $X^+_{\Gamma}$ denote the set of weights in $\Gamma$. Then we define a map \[ \phi = \phi_{\Gamma}: X^+_{\Gamma} \to  \{ \text{plots of rank } i \} \] by sending $\lambda$ to the plot of the weight of the irreducible representation $\phi_n^i(L(\lambda))$. Then $\phi_{\Gamma}$ is a bijection. Each plot has defining segments and sectors, and by transfer with $\phi_{\Gamma}$ this defines the segments and sectors of a given weight diagram in $X_{\Gamma}^+$.

\medskip{\it Shifting $\times$ and $\circ$}.  
We now quote from \cite{Brundan-Stroppel-4}, lemma 2.4

\begin{lem}\label{translation}
Let $\lambda \in X^+(n)$ and $i \in \mathbb Z$.
For symbols $x,y \in \{\circ,\wedge,\vee,\times\}$
we write $\lambda_{xy}$ for the diagram obtained from $\lambda$
with the $i$th and $(i+1)$th vertices
relabeled by $x$ and $y$, respectively.
\begin{itemize}
\item[\rm(i)]
If $\lambda = \lambda_{{\vee}\times}$ then $E_i L(\lambda) \cong L(\lambda_{\times {\vee}})$. If $\lambda = \lambda_{\times \vee}$ then $F_i L(\lambda) \cong L(\lambda_{{\vee} \times})$.
\item[\rm(ii)]
If $\lambda = \lambda_{{\wedge}\times}$ then $E_i L(\lambda) \cong L(\lambda_{\times {\wedge}})$. If $\lambda = \lambda_{\times \wedge}$ then $F_i L(\lambda) \cong L(\lambda_{\wedge \times})$.
\item[\rm(iii)]
If $\lambda = \lambda_{{\vee}\circ}$ then $F_i L(\lambda) \cong L(\lambda_{\circ {\vee}})$. If $\lambda = \lambda_{\circ \vee}$ then $E_i L(\lambda) \cong L(\lambda_{{\vee} \circ})$.
\item[\rm(iv)]
If $\lambda = \lambda_{{\wedge}\circ}$ then $F_i L(\lambda) \cong L(\lambda_{\circ {\wedge}})$. If $\lambda = \lambda_{\circ \wedge}$ then $E_i L(\lambda) \cong L(\lambda_{\wedge \circ})$.
\item[\rm(v)]
If $\lambda = \lambda_{{\times}\circ}$ then:
$F_i L(\lambda)$ has irreducible
socle and head both isomorphic to $L(\lambda_{\vee\wedge})$, and all other composition
factors are of the form $L(\mu)$ for $\mu \in\lambda$
such that
$\mu = \mu_{\vee \wedge}$,
$\mu = \mu_{\wedge \vee}$ or
$\mu = \mu_{\wedge \vee}$.
Likewise for $\lambda = \lambda_{\circ \times}$ and $E_i L(\lambda)$.
\item[\rm(vi)]
If $\lambda = \lambda_{{\vee\wedge}}$ then $F_i L(\lambda) \cong L(\lambda_{\circ{\times}})$.
\end{itemize}
\end{lem}

\medskip
For a pair of neighbouring vertices $(i,i+1)$ in the weight diagram of $\lambda = \lambda_{\vee\times }$, labelled by $( \vee \times)$,  we get $$E_i L(\lambda_{\vee\times}) = L(\lambda_{ \times\vee})$$ from \ref{translation}.1. In other words, the functor replaces the irreducible representation
of weight $\lambda_{\times \vee}$ by the irreducible representation of weight $\lambda_{\vee \times}$, which has the same weight diagram as $\lambda_{\times \vee}$, except that
the positions of $\times$ and $\vee$ are interchanged. Note that
$$  \phi(\lambda_{\vee\times }) \ =\  \phi(\lambda_{\times\vee }) \ ,$$
but $L=L(\lambda_{\vee\times})$ and $L^{up}= L(\lambda_{\times\vee})$
lie in different blocks.

\begin{lem}
Suppose for the representation
$L=L(\lambda_{\vee\times})$ in $\calR_n^{i}$ the assertion of theorem \ref{mainthm} holds. Then it also holds for
the representation $L^{up}= L(\lambda_{\times\vee})$.
\end{lem}

\medskip{\it Proof}. By assumption we have a commutative diagram
$$    \xymatrix{ L \ar[rr]^p\ar[dd]   & &  \lambda \ar[dd] \cr
& & \cr 
 DS(L) \ar[rr]^p   & &  \lambda'  \cr}
$$
We have to show that we have the same diagram for $L^{up}$ instead
of $L$. Let $S_\nu$ denote the sectors of the plot $\lambda=\phi(\lambda_{\vee\times})$
and let $S_j$ denote the sector containing the integer $p(i)$. 
Then $DS(L)$ is a direct sum of irreducible representations $L_\nu$, whose
sector structure either is obtained by replacing one of the sectors $S_\nu, \nu\neq j$
by $\partial S_\nu$, and there is the unique irreducible summand $L_j$
whose sector structure either is obtained by replacing the sectors $S_j$
by $\partial S_j$. We would like to show that $DS(L^{up})$ can be similarly
described in terms of the sector structure of $L^{up}$. 
The sectors of $L^{up}$
literally coincide with the $S_\nu$ for $\nu\neq j$, and for $\nu=j$ the remaining
sector of $L^{up}$ is obtained from the sector $S_j$ by transposing the positions
at the labels $i,i+1$ (within this sector). Hence to show our claim, it remains to
show that $DS(L^{up})$ is isomorphic to a direct sum of irreducible representations
$L^{up}_\nu$ with the sector structures such that $L^{up}_\nu$ is obtained
from $L_\nu$ by applying the functor $E_i$ (i.e. replacing the positions of $\vee$ and $\times$ at the labels $i,i+1$). Indeed, the derivative $\partial$ for sectors commutes with
the interchange of labels at $i,i+1$ in our situation (the sign rule is obviously preserved).
Hence it remains to show  
$$E_i(DS(L(\lambda_{\vee\times}))) = DS(E_i (L(\lambda_{\vee\times}))) \ .$$ 
But this assertion follows by an argument similarly to the one
used for the proof of lemma \ref{tildeA}. \qed

\medskip
 Likewise by lemma \ref{translation} one can show  
 
\begin{lem}\label{shift-1}
Suppose for the representation
$L=L(\lambda_{\vee\circ})$ in $\calR_n^{i}$ the assertion of theorem \ref{mainthm} holds. Then it also holds for
the representation $L^{up}= L(\lambda_{\circ\vee})$.
\end{lem}
 
 \begin{lem}\label{shift-2}
Suppose the main theorem holds for the representation
$L=L(\lambda_{\wedge\times})$ in $\calR_n^{i}$. Then it also holds for
the representation $L^{up}= L(\lambda_{\times\wedge})$.
\end{lem}

\begin{lem}\label{shift-3}
Suppose the main theorem holds for the representation
$L=L(\lambda_{\wedge\circ})$ in $\calR_n^{i}$. Then it also holds for
the representation $L^{up}= L(\lambda_{\circ\wedge})$.
\end{lem}



\section{Inductive Control over $DS$}\label{sec:inductive}

\medskip\noindent

We prove now the main theorem under the assumption that there exist objects $\A$ with certain nice properties. Under these assumptions we give an inductive proof of theorem \ref{mainthm} using the proposition \ref{3} below. We verify in section \ref{sec:moves} that certain objects $F_i(L_{\times \circ})$ verify these conditions.

\medskip

First recall that for $\varepsilon \in \{ \pm 1\}$  the full abelian subcategories
$\calR_n(\varepsilon)$ of $\calR_n$ consist of all objects whose irreducible
constituents $X$ have sign $\varepsilon(X)=\varepsilon$. We quote from section \ref{sec:main}
the following

\begin{prop} \label{ext-0}
The categories $\calR_n(\varepsilon)$ are semisimple abelian categories.
\end{prop}

\medskip
{\it Definition}. An object $M$ in $\calR_n$ is called {\it semi-pure} (of sign $\varepsilon$), if its socle is in the category $\calR_n(\varepsilon)$. Every subobject of a semi-pure object is semi-pure.
For semi-pure objects $M$ the second layer of the lower Loewy series (i.e. the socle of
$M/socle(M)$) is in $\calR_n(-\varepsilon)$ by the last proposition. Hence by induction,  
the $i$-th layer of the lower Loewy filtration is in $\calR_n((-1)^{i-1}\varepsilon)$. Hence all
layers of the lower Loewy filtration are semi-pure. The last layer $top(M)$ of the lower Loewy series is semisimple. Since
$cosocle(M) \cong cosocle(M)^* \cong socle(M^*)$ 
this easily implies

\begin{lem} \label{purity}
For semi-pure $*$-selfdual indecomposable objects $M$ in $\calR_n$ of Loewy length $\leq 3$ the lower and the upper Loewy series coincide.
\end{lem}

\medskip\noindent
We now formulate certain {\it axioms} for an  
object $\A$ of $\calR_n$. Along with the results of section \ref{sec:loewy-length} we will see in section \ref{sec:moves} that the translation functors $F_i(L^{\times \circ})$ verify these conditions.
\begin{enumerate}
\item $\A \in \calR_n$ is indecomposable with Loewy structure $(L,A,L)$.
\item  $\A$ is $*$-selfdual.
\item $L\in \calR_n(\varepsilon) $ is irreducible and satisfies theorem \ref{mainthm} with $A\in \calR_n(-\varepsilon)$.
\item  $\tilde \A:=DS(\A) = \tilde \A^+ \oplus \Pi(\tilde \A^-)$ 
is the direct sum of $\tilde \A^+ :=H^+(\A)$ and $\tilde \A^- = H^-(\A)$
such that $\tilde \A^+= \bigoplus_{\nu} \tilde \A_{\nu}$ and $\tilde \A^-= \bigoplus_{\overline\nu} \tilde \A_{\overline\nu}$ with indecomposable objects 
$\tilde\A_\nu \in \calR_{n-1}$ of Loewy structure $\tilde\A_\nu = (\tilde L_\nu,\tilde A_\nu,\tilde L_\nu)$ resp. $\tilde\A_{\overline\nu}, \in \calR_{n-1}$ of Loewy structure $\tilde\A_{\overline\nu} = (\tilde L_{\overline\nu},\tilde A_{\overline\nu},\tilde L_{\overline\nu})$.
\item All $\tilde L_\nu$ and $\tilde L_{\overline\nu}$ are irreducible so that  $\tilde L_\nu \in \calR_{n-1}(\varepsilon)$ and 
$\tilde L_{\overline\nu} \in \calR_{n-1}(-\varepsilon)$; furthermore $\tilde A_\nu \in \calR_{n-1}(-\varepsilon)$ and $\tilde A_{\overline\nu} \in \calR_{n-1}(\varepsilon)$ 
\item For each $\mu=\nu$ (resp. $\mu=\overline\nu$) there exist irreducible {\it detecting} objects $$A'_\mu \subseteq \tilde A_\mu \ ,$$ also contained in $H^+(A)$ (resp. in $H^-(A)$), such that $$Hom_{\calR_{n-1}}(A'_\mu,H^\pm(L))=0 \ \ \text{ and } \ \
Hom_{\calR_{n-1}}(A'_\mu,\bigoplus _{\rho\neq \nu} \tilde A_\rho)=0\ .$$
\end{enumerate}

\medskip
{\bf Remark}. For $*$-selfdual indecomposable objects as above the layers (graded pieces) of the upper and lower Loewy filtrations coincide, since otherwise proposition \ref{ext-0} would give a contradiction.  In the situation above we assume that $\A$ is $*$-selfdual of Loewy length 3 with socle $socle(\A) \cong L$ and $cosocle(\A) \cong socle(\A)^* \cong L^* \cong L$ and middle layer $A$.

\medskip
{\bf Remark}. 
For the later applications we notice  that we will construct the detecting objects $A'_\mu$ in $H^+(A^{down})$ (resp. $H^-(A^{down}))$
where $A^{down}$ will be an accessible summand of $A$. By induction we later will also know that these submodules $A'_\mu$ therefore already satisfy theorem \ref{mainthm}.  Hence it suffices to check the properties
$A'_\mu \subseteq \tilde A_\mu$ and $A'_\mu \subset H^{\pm}(A)$, since these already imply by the main theorem (valid for summands of $A^{down}$) the stronger assertion made in the axiom telling whether $A'_\mu$ appears in $H^+(A)$
or $H^-(A)$. Notice $A'_\mu \subseteq \tilde A_\mu$ and  $\tilde A_\mu \in \calR_{n-1}(\mp \varepsilon)$ depending on $\mu=\nu$ resp. $\overline \nu$. On the other hand
$A^{down} \subset A\in \calR_n(-\varepsilon)$. Hence, if the main theorem is valid for $A^{down}$,
we get $A'_\mu \in H^+(A)$ for $\mu=\nu$ and $A'_\mu \in H^-(A)$ for $\mu=\overline\nu$.

\begin{prop} \label{3} Under the assumptions on $\A$ from above the $H^\pm(A)$ will be semisimple objects in   $\calR_{n-1}(\mp \varepsilon)$.
\end{prop}

We will prove the key proposition \ref{3} below after listing some of its consequences.

\medskip
{\it The ring homomorphism $d$}. As an
element of the Grothendieck group $K_0(\calR_{n-1})$  we define for a module $M \in \calR_n$  $$d(M)= H^+(M) - H^-(M)\ .$$ Notice 
$d$ is additive by lemma \ref{hex}. 
Notice $$K_0(T_n) = K_0(\calR_n) \oplus K_0(\Pi \calR_n) = K_0(\calR_n) \otimes
(\mathbb Z \oplus \mathbb Z\cdot \Pi) \ .$$
We have a commutative diagram
$$  \xymatrix{  K_0(T_n) \ar[d]_{DS} \ar[r] &    K_0(\calR_n) \ar[d]^d \cr
 K_0(T_{n-1})  \ar[r] &    K_0(\calR_{n-1}) \cr} $$
 where the horizontal maps are surjective ring homomorphisms defined by
$\Pi \mapsto -1$. 
Since $DS$ induces a ring homomorphism, it is easy to see that $d$ defines
a ring homomorphism.

\medskip
The assertion of the last proposition implies that $H^+(A)$ and $H^-(A)$ have no common constituents in $\calR_{n-1}$ and that they are semisimple. 
Therefore $d(A)= H^+(A) - H^-(A) \in K_0(\calR_{n-1})$ uniquely determines $H^{\pm}(A)$ up to an isomorphism. By the additivity of $d$ and $d(\A) = \tilde\A$ we get  $2d(L) + d(A)= 2\tilde L + \tilde A$ in $K_0(\calR_{n-1})$. Hence 

\begin{cor}\label{determine-H+} $H^+(A)\in \calR_{n-1}(-\varepsilon)$ and $H^-(A)\in \calR_{n-1}(\varepsilon) $ are uniquely determined by the following formula in $K_0(\calR_{n-1})$ 
\[ H^+(A) - H^-(A) = d(A)= \tilde A  + 2(\tilde L - d(L)) \ .\] 
\end{cor}

We later apply this in situations where $\tilde L - d(L) = (-1)^{i+n}L^{aux}$ holds by lemma \ref{tildeL} and \ref{tildeLII} and $A' - \tilde A = 2(-1)^{i+n}L^{aux}$ holds
by lemma \ref{commute} and \ref{comII}, for some object $L^{aux}$. Here $A'$ denotes the normalized derivative
 of $A$, introduced in section \ref{signs}, defining a homomorphism
\[ {}': K_0(\calR_n) \to K_0(\calR_{n-1}) \ .\] 
Hence the last corollary implies the following theorem which 
repeatedly applied proves theorem \ref{mainthm}  by induction.

\begin{thm} \label{derive}
Under the axioms on $\A$ from above 
$d(A) = H^+(A) - H^-(A)$ is the derivative $A'$ of $A$. 
\end{thm}

{\it Proof of the proposition \ref{3}}.  Step 1).
Assumption (4) implies $H^+(\A)=\tilde \A^+$ and $H^-(\A)=\tilde \A^-$ in $\calR_{n-1}$.

\medskip Step 2). 
Axiom (1) on the Loewy structure of $\A$ therefore gives exact hexagons 
in $\calR_{n-1}$ for  $K:= Ker(\A \to L)$ using $K/L \cong A$:
$$ \xymatrix@+0.3cm{ H^+(K) \ar[r] & \tilde\A^+ \ar[r]^-{H^+(p)} & H^+(L) \ar[d]\cr
H^-(L) \ar[u]^\delta  & \tilde\A^- \ar[l]^{H^-(p)} &  H^-(K) \ar[l] }\  \quad
 \xymatrix@R+0.3cm{ H^+(K) \ar[r] & H^+(A) \ar[r] & H^-(L) \ar[d]\cr
H^+(L) \ar[u]^-{H^+(j)}  & H^-(A) \ar[l] &  H^-(K) \ar[l] }\ $$ 

Step 3) Assumption (3),(4),(5) on the Loewy structure of the $\tilde\A_\nu$ and $H^{\pm}(L)$ imply the following factorization property for $\tilde\A^+$ (and then similarly also for $\tilde\A^-$) 
$$  \xymatrix@-0.5cm{ \tilde \A^+= \bigoplus_{\nu} \tilde\A_\nu \ar[dr]_{\oplus q_\nu} \ar[rr]^{H^+(p)} & & H^+(L) \cr & \bigoplus_{\nu} \tilde L_\nu \ar@{.>}[ur]^{\exists !}_{\oplus p_\nu} &  } \ $$   

Step 4) Let $\Sigma$ be the set of all $\nu$ such that $p_\nu=0$. (Similarly
let $\overline\Sigma$ be the  set of all $\overline\nu$ such that $p_{\overline\nu}=0$).
Then we obtain exact sequences
$$ \xymatrix@R-0.2cm@C-0.1cm{  & & H^+(L)\ar[d]^{H^+(j)} & & \cr  
0\to \bigoplus_{\overline\nu\notin \overline\Sigma}\tilde L_{\overline\nu} \ar[r] & H^-(L)\ar[r]^\delta & H^+(K) \ar[d] \ar[r] & \bigoplus_{\nu\in\Sigma} \tilde\A_\nu \oplus \bigoplus_{\nu\notin\Sigma} \tilde K_\nu \ar[r] & 0  \cr
 & A'_\nu \ar@{.>}[ur]\ar[dr]^0 \ar@{^{(}->}[r] & H^+(A) \ar[d] & & \cr
   & & H^-(L) & & \cr } $$

\medskip   
Step 5) The detecting object $A'_\nu \hookrightarrow H^+(A)$ has trivial
image in $H^-(L)$ by axiom (6), hence can be viewed as a quotient object of $H^+(K)$. Again
by axiom (6) we can then view $A'_\nu$ as a nontrivial quotient object
of $$  H^+(K)/(I + \delta(H^-(L))) $$
where $ I := H^+(j)(H^+(L))$
is the image of $H^+(L)$ in $H^+(K)$.

\medskip
Step 6) The cosocle of $\bigoplus_{\nu\in\Sigma} \tilde\A_\nu \oplus \bigoplus_{\nu\notin\Sigma} \tilde K_\nu $ is $\bigoplus_{\nu\in\Sigma} \tilde L_\nu \oplus \bigoplus_{\nu\notin\Sigma} \tilde A_\nu $ by assumption (4) on the Loewy structure of $\tilde\A_\nu, \tilde K_\nu$.

\medskip
Step 7) The simple quotient object $A'_\nu$ of $H^+(K)$ can be viewed as a nontrivial
quotient object of the cosocle of $H^+(K)$ by step 5). We have 
an exact sequence
$$  H^-(L)/\bigoplus_{\overline\nu\notin \overline\Sigma}\tilde L_{\overline\nu} \to
cosocle(H^+(K)) \to cosocle( \bigoplus_{\nu\in\Sigma} \tilde\A_\nu \oplus \bigoplus_{\nu\notin\Sigma} \tilde K_\nu) \to 0$$
and we can view $A'_\nu$ as a nontrivial quotient object of  
$$ cosocle( \bigoplus_{\nu\in\Sigma} \tilde\A_\nu \oplus \bigoplus_{\nu\notin\Sigma} \tilde K_\nu)
 \ =\ \bigoplus_{\nu\in\Sigma} \tilde L_\nu \oplus \bigoplus_{\nu\notin\Sigma} \tilde A_\nu $$
by step 5) and 6). Notice that we can consider $A'_\nu$ for arbitrary $\nu$. For
$\nu\in \Sigma$ the last assertion contradicts axiom (6):  
$$  Hom_{\calR_{n-1}}(\bigoplus_{\mu\in\Sigma} \tilde L_\mu \oplus \bigoplus_{\mu\notin\Sigma} \tilde A_\mu, A'_\nu) \ = \ 0 \ .$$  
This contradiction forces
$$  \Sigma = \emptyset  \ \ \text{ and similarly } \ \ \overline\Sigma = \emptyset \ ,$$ 
so we obtain two exact sequences
$$ \xymatrix@C+0.5cm{ 0 \ar[r] & \bigoplus_{\overline\nu}\tilde L_{\overline\nu} \ar[r] & H^-(L) \ar[r] & H^+(K) \ar[r] & \bigoplus_\nu \tilde K_\nu \to 0} $$
$$ \xymatrix@C+0.5cm{ 0 \ar[r] & \bigoplus_{\nu}\tilde L_{\nu} \ar[r] & H^+(L) \ar[r] & H^-(K) \ar[r] & \bigoplus_{\overline\nu} \tilde K_{\overline\nu} \to 0} $$
Step 8) The last step 7) proves that
$$  \text{\it $H^+(p)$ is injective on the cosocle $\bigoplus_\nu \tilde L_\nu$ of $H^+(\A)$} \ .$$
Let  $i: L\hookrightarrow \A$ be the composition 
 of $j: L\hookrightarrow K$ and the inclusion $K\hookrightarrow \A$.
 Then $i: L \hookrightarrow \A$ is the $*$-dual of the projection $p: \A \twoheadrightarrow L$
by the axiom (2). Hence by $*$-duality we get from the previous assertion on $H^+(p)$
the following assertion 
$$ \text{ \it  $H^+(i)$ surjects onto the socle $\bigoplus_\nu \tilde L_\nu$ of $H^+(\A)$} \ .$$ 
Now considering
$$ \xymatrix{  L \ar@{^{(}->}[d]_j \ar@{^{(}->}[dr]^i &     \cr
K \ar@{^{(}->}[r] &  \A } \quad
\xymatrix{ \cr \text{and}  } \quad
 \xymatrix{  I \ar@{^{(}->}[d] \ar@{->>}[dr]   &  \cr
socle(H^+(K)) \ar[r] & socle(\tilde \A^+) \cong \bigoplus_\nu \tilde L_\nu } $$
we see that $\bigoplus_\nu \tilde L_\nu$ can also be embedded into the semisimple 
$I$ as a submodule 
$ \bigoplus_\nu \tilde L_\nu \ \hookrightarrow\  I  $.

\medskip Step 9)  Recall the following diagram 
$$ \xymatrix@R-0.4cm{ & & & I \ar[d] & & \cr 
& \bigoplus_{\nu}\tilde L_{\nu} \ar@{^{(}->}[r]  & H^-(L) \ar[r]^\delta & H^+(K) \ar[d]\ar[r]^\pi  & \bigoplus_\nu \tilde K_\nu \ar[r] & 0 \cr
& & & H^+(A) & & }$$
Since $I$ is in $\calR_{n-1}(\varepsilon)$ and  
$H^-(L) \in \calR_{n-1}(-\varepsilon)$ by our axioms, we also have
$$ \fbox{$ \delta(H^-(L)) \ \cap \ I \ = \ \{ 0\} $} \ .$$
Hence the composite of $\pi$ and the inclusion $I \hookrightarrow H^+(K)$ 
maps the semisimple module $I$
injectively into the socle of $\bigoplus_\nu \tilde K_\nu$. Since
$socle(\bigoplus_\nu \tilde K_\nu) = \bigoplus_\nu \tilde L_\nu$ and since $I$
contains $\bigoplus_\nu \tilde L_\nu $ as a submodule, this implies that
$$ \xymatrix{  \pi: I \ \ar[rr]^-\sim & & \ socle( \bigoplus_\nu \tilde K_\nu) \cong \bigoplus_\nu \tilde L_\nu } \ $$
is an isomorphism. Notice $(\bigoplus_\nu \tilde K_\nu)/socle( \bigoplus_\nu \tilde K_\nu) 
\cong \bigoplus_\nu \tilde A_\nu$.  

\medskip
Step 10) The last isomorphism of step 9) gives the exact sequence
$$ 0 \to \bigoplus_{\overline\nu} \tilde L_\nu \to H^-(L) \to \Bigl( H^+(K)/I \Bigr) \to \bigoplus_\nu \tilde A_\nu \to 0 \ .$$ 
By our assumptions $H^-(L)$ is in $\calR_{n-1}(-\varepsilon)$, and hence semisimple. 
Furthermore all  $\tilde A_\nu$ are semisimple and contained in 
$\calR_{n-1}(-\varepsilon)$. Hence by proposition \ref{ext-0} 
$H^+(K)/I$ is semisimple and contained in $\calR_{n-1}(-\varepsilon)$.

\medskip
Step 11). By step 10) and the exact hexagon 
$$ \xymatrix@R-0.7cm{ & H^+(K) \ar[r] & H^+(A) \ar[r] & H^-(L) \ar[dd]\cr
I \ar@{^{(}->}[ur]  & & &  \cr
& H^+(L) \ar@{->>}[ul] \ar[uu]_{H^+(j)}  & H^+(A) \ar[l] &  H^-(K) \ar[l] }\ .$$ 
$H^+(A)$ defines an extension of the semisimple module
$H^+(K)/I $ by a submodule of $H^-(L)$
$$ 0 \to H^+(K)/I \to H^+(A) \to Ker\Bigl(H^-(L)\to H^-(K)\Bigr) \to 0 \ .$$
Since $H^+(K)/I$ and $Ker(H^-(L)\to H^-(K))$ are both in $\calR_{n-1}(-\varepsilon)$, 
the proposition \ref{ext-0} implies that
$$ H^+(A) \ \cong \  
\Bigl( H^+(K)/I \Bigr) \ \oplus \  Ker\Bigl(H^-(j): H^-(L)\to H^-(K)\Bigr) \ $$
is semisimple and contained in $\calR_{n-1}(-\varepsilon)$. The first summand has been computed above.
Similarly then 
$$ H^-(A) \ \cong \  
\Bigl( H^-(K)/\overline I \Bigr) \ \oplus \  Ker\Bigl(H^+(j): H^+(L)\to H^+(K)\Bigr) \ $$
is semisimple and contained in  $\calR_{n-1}(\varepsilon)$.
\qed

\bigskip
{\bf Example}. Recall the indecomposable
$*$-selfdual objects $\A_{S^i}$ in $\calR_{n}, n \geq 2$ for $i=1,2,...$ with
Loewy structure $(L,A,L)$ where $L=S^{i-1}$ and 
$$  A = S^i \oplus S^{i-2} \oplus \delta_n^i \cdot Ber_n^{-1} \ .$$
Concerning the notations: $\delta_n^i$ denotes Kronecker's delta 
and $S^{-1}=0$. The conditions (1)-(5) are satisfied for $\varepsilon =(-1)^{i-1}$
and $A'=S^{i-2}$. Indeed condition (5) follows, since by induction on $i$
one can already assume that $H^-(L)=H^-(S^{i-1})$ is $Ber^{-1}$ or zero
and that $H^+(S^{i-2})$ contains $S^{i-2}$. Then by induction on $i$
the computation of $H^{\pm}(A)$ in terms of $\tilde A, \tilde L, H^\pm(L)$
from above easily gives as in section \ref{sec:strategy} the following result 

\begin{prop} Suppose $n\geq 2$. Then for 
the  functor $DS: \calR_n \to T_{n-1}$ of Duflo-Serganova we
obtain $DS(Ber_{n}) = \Pi(Ber_{n-1})$ and 
\begin{enumerate}
\item $DS(S^i) = S^i\ $ for $i<  n-1$, 
\item $DS(S^{i}) = S^{i} \oplus \Pi^{n-1-i} Ber^{-1}$  for $i\geq n-1$.
\end{enumerate}
\end{prop}

\medskip\noindent



\section{Moves}\label{sec:moves}

\medskip
We verify now the conditions on the indecomposable objects $\A$ of section \ref{sec:inductive} for the translation functors $F_i(L(\lambda_{\times \circ}))$. Additionally we verify the commutation rules in and after corollary \ref{determine-H+}. Instead of working directly with the irreducible representation $L$ we use the associated plot as in section \ref{derivat} and \ref{sec:loewy-length}. Recall that a plot $\lambda$ is a map $\lambda: \mathbb Z \to \{\boxplus,\boxminus\}$. We also use the notation $\boxplus_i$ to indicate that the $\lambda(i) = \boxplus$ and likewise for $\boxminus$. For an overview of the algorithms I and II used in this section see section \ref{sec:strategy}.

\medskip
 Let $L=(I,K)$ for $I=[a,b]$ be a segment
with sectors $S_1,..,S_k$ from left to right. Suppose $S_j =[i,i+1]$
is a sector of rank 1. Then the segment may be visualized as
$$  L \ = \ (S_1 \cdots S_{j-1} [\boxplus_i,\boxminus_{i+1}] S_{j+1} \cdots S_k) \ .$$
We define the {\it upward move} of the segment $L$ as the plot defined by
the {\it two segments} with intervals $[a,i-1]$ and $[i+1,b+1]$
$$ L^{up} \ =\ (S_1\cdots S_{j-1})\ \boxminus_i \ (\int(S_{j+1}\cdots S_k)) \ .$$
Similarly we define
the {\it downward move} of the segment $L$ as the plot defined by
the {\it two segments} with intervals $[a-1,i]$ and $[i+2,b]$
$$ L^{down} \ = \ (\int(S_1\cdots S_{j-1}))\ \boxminus_{i+1}\  (S_{j+1}\cdots S_k) \ .$$
Furthermore for $r\neq j$ we define additional $r$-th {\it internal} lower resp. upper
{\it downward moves} $L^{down}_r$ by the plots associated\footnote{For $r=j-1$ or $r=j+1$ the inner integral over the
empty sector is understood to give the sector $([i-1,i],\{i-1\})$ respectively $([i+1,i+2],\{i+1\})$.} to the {\it single segments}
$$ (S_1 \cdots S_{r-1} \int\bigl(S'_r \int (S_{r+1} \cdots S_{j-1})\bigr)\ S_{j+1} \cdots S_k) $$
for each $1\leq r \leq  j-1$ respectively
$$ (S_1 \cdots S_{j-1} \int\bigl(\int(S_{j+1} \cdots S_{r-1}) \ S'_r\bigr)\ S_{r+1} \cdots S_k) $$
for each $j+1\leq r \leq k$. Explaining the notion 'internal',
notice that the segments defined by these 
internal downward moves have the same underlying interval $I=[a,b]$ as the segment $L$
we started from.
We remark that the last formulas do remind on partial integration.
Formally by setting $L^{down}_r := L^{down}$ for $r=j$, we altogether obtain $k$ downward moves
and one upward move. All
these moves preserve the rank.

\medskip
The plot $L$ has a sector $[i,i+1]$ of rank 1. The {\it auxiliary plot} $L^{aux}$ attached to $L$ (and $[i,i+1]$) is the plot of rank
$r(L)-1$ defined by
two segments with intervals $[a,i-1]$ and $[i+2,b]$
$$ L^{aux}\ = \ (S_1\cdots S_{j-1})\ \boxminus_i \boxminus_{i+1} (S_{j+1}\cdots S_k) $$
and we also consider 
$$ L^{\times\circ} = (S_1\cdots S_{j-1})\ \times_i \circ_{i+1}\ (S_{j+1}\cdots S_k) \ .$$

\medskip
{\bf Algorithm I} (lowering sectors). 
For a plot with $k$ sectors $S_\nu$ with ranks $r_\nu=r(S_\nu) \geq 0$ and the distances
$d_\nu \geq 0$ for $\nu=1,....k$ (from left to right) we formally define $r_{k+1} = r_{k+1} = ... =0$ and $d_k=d_{k+1}=...= 0$. We  can then compare different plots with respect to the
lexicographic ordering of the sequences
$$  (-r_1,d_1,-r_2,d_2,..... ) \ .$$
Within the set of plots of fixed rank say $n$, the minimum with respect to this
ordering is attained if $r_1=n$, i.e. if there exists only one sector.

\medskip
{\it Algorithm I will be applied to given plots, say $\lambda_{\wedge\vee}$, 
with {\it more than one segment}}.
The upshot is: In this situation
one can always find a lexicographic smaller plot $L$ so that
the given plot is of the form $\lambda_{\wedge\vee} = L^{up}$ and such that  $L$ and all plots obtained by the moves
$L^{down}_r$ of $L$ are strictly smaller than the starting plot $L^{up}$. 
Algorithm I is used for induction arguments to reduce certain statements
(e.g. theorem \ref{mainthm}) to the case of plots with 1 segment.

\medskip
{\it Definition of $L$}. For a given plots say $\lambda_{\wedge\vee}$, 
with more than one segment, $d_\nu>0$ holds for some integer $\nu$. So choose $j$ so that the distances $dist(S_1,S_2)=...=dist(S_{j-2},S_{j-1}) =0$ for the sectors $S_1,..,S_{j-1}$ of $\lambda_{\wedge\vee}$ and $dist(S_{j-2},{j-1}) > 0$. We temporarily write $S$ for the next sector $S$ of $\lambda_{\wedge\vee}$. Interpret $S=\int(S_{j+1} \cdots S_k)$ for some sectors $S_j,...,S_k$. This is possible, but keep in mind that $S_j, ... , S_k$ are not sectors of $\lambda_{\wedge\vee}$ but will be sectors of $L$, and this explains the notation. Indeed, 
for $i+1=min(S)$, we define $L$ to be 
$$ L \ = \  (S_1 \cdots S_{j-1}) ... d_{j-1} .... (S_jS_{j+1}\cdots S_k) ... d_k ....   $$
with $S_{j}$ of rank 1 at the positions $[i,i+1]$. 
To  simplify notations we do not write further sectors to the right, since the sectors of $\lambda_{\wedge\vee}$ to the
right of $S$ will not play an essential role in the following. Indeed, they will
appear verbatim in the sector structure of $L$ up to some distance shifts at the 
following positions
$$   dist(S_{j-1},S) = 1+ d_{j-1} \quad , \quad dist(S,\text{next sector}) = d_k -1 \ .$$
Concerning the lexicographic ordering $$dist(S_{j-1},S_j) =d_{j-1} < dist(S_{j-1},S) = 1+ d_{j-1}$$
shows that $L$ is smaller than $L^{up}= \lambda_{\wedge\vee}$. We leave it to the reader
to check that also all $L^{down}_j$ are smaller than $L^{up}= \lambda_{\wedge\vee}$.
Notice, here we apply the moves as in the preceding paragraph with the notable exceptions that
\begin{enumerate}
\item  There may be further sectors beyond $S_k$. These are just appended, and do not
define new moves.
\item  If $d_{j-1} \geq 1$ the sector $S_{j-1}$ has distance $>0$ to the sector 
$S_j$ and therefore does not define downward moves, so that only the doenward moves $L^{down}_r$ for
$r=j,...,k$ are relevant.   
\end{enumerate}
In the later discussion we always display the {\it more complicated} case where $d_{j-1}=0$
(without further mentioning). For the case $d_{j-1}>0$ one can simply \lq{omit}\rq\ $S_1,...,S_{j-1}$, by just appending them 
in the same way as we agreed to \lq{omit}\rq\ sectors to the right of $S_k$.

\medskip
{\bf Construction of detecting objects for algorithm I}.
Fix $L= L(\lambda_{\wedge\vee})$ with the sector $[i,i+1]$. 
Then $L$ is determined by its sectors. For the construction
of detecting objects we are only interested in down moves.
In the following it therefore suffices only to keep track of the sectors
below $[i,i+1]$ in the segment containing the sector $[i,i+1]$. 
Notice that $L$ is a union of the sector $[i,i+1]$ and, say $s$, other sectors $S_\nu$.
Let $S_1,..,S_{j-1}$ denote the sectors below $[i,i+1]$ in the segment
of $[i,i+1]$. Hence $L$ is
$$  \boxminus S_1 \cdots S_{j-1} [\boxplus_i \boxminus_{i+1}] $$
and the union of other disjoint sectors $S_\nu$ for $j+1 \leq \nu \leq s$.
Then $L^{\times \circ}$ is
$$  \boxminus S_1 \cdots S_{j-1} [\times_i \circ_{i+1}] $$
and the union of other disjoint sectors $S_\nu$ for $j+1 \leq \nu \leq s$.
We define $\A = F_i L^{\times \circ}$. Then $\A$ is $*$-self dual of Loewy length 3
with socle and cosocle $L$. The term $A$ in the middle is semisimple
and the weights of its irreducible summands are given by $L^{up}$ and 
the $k$ down moves of $L$ according to section \ref{sec:loewy-length}.

\medskip
To determine $\tilde \A = DS(\A)$ we use induction and 
lemma \ref{tildeA}. This implies that $\tilde \A$ is the direct sum
of $\Pi^{m_{\nu}} \tilde \A_\nu$ for indecomposable objects $\tilde \A_\nu$ in $\calR_{n-1}$, which uniquely correspond to the
irreducible summands of $DS(L^{\times \circ})$. However these correspond to
the irreducible
summands $\tilde L_\nu$ of $DS(L^{aux})$. 
Again by induction (now induction on the degree
of atypicity) the summands of $DS(L^{\times\circ})$ respectively $DS(L^{aux})$ are already known to be given by the derivative of $\lambda_{aux}$. These facts imply
the next

\medskip
\begin{lem} \label{tildeL} We have 
$$  \tilde \A = \bigoplus_{\mu=1}^s \tilde \Pi^{m_{\mu}} \A_\mu$$
where each $\tilde \A_\mu \in\calR_{n-1}$ has Loewy lenght 3 with irreducible socle
and cosocle $\tilde L_\mu$ defined by the $s$ plots for $\mu=1,...,s$
$$ [\boxplus_i \boxminus_{i+1}] \cup S'_\mu \cup \bigcup_{\nu\neq \mu} S_\nu \ .$$
In particular, for $m_{aux}$ (which is congruent to $ i+n-1$ modulo 2), we get\footnote{assuming that theorem \ref{mainthm} holds for $L$, say by induction assumption.}
$$ \fbox{$ DS(L) \cong \Pi^{m_{aux}} L^{aux} \ \oplus \ \bigoplus_{\mu=1}^s \tilde \Pi^{m_{\mu}} L_\mu $} \ .$$ 
Hence in $K_0(\calR_{n-1})$
$$ \fbox{$ d(L)\ =\ L' \ =\ \tilde L \ + \ (-1)^{i+n-1} \cdot L^{aux} $}  \ .$$ 
\end{lem}

\medskip
Now each $\tilde A_\mu$ is determined from $\tilde L_\mu$
by applying certain upward and the downward moves starting
from $\tilde L_\mu$. 

\medskip
We indicate that the segment of $\tilde L_\mu$ containing 
$[\boxplus_i \boxminus_{i+1}]$ has less than $r$ sectors,
if $1 \leq \mu \leq j-1$. Indeed the union  of the sectors of $\tilde L_\mu$
in the segment of $[\boxplus_i \boxminus_{i+1}]$ 
is
$$ ... \boxminus S_{\mu+1} \cdots S_{j-1} [\boxplus_i \boxminus_{i+1}] S_{j+1} \cdots S_r  \boxminus ... \ $$
for $\mu \leq j-2$ and by $[\boxplus_i \boxminus_{i+1}] S_{j+1} \cdots S_r  \boxminus ... $
for $\mu= j-1$. We are now able to define the {\it detecting objects} $A'_\mu \subseteq
\tilde A_\mu$ 
for $\mu=1,...,s$ by $\tilde L_\mu^{down}$, given by induction as follows 
\begin{enumerate}
\item $ (\int(S_1 \cdots S_{j-1})) \boxminus_{i+1}  \cup (S'_\mu) \cup \bigcup_{j-1 <\mu\neq \ell} S_\ell $ for $\mu \notin \{1,...,j-1\}$,
\item $S_1 \cdots S'_\mu (\int (S_{\mu +1} \cdots S_{j-1})) \boxminus_{i+1} S_{j+1} \cdots S_k 
\cup \bigcup_{k <\ell} S_\ell $ for $\mu \leq j-2$, 
\item $S_1 \cdots S_{j-2} \boxminus S'_{j-1} (\boxplus \boxminus_i) \boxminus_{i+1} S_{j+1} \cdots S_k 
\cup \bigcup_{k <\ell} S_\ell $ for $\mu= j-1$.
\end{enumerate}
It is therefore clear that the detecting object is different from all objects 
in $DS(L)$, which by induction are known to be given by the derivative of $L$.
Furthermore $A'_\mu \subseteq \tilde A_\mu$. It requires some easy but tedious inspection
to see that $A'_\mu$ is not contained in $\tilde A_\nu$ for $\nu\neq \mu$.
Hence to see that the $A'_{\mu}$ are detecting objects, it suffices to show the next

\begin{lem}
The objects $A'_\mu$ are contained in $DS(A)$. If $L^{up}$ is stable, then $L$ is stable
and $A'_\mu \subset H^+(A)\oplus H^-(A)$ for all $\mu$.
\end{lem} 

\medskip
{\it Proof}. Recall $$A \cong A^{up} \oplus A^{down} $$
for $ A^{up} := L^{up}$ and  $A^{down} := \bigoplus_{i=1}^k L^{down}_r$.

\bigskip\noindent
We do not know how to compute $DS(A^{up})$. However by induction we already know
that the derivative computes $DS(A^{down})$. In $A^{down}\subset A$ we have the following
objects $A_\mu$
\begin{enumerate}
\item $ (\int(S_1 \cdots S_{j-1})) \boxminus_{i+1} \cup \bigcup_{j-1 <\ell} S_\ell $ for $\mu \notin \{1,...,j-1\}$,
\item $S_1 \cdots \int( S'_\mu \ \int (S_{\mu +1} \cdots S_{j-1})) S_{j+1} \cdots S_k 
\cup \bigcup_{k <\ell} S_\ell $ for $\mu \leq j-2$,
\item $S_1 \cdots S_{j-2} (\boxplus S'_{j-1} \boxplus \boxminus_i \boxminus_{i+1}) S_{j+1} \cdots S_k 
\cup \bigcup_{k <\ell} S_\ell $ for $\mu= j-1$.
\end{enumerate} 
Their derivative $DS(A_\mu)$ contains 
\begin{enumerate}
\item $ (\int(S_1 \cdots S_{j-1}) \boxminus_{i+1}) \cup (S'_\mu) \cup \bigcup_{j-1 <\ell\neq \mu} S_\ell $ for $\mu \notin \{1,...,j-1\}$,
\item $S_1 \cdots S'_\mu (\int (S_{\mu +1} \cdots S_{j-1})) \boxminus_{i+1} S_{j+1} \cdots S_k 
\cup \bigcup_{k <\ell} S_\ell $ for $\mu \leq j-2$,  
\item $S_1 \cdots S_{j-2} \boxminus S'_{j-1} (\boxplus \boxminus_i) \boxminus_{i+1} S_{j+1} \cdots S_k 
\cup \bigcup_{k <\ell} S_\ell $ for $\mu= j-1$.
\end{enumerate}
This proves $A'_{\mu} \subset DS(A_\mu)$ and
hence our claim. \qed

\medskip
{\bf Commutation rule for algorithm I}. Now we discuss how moves commute with
differentiation for a given $L$ as above. It is rather obvious from the
definitions that 
for this we can restrict ourselves to the situation where $L$ is
the single segment 
$$  L = (S_1 \cdots S_{j-1} \boxplus_i \boxminus_{i+1} S_{j+1} \cdots S_k) \ .$$
So let us assume this for simplicity of exposition.

\medskip
1) {\it Computation of $\tilde A$}. Taking first the derivative we obtain $L^{aux}$ and $(k-1)$ plots $\tilde L_\mu$ of the form
$$ S_1 \cdots S'_\mu (S_{\mu+1} \cdots S_{j-1} \boxplus_i \boxminus_{i+1} S_{j+1} \cdots S_k) \ $$
(lower group where $\mu \leq j-1$) respectively
$$ (S_1 \cdots S_{j-1} \boxplus_i \boxminus_{i+1} S_{j+1} \cdots S_{\mu -1}) S'_\mu \cdots S_k \ $$
(upper group where $\mu \geq j+1$). The sign in the Grothendieck group attached to these is $(-1)^{a_\mu + n-1} =(-1)^{i+n-1}$ for $S_\mu = [a_\mu,b_\mu]$.

\medskip
Notice that
$L^{aux}$ does not define any moves. 
The segment containing $\boxplus_i \boxminus_{i+1}$ (indicated by the brackets)
defines the possible moves of each of these derived plots $\tilde L_\mu$. 
These are e.g. in the {\it lower group case}  the upward move
$$ \fbox{$ S_1 \cdots S'_\mu S_{\mu +1} \cdots S_{j-1} \boxminus_i \int (S_{j+1} \cdots S_k) $} \ $$
and the downward move
$$ \fbox{$ S_1 \cdots S'_\mu \int(S_{\mu +1} \cdots S_{j-1}) \boxminus_{i+1} S_{j+1} \cdots S_k$} \ $$
and the internal upper/lower downward moves 
$$ \fbox{$ S_1 \cdots S'_\mu S_{\mu+1} \cdots S_{j-1} \int(\int ( S_{j+1} \cdots S_{r-1}) \ S'_r\ ) 
S_{r+1} \cdots S_k $} \ $$
$$ \fbox{$ S_1 \cdots S'_\mu S_{\mu+1} \dots S_{r-1} \int (S'_r \int (S_{r+1} \cdots S_{j-1}))\  S_{j+1} \cdots S_k $} \ $$

\medskip
2) {\it Computation of the derivative $A'$}. Now we revert the situation and first consider the moves of $L$, the upward
and downward moves 
$$ L^{up} \ =\ S_1\cdots S_{j-1}\ \boxminus_i \ \int(S_{j+1}\cdots S_k) \ ,$$
$$ L^{down} \ = \ \int(S_1\cdots S_{j-1})\ \boxminus_{i+1}\  S_{j+1}\cdots S_k \ ,$$
and the internal downward moves (for lower sectors)
$$ S_1 \cdots S_{r-1} \int\bigl(S'_r \int (S_{r+1} \cdots S_{j-1})\bigr)\ S_{j+1} \cdots S_k $$
respectively (for upper sectors)
$$ S_1 \cdots S_{j-1} \int\bigl(\int(S_{j+1} \cdots S_{r-1}) \ S'_r\bigr)\ S_{r+1} \cdots S_k \ .$$ 
If we differentiate $L^{up}$, we get the plots of the form
$$ \fbox{$ S_1\cdots S'_\mu \cdot S_{j-1}\ \boxminus_i \ \int(S_{j+1}\cdots S_k) $} $$
with sign $(-1)^{s_\mu + n-1} = (-1)^{i +n-1}$ and similarly
$$ L^{aux} = S_1\cdots S_{j-1}\ \boxminus_i \ \boxminus_{i+1} S_{j+1}\cdots S_k $$
with sign $(-1)^{i+n}$.
If we differentiate $L^{down}$, we get 
$ \int(S_1\cdots S_{j-1})\ \boxminus_{i+1}\  S_{j+1}\cdots S'_\mu  \cdots S_k $
and similarly 
$$ L^{aux} = S_1\cdots S_{j-1}\ \boxminus_i \ \boxminus_{i+1} S_{j+1}\cdots S_k \ $$
with sign $(-1)^{i+n}$. 
If we derive the plots defined by the internal moves ({\it lower group}, where we derive at  
$\nu \leq j-1$) 
we get the plots of the form 
$$ \fbox{$ S_1 \cdots S_{r-1} S'_r \int (S_{r+1} \cdots S_{j-1})\boxminus_{i+1} \ S_{j+1} \cdots S_k $} $$ with sign $(-1)^{s_r + n-1} = (-1)^{i+n-1}$ 
together with
$$ \fbox{$ S_1 \cdots S'_\mu S_{\mu+1} \cdots S_{j-1} \int(\int ( S_{j+1} \cdots S_{r-1}) \ S'_r\ ) 
S_{r+1} \cdots S_k $} \ $$
$$ \fbox{$ S_1 \cdots S'_\mu S_{\mu +1} \cdots S_{r-1} \int\bigl(S'_r \int (S_{r+1} \cdots S_{j-1})\bigr)\ S_{j+1} \cdots S_k $} $$
of sign $(-1)^{i+n-1}$ respectively similar terms for the upper group, where we differentiate at $\mu \geq j+1$.
Altogether, besides two additional signed plots of the form $L^{aux}$, these give precisely the plots
obtained before. This implies

\begin{lem} \label{commute}The differential of the moves of $L$ gives the term $2 \cdot L^{aux}$
plus the moves of the differential of $L$, i.e.
$$   \fbox{$ {A'} \ =\ \tilde A \ +\  2(-1)^{i+n} \cdot L^{aux} $} \  $$
holds in $K_0(\calR_{n-1})$.
\end{lem}

\medskip
{\bf Algorithm II} (melting sectors). Suppose $\lambda_{\wedge\vee}$
is a plot  with a single segment and at least two sectors. $[a,i]$ and
point $i$ respectively left boundary point $i+1$ of the segment defining 
the plot $\lambda_{\wedge\vee}$.
In algorithm II we melt the first two adjacent sectors $[a,i]$ and $[i+1,b]$ 
together into a single sector $S^{melt}$ to obtain
a new plot $\lambda_{\vee\wedge}$ so that $$supp(\lambda_{\wedge\vee}) - \{i+1\} = supp(\lambda_{\vee\wedge}) - \{i\} \ .$$
This new plot $\lambda_{\wedge\vee}$ again has a unique segment 
with the same underlying interval as  
the plot $\lambda_{\wedge\vee}$. But the sector structure is different, since 
the number of sectors decreases by one.

\medskip
Notice, opposed to algorithm I, the interval $[i,i+1]$ does not define
a sector of the original plot $\lambda_{\wedge\vee}$. However $[i,i+1]$
defines a sector
of the 'internal' plot $$L_{int} := \partial(S^{melt}) \ ,$$ with sector structure say
$$  L_{int} = S_1 \cdots S_{j-1} [\boxplus_i \boxminus_{i+1}] S_{j+1} \cdots S_k \ ,$$
so that
$$ \lambda_{\vee\wedge} = (\int L_{int})  \ \text{ other sectors} \quad , \quad
 \lambda_{\times\circ} :=  \ (\int (L_{int})^{\times\circ} )  \ \text{ other sectors} $$
We similarly define 
for $r=1,...,k$ and $r\neq j$ the plots
$$  \lambda^{down}_r \ := \  (\int (L_{int})^{down}_r )\ \text{ other sectors} \ .$$
Finally $\  \lambda_{\wedge\vee} =  (\int (L_{int})^{up} )\text{other sectors }$, 
which is the plot we started from. Since $\int (L_{int})^{up}$ has two sectors,
all the plots $\lambda^{down}_r$ for $1 \leq r\neq j \leq k$ have less sectors than the 
plot $\lambda_{\wedge\vee}$. Indeed, 
the plots $\int (L_{int})^{down}_r$ are irreducible as an easy consequence 
of the integral criterion.

\medskip
{\bf Construction of detecting objects for algorithm II}.
Fixing $\lambda_{\wedge\vee}$ as above, 
$$  \A = F_i(L(\lambda_{\times\circ})) \ = \ (L,A,L) $$
defines a $*$-self dual object in $\calR_n$ of Loewy length 3
with socle and cosocle $L$,
where 
$$  L= L(\lambda_{\vee\wedge}) $$
and $A=A^{up} \oplus A^{down}$ for
$ A^{up} = L(\lambda_{\wedge\vee}) $ and $
A^{down} = \bigoplus_{r\neq j,1}^k L(\lambda_r^{down})$. 

\medskip
To determine $DS(\A) = \bigoplus_\mu \Pi^{m_{\mu}} \tilde\A_\mu$ we use induction and 
lemma \ref{tildeA}. This implies that $\tilde \A$ is the direct sum
of $\tilde \Pi^{m_{\mu}} \A_\mu$ for indecomposable objects $\tilde \A_\mu$ in $\calR_{n-1}$, which uniquely correspond to the irreducible
summands $\tilde L_\mu$ of $DS(L^{\times\circ})$. But by induction (now induction on the degree 
of atypicity !) the irreducible summands of $DS(L^{\times\circ})$, that determine the 
irreducible modules $\tilde L_\mu$, can be computed by the derivative of $\lambda_{aux}$. Since in the present situation
replacing $i,i+1$ by $\times,\circ$ commutes with the derivative, these facts imply
the next

\medskip
\begin{lem} \label{tildeLII} If $L$ has $s$ sectors, for the melting algorithm we have 
$$  \tilde \A = \bigoplus_{\mu=1}^s \Pi^{m_{\mu}} \tilde \A_\mu $$
where each $\tilde \A_\mu \in\calR_{n-1}$ has Loewy length 3 of Loewy structure $(\tilde L_\mu, \tilde A_\mu, \tilde L_\mu)$ with irreducible socle
and cosocle $\tilde L_\mu$. 
 For the various summands, for varying $\mu$, up to the shift $m_\mu$ the socles $L_\mu$ are defined by the $s-1$ different plots plots arising from the derivative 
$$   (\int L_{int})\
\text{(other sectors)}' $$ together with the plot $$  L_{int}
\ \text{(other sectors)} \ .$$
In particular\footnote{assuming that theorem \ref{mainthm} holds for $L$ and $L^{\times\circ}$, say by induction assumption.},
if $L$ has $s$ sectors,
$ DS(L) \cong  \bigoplus_{\mu=1}^s \Pi^{m_{\mu}} \tilde L_\mu $.
This gives in $K_0(\calR_{n-1})$ the formula
$$ \fbox{$ d(L)\ =\ L' \ =\ \tilde L  $}  \ .$$ 
\end{lem}

\begin{cor} \label{cohomII}
In the situation of the last lemma the morphisms $$H^i(p): H^i(\A) \to H^i(L)$$ are surjective for all $i\in \mathbb Z$.
\end{cor}

\medskip
{\it Proof}. We already know that $H^{\pm}(p) : H^{\pm}(\A) \to H^{\pm}(L)$
induces injective maps on the cosocle of $H^{\pm}(\A)$. By lemma \ref{tildeLII}
therefore these induced maps are bijections between the cosocle of $H^{\pm}(\A)$
and $H^{\pm}(L)$. In particular  the morphisms  $H^{\pm}(p) : H^{\pm}(\A) \to H^{\pm}(L)$
are surjective. This implies the assertion. \qed 

\medskip
This being said  
note that $2d(L)+d(A) = d(\A) =d(\tilde \A) = 2\tilde L + \tilde A$
together with the assertion $d(L) = \tilde L$ from the lemma \ref{tildeLII} above 
implies $d(A) = \tilde A$. Any $\tilde L_\mu$ defines a nontrivial
term $\tilde A_\mu$. We claim that any irreducible summand 
$A'_\mu \subset \tilde A_\mu$ is a detecting object now. 
Indeed any summand $A'_\mu $ of $\tilde A$ appears in $H^+(A)$
by the formula $d(A) = H^+(A) - H^-(A) = \tilde A$. Checking the possible
moves that define the constituents of $\tilde A_\mu$ from $\tilde L_\mu$
it is clear that $A'_\mu$ is not a constituent of any $\tilde A_\nu$ for $\nu\neq \mu$.  
Hence

\begin{lem} Detecting objects $A'_\nu$ exist for algorithm II.
\end{lem}

\medskip
{\bf Commutation rule for algorithm II}. Now we discuss how moves commute with
differentiation for a given $L$ as above. It is rather obvious from the
definitions that we can restrict ourselves for this to the situation where the segment of
plot $\lambda_{\wedge\vee}$ has only two sectors. In other words we claim that
we can assume without restriction of generality that the terms 'other factors' does not appear, so that $s=2$ holds in the last lemma \ref{tildeLII}. The reason for this is, that moves for $L_{int} \ (others\ sectors) $ are the same as for $\int L_{int}\ (others\ sectors)'$, since by \cite{Brundan-Stroppel-2} the relevant moves are moves 'within' the sector $\int L_{int}$. Hence for the proof of the next 
lemma we can assume that $L=\int L_{int}$ has a unique sector so that
$d(L)=\tilde L=L'= L_{int}$, which has a single segment.

\begin{lem}  \label{comII} 
The differential of the moves of $L$ gives the moves of the differential of $L$, i.e.
$$   \fbox{$ {A'} \ =\ \tilde A \ $} \  $$
holds in $K_0(\calR_{n-1})$.
\end{lem}

\medskip{\it Proof}. Without restriction of generality we can assume
that the plot $\lambda_{\vee\wedge}$, we are starting with,  is a segment with only two sectors,
so that $L= \int L_{int}$. Let the single segment of $L_{int}$ have the form
$$ S_1 \cdots S_{j-1} [\boxplus_i \boxminus_{i+1}] S_{j+1} \cdots S_k \ $$
with $k$ sectors $S_1,...,S_k$ where the underlying interval of $S_j$ is $[i,i+1]$.

\medskip
1) {\it Computation of $\tilde A$}. According to \cite{Brundan-Stroppel-2}, \cite{Weissauer-gl} the constituents of $\tilde\A = \bigoplus_{\mu=1}^s \Pi^{m_{\mu}} \tilde\A_\mu$ are obtained from the socle module $\tilde L_\mu$ of $\tilde\A_\mu$
by moves. The last lemma shows that $\tilde L = \bigoplus_\mu \tilde L_\mu$ is 
the derivative $L' = (\int L_{int})' = L_{int}$ of $L$ up to a shift determined by the sign factor $(-1)^{a + n-1}$. Since $s=1$ by assumption, $\tilde \A= (\tilde L, \tilde A, \tilde L)$ is an indecomposable module
with socle
$$ \tilde L  \ = \ L_{int} \ .$$
Up to a parity shift by $m=a+n-1$, the module $\tilde A$ therefore is the direct sum 
$$ \tilde A \ =\ (L_{int})^{up} \oplus\  \bigoplus_{r=1}^k  \ (L_{int})^{down}_r \ $$
of the irreducible modules
obtained from $\tilde L = L_{int}$ by the unique upward and the $k$ downward moves.
Notice $(L_{int})^{down}_j = (L_{int})^{down}$ is the 'nonencapsulated'
downward move in the notions of \cite{Weissauer-gl}. Here it occurs, since $[i,i+1]$ is one of the sectors
of the $\tilde L$.

\medskip
2) {\it Computation of the derivative $A'$}. Now we revert the situation and first consider $\A =(L,A,L)$ and the moves of $L=\int L_{int}$ that determine the irreducible summands of $A$. Indeed
$$  A = \ \Bigl(\int L_{int}\Bigr)^{up} \oplus \bigoplus_{r\neq j, r=1}^k  \Bigl(\int L_{int}\Bigr)^{down}_r \ $$
holds for the irreducible modules
obtained from $L = \int L_{int}$ by the upward move $L^{up}$ and the $k-1$ internal inner/upper downward moves $(L_{int})^{down}_r$ for $r\neq j$. Notice that $(L_{int})^{down}_j = (L_{int})^{down}$, as opposed to the situation above, this time does not
appear as a move, since we are in the 'encapsulated' case in the notions  \cite{Weissauer-gl} where $[i,i+1]$ is not a sector of $L$ (but
only an internal sector of $L$). 

\medskip
The formulas above imply that $A' = (A^{up})' \oplus (A^{down})'$ is a direct sum  of the two irreducible summands $$ (A^{up})' = B_1 \oplus B_2 \ ,$$
coming from
$ ((\int L_{int})^{up})' = (L^{up})' = L(\lambda_{\wedge\vee})' $ for
$\lambda_{\wedge\vee}=[a,i][i+1,b]$ with derivative $(-1)^{a+n-1}( \partial([a,i]) \cup [i+1,b] ) + (-1)^{i+n-1} ([a,i] \cup \partial([i+1,b]))$,
and the $k-1$ irreducible summands $(A^{down}_r)'$ of $(A^{down})'$ given by 
$$ (A^{down}_r)' \ = \  \Bigl(\Bigl(\int L_{int}\Bigr)^{down}_r\Bigr)' \ .$$
This gives $2+(k-1)=k+1$ irreducible factors in $\tilde A$,
and all signs coincide by $(-1)^{i+n-1} = (-1)^{a+n-1}$.

\medskip
{\it The comparison}. Since all signs are $(-1)^{a+n-1}$ for both computations, we can ignore the parity shift. Then observe that $ (\int L_{int})^{down}_r = \int (L_{int})^{down}_r)$
holds for $r\neq j$, hence $(A^{down}_r)' = ((\int L_{int})^{down}_r)' =  L_{int})^{down}_r $
for $r\neq j$.
So it remains to compare the two remaining summands
$$ B_1 \quad \ \ , \ \  \quad B_2$$ of $A'$ and the two remaining summands  
$$ (L_{int})^{up} \quad , \quad  (L_{int})^{down}_j $$ 
of $\tilde A$. 
The latter correspond to the plots $S_1... S_{j-1} \boxminus_i \int (S_{j+1} ... S_k)$, giving the upward move, 
resp. $\int (S_1 ... S_{j-1}) \boxminus_{i+1} S_{j+1} \cdots S_k$, giving the downward move. Obviously  these two
define the plots $\partial([a,i]) \cup [i+1,b]$ respectively $[a,i] \cup \partial([i+1,b]) $
defining the two summands $B_1$ and $B_2$. \qed

\bigskip



\part{Consequences of the Main Theorem}
\medskip
We describe some applications of the main theorem. The main result is the computation of the $\Z$-grading of $DS(L)$ for any irreducible representation $L$ in sections \ref{kohl} - \ref{koh3}. This result is based on the main theorem and its proof. For that we first need a description of the dual of an irreducible representation in section \ref{duals}. In the later sections \ref{kac-module-of-one} - \ref{hooks} we obtain various results about the cohomology of maximally atypical indecomposable representations.

\medskip

\section{Tannaka Duals}\label{duals}

\medskip
Let $\lambda$ be an atypical weight, and  $L(\lambda)$
the associated irreducible representation. Note that $(Ber^k \otimes L(\lambda))^{\vee} = Ber^{-k} \otimes L(\lambda)^{\vee}$. We use the description of the duals obtained in \cite{Heidersdorf-mixed-tensors}. Note that $L(\lambda)=socle(P(\lambda))= cosocle(P(\lambda))$, since projective modules are $*$-self dual. Hence $L(\lambda)^\vee = socle(P(\lambda)^{\vee})$, so it suffices to compute the socle of $P(\lambda)^\vee$. Now $P(\lambda) = R(\lambda^L,\lambda^R)$ for the bipartition $(\lambda^L,\lambda^R) = \theta^{-1}(\lambda)$ satisfying $k(\lambda^L; \lambda^R) = n$ by \cite{Heidersdorf-mixed-tensors}. The dual of any mixed tensor is $ R(\lambda^L,\lambda^R)^{\vee} = R(\lambda^R,\lambda^L)$, hence we simply have to calculate the socle of $R(\lambda^R,\lambda^L)$. 

\medskip

If $k(\lambda^L,\lambda^R) = n$, the description of the map $\theta$ is easy: Calculate the weight diagram of $(\lambda^L, \lambda^R)$ as in section \ref{stable0} and write down its labeled cup diagram. Then turn all $\vee$'s which are not part of a cup into $\wedge$'s and leave all other symbols unchanged. The resulting diagram is the weight diagram of $socle(P(\lambda))$. Hence in order to calculate the dual of $L(\lambda)$ we simply have to understand the effect of changing $(\lambda^L,\lambda^R)$ to $(\lambda^R,\lambda^L)$ on the weight diagram. Recall from section \ref{stable0} that \[ I_{\wedge}(\lambda)  := \{ \lambda_1^L, \lambda_2^L - 1, \lambda_3^L - 2, \ldots \} \quad \text{and}\quad  I_{\vee}(\lambda)  := \{1  -\lambda_1^R, 2 - \lambda_2^R, \ldots \}\ . \]

If $\lambda_i^L - (i-1) =s$, then $i - \lambda_i^L = i - s - i + 1 = 1 -s$ and likewise for $\lambda_j^R$. Hence interchanging $\lambda^L$ and $\lambda^R$ means reflecting the symbols $s \mapsto 1-s$ and swapping $\vee$'s with $\wedge$'s. If the vertex $s$ is labelled by a $\times$, then there exist $i,j$ such that $\lambda_i^L - (i-1) = j - \lambda^R_j = s$. But then $\lambda_j^R - (j-1) = i - \lambda_i^L  = 1-s$ and we obtain a $\times$ at the vertex $1-s$. We argue in the same way for the $\circ$'s. If $(s,s+r)$ is labelled by $(\vee, \wedge)$-pair such that we have a cup connecting $s$ and $s+r$, we obtain a $(\vee, \wedge)$-pair at $(1-s-r,1-s)$ which is connected by a cup. To obtain the highest weight $\theta(\lambda^R,\lambda^L)$ the $\vee$'s not in cups get flipped to $\wedge$'s.

\begin{prop} \label{irreducible-dual} The weight diagram of the dual of an irreducible representation $L$ is obtained from the weight diagram of $L$ as follows: Interchange all $\vee \wedge$-pairs in cups, then apply the reflection $s \mapsto 1-s$ to each symbol.
\end{prop} 

It is easy to see that this description is valid for $m \geq n$ if we use the reflection $1 - \delta -s$ instead of $1-s$ where $\delta = m-n$.

\medskip \textit{The maximal atypical case}. We describe the dual in the language of plots. We assume here that $L(\lambda)$ is maximally atypical, but we can reduce the general case to this one, using the map $\phi$ from section \ref{sec:loewy-length} and lemma \ref{phi-dual}. Let $\lambda$ denote the unique plot corresponding to the weight $\lambda$.
Let $\lambda(s) = \prod_i \lambda_i(s)$ be its prime factorization. For each prime factor $\lambda_i(s)=(I,K)$ with segment $I$ and support $K$ we define $ \lambda_i^c(s):= (I,K^c)$,
where $K^c = I\!-\!K$ denotes the complement of $K$ in $I$.
Then put $$ \lambda^c(s) := \prod_i \lambda_i^c(s)  \ .$$

The previous description of the duals implies the next proposition.

\begin{prop}\label{dual-plot}  The Tannaka dual representation $\lambda^{\vee}$ of a maximal atypical representation $\lambda$ is 
given by the plot \[ \lambda^\vee(s) = \lambda^c(1-s).\]
\end{prop}

\medskip {\bf Example 1}. Suppose $\lambda = [0,\lambda_2,\ldots,\lambda_n]$ holds with $0> \lambda_2$ and $\lambda_i > \lambda_{i+1} $ for $2 \leq i \leq n-1$. Then $\lambda^{\vee} = [n - \lambda_n - 1, n - \lambda_{n-1} - 1, \ldots, n- \lambda_2 - 1,n-1]$.

\medskip

Dualising is compatible with the normalized block equivalence $\phi_n^i$ of section \ref{signs}.

\begin{lem} \label{phi-dual} For irreducible $i$-atypical $L$ we have $\phi_n^i (L^{\vee}) = \phi_n^i(L)^{\vee}$.
\end{lem}

\medskip{\it Proof}.  If $L$ is $i$-atypical, then $\tilde{\phi}_n^i$ preserves the distances between the sectors, hence $\tilde{\phi}_n^i(L)^{\vee} = Ber^{\ldots} \otimes \tilde{\phi}_n^i(L^{\vee})$. Since we remove $2(n-i)$ symbols from the weight diagram of $L$, we obtain the shift \[ \tilde{\phi}_n^i(L^{\vee}) = Ber^{-2(n-i)} \otimes \tilde{\phi}_n^i(L)^{\vee}.\] Now we calculate for the normalised block equivalence \begin{align*} \phi_n^i(L^{\vee}) & = Ber^{n-k} \tilde{\phi}_n^i(L^{\vee}) = Ber^{n-k} Ber^{-2(n-k)} \tilde{\phi}_n^i(L)^{\vee} \\ & = Ber^{-n+k} \otimes \tilde{\phi}_n^i(L)^{\vee} = (Ber^{n-k})^{\vee} \otimes \phi(L)^{\vee} = (Ber^{n-k} \otimes \tilde{\phi_n^i}(L))^{\vee} \\ & = \phi_n^i(L)^{\vee}.  \end{align*} \qed

\begin{lem} \label{basic}
For maximal atypical irreducible $L=[\lambda_1,...,\lambda_n]$
such that $\lambda_n=0$ the following assertions are equivalent.
\begin{enumerate}
\item $L^\vee \cong [\rho_1,...,\rho_n]$ holds such that $\rho_n\geq 0$.
\item $L$ is basic, i.e. $\lambda_1 \geq \ldots \lambda_n \geq 0$ and $\lambda_i \leq n-i$ holds for all $i=1,...,n$.
\item $\lambda_1 \leq n-1$ and $L^\vee \cong [\lambda^*_1,...,\lambda^*_n]$ holds
for the transposed partition $\lambda^*=(\lambda^*_1,...,\lambda^*_n)$ of the 
partition $\lambda=(\lambda_1,...,\lambda_n)$. 
\end{enumerate}
\end{lem}

\medskip
{\bf Remark}. The number of {\it basic} maximal atypical weights in $X^+(n)$ is equal to the Catalan number $C_n$.

\medskip
{\it Proof}. i) implies ii): If $\rho_n = 0$ the leftmost $\vee$ in the weight diagram of $[\rho]$ is at position $-n+1$. Then the smallest $\wedge$ bound in a cup is at a position $\leq 1$ and $\geq 1-n$. After the change $(I,K) \to (I,I-K)$ and the reflection $s \mapsto 1-s$ this means that the rightmost $\vee$ in $[\rho]^{\vee}$ is at position $\leq n-1$ and $\geq 0$ which is equivalent to $0 \leq \lambda_1 \leq n-1$. Likewise the $i$-th leftmost $\wedge$ bound in a cup is at a position $ \geq -n + i + 1$ and $\leq n$. It will give the $i$-th largest $\vee$ in the weight diagram of $[\lambda]$. After the change $(I,K) \mapsto (I,I-K)$ and the reflection the $i$-th largest $\vee$ is at a position $\leq n- 2i +1$ which is equivalent to $\lambda_i \leq n-i$. ii) implies i): If $\lambda$ is basic the largest $\vee$ is at position $\leq n-1$, hence the largest $\wedge$ bound in a cup is at position $\leq n$. It gives the smallest $\vee$ of $[\lambda]^{\vee}$. Hence the smallest $\vee$ of $[\lambda]^{\vee}$ is 
at a position $\geq 1-n$ which is equivalent to $\lambda_n^{\vee} \geq 0$.

ii) implies iii): If $\lambda$ is basic, the $2n$ vertices in cups form the intervall $J := [-n+1,n]$ of length $2n$. If $J_{\vee}$ is the subset of vertices labelled by $\vee$, the subset $J \setminus J_{\vee}$ is the subset of vertices labelled by $\wedge$. The intervall $J$ is preserved by the reflection $s \mapsto 1-s$. If $\lambda$ is basic, so is $\lambda^*$. We use the following notation: If\[ \lambda_1 = \ldots = \lambda_{s_1} > \lambda_{s_1 + 1} = \ \ldots \  = \lambda_{s_2} > \lambda_{s_2 + 1} = \ldots = \lambda_{s_r } > \lambda_{s_r + 1} = 0 \] put $\delta_1 = s_1$ and $\delta_i  = s_i - s_{i-1}$ and $\Delta_i = \lambda_{s_i}- \lambda_{s_i + 1}$: Likewise for $\lambda^*$ with $\delta_i^*$ and $\Delta_i^*$. Then \[ \delta_i = \Delta_{i-r}^*, \ \Delta_i = \delta_{i-r}^*.\] Then the weight diagram of $ [\lambda^*]$ looks, starting from $n$ and going to the left \[  \ldots \overbrace{\vee \ldots \vee}^{\delta_3^*} \overbrace{\wedge \ldots \wedge}^{\Delta_2^*}  \overbrace{\vee \ldots \vee}^{\delta_2^*}  \overbrace{\wedge \ldots \wedge}^{\Delta_1^*}  \overbrace{\vee \ldots \vee}^{\delta_1^*} \wedge \ldots \wedge\] and the weight diagram of $[\lambda]$ looks, starting from $-n +1$ and going to the right like\[ \vee \ldots \vee \overbrace{\wedge \ldots \wedge}^{\Delta_r = \delta_1^*}  \overbrace{\vee \ldots \vee}^{\delta_r = \Delta_1^*}  \overbrace{\wedge \ldots \wedge}^{\Delta_{r-1} = \delta_2^*}  \overbrace{\vee \ldots \vee}^{\delta_{r-1} = \Delta_2^*}  \ldots.\] The two weight diagrams are mirror images of each other and the rule for the $\vee$'s in cups in one is the same as the rule for the $\wedge$'s in the cups of the other. Hence after the change $(I,K) \mapsto (I,I-K)$ and the reflection $s \mapsto 1-s$ the two weight diagrams agree.
iii) implies i): trivial. \qed

\medskip {\bf Example 3}. Duals in the ${\calR}_3$-case. If $a>b>0$, then  $[a,b,0]^{\vee} = [2,2-b,2-a] = Ber^{2-a} [a,a-b,0]$. If $a \geq 1$ then $[a,a,0]^{\vee} = [2,1-a,1-a] = Ber^{1-a}[a+1,0,0] = Ber^{1-a} S^{a+1}$.

\medskip

{\it A better description.} If $L=L(\lambda)$ is an irreducible maximal atypical representation in $\calR_n$, its weight $\lambda$ is uniquely determined by its plot.
Let $S_1...S_2... S_k$ denote the segments of this plot. Each segment $S_\nu$ has
even cardinality $2r(S_\nu)$, and can be identified up to a translation with a unique basic
weight of rank $r(S_\nu) = r_{\mu}$ and a partition in the sense of lemma \ref{basic}.
For the rest of this section we denote the segment of rank $r(S_\nu)$
attached to the dual partition by $S_\nu^*$, hoping that this will not be confused 
with the contravariant functor $*$. Using this notation, Tannaka duality
maps the plot $S_1 .. S_2 ... S_k$ to the plot $S_k^* ... S^*_2 .. S_1^*$
so that the distances $d_i$ between $S_i$ and $S_{i+1}$ coincide
with the distances between $S^*_{i+1}$ and $S^*_i$. This follows from
 proposition \ref{dual-plot} and determines
the Tannaka dual $L^\vee$ of $L$ up to a Berezin twist. 

\medskip

{\it The dual forest}. If we identify the basic plots with rooted trees $S_i \leftrightarrow \calT_i$, we can describe a weight by a \textit{spaced} forest \[ \calF = (d_0, \calT_1, d_1, \calT_2, \ldots, d_{k-1}, \calT_k)\] where $d_0 = \lambda_n$. We describe the dual in this language. 

\medskip

{\it Grafting}. Given a planar forest $\calF = \calT_1 \ldots \calT_n$ of planar rooted trees, we can introduce a new $n$-ary root and graft the trees $\calT_i$ onto this root. This new tree is called the grafting product $\vee(\calT_1 \ldots \calT_n)$ of $\calT_1 \ldots \calT_n$. The grafting product of the trees in a spaced forest is obtained by forgetting the distances and simply taking the grafting product of the trees.

\medskip {\bf Example}. Consider the forest of two rooted planar trees
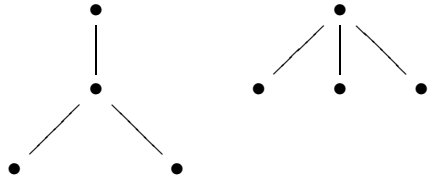
\[ \xymatrix@R-0.2cm@C-0.2cm{ & \bullet \ar@{-}[d] & & & \bullet \ar@{-}[dl] \ar@{-}[d] \ar@{-}[dr] & \\ & \bullet \ar@{-}[dl] \ar@{-}[dr] & & \bullet & \bullet & \bullet  \\ \bullet & & \bullet & & & } \] Grafting this planar forest gives the forest with the single tree

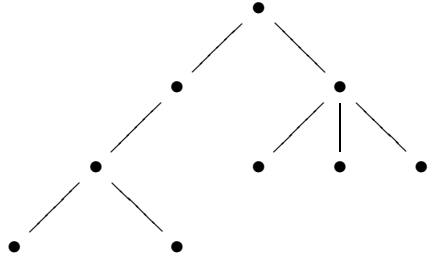
\[ \xymatrix@R-0.2cm@C-0.2cm{ & & & \bullet  \ar@{-}[dr] \ar@{-}[dl] & & \\ & & \bullet \ar@{-}[dl] &  & \bullet \ar@{-}[dl] \ar@{-}[d] \ar@{-}[dr] & \\ & \bullet  \ar@{-}[dl] \ar@{-}[dr] & & \bullet & \bullet & \bullet \\ \bullet & & \bullet & & &} \] 

{\it Mirror tree}. If $\calT$ is a planar rooted tree, then the mirror image $\calT^*$ of $\calT$ along the root axis is recursively defined as follows: Put $(\vee(\emptyset))^* = \vee(\emptyset), \ \emptyset^* = \emptyset$ where $\emptyset$ is the empty tree and extend via \[ (\vee(\calT_1 \ldots \calT_n))^* = \vee(\calT_n^* \ldots \calT_1^*).\]

\medskip {\bf Example}. The mirror image of the grafted planar tree above is 

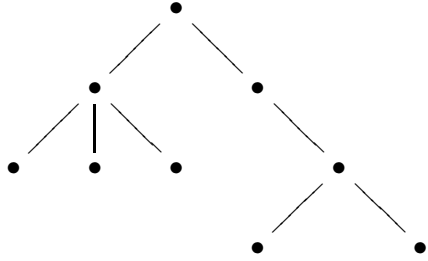
\[ \reflectbox{ \xymatrix@R-0.2cm@C-0.2cm{ & & & \bullet  \ar@{-}[dr] \ar@{-}[dl] & & \\ & & \bullet \ar@{-}[dl] &  & \bullet \ar@{-}[dl] \ar@{-}[d] \ar@{-}[dr] & \\ & \bullet  \ar@{-}[dl] \ar@{-}[dr] & & \bullet & \bullet & \bullet \\ \bullet &  & \bullet & & &} } \] 

\begin{lem} The weight of the dual representation corresponds to the spaced forest \[ \calF^{\vee} = (d_0^*, \calT_k^*, d_1^*, \calT_{k-1}^*, d_2^*, \ldots, d_{k-1}^*, \calT_1^*)\] where $d_i^* := d_{k-i}$ for $i=1,\ldots,k-1$ and $d_0^*  =   - d_0 - d_1 - \ldots - d_{k-1}$ and $\calT_i^{*}$ denotes the mirror image (along the root axis) of the planar tree $\calT_i$.
\end{lem}

{\it Proof}. The claim about the distances $d_1^*,\ldots,d_{k-1}^*$ follows from the description of the dual plots. We first prove the claim about $d_0^*$. Now $d_0^* = (1- b) + n-1$ where $b$ is the last point of the rightmost sector with rank $r_k$ \[ b = \lambda_1 + (2r_k - 1).\] Hence $d_0^* =-b + n=  -\lambda_1 - 2 r_k + n + 1$. 
Now use that $\lambda_1 = (\lambda_n - n+ 1) + 2 r_1 + \ldots + 2r_{k-1} + d_1 + \ldots + d_{k-1}$, hence \begin{align*} d_0^* & = [- \lambda_n + n - 1 - 2(r_1 + \ldots + r_{k-1}) - (d_1+ \ldots + d_{k-1})] - 2r_k + n + 1 \\ & = - \lambda_n  - (d_1+ \ldots + d_{k-1}) = -d_0 - d_1 … - d_{k-1}.\end{align*} 
It remains to prove that if $S_i$ corresponds to $\calT_i$, then the dual plot $S_i^*$ corresponds to the mirrored tree $\calT_i^*$. We induct on the rank of the sector. The case $r_k =1$ is clear. If $[a,b]$ is a sector, then $\lambda(a) = \boxplus$ and $\lambda(b) = \boxminus$. According to proposition \ref{dual-plot} the dual plot is obtained by first exchanging $\boxplus$ and $\boxminus$ and then reflecting $s \mapsto 1-s$. Hence the dual plot of $S_i^*$ is obtained (ignoring distances) by keeping the outer labels $\boxplus$ and $\boxminus$ of the sector and dualising the plot of the inner segment $[a+1,b-1]$. This corresponds to keeping the root of the tree $\calT_i$ and calculating the dual of the forest of the inner trees obtained from $\calT_i$ by removing the root of $\calT_i$. More precisely: The tree to the dual plot $S_i^*$  is obtained by taking the grafting product of the inner subtrees corresponding to the dual of the plot of the inner segment. The interval $[a+1,b-1]$ is a segment consisting of sectors $\tilde{S}_1\ldots \tilde{S_l}$ corresponding to the trees $\tilde{\calT}_1 \ldots\tilde{\calT}_l$. Dualising the inner segment yields by induction the forest $\tilde{\calT}_l^* \ldots \tilde{\calT}_1^*$ since the ranks of inner sectors are smaller than the rank of $S_i$. Hence the tree corresponding to $S_i^*$ is obtained by grafting the forest $\tilde{\calT}_l^* \ldots \tilde{\calT}_1^*$. This is just the definition of the mirror image of $\calT_i$.\qed

\medskip {\bf Example}. Consider the irreducible representation $[11,9,9,5,3,3,3]$ in $\calR_7$. It has sector structure $S_1 = [-3,4]$, $S_2 = [7,10]$ and $S_3 = [11,12]$ with distances $d_0 = 3$, $d_1 = 2$ and $d_2 = 0$. The associated spaced forest is  
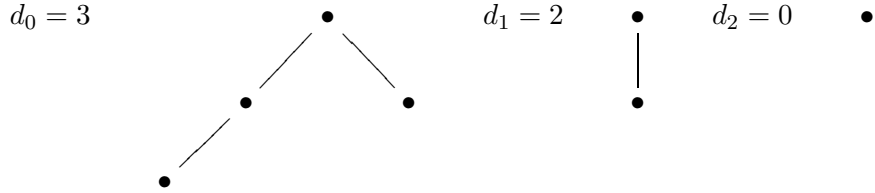
\[ \xymatrix@R-0.2cm@C-0.2cm{ d_0 = 3 & & & \bullet \ar@{-}[dl] \ar@{-}[dr] & & d_1 = 2 & \bullet \ar@{-}[d] & d_2 = 0 & \bullet \\ & & \bullet \ar@{-}[dl] & & \bullet & &  \bullet & & \\ & \bullet  & & & & & & & } \] 

The dual is the representation $[1,1,0,0,-4,-4,-5]$ with sectors (from left to right) $S_3^* = [-11,-10]$, $S_2^* = [-9,-6]$ and $S_1^* = [-3,4]$ with associated spaced forest 
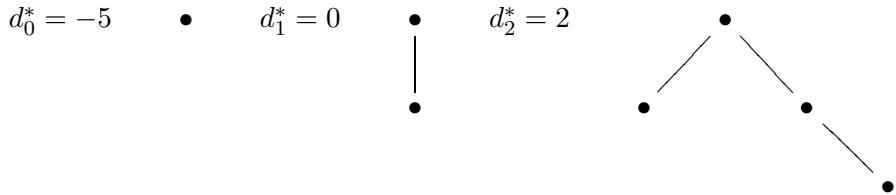
\[ \xymatrix@R-0.2cm@C-0.2cm{ d_0^* = -5 & \bullet & d_1^* = 0 & \bullet \ar@{-}[d] & d_2^* = 2 & & \bullet \ar@{-}[dl]  \ar@{-}[dr] & & \\  & & & \bullet & &  \bullet & & \bullet \ar@{-}[dr] &  \\  & & & & & & & & \bullet } \]



\section{Cohomology I} \label{kohl}

\medskip
In corollary \ref{cohomII} we have seen that in the situation of the melting algorithm
one obtains surjective maps $H^i(p): H^i(\A) \to H^i(L)$ for all $i\in\mathbb Z$. For $K=Ker(p: \A \to L)$
we therefore get exact sequences
$$   0 \to H^i(K) \to H^i(\A) \to H^i(L) \to 0 $$
for all integers $i$. Hence, if in addition $H^i(\A)=0$ and $H^i(L)=0$ vanish for all $i\neq 0$,
then $H^i(K)=0$ holds for all $i\neq 0$. Then $K/L \cong A$ implies $H^i(A)=0$
for $i\neq -1,0$. Suppose, the same conditions are satisfied for $\A^\vee$ as well. Then also $H^i(A^\vee)=0$ holds for $i\neq -1,0$. Then, by duality
$H^i(A)^\vee \cong H^{-i}(A^\vee)$, the cohomology modules $H^i(A)$ vanish for $i\neq 0$.
This proves
 
\begin{prop} \label{cohomvan}
For  irreducible basic modules $V=[\lambda_1,...,\lambda_{n-1},0]$ in $\calR_n$ the cohomology modules $H^i(V)$
vanish for all $i\neq 0$.
\end{prop}

\medskip
{\it Proof}. We use induction with respect to the degree $p = p(\lambda)=\sum_i \lambda_i$, where
$\lambda_i$ for $i=1,..,n$ denote the coefficients of the weight vector. 
By induction assume the assertion holds for all irreducible basic modules of degree $<p$. For $V$ of degree $p$ by the melting algorithm
there exists an irreducible basic module $L$ of degree $p-1$ and $\A$ with layer structure
$(L,A,L)$ such that $A=V \oplus A'$, where $A'$ is a direct sum of irreducible basic modules of degree $<p$. Since $H^i(L)=0$ for $i\neq 0$,
$H^i(\A)=0$ for $i\neq 0$ now follows from lemma \ref{tildeLII}. The same applies for the dual
modules $\A^\vee$ and $L^\vee$. Indeed the dual module of a basic irreducible module is  basic irreducible again with the same
degree $<p-1$ (lemma \ref{basic}) using $\sum_i \lambda_i = \sum_i \lambda^*_i$. Hence the remarks preceding proposition
\ref{cohomvan} imply $H^i(A)=0$ for $i\neq 0$. Since $V$ is a direct summand of $A$, this proves our assertion. \qed

\bigskip



\section{Cohomology II}\label{kohl-2}

\medskip

We calculate the $\Z$-grading of $DS(L)$ for maximal atypical irreducible $L$. The case of general $L$ is treated in section \ref{koh3}.

\begin{prop}\label{hproof}
For maximal atypical irreducible $L(\lambda)$ in $\calR_n$ with
weight $\lambda$, normalized so that $\lambda_n=0$, suppose
$\lambda$ has sectors $S_1,.., S_i,..,S_k$ (from left to right). Then the constituents
$L(\lambda_i)$ of $DS(L(\lambda))$ for $i=1,...,k$ have sectors $S_1,.., \partial S_i , .. S_k$, and the cohomology
of $L(\lambda)$ can be expressed in terms of the added distances $\delta_1,...,\delta_{k}$ between these
sectors as follows:
$$  \fbox{$ H^{\bullet}(L(\lambda)) \ =  \ \bigoplus_{i=1}^k \ L(\lambda_i)\langle -\delta_i \rangle $} \ .$$
\end{prop}

\medskip
\textbf{Example.} We know by the main theorem that $DS([6,4,4,1]) = \Pi [3,3,0] \oplus \Pi[6,4,0] \oplus \Pi[6,4,4]$. The proposition above tells us the $\Z$-grading using \[DS(V)
 = \bigoplus_{\ell \in \Z} \ \Pi^\ell (H^\ell(V)).\] In this example $d_0 = 1, \ d_1 = 2$ and $d_2 = 0$. The summand $L(\lambda_i)$ is obtained by differentiating the $i$-th sector in the plot associated to $\lambda$, hence $L(\lambda_1) = [6,4,4], \ L(\lambda_2) = [6,4,0]$ and $L(\lambda_3) = [3,3,0]$. We obtain \[ H^{\bullet}([6,4,4,1]) = [6,4,4]\langle-1\rangle \oplus [6,4,0]\langle-3\rangle \oplus [3,3,0]\langle-3\rangle.\] 

\medskip
{\it Proof}. 
In the special case where all distances vanish $d_1=\cdots =d_k=0$, i.e. in case where the plot of $\lambda$ has only one segment, the assertion of the proposition has been shown
in proposition \ref{cohomvan}. We then prove the general case of nonvanishing distances
by induction with respect to ($n$ and) the lexicographic ordering
used for algorithm I. This means: We prove proposition \ref{hproof} recursively for $L^{up}$,
thereby assuming that we already know the cohomology degrees
of $L^{down}_r$, $L$ and $L^{aux}$ (using the notations of algorithm I). 
First recall the notations used for algorithm I:
$$ L = (S_1\cdots S_{j-1}) \leftarrow \text{distance } d_j \to (S_jS_{j+1} \cdots S_k) \leftarrow \text{ distance }  d_k \ \to ... $$
$$ L^{up} =  (S_1\cdots S_{j-1}) \leftarrow \text{dist. } (d_j + 1) \to \int(S_{j+1} \cdots S_k) \leftarrow  \text{dist. }  (d_k-1) \to ... $$
for a sector $S_j$ with $r(S_j)=1$ supported at $i\in \mathbb Z$.
Recall $\A = (L,A,L)$ with $$ A \ = \ L^{up}\ \oplus\ \bigoplus_{r=1}^k L^{down}_r \ .$$ Furthermore $DS(L)= \Pi^{m_{aux}} L^{aux} \oplus
\bigoplus_{\mu=1}^s \Pi^{m_{\mu}} \tilde L_\mu$ for $DS(\A)= \bigoplus_{\mu=1}^s \tilde\A_\mu$
and $\tilde\A_\mu=(\tilde L_\mu,\tilde A_\mu,\tilde L_\mu)$ such that
the derivative $d(A)$ of $A$ is \[ d(A) = \tilde A + 2(-1)^{i+n-1} L^{aux}\] in $K_0(\calR_n)$. 
Obviously $DS(L^{up})$ has the summands
$$ \bigoplus_{\nu=1}^{j-1} (S_1 \cdots \partial S_\nu \cdots S_{j-1}) ... (d_{j-1}+1) ... S 
... (d_k -1) ... S_{k+1} \cdots) $$
$$ \bigoplus_{\nu >k} (S_1 \cdots  S_{j-1}) ... (d_{j-1}+1) ... S ... (d_k -1) ...S_{k+1} \cdots \partial S_\nu \cdots ) $$
and
$$ L^{aux} = 
(S_1 \cdots  S_{j-1}) ... (d_{j-1}+2) ... (S_{j+1} \cdots S_k) ... (d_k) ... S_{k+1} \cdots  \ .$$
This immediately implies the next

\begin{lem} \label{lastl} The following holds
\begin{enumerate}
\item $DS(L^{up}) \subseteq L^{aux} \oplus DS(L)^{up}$.
\item None of the summands of $DS(L^{up})$ different from $L^{aux}$
is contained in $DS(L)$.  
\item $L^{aux}$ is not a
summand of   $\bigoplus_\mu \tilde A_\mu$.
\end{enumerate}
\end{lem}

\medskip
{\it Proof}. The last assertion holds, since
the constituents of $\bigoplus_\mu \tilde A_\mu$ are obtained from $\tilde L_\mu$
by moves. It can be checked that $L^{aux}$ can not be realized
in this way.  \qed

\medskip
{\it The $d_{j-1}\pm 1$ alternative}. 
By the induction assumption $H^\bullet(L)$ contains $L^{aux}$
with multiplicity 1, and $L^{aux}$ appears in cohomology at the degree
$d_{j-1}$. To determine
$  H^i(L^{up}) \subseteq H^i(A) $ we may use step 11) of the proof of theorem \ref{3}. It 
 easily implies by a small modification of the arguments that 
$$ H^i(A)  = 
\bigoplus_{m_\mu =-i} \tilde A_\mu \oplus H^{i-1}(L)/H^{i-1}(\bigoplus_\mu\tilde L_\mu) \oplus Kern(H^{i+1}(L) \to H^{i+1}(K)) \ .$$ 
Since $ H^{\bullet}(L)/(\bigoplus_{\mu}\tilde L_\mu) \cong  L^{aux}$ by lemma \ref{tildeL}
and since 
\[ H^{i-1}(L)/(\bigoplus_{\mu=1-i}\tilde L_\mu) \ \cong \ L^{aux}\ \] 
for $i-1=d_{j-1}$ by the induction assumption, we get
$$ Kern(H^{\bullet}(L) \to H^{\bullet}(K)) =  L^{aux} \ ,$$
and this implies $$ Kern(H^{i+1}(L) \to H^{i+1}(K)) =  L^{aux} \ $$
for $i+1=d_{j-1}$. In other words  $DS(A) = \tilde A + 2\cdot L^{aux}$
and the two copies of $L^{aux}$ occur in the two possible cohomology degrees
$$ d_{j-1} \pm 1\ .$$
 
\medskip
{\it Continuation of the proof for proposition \ref{hproof}}. 
By lemma \ref{lastl}   the cohomology degree of the constituents
of $H^\bullet(L^{up})$ that appear in \[ \tilde A_\mu \subseteq \bigoplus_{m_\mu =-i} \tilde A_\mu \] can be immediately read of from the degrees $m_\mu$, i.e. from the cohomology degrees
of $\tilde L_\mu$ in $H^\bullet(L)$. These degrees are known by the induction assumption. This
easily proves proposition \ref{hproof} for all constituents $L(\lambda_i)$ of $H^\bullet(L^{up})$
that are not isomorphic to $L^{aux}$. Indeed, according to our claim the cohomological degrees 
for the constituents $L(\lambda_i)\not\cong L^{aux}$ of $H^\bullet(L^{up})$ are given by
$$ 0,\cdots , 0, d_{j-1}+1,d_{j-1}+d_k, \cdots     \  ,    $$
and the summand $L^{aux}$ should occur in degree $d_{j-1}+1$.
The cohomology of $H^\bullet(L)$ on the other hand is concentrated
in the degrees 
$$ 0,\cdots , 0, d_{j-1},d_{j-1}+d_k, \cdots                         \ $$
with the summand $L^{aux}$ corresponding to degree $d_{j-1}$.
All summands $\not\cong L^{aux}$ precisely match, so this proves proposition \ref{hproof}
for all constituents of $H^\bullet(L^{up})$ except for $L^{aux}$.

\medskip
It remains to determine the cohomology degree of $L^{aux} \subseteq H^\bullet(L^{up})$.
As already explained, the summand $L^{aux}$  occurs in degree $d_{j-1}-1$ or $d_{j-1}+1$. So to show that $L^{aux}$ 
occurs in $H^\bullet(L^{up})$ for degree $d_{j-1}+1$,  it now suffices by the $d_{j-1}\pm 1$ alternative to show that $L^{aux}$ occurs in $H^\nu(\bigoplus_r L^{down}_r)$ in the degree $\nu=d_{j-1}-1$. 
Indeed $L^{aux}$ appears 
in $DS(L^{down}) = \bigoplus_\nu H^\nu(L^{down})$ for $L^{down} := L^{down}_j$. 
This follows from the structure of the sectors of
$$ L^{down} = \int(S_1\cdots S_{j-1}) \boxminus... (d_{j-1}-1) ...  \boxminus_{i+1} S_{j+1} \cdots S_k \ $$
and the induction assumption. It gives the  degree $d_{j-1}-1$, for $d_{j-1}\geq 1$, respectively in degree $d_{j-1}-1=-1$, for $d_{j-1}=0$, for the summand $L^{aux}$ in $H^\bullet(L^{down})$.
Hence $$ \fbox{$ L^{aux} \subseteq H^{d_{j-1} +1}(L^{up}) $}\ ,$$
which completes the proof of proposition \ref{hproof}. 
\qed 

\bigskip\noindent



\section{Cohomology III}\label{koh3}

The cohomology of an $i$-atypical $L$ can be calculated in the same way using the normalised block equivalence $\phi_n^i$ of section \ref{signs}. We call an irreducible module $L$ of atypicality $i$ $\phi$-basic if $\phi_n^i (L)$ is basic in $\calR_i$. These will replace the basic modules in the proof of proposition \ref{cohomvan}. The unique mixed tensor in a block of atypicality $i$ replaces the trivial representation.

\begin{prop} \label{cohomvan-2}
For  irreducible $\phi$-basic modules $V$ in $\calR_n$ the cohomology modules $H^i(V)$
vanish for all $i\neq 0$.
\end{prop}

\medskip
{\it Proof}. The remarks preceeding proposition \ref{cohomvan} are valid. By lemma \ref{mixed-tensor-derivative} the cohomology of the mixed tensor $L(\lambda)$ is concentrated in one degree, and by lemma \ref{stable} this degree is zero since $\lambda_n = 0$. Since $\phi_n^i(L(\lambda)) = \one$, we induct as in the proof of \ref{cohomvan} on the sum $p=\sum_i \lambda_i$ of the coefficients of $\phi_n^i(L)$. The rest of the proof works verbatim. Note that the dual of a $\phi$-basic module is $\phi$-basic again of the same degree using $\phi_n^i(L)^{\vee} = \phi_n^i(L^{\vee})$ of lemma \ref{phi-dual} and lemma \ref{basic}. \qed

\medskip
We can now copy the proof of proposition \ref{hproof} to obtain the next statement. Here the added distances $\delta_i$ are the distances in the plot $\phi(\lambda)$ associated to $\lambda$ in section \ref{sec:loewy-length}.

\begin{prop}\label{hproof-2}
For irreducible $L(\lambda)$ in $\calR_n$ with
weight $\lambda$, normalized so that $\phi_n^i(L(\lambda)) = [\lambda_1^{\phi}, \ldots, \lambda_i^{\phi}]$ satisfies $\lambda_i^{\phi}=0$, suppose $\lambda$ has sectors $S_1,.., S_j,..,S_k$ (from left to right). Then the constituents
$L(\lambda_j)$ of $DS(L(\lambda))$ for $j=1,...,k$ have sectors $S_1,.., \partial S_j , .. S_k$, and the cohomology
of $L(\lambda)$ can be expressed in terms of the added distances $\delta_1,...,\delta_{k}$ between these
sectors as follows:
$$  \fbox{$ H^{\bullet}(L(\lambda)) \ =  \ \bigoplus_{j=1}^k \ L(\lambda_j)\langle -\delta_j \rangle $} \ .$$
\end{prop}

\medskip\noindent



\section{The forest formula} \label{sec:forest}


\medskip
Recall the functor $DS_{n,0}: T_n \to T_0=svec_k$ with its decomposition 
$DS_{n,0}(V) = \bigoplus_{\ell\in\mathbb Z} D_{n,0}^\ell(V)[-\ell]$ for objects $V$ in $T_n$ and
objects $D_{n,0}^\ell(V)$ in $svec_k$.
For $V\in T_n$ we define the Laurent polynomial
$$ \omega(V,t) = \sum_{\ell\in\mathbb Z} sdim(D_{n,0}^\ell(V)) \cdot t^\ell \ $$
 as the Hilbert polynomial of the graded module $DS^\bullet_{n,0}(V)= \bigoplus_{\ell\in\mathbb Z} D_{n,0}^\ell(V)$.  Since $sdim(W[-\ell])=(-1)^\ell sdim(W)$ and $V= \bigoplus D_{n,0}^\ell(V)[-\ell]$ holds, the formula  $$sdim(V) = \omega(V,-1)$$ follows. For $V= Ber_n^i$ 
$$    \omega(Ber_n^i,t) \ = \ t^{ni}  \ .$$
Indeed, $H^\ell(Ber_n^i) =0$ for $\ell\neq i$ 
and $H^\ell(Ber_n^i) =Ber_{n-1}^i$ for $\ell=i$ implies $DS(Ber_n^i) = Ber_{n-1}^i[-i]$.
If we apply this formula $n$-times and consider $B_0={\bf 1}$, we obtain 
$DS_{n,0}(Ber_n^i)=DS^n(Ber_n^i) ={\bf 1}[-ni]$  from the fact that $DS_{n,0}(L)=DS^n(L)$ holds for simple objects $L$. This implies $DS_{n,0}^{ni}(Ber_n^i)=\bf 1$ and that $DS_{n,0}^{\ell}(Ber_n^i)$ is zero otherwise.

\medskip
Since $DS_{n,0}$ is a tensor functor, 
$\omega(M\otimes L,t)= \omega(M,t)\omega(L,t)$ holds. Hence
$$ \omega(Ber_n^i\otimes L,t) \ =  \ t^{ni} \cdot \omega(L,t) \ .$$
Similar as in the proof of  lemma \ref{-ell} one shows $$\omega(V^\vee,t) = \omega(V,t^{-1})\ .$$

\medskip
Let now $L=L(\lambda)$ be a maximal atypical irreducible representation in $\calR_n$. Associated to its plot
$\lambda$ we have the basic plot $\lambda_{basic}$ and the numbers $d_0,\ldots,d_{k-1}$.
Furthermore, let $S_1S_2\cdots S_k$ be the sector structure of $\lambda_{basic}$. 
For the degrees $r_i=r(S_i)$ we define the number 
$$ D(\lambda) \ =\ \sum_{i=1}^k r_i \sum_{0 \leq j<i} d_j = \sum_{i=1}^k r_i \delta_i\ ,$$
Recall 
$\delta_i=\sum_{\nu=0}^{i-1} d_{\nu}$ implies 
$\delta_1\leq \delta_2 \leq \cdots \leq \delta_k$ and $\delta_i\in \mathbb Z^k$.
Consider the vector $D$ with coordinates $\delta_1, \ldots,\delta_k$. 
Together with $\lambda_{basic}$ the knowledge of $D$ determines $\lambda$.
For simplicity, we express this by writing $\lambda = D \times \lambda_{basic}$ in the following argument. With this notation, our proposition \ref{hproof} gives for  $DS(L)$ the following element in the Grothendieck group $K_0(\calR_{n-1})\otimes
k[t]$ 
$$ DS\bigl( \begin{pmatrix} \delta_1  \cr . \cr \delta_{i-1} \cr \delta_i \cr
\delta_{i+1} \cr . \cr \delta_k  \end{pmatrix} \times (S_1\cdots S_k)_{basic}\Bigr) \ = \
\sum_{i=1}^k   t^{\delta_i} \cdot 
\begin{pmatrix} \delta_1 -1 \cr . \cr \delta_{i-1} - 1\cr \delta_i \cr
\delta_{i+1}+1 \cr . \cr \delta_k +1 \end{pmatrix} \times (S_1 \cdots \partial S_i \cdots S_k)_{basic} \ $$
where formally (and without loss of information) we  replace the shifts $[-\nu]$ by $t^\nu$. In the following, we refer to
this formula as the {\it key formula}.
Now $\partial S_i$ may introduce new sectors in $(S_1 \cdots \partial S_i \cdots S_k)_{basic}$. So if we want  to treat everything on an equal footing, we better
count each sectors $S_i$ with the multiplicity $r_i$. This amounts to consider instead of the
vector $D$ the new refined vector $\delta$ in $\mathbb Z^n$ with the coordinates
\[ (\underbrace{\delta_1, \ldots,\delta_1}_{r_1}, \ldots, \underbrace{\delta_k, \ldots, \delta_k}_{r_k}).\] Then the number $D(\lambda)$ defined above is just the sum of the coordinates of this vector. 
With this new vector we have an analogous formula expressing $DS$ as above, where for the $i$-th summand on the right side one of the entries $\delta_i$ of $\delta$ has to be removed to obtain a vector in $\mathbb Z^{n-1}$. The right side is now of the correct form to enable the application of the formula for $DS$ to the right side again. Inductively, after $n$ steps this gives a complicated expression
with at most $n!$ summands. The number of summands depends only on $\lambda_{basic}$. Since the additional monomial term in $t$ obtained from each derivative is of the 
form $t^{\delta_\nu \pm s_\nu}$, for some shift factors $s_\nu$ not depending on $\delta$, and since in each summand all coordinate entries of $\delta$
will be finitely successively deleted after $n$ times applying $DS$, this vector disappears 
and each of these summands has the form
$$   t^{\sum_{i=1}^k r_i\delta_i} \cdot P(t) \times \emptyset $$
for a certain Laurent polynomials $P(t)$ that depends on the specific summand and on $\lambda_{basic}$, but that does not depend
on the coefficients $\delta_1,...,\delta_k$. If we compare with the case $\delta_1\!=\! ...\! = \! \delta_k\! =\! 0$,
we therefore obtain the following {\it translation formula}: 
$$   \omega(L(\lambda),t) \ = \ t^{D(\lambda)} \cdot \omega(L(\lambda_{basic}),t) \ .$$
This being said, we use that the basic plots of rank $n$ are in 1-1 correspondence with
planar forests $\calF$ with $n$ nodes $x\in \calF$ as in sections \ref{sec:main} \ref{duals} and \cite{Weissauer-gl}.
For a planar forest, let $\#{\calF}$ denote the number of its nodes.
We visual each of the trees in a plain forest top down, i.e. with their root on the top of the tree. 
Then, for each node $x\in \calF$ let
 $\calF(x)$ denote the subtree of  the tree containing $x$ with all nodes removed that are not below the node $x$. In this way the node $x$ becomes the root of the tree $\calF(x)$ by definition. 
For a forest $\calF$ we recursively define  the quantum forest factorials
$$   [\calF]_t ! = \prod_{x \in \calF} [\#\calF(x)]_t  \ \in\ \mathbb Z[t] $$ 
using the following abbreviations: For the number $m=\#\calF(x)$ of nodes in $\calF(x)$ 
we define the quantum numbers
$$   [m]_t \ := \ \frac{t^m - t^{-m}}{t - t^{-1}} \ .$$
Clearly $[m]_t = (-1)^{m-1}[m]_{-t} = [m]_{t^{-1}}$. Obviously
the tree factorial $\calT !$ of section \ref{sec:main} equals $[\calT]_1!$. 
Example: For the forest $\calF$ that contains only one linear tree, the forest factorial $[\calF]_t !$  
specializes to the quantum factorial $[n]_t! = \prod_{m=1}^n [m]_t$. 
For a planar forest $\calF$, given as the union of trees $\calT_i$ for $i=1,\ldots,k$
with $r_i$ nodes respectively,
one has $[\calF]_t ! = \prod_{i=1}^k [\calT_{i}]_t ! $ and hence
$$(*) \quad \quad \frac{[\#\calF]_t!}{[\calF]_t !} = \frac{ [\sum_i r_i]_t! }{ [r_1]_t ! \cdots [r_k]_t!} \cdot \prod_{i=1}^k 
\frac{ [\#\calT_i]_t ! }{[\calT_i]_t ! }\ .$$
Observe, for a tree $\calT$ the value $\frac{[\#\calT]_t!}{[\calT]_t !}$ does not change
under grafting, i.e. replacing $\calT$ by a new tree with $\#\calT +1$ nodes by putting a new root on top. Similar $\frac{[\#\calF]_t!}{[\calF]_t !}$ does not change under the grafting of the planar forest ${\calF}$, that replaces $\calF$ by a forest with a single tree with $\#\calT +1$ nodes obtain by putting a new root on top of all trees connected to the old roots of the trees of $\calF$.


\begin{lem} \label{thm:forest-formula} For irreducible maximal atypical 
representations $L=L(\lambda)$ in $\calR_n$
we have the forest formula
$$   \fbox{$ \omega(L,t) \ =\  t^{D(\lambda)} \cdot \frac{ [n]_t !}{[\lambda_{basic}]_t !} $}\ $$
where $\lambda_{basic}$ is viewed as the planar forest associated  to $L$. 
\end{lem}

\medskip
{\it Proof}.  From the translation formula we may assume $\lambda=\lambda_{basic}$. Let us first consider the simple case of basic representations $L$, where
all sectors $S_i$ for $i=1,..,k$ are intervals $I_i=[a_i,a_i +2r_i-1]$
where the support of the plot $S_i$ is $[a_i,...,a_i+r_i -1]$.
The corresponding $\omega(L,t)$ then only depends
on the ranks $r_1,...,r_k$ of the sectors $S_1,...,S_k$, hence will be denoted $\omega_{r_1,...,r_k}(t)$ in the following. From the key formula then, for general basic
$\lambda$ with sector structure $S_1 \cdots S_k$ and $r_i=r(S_i)$, similarly to the translation formula we easily obtain
the following generalized {\it Leibniz formula} 
$$  \omega(L(\lambda),t ) = \omega_{r_1,...,r_k}(t) \cdot \prod_{i=1}^k \ \omega(L(S_i),t) \ ,$$
where $L(S_i)$ denotes the irreducible basic maximal atypical 
representation in $\calR_{r_i}$ whose plot is $S_i$ (up to a translation on the number line).
Now each of these $S_i$ has a unique sector. For basic plots $S$ 
with a unique sector (like the $S_i$) the key formula obviously implies the following {\it grafting formula}
$$ \omega(L(S_i),t) \ = \ \omega(L(\partial S),t) \ ,$$
where $L(\partial S)$ denotes the unique maximal atypical basic representation  in $\calR_{r(S) -1}$
whose plot is $\partial S$. So the forest attached to $L(S)$ is obtained by grafting the forest
of $\partial S$.

\medskip
It is clear that inductively the translation formula, the generalized Leibniz formula and the grafting
formula determine the Laurent polynomials $\omega(L,t)$ for irreducible maximal atypical $L \in {\calR}_n$ uniquely. Hence for the proof it suffices  that the expression on the right side of the identity stated in lemma \ref{thm:forest-formula} satisfies the analogous formulas and that it holds for $n=1$. Indeed, our assertion is obvious for $n=1$. 
The translation formula and the grafting formula for the right side are also obvious. To check the generalized Leibniz formula, by the formula (*) from above it suffices to prove 
$$ \omega_{r_1,...,r_k}(t) = \frac{ [\sum_i r_i]_t !}{ \prod_{i=1}^k [r_i]_t! } \ .$$
To this end it is helpful
that for $a=r_1+\cdots + r_i$ and $b=r_{i+1} + \cdots + r_k$, by the key formula, also the following version
of the generalized Leibniz formula holds
$$  \omega_{r_1,...,r_k}(t) \ = \ \omega_{a,b}(t) \omega_{r_1,...,r_i}(t) \omega_{r_{i+1},...,r_k}(t) \ .$$ 
So it suffices to verify that $\omega_{a,b}(t) [a]_t! [b]_t! = [a+b]_t !$, which
finally is proved by induction on $n=a+b$.  For this notice that the key formula immediately
implies the following generalized {\it Pascal rule} 
$$\omega_{a,b}(t) = t^b\cdot \omega_{a-1,b}(t) + t^{-a}\cdot \omega_{a,b-1}(t)   $$
for the generalized binomial coefficients $\omega_{a,b}(t)$. Indeed, the derivative of the two sectors $S_1S_2$ give $\partial S_1 S_2$ with
$d_0=0$, $d_1=1$ respectively $S_1 \partial S_2$ with $d_0=-1$, $d_1=1$. 
Hence $D(\partial S_1 S_2)= 0\cdot (a-1) + 1 \cdot b=b$ and 
$D(S_1 \partial  S_2)= -a + 0\cdot (b-1) =-a$. Hence using the induction assumption, we already know 
$\omega_{a-1,b}(t) [a-1]_t! [b]_t! = [a+b-1]_t !$ and $\omega_{a,b-1}(t) [a]_t! [b-1]_t! = [a+b-1]_t !$. 
Hence the proof of the induction step finally amounts  for the quantum numbers $[m]_t$ to the following generalized {\it additivity} 
 $$ [a+b]_t = t^b \cdot [a]_t + t^{-a} \cdot [b]_t \ $$ 
 that is easily verified. This completes the proof. \qed

\bigskip
Since $[n]_t!$ and $[\lambda_{basic}]_t!$ are products of certain Laurent polynomials
$[m]_t$ for integers $m$, the forest formula implies $\omega(L,-t) = \pm\ \omega(L,t)$.
The forest formula also gives $\omega(L,1) >0$. Since $\omega(L,-1)=sdim(L)$, hence  $\omega(L,-t) = sign(sdim(L))\cdot \omega(L,t)$. Recall that for irreducible $L$ in $\calR_n$ we defined a sign $\varepsilon(L)$ and that the sign of $sdim(L)$ is $\varepsilon(L)$,
as shown in section \ref{sec:main} and also in \cite{Weissauer-gl}. Finally, the forest formula also implies 
$\omega(L,t^{-1}) /\omega(L,t)= t^{-2D(\lambda)}$ for $L=L(\lambda)$. Hence we obtain

\begin{lem}\label{ev} For irreducible (maximal atypical)
$L=L(\lambda)$ in $\calR_n$ one has the formulas $\omega(L^\vee,t)=\omega(L,t^{-1}) = t^{-2D(\lambda)} \omega(L,t)$ and $$\omega(L,-t) = \varepsilon(L) \cdot \omega(L,t)\ .$$
\end{lem}

\medskip
{\bf Example}. For $S^{n-1+d}$ in $\calR_n$ and for integers $d\geq 0$ 
$$  \omega(S^{n-1+d},t) =  t^{d-n+1}  + t^{d-n+3} + \cdots + t^{d+n-1} = t^d\cdot \omega(S^{n-1},t) \ .$$

\begin{lem}\label{maxde} For irreducible max. atypical representations
$L=L(\lambda)$ in $\calR_n$  the Laurent polynomial $\omega(L,t)$
has degree $p(\lambda)=\sum_{i=1}^n \lambda_i $ in the sense that 
$$ \omega(L,t)  = t^{p(\lambda)} + \sum_{\ell<p(\lambda)} a_\ell \cdot t^\ell  \ .$$  
\end{lem}

{\it Proof}. 
This follows from the  key formula. 
Indeed its $i$-th summand gives rise to shifts by $\delta_{i+1} +1, ... , \delta_k +1$. To determine the highest cohomology degree of $L(\lambda)$ one has to look for the maximal contributions from all these shifts. Each time we apply $DS$, the maximal contribution is obtained from the first summand  $i=1$. Hence the highest $t$-power arises from the first summands of the key formula each times we apply $DS$ ($n$-times), 
in other words by applying the derivative $\partial$ each time to the leftmost sector. In particular the highest
$t$-power of $\omega(L(\lambda),t)$ is $t^{\delta_1}$ times the highest $t$-power of $\omega(\overline L,t)$ for the representation $\overline L\in T_{n-1}$ associated to the plot $(\partial S_1)S_2\cdots S_k$ with the new vector $\delta=(\delta_1,...,\delta_1, \delta_2 +1,\cdots, \delta_k+1)$ with one copy of $\delta_1$ deleted. Now it is not hard to see, by unraveling the weight associated to this spaced forest, that the associated representation $\overline L$ is the highest weight module $L(\overline \lambda)$ in the sense
of Lemma \ref{stable}. Therefore, the highest $t$-power  of $\omega(L,t)$ is $t^{\lambda_n}$ times the
highest $t$-power of $\omega(\overline L,t)$. In other words $deg_t(\omega(L,t))= \lambda_n + deg_t(\omega(\overline L,t))$. By induction on $n$ hence
$deg_t(\omega(L,t))= \lambda_n + p(\overline \lambda) = p(\lambda)$.
\qed

\medskip
Therefore the forest formula implies 

\begin{cor}\label{corforest}
For irreducible maximal atypical representations
$L=L(\lambda)$ in $\calR_n$ one has the formula $D(\lambda)= p(\lambda) - p(\lambda_{basic})$.
\end{cor}

\medskip
By lemma \ref{maxde} and \ref{ev} furthermore

\begin{cor} \label{lem:top-degree} For irreducible maximal atypical representations
$L=L(\lambda)$ in $\calR_n$ 
$$ \omega(L,t)  = t^{p(\lambda)} + \sum_{q(\lambda) <\ell<p(\lambda)} a_\ell \cdot t^\ell \ \ + \ t^{q(\lambda)}  \ $$  
holds for
$p(\lambda) - q(\lambda) = p(\lambda) + p(\lambda^\vee) =
2\cdot p(\lambda_{basic})$.
\end{cor}

\medskip
{\it Proof}.  From the forest formula and lemma \ref{ev} we obtain
$$ q(\lambda)= D(\lambda) + q(\lambda_{basic}) = D(\lambda) - p(\lambda_{basic})\ .$$
Hence $p(\lambda) - q(\lambda) =
2 p(\lambda_{basic})$. Since $\omega(L^\vee,t)= \omega(L,t^{-1})$, we obtain $p(\lambda^\vee)= -q(\lambda)= - D(\lambda) + p(\lambda_{basic})$. Combined with corollary \ref{corforest} this last formula gives $ p(\lambda) + p(\lambda^\vee) = 2 p(\lambda_{basic})$. \qed





\section{$I$-module structure on the cohomology $H^\bullet_{DS_n}$} \label{sec:chevalley-eilenberg}

In this section we show that the cohomology of the operator $DS_{n,0}$ is a graded module under the invariant algebra $I = \Lambda^{\bullet}(\mathfrak{p}_{-1})^H$ defined below. As an application we compute the cohomology and the Hilbert polynomial of a maximal atypical Kac module $V(\lambda)$ for the operator $DS_{n,0}$. 
 We also show that the projection of $V(\lambda)$ to $L(\lambda)$ induces a map on the $DS_{n,0}$-cohomology which vanishes except in the top degree $p(\lambda)$. Note that it does not make sense to consider the Hilbert polynomial $\omega(V,t)$ for $V = V(\lambda)$ and the Dirac operator $D$ since any Kac module is in the kernel of $H_D$.

\medskip\noindent
The tensor functor associated to an element $x$ in $X = \{x \in \g_1 \ | \ [x,x] = 0\}$ only depends by \cite{Duflo-Serganova} on the $G_0$-orbit on $X$. We therefore work in this section with the operator $DS_n$ associated to the action of the element \[ \mathbb{D} = \begin{pmatrix} 0 & id_n \\ 0 & 0 \end{pmatrix} \] which is clearly in the same $G_0$-orbit as our usual choice of $x \in \g_1$ with $1$'s in the anti-diagonal. It defines a (graded) tensor functor $DS_n\!:\! T_n \to T_0$ which is isomorphic to $DS_{n,0}$. 


\medskip\noindent
{\it Notations and conventions}. \begin{itemize} 
\item Let $H$ denote $Gl(n)$ diagonally embedded into $G_0\! =\! Gl(n) \times Gl(n)$ via $g \!\mapsto\! diag(g,g) \in G_0$. Then $Lie(H) \cong \mathfrak{gl}(n)$.

\item We consider the subalgebra $\p \subset \g$ with grading $\p = \p_{-} \oplus \p_0 \oplus \p_+$ \begin{align*} \p_{-} & = \{ \begin{pmatrix} 0 & 0 \\ x & 0 \end{pmatrix} \ | \ x \in \mathfrak{gl}(n) \} \\  \p_0 & = \{ \begin{pmatrix} x & 0 \\ 0 & x \end{pmatrix} \ | \ x \in \mathfrak{gl}(n) \} = Lie(H) \\  \p_+ & = \{ \begin{pmatrix} 0 & id_n \\ 0 & 0 \end{pmatrix}  \}. \end{align*} 

\item Recall that the restriction of the Kac module $V(\lambda)$ to $\frak p$ is given by \[ V(\lambda)\vert_{\frak p} = \Lambda^{\bullet}(\p_{-}) \otimes L_0(\lambda) \] where $L_0(\lambda)$ is the irreducible $\g_0$-module $L_0(\lambda)$ trivially extended to the parabolic subalgebra of upper triangular block matrices.

\item We write in this section $\rho^{\vee} \boxtimes \rho$ for the irreducible representation $L_0(\lambda)$ of $\g_0$ which is given by the external tensor product of the irreducible $\gl(n)$-representation $\rho^{\vee}$ of weight $(\lambda_1,\ldots,\lambda_n)$ with its dual of weight $(-\lambda_n, \ldots, - \lambda_1)$. If viewed as a representation of $H \subset G_0$ this becomes $\rho^{\vee} \otimes \rho \cong End(\rho)$. In this notation $V(\lambda) = V(\rho^{\vee} \boxtimes \rho)$ and $L(\lambda) = L(\rho^{\vee} \boxtimes \rho)$. 

\item The tensor product $\rho^{\vee} \boxtimes \rho$ contains the trivial representation with multiplicity 1. We call a vector in this subspace an $H$-\textit{spherical} vector. In this sense $L_0(\lambda)$ has an $H$-spherical vector if and only if $\lambda$ is maximal atypical.

\end{itemize}

The action of the generator $\mathbb{D}$ of $\frak p_+$ on $\frak p$-modules induces the operator denoted
$DS_n$ in the following. In particular $\mathbb{D}$ acts on ${\frak p}$ and the ideal ${\frak p}_- \oplus {\frak p}_0$ via the adjoint representation. Notice $Lie(H)$ acts on ${\frak p}_0 \cong Lie(H)$ by the adjoint representation of $Lie(H)$
and on ${\frak p}_{-}$ by the adjoint action of ${\frak p}$ such that the map
$$    \begin{pmatrix} 0 & 0  \\ x & 0 \end{pmatrix}  \ \mapsto \  \begin{pmatrix} x & 0 \\ 0 & x \end{pmatrix} $$ 
is a $Lie(H)$-linear isomorphism ${\frak p}_{-} \cong {\frak p}_0$ inducing a canonical  identification 
$\Lambda^\bullet({\frak p}_{-}) \cong \Lambda^\bullet({\frak p}_0)$ of $H$-modules.

\bigskip\noindent
The universal enveloping algebra $U({\frak p})$ of $\frak p$ contains the universal enveloping algebras $U({\frak p}_{-}\oplus {\frak p}_0)$ and $U({\frak p}_{-})$ as subalgebras.
For $\theta\in \Lambda^\bullet({\frak p}_{-}) \cong U({\frak p}_{-})$ the supercommutator $[ \mathbb{D}, \theta]$
is contained in $U({\frak p}_{-}\oplus {\frak p}_0)$. For a basis $x_{ij}$ of $Lie(H)$ it has the
form
$  [ \mathbb{D}, \theta] \ = \ \theta_0 + \theta_1$ with $\theta_1= \sum_{i,j=1}^n \theta_{ij}\cdot  x_{ij} $
for uniquely defined elements $\theta_0, \theta_{ij} \in \Lambda^\bullet({\frak p}_{-}) \cong U({\frak p}_{-})$.

\medskip\noindent
We now consider $V(\lambda)$ as a $\p$-module \[ V(\rho^{\vee} \boxtimes \rho)|_{\p} = \Lambda^{\bullet}(\p_{-}) \otimes End(\rho).\]

\begin{lem} \label{lem:CE} The operator $DS_n$ induces the Lie algebra homology differential $\delta$ on the Chevalley-Eilenberg complex $\Lambda^{\bullet}(\p_{-}) \otimes End(\rho)$.
\end{lem}

We shortly recall the definition. For a Lie algebra $\g$ and a $\g$-module $V$ we consider the complex with $p$-th entry $V_p(\g,V) = \Lambda^p(\g) \otimes V$ and differential $\delta: V_{p}(\g,V) \to V_{p-1}(\g,V)$ given for $p \geq 2$ by $\delta(x_1 \wedge \ldots \wedge x_p \otimes v) = \theta_0\otimes v  + \theta_1(v) $ for $x_1,...,x_p\in \frak g$ where \begin{align*}  \theta_0\otimes v & = \Bigl(\sum_{\mu < \nu} (-1)^{\mu + \nu} [x_{\mu},x_{\nu}] \otimes x_1 \wedge \ldots \wedge \hat{x}_{\mu} \wedge \ldots \wedge \hat{x}_{\nu} \wedge \ldots \wedge x_p\Bigr) \otimes v \\
\theta_1(v) & = \sum_{\nu=1}^p \Bigl(  (-1)^{\nu+1} x_1 \wedge \ldots \wedge \hat{x}_{\nu} \wedge \ldots \wedge x_p\Bigr) \otimes x_{\nu}(v) \ .\end{align*}

\medskip\noindent
{\it Proof}. $\mathbb{D}$ acts on an element $x_1 \wedge \ldots \wedge x_r\otimes \varphi$  in $\Lambda^{\bullet}(\p_{-}) \otimes End(\rho)$ for $x_1,...,x_r\in {\frak p}_{-}$ and $\varphi\in End(\rho)$ as \begin{align*} & \mathbb{D}(x_1 \wedge \ldots \wedge x_r \otimes \varphi) \\& = \sum \pm x_1 \wedge \ldots \mathbb{D}(x_i) \wedge \ldots x_r \otimes \varphi \pm \sum x_1 \wedge \ldots \wedge x_r \otimes \mathbb{D}(\varphi) \end{align*} with $\mathbb{D}(x_i) \in \p_0$. The second sum vanishes since $\mathbb{D}(\varphi) = 0$ for all $\varphi \in End(\rho)$ by definition of the Kac module. We now evaluate the first sum. $\mathbb{D}$ acts on an element in $\p_{-1}$ by the supercommutator \[ \left[\begin{pmatrix} 0 & id_n \\ 0 & 0 \end{pmatrix},\begin{pmatrix} 0 & 0  \\ x & 0 \end{pmatrix} \right] \ =\  \begin{pmatrix} x & 0 \\ 0 & x \end{pmatrix} \ \in \p_0.\] Therefore $\mathbb{D}$ acts on an element in $V(\lambda)$ as  \begin{align*} & \mathbb{D}(x_1 \wedge \ldots \wedge x_r \otimes \varphi) \\ & = [\mathbb{D},x_1]x_2 \wedge \ldots \wedge x_r \otimes \varphi - x_1 \wedge [\mathbb{D},x_2]x_3 \wedge \ldots \wedge x_r \otimes \varphi \ldots \\ & + (-1)^{r+1} x_1 \wedge \ldots [\mathbb{D},x_r]\otimes \varphi \\ & = \begin{pmatrix} x_1 & 0 \\ 0 & x_1 \end{pmatrix}  (x_2 \wedge \ldots \wedge x_r \otimes \varphi) - x_1 \wedge \begin{pmatrix} x_2 & 0 \\ 0 & x_2 \end{pmatrix}  (x_3 \wedge \ldots \wedge x_r \otimes \varphi) + \ldots \\ & + (-1)^{r+1} x_1 \wedge \ldots \wedge x_{r-1}\otimes \begin{pmatrix} x_r & 0 \\ 0 & x_r \end{pmatrix}  (\varphi) \end{align*} where the derivations $[\mathbb{D},x_{\nu}] \in \p_0$ act on all  terms to the right. The $\theta_1$-term arises from the action of  $[\mathbb{D},x_{\nu}] $ on the last term $\varphi$ to the right, the remaining terms  lead to a sum with the $\sum_{\mu < \nu}$-condition defining the $\theta_0$-term.  \qed

\medskip
Viewing $\theta:= x_1\wedge \cdots \wedge x_r$ as an element in the universal enveloping algebra $U({\frak p}_{-})$
of ${\frak p}_{-}$,  the super commutator $[\mathbb{D}, \theta]$ in the universal enveloping algebra $U({\frak p})$ of $\frak p$ is $[\mathbb{D}, \theta] = \theta_0 + \theta_1$ with $\theta_1$ in the universal enveloping algebra $U({\frak p}_{-}\oplus {\frak p}_0)$ of ${\frak p}_{-}\oplus {\frak p}_0$
and $\theta_0$ in the universal enveloping algebra $U({\frak p}_{-})$ of ${\frak p}_{-}$ as defined above, but viewed as element in the universal enveloping algebra of $\frak p$. Furthermore $\theta_1$ annihilates $H$-invariant vectors in any ${\frak p}$-module. Stated in this form, the assertion obviously holds for  arbitrary elements $\theta$ in the universal enveloping algebra of ${\frak p}_{-}$.
%

\medskip\noindent
{\it $H^{\bullet}(V(\lambda))$ and the theorem of Hopf}. Lemma \ref{lem:CE} identifies $H_{DS_n}^{\bullet}(V(\lambda))$ with the Lie algebra homology ring $H_{\bullet}(\p,End(\rho))$. We recall some facts about Lie algebra (co)homology. Note that $H_{\bullet}(\g) = H^{\bullet}(\g)$.

\medskip\noindent
Let $\g$ be a reductive Lie algebra and $\Lambda^{\bullet}(\g)^{\g}$ the space of invariants under the adjoint action of $\g$. It has the structure of a graded super Hopf algebra. Let $P(\g)$ denote the space of primitive elements, i.e. \[ P(\g) = \{ x \in \Lambda^{\bullet}(\g)^{\g} \ | \ \Delta(x) = x \otimes 1 +  1 \otimes x \} \] where $\Delta$ denotes the comultiplication. Define a grading on $P(\g)$ by requiring that the inclusion $P(\g) \to \Lambda^{\bullet}(P(\g))$ preserves degrees.

\begin{thm} \cite[Theorem 10.2, Corollary 10.2, Corollary 10.3]{Meinrenken} (Hopf-Koszul-Samelson) \label{thm:hopf}
\begin{enumerate}
\item The inclusion of $P(\g)$ in $\Lambda^{\bullet}(\g)^{\g}$ extends to an isomorphism of graded super Hopf algebras $\Lambda^{\bullet}(P(\g)) \cong \Lambda^{\bullet}(\g)^{\g}$.
\item There is an isomorphism $H^{\bullet}(\g) \cong \Lambda^{\bullet}(P(\g))$ of graded super Hopf algebras, i.e. the cohomology ring is an exterior algebra over the primitive elements. In particular the elements in $\Lambda^{\bullet}(P(\g))$ are closed.
\item The space of primitive elements has dimension $rank(\g)$. For $\mathfrak{gl}(n)$ the basis elements $f_1,f_3, .. , f_{2n-1}\in P(\g)$ have degree $1, 3, .. , 2n-1$.
\end{enumerate}
\end{thm}

We now apply this theorem for the Lie algebra $\g$ of $H$ and the $H$-invariant ring $I$ in the universal enveloping algebra 
of ${\frak p}_-$ using the following identifications
$I \cong \Lambda^\bullet({\frak p}_{-})^{H}\cong  \Lambda^{\bullet}({\frak p}_0)^{H}\cong  \Lambda^{\bullet}(\g)^{\g} \cong V(\one)^H$ for the invariant ring
$$ I:=  U({\frak p}_{-})^{H} \ .$$  From theorem \ref{thm:hopf} we obtain 
the following corollary.

\begin{cor}
The cohomology $H_{DS_n}^{\bullet}(V(1))$ is isomorphic to $I\cong V(\one)^H$ and $I$ has the structure of a supercommutative polynomial ring $\mathbb C\{ f_1,..,f_{2n-1}\}$ generated by elements $f_\nu$ in the degrees $1-2\nu$
for $\nu=1,..,n$. In particular $$\omega_{DS_n}(t) = \prod_{\nu=1}^n (1 + t^{1-2\nu})\ .$$
\end{cor}

\begin{lem} \label{HOPF} For any $\p$-module $V$, the cohomology group $H^{\bullet}_{DS_n}(V^H)$ is a graded $I$-module.
\end{lem}

{\it Proof}. By theorem \ref{thm:hopf} we have $\mathbb{D}(v)=\theta_0(v) + \theta_1(v)=0$ for every element $v\in V(1)^H$.
Since $\theta_1(v)=0$ holds for $H$-invariant vectors, we get $\theta_0(v)=0$ and hence $\theta_0=0$
holds in $U({\frak p}_{-})$. This implies $[\mathbb{D},\theta] = \theta_1$ for all $\theta\in U({\frak p}_{-})^H \cong \Lambda^\bullet({\frak p}_{-})^H = I$.
For any  $P\in I$, hence $[\mathbb{D},P] = P_1 \in U({\frak p})$ annihilates  $H$-invariant vectors. For any finite dimensional
algebraic ${\frak p}$-module $V$, $V$ in particular is an $U({\frak p}_{-})$-module and the subspace $V^H$ obviously is
an $I$-module.  Since $\mathbb{D}$ commutes with $H$, we obtain
a linear map $\mathbb{D}: V^H \to V^H$. For $P\in I$ and $v\in V^H$ the formulas $\mathbb{D}(P v)= [\mathbb{D},P ]v + P \mathbb{D}(v)$ and $[\mathbb{D},P] = P_1$
and $P_1v=0$ imply $\mathbb{D}(Pv)=P\mathbb{D}(v)$ for all $v\in V^H$. Hence the subspace of $\mathbb{D}$-coboundaries resp. of $\mathbb{D}$-closed elements in $V^H$ are both $I$-modules. \qed


\begin{lem} \label{lem:cohom-invar0}  For finite-dimensional $\mathfrak{gl}(n\vert  n)$-modules $M$ the following holds:
$H_{DS_n}^{\bullet} (M) \cong H_{DS_n}^{\bullet} (M^H)$.
\end{lem}

{\it Proof}. $H$ commutes with $\mathbb{D}$ and operates therefore on the cohomology $H_{DS_n}^{\bullet} (M)$. Since $H$ is reductive, a finite-dimensional representation of $H$ is trivial if and only if its restriction to a Cartan subgroup is trivial. We therefore show that the diagonal torus $T \subset H$ acts trivially on the cohomology. By the Leray spectral sequence \[ DS_{n,n-1} \circ DS_{n-1,n-2} \circ \ldots \circ DS_{1,0} \Longrightarrow DS_{n,0} = DS_n.\] By section \ref{DDirac} $DS_{n,n-1}$ is invariant under \[ H_{n,n-1} = \begin{pmatrix} 0 &  & & &\\ & \ldots & & &  \\ & & 0 & \\ & & & 1 \end{pmatrix}, \] $DS_{n-1,n-2}$ is invariant under \[ H_{n-1,n-2} = \begin{pmatrix} 0 &  & & &\\ & \ldots & & &  \\ & & 1 & \\ & & & 0 \end{pmatrix} \] and so on. Hence $H_{DS_n}^{\nu}(M)$ has a filtration which is respected by $T$ such that $T$ acts trivially on the graded pieces. Since $T$ acts in a semisimple way, this implies that the operation of $T$, and therefore of $H$, is trivial.  \qed

\begin{prop} \label{lem:cohom-invar} For $M \in \mathcal{R}_n$ the cohomology $H_{DS_n}^\bullet(M)$ is a graded $I$-module for the graded polynomial ring $I$. For morphisms $f: M \to M'$ in $\mathcal{R}_n$ the induced map
$H^\bullet_{DS_n}(M) \to H^\bullet_{DS_n}(M')$ is graded $I$-linear. 
\end{prop}

{\it Proof}. This follows from the lemmas \ref{HOPF}  and \ref{lem:cohom-invar0}. 
\qed

\begin{lem} Let $\rho$ be an irreducible representation of $Gl(n)$. Then the map  \[ \varphi: V(\one)|_{\p} \to V(\rho^{\vee} \boxtimes \rho)|_{\p}\ \ , \ \ v \otimes 1 \mapsto v \otimes id_{\rho} \] is a $\p$-linear inclusion.
\end{lem}

{\it Proof}. We have $x_{\nu} (\varphi) = 0$ for all $\nu$ if $\varphi \in End_H(\rho) = \mathbb{C} id_{\rho}$. \qed

\medskip\noindent
{\it Remark.} Every maximal atypical $\g$-module $V$, when restricted to $\g_0$, has the form $V|_{\g_0} = \bigoplus \rho_{\nu} \boxtimes \rho_{\mu}$ with $deg(\rho_{\nu}) = deg(\rho_{\mu})$. This degree makes $V$ into a graded $\p$-module. For $V(\one)$ we obtain the degree defined previously. For $V(\rho^{\vee} \boxtimes \rho)$, $id_{\rho}$ has degree $deg(\rho)$. Therefore $\varphi$ shifts the degrees by $deg(\rho)$. The degree $deg(\rho)$  coincides with $p(\lambda)$ for $\rho$ with highest weight $\lambda = (\lambda_1,\ldots,\lambda_n)$.

\begin{lem} \label{lem:kac=kac} The induced morphism \[ \xymatrix{ H^{\bullet}_{DS_n}(V(\one)) \ \ar[r]^-{H_{DS_n}^{\bullet}(\varphi)} & \ H_{DS_n}^{\bullet + deg(\rho)} (V(\rho^{\vee} \boxtimes \rho)) } \] is a graded isomorphism on the cohomology. Hence $ H_{DS_n}^{\bullet + deg(\rho)} (V(\rho^{\vee} \boxtimes \rho))  $ is the free $I$-module of rank one
generated by the top cohomology.
\end{lem}

{\it Proof}. As a $\p$-module \[ V(\rho^{\vee} \boxtimes \rho)|_{\p} \cong \varphi(V(\one)) \oplus (\Lambda^{\bullet}(\p_-) \otimes End^0(\rho)) \] where \[ End^0(\rho)  = \{ \phi \in End(\rho) \ | \ Tr(\phi) = 0 \}. \] Since $D \in \p$, we obtain \[ H_{DS_n}^{\bullet} (V(\rho^{\vee} \boxtimes \rho))\ =\ H_{DS_n}^{\bullet} (\varphi( V(\one))) \ \oplus\  H^{\bullet}_{DS_n} (\Lambda^{\bullet}(\p_-) \otimes End^0(\rho)). \] By \cite{Hochschild-Serre} the Lie algebra cohomology for reductive $H$ with coefficients in a representation $W$ is trivial except for the trivial representation \[ H^{\bullet}(Lie(H),W) \ \cong\ H^{\bullet}(Lie(H),\one) \otimes W^H.\] Since the Lie algebra cohomology is dual to the homology, this shows \[  H^{\bullet}_{DS_n} (\Lambda^{\bullet}(\p_-) \otimes End^0(\rho)) = 0. \] Since $\delta = ad_{\mathbb{D}}$ commutes with $\varphi$, we get \[ \xymatrix{ H^{\bullet}_{DS_n}(V(\one)) \ \ar[r]_-{H^{\bullet}_{DS_n}(\varphi)}^-{\cong} & \ H^{\bullet}_{DS_n}(\Lambda^{\bullet}(\p_-) \otimes id_{\rho}) } \] up to the degree shift with $deg(\rho)$. \qed

\begin{cor} For the Hilbert polynomial of $V(\lambda)$ relative to $D$ we obtain \[ \omega_{DS_n}(V(\lambda),t) = t^{p(\lambda)}\cdot  \omega_{DS_n} (V(\one),t) =  t^{p(\lambda)}\cdot  \prod_{\nu=1}^n (1 + t^{1-2\nu}) .\]
\end{cor}

%

\begin{thm}\label{cohom-proj} Let $L(\lambda)$ be an  irreducible and maximal atypical
representation and
$pr: V(\lambda) \to L(\lambda)$ be a projection onto the top. Then
the induced homomorphism \[ H^{\nu}_{DS_n}(pr): H_{DS_n}^{\nu}(V(\lambda)) \to H_{DS_n}^{\nu}(L(\lambda)) \] is zero in degrees $\nu < p(\lambda)$ and an isomorphism for $\nu = p(\lambda)$. 
\end{thm}

{\it Proof}. As a graded $I$-module $H_{DS_n}^{\nu}(V(\rho^{\vee} \boxtimes \rho))$ is 
the free $I$-module generated by the cohomology in the top degree.
To prove our claim it suffices that the primitive elements $f_1, f_3,\ldots, f_{2n}-1\in I$
act trivially on $H_{DS_n}^{\nu}(L(\rho^{\vee} \boxtimes \rho)) $. This follows from the discussion in section \ref{sec:forest}, lemma \ref{ev}, which shows  for $\nu = 1,\ldots,n$ \[ H_{DS_n}^{deg(\rho) - 2\nu +1} (L(\rho^{\vee} \boxtimes \rho)) = 0 \ .\]
\qed

\section{Primitive elements of $H^{\bullet}_{DS_n}(V(\one))$} \label{sec:primitive}

We will now describe the primitive elements of $H^{\bullet}_{DS_n}(V(\one))$ in terms of the representation
theory of the superlinear group $Gl(n\vert n)$. The radical filtration on $V(\one)$ defines a decreasing filtration $F_i$ of $V(\one)$. The $H$-invariants $F_i^H$ coincide with the powers $(I^+)^i$ of the augmentation ideal $I^+$ of the invariant ring $I = V(\one)^H$. In this way monomials of degree $i$ in the primitive generators $f_1,...,f_{2n-1}$ can be identified with the generators of the  cohomology $H^\bullet(F_i^H/F_{i+1}^H)$.

\medskip\noindent
{\it The Murnagan-Nakayama rule}. Let $\lambda=(\lambda_1,...,\lambda_r)$ with $\lambda_1\geq \lambda_2\geq \cdots \geq
\lambda_r$ be a partion of degree $n=deg(\lambda)$. For partitions $\nu$ and $\mu$ 
of $m$ and $n-m$ let $c_{\mu\nu}^\lambda$ denote the {\it Littlewood-Richardson}
coefficient. Assume that $\nu$ is a {\it hook}, i.e a partition of type $\nu_1=r$, $\nu_2=\cdots =\nu_{m-r+1}=1$
and $\nu_i=0$ for $i>m-r+1$. Recall that a hook is a special case  of a {\it rim hook} (also called {\it skew hook}).
We say $\nu$ is a {\it symmetric hook} if $m=2r-1$.
According to \cite[Section 4.10]{Sagan} we have

\begin{prop} \label{murnaghan-nakayama} Suppose $\nu$ is a hook. Then
 $c_{\mu\nu}^\lambda=0$ unless the Young diagram of $\mu$ is contained in the Young diagram of $\lambda$
and the complement $\lambda/\mu$ is a union of $k$ edgewise connected rim hooks. If this is the case, then
$$  c_{\mu\nu}^\lambda \ = \ {k-1 \choose c - r } $$
where $r=\nu_1$ and $c$ is the number of rows spread by the rim hooks contained in $\lambda/\mu$.
\end{prop}
 
We remark that in \cite{Sagan} $c$ denotes the number of columns instead of rows, since Young diagrams in \cite{Sagan}
are written top down instead of being written from left to right, as with our conventions.

\begin{cor} \label{cor:LRR} Suppose $\lambda$ and $\nu$ are symmetric hooks.
Then $c_{\mu\nu}^\lambda\!=\!0$ unless $\mu\!=\!\nu$ or $\nu\!=\!0$.
\end{cor}

{\it Proof}. Suppose $c_{\mu\nu}^\lambda\neq 0$. By the proposition the edgewise connected components of $\lambda/\mu$ are rim hooks, hence $\#(\lambda/\mu) \leq 2$. Since $deg(\nu)$ is odd for symmetric hooks and $ \#(\lambda/\mu) =deg(\nu)$, we may assume without restriction of generality that $\#(\lambda/\mu) =1$. But this gives a contradiction since
$deg(\lambda) - deg(\mu) =1$ would  be the difference of two odd numbers. \qed

\medskip\noindent
Let $\rho^\vee$ denote the dual representation of $\rho$.
Suppose $\lambda$ is a partition of $n$ and $\lambda^*$ is the dual partition of $n$,
then define  $(\rho_\lambda)^* := \rho_{\lambda^*}$ for the representations $\rho=\rho_\lambda$
of $GL(n)$ with highest weight $\lambda$.  

\begin{cor} Suppose that $\nu$
is a symmetric hook of degree $2r-1$  and suppose $k=1$ (in the notation of proposition \ref{murnaghan-nakayama}). Suppose the rim hook $\lambda/\mu$ reaches from $(i, \lambda_i)$
to $(j,\lambda_j)$ where $i>j$. Then   $c_{\mu\nu}^\lambda\!=\!0$ hold unless $\lambda_i - \lambda_j \! =\! i - j\!=\! r$.
\end{cor}

{\it Proof}. Since $k=1$, ${k-1 \choose c - r }\neq 0$ if and only if $\lambda_i-\lambda_j=c=r$.
Since $\nu$ is a rim hook, we have $2r-1=deg(\nu)= (\lambda_i-\lambda_j)+ (j-i) - 1$.  Hence
$\lambda_i-\lambda_j=r$ implies $j-i=r$. \qed

\medskip\noindent
{\it The Lie superalgebra $\mathfrak{gl}(n\vert n)$ and primitive elements of $\mathfrak{gl}(n)$}. The following proposition is a well-known consequence of the dual Cauchy identity.

\begin{prop} \cite[Theorem B.17]{Bump-Schilling} The space of matrices  $M_{nn}(k)$ is a $Gl(n,k)\times Gl(n,k)$-module in a natural way by left and right multiplication, hence also the Gra\ss mann algebra $\Lambda:=\Lambda^\bullet(M_{n}(k))$. As a representation
of $Gl(n,k)\times Gl(n,k)$ we have
$$ \Lambda^\bullet(M_{nn}(k)) \ \cong \ \bigoplus_{\rho}
  \rho^\vee \boxtimes \rho^*$$
  where $\rho=\rho_\lambda$ runs over all partitions in \[ P(n,n)=  \{ \lambda \in {\bf Z}^n\ \vert \ n \geq \lambda_1 \geq \lambda_2 \geq ... \geq \lambda_n \geq 0 \} \ .\]
\end{prop}

\medskip\noindent
Warning: The degree $deg(\rho^\vee)$ is the the negative of the degree in the Gra\ss mann algebra $\Lambda$!

\begin{cor} Let $H=GL(n,k)$ be embedded diagonally. Then \[ I:=\Lambda^\bullet(M_{nn}(k))^H \cong  \bigoplus_{\rho}
  (\rho^\vee \boxtimes \rho^*)^H \ ,  \] and $ (\rho^\vee \boxtimes \rho^*)^H\neq 0$ if and only if $\rho=\rho_\lambda$ 
  for
  a symmetric Young diagram  $\lambda=\lambda^*$. There exist $2^n$ symmetric Young diagrams
  with $\lambda=(\lambda_1,...,\lambda_n)$ and  $\lambda_1\leq n$.
\end{cor}

\bigskip\noindent  
The  space $I$ is an algebra with respect to the wedge product. The subspace $I^+\subseteq I$ of elements of degree $\geq 1$ is an ideal (the augmentation ideal).

\begin{prop} \cite[Proposition 10.11]{Meinrenken} $I^+$ decomposes as 
\[  I^+ \ = \ P(H) \oplus (I^+)^2 \ .\]
\end{prop}

\begin{cor} \label{prim-hooks} With summation over all $\rho\!=\!\rho_\lambda$ for symmetric hook diagrams
$\lambda$ of degrees $deg(\lambda)\!=\!1,3,5,....,2n-1$ the space of primitive elements 
is
$$ P(H) \ = \ \bigoplus_{\rho} \ (\rho^\vee\boxtimes \rho^*)^H \ .$$
\end{cor}

\noindent
{\it Proof}. This follows from the fact that for hook diagrams $\lambda$
the space $\rho_\lambda $ cannot be a constituent of $\rho_\mu \otimes\rho_\nu$ 
for $\mu=\mu^*$ and $\nu=\nu^*$ where $(\rho_\mu^\vee \boxtimes \rho_\mu^*)^H$ and $(\rho_\nu^\vee \boxtimes \rho^*_\nu)^H$ are constituents of $I^+$. Hence $(\rho_\lambda^\vee \boxtimes \rho^*_\lambda)$ cannot be contained
in $(I^+)^2$. \qed 

\medskip\noindent
{\it The index}. The selftransposed weights $\lambda_{(i)}=(i,..,i,0,..,0)$ for $i=0,..,n$ in $P(n,n)$ are called the {\it basic selftransposed weights}. The index $ind(\lambda)$ of a selftransposed $\lambda$ in $P(n,n)$
is the maximal index $i$ of a basic selftransposed $\lambda_i$ whose Young diagram is 
contained in the Young diagram of $\lambda$. The index
of $\lambda$ is the unique $i$ between $1$ and $n$ such that $\lambda_i \geq i$ and 
$\lambda_{i+1} \leq i$. We denote by $P_i(n,n)$ the set of all weights in $P(n,n)$ with index $i$.

\begin{prop} Using that $\Lambda \cong V(\one)$, the canonical filtration
defined by the radical filtration of $V(\one)$ in the category of $\mathfrak{gl}(n\vert n)$-modules
gives a filtration $F_i$ on $\Lambda$ such that
$$   F_i   \ = \ \bigoplus_{\rho} \ (\rho^\vee \boxtimes \rho^*) $$
for all $\rho=\rho_\lambda$ running over all partitions $\lambda$ containing the partition $(i^i)$ of degree $i^2$. 
\end{prop}

Before the proof we recall that $V(\one)$ has a decreasing filtration (the radical filtration) of $Gl(n\vert n)$-subrepresentations 
with $n+1$ irreducible graded pieces $L_i$ such that $L_0=k$ is the maximal
irreducible quotient representation. The highest weights of the $L_i$ can be computed
from \cite[Theorem 5.2]{Brundan-Stroppel-1} to be the duals
$$   \lambda_{(i)}^\vee = (0,\cdots,0,-i,...,-i)   \quad  , \quad  \mbox{ for }\  i=0,...,n  \  $$
of the basic selftransposed weights $\lambda_i$ in $P(n,n)$.

\medskip\noindent
{\it Proof}. We need to show that the representation $L_i$, considered as a representation of $G \subset Gl(n\vert n)$, decomposes into a direct sum over the duals of all irreducible representations $\rho(\lambda)\boxtimes \rho(\lambda^*)$ for which $$\lambda \in P_i(n,n) \ .$$ Consider the decomposition of $L_i^\vee$ under $G=Gl(n)\times Gl(n)$. Let $\lambda=(\lambda_1,..,\lambda_n)$ be a corresponding highest weight of $G$ in $L_i$. We then claim $ind(\lambda)=i$.
Obviously 
$$  \lambda \geq \lambda_{(i)} = (i,...,i,0,...,0)\ .$$
On the other hand we have  $$ V(\one)^\vee \cong (det^n\boxtimes det^n) \otimes V(\one) \ $$
since the dual of a Kac module is a Kac module. Hence, since the order of the socle layers in the dual Kac module is reversed 
and since the Loewy length of $V(\one)$ is $n+1$ \cite[Theorem 5.2]{Brundan-Stroppel-1}, this implies
$$ L_i^\vee \cong (det^n\boxtimes det^n) \otimes L_{n-i} \ .$$
This in turn implies
$$ \lambda \leq (n,...,n) + \lambda_{(n-i)}^\vee = (n,...,n,i,...,i) \ $$
with $i$ copies of $n$ and $n-i$ copies of $i$.
Both estimates together force
$ \lambda_i \geq i$ and $\lambda_{i+1}\leq i$, hence $ind(\lambda)=i$.
This proves our claim. Since any $\lambda$ in $P(n,n)$ appears in one of the $L_i^\vee$,
$L_i^\vee$ then consists precisely of the $G$-constituents $\rho(\lambda)\boxtimes \rho(\lambda^*)$
for $\lambda$ in $P_i(n,n)$. \qed

 
\begin{cor} $F_i^H = (I^+)^i$, hence we can identify 
monomials of degree $i$ in the primitive generators $f_1,...,f_{2n-1}$
with the generators of the  cohomology $H^\bullet(F_i^H/F_{i+1}^H)$.
\end{cor}

We introduce the notation $Prim_i \subseteq I^+$ for the space that is spanned by monomials in the primitive
elements $f_{2\nu-1}$ with exactly $i$ factors. In this notation $Prim_1 = P(H)$.

\medskip\noindent
{\it Proof}. By corollary \ref{prim-hooks} $Prim_1$ is a complement to $F_2^H$. We now show $Prim_i \cap F_{i+1}^H = 0$ by induction on $i$.
 Using the induction assumption and $Prim_i = Prim_{i-1} \cdot Prim_1$, the space $Prim_i$ only gives rise to Young diagrams $\lambda$ that occur in the tensor product of some $\mu \in P_{j}(n,n)$ for $j< i$ and a symmetric hook $\nu\in P_1(n,n)$. By proposition \ref{murnaghan-nakayama}, $\mu$ is obtained from $\lambda$ by removing a (possibly disconnected) rim hook. If $\lambda\in P_k(n,n)$, this implies $j=k$ or $j=k-1$ and hence $k\leq j+1 \leq i$.
This proves $Prim_i \cap F_{i+1}^H = 0$ since all selftransposed weights $\lambda$ in $F_{i+1}^H$
are contained in $P_k(n,n)$ for $k \geq i+1$. 

\medskip\noindent
This implies $Prim_\nu \cap F_{j}^H = 0$ for $\nu<j$. Since $I^+$ is the direct sum of $ \bigoplus_{\nu =0}^i Prim_\nu$ and $(I^+)^{i+1}$,
this implies that $F_{i+1}^H$ is in the complement of $\bigoplus_{\nu =0}^i Prim_\nu$ and therefore $F_{i+1}^H \subseteq (I^+)^{i+1}$. There are ${n \choose i}$ selftransposed weights  $\lambda\in P_i(n,n)$ and all of them
occur in $V(\one)$ with multiplicity one. On the other hand, 
the space $Prim_i \subseteq I^+$ that is spanned by monomials in the primitive
elements $f_{2\nu-1}$ with exactly $i$ factors also has dimension ${n \choose i}$. Since the dimensions agree, this implies $(I^+)^{i+1} = F_{i+1}^H$ and $Prim_i$ is therefore represented by $(F_i/F_{i+1})^H \cong H^\bullet_{DS_n}(F_i/F_{i+1})$. \qed

\section{Kac module of $\one$}\label{kac-module-of-one}

\medskip\noindent

\textit{Overview.} We now study the effect of $DS$ on indecomposable modules in the remaining sections \ref{kac-module-of-one} - \ref{hooks}. The easiest examples are perhaps the extensions of two irreducible modules, and we focus here on the case of extensions of the trivial representation $\one$ by another irreducible module. Our main result in these sections is corollary \ref{splitting1}, saying that a representation $Z$, such that the projection onto its cosocle $\one$ induces a surjection $\omega: \omega(Z) \to \omega(\one)$, is equal to the trivial representation. Such a representation $Z$ contains extensions of the trivial representation $\one$ with other irreducible representations. We show in the resum\'{e} of section \ref{strictmorphisms} that if $Z$ is not irreducible, we obtain extensions $V$ of $\one$ by an irreducible representation $S$ such that the induced morphism $\omega(V) \to \omega(\one)$ is surjective. Since the dimension of $Ext^1(S,\one)$ is at most one-dimensional, any two such extensions are isomorphic. Hence we can study them by realizing them as quotients of modules whose cohomology is sufficiently understood. A typical example occurs in the current section \ref{kac-module-of-one}: The Kac module $V(\one)$ of $\one$ contains an extension of $\one$ with the irreducible representation $[0,\ldots,0,-1]$; and by considering the cohomology of the Kac module we are able to compute the cohomology of this extension and its dual in lemma \ref{Ia} and lemma \ref{aKac}. We also show in corollary \ref{aKac2} that $\omega^0(V) = 0$ for the extension $V$ between $\one$ and $Ber \otimes S^{n-1}$. The other $n$ nontrivial extensions of $\one$ (listed in lemma \ref{ext}) are studied in section \ref{hooks}. We realize these extensions as a quotient of the mixed tensor $R(n^n)$ studied in section \ref{sec:n^n}. The key proposition \ref{trivialextension} shows that for any of our nontrivial extensions $V$ the zero degree part $\omega^o$ of the induced map $\omega(q_V): \omega(V) \to \omega(\one)$ vanishes, a contradiction our analysis in the resum\'{e} of section \ref{strictmorphisms}, hence $Z \simeq \one$.

\medskip
The constituents of the Kac module
$V({\bf 1}) \in \calR_n$ are \cite{Brundan-Stroppel-1}, thm. 5.2, $$L_a = Ber^{-a} \otimes [a,...,a,0,...,0] \ \text{ for } \ a=0,...,n \ ,$$ where the last entry of $a$ is at the position $i=n-a$.  Therefore $Ber^a \otimes L_a$ is basic
and therefore has cohomology concentrated in degree zero, hence
the cohomology of $L_a$ is concentrated in degree $-a$ and
$$  H^{-a}(L_a) \ \cong \ I_a \ \oplus \ I_{a-1} \quad , \quad a=0, 1, \ ...\ , n $$
where $I_{-1}:=I_n:=0$ and 
$$ I_a \ := \ Ber^{-a-1} \otimes [a+1,...,a+1,0,...,0] $$
(with $n-a-1$ entries $a+1$ and $a$ entries $0$). Notice $I_1^\vee \cong Ber \otimes S^{n-1}$
and $I_0={\bf 1}$, $I_1=[0,..,0,-2], ... , I_{n-1}=Ber^{-n}$. 
For the cyclic quotient $Q_a$ of $V({\bf 1})$ with socle $L_a$ this implies inductively 

\begin{lem}\label{Qa}
The natural quotient map $Q_a \to {\bf 1}$ induces an isomorphism 
$H^0(Q_a)\cong H^0({\bf 1}) \cong {\bf 1}$ and
$$  H^{-\nu}(Q_a) =
\begin{cases}
I_\nu \oplus I_{\nu -1}     & \text{$\nu=0,...,a$}, \\
   0   & \text{otherwise}.
\end{cases} \ .$$ 
\end{lem}

Notice $Q_a=V({\bf 1})$ for $a=n$ and $Q_a={\bf 1}$ for $a=0$. 
Similar as in the proof of the last lemma, for $K_a=Ker(V({\bf 1})\to Q_a))$ and $a \leq n-1$ 
we obtain exact sequences
$$  0 \to H^\bullet(K_a) \to H^\bullet(V({\bf 1}) ) \to H^\bullet(Q_a ) \to 0 \ .$$
Indeed,  the cohomology of $H^\bullet(K_a)$ is concentrated in degrees $\leq -a-1$,
whereas the cohomology of  $H^\bullet(Q_a )$ is concentrated in degrees $\geq -a$.
We can view these as short exact sequences of homology complexes
$$  0 \to (H^\bullet(K_a),\overline\partial) \to (H^\bullet(V({\bf 1}) ),\overline\partial) \to (H^\bullet(Q_a ),\overline\partial) \to 0 \ .$$
The long exact homology sequence for the $H_{\overline\partial}$-homology together with $H_{\overline\partial}(H^{\nu}(V)) = H_D^{\nu}(V)$ (lemma \ref{abutment}) implies
$$     \xymatrix@-0.3cm{  H_D^{-\nu}(K_a)  \ar[r] &  H_D^{-\nu}(V({\bf 1})) \ar[r] & H_D^{-\nu}(Q_a) \ar[r]^-\delta &
H_D^{-\nu-1}(K_a) \ar[r] & H_D^{-\nu-1}(V({\bf 1}))}  \ $$
and $H^\nu_D(V({\bf 1}))=0$ for all $\nu$ hence gives
$H^ \nu_D(Q_a) \cong H_D^{\nu -1}(K_a)$. Now $H_D^{\nu}(Q_a)$ vanishes unless $\nu \geq -a$ by
lemma \ref{Qa}. The right hand side $H^{\nu - 1}(K_a)$ is concentrated in degrees $\nu \leq -a$. Hence the long exact homology sequence has
at most one nonvanishing connecting morphism $\delta$, namely $\delta: H_D^{-a}(Q_a) \to 
H_D^{-a-1}(K_a)$ in degree $-a$. Hence 
$H_D^\nu(Q_a)=0$ for
$\nu\neq {-a}$. Since there is a unique common irreducible module $I_{a}$ in the cohomology
$H^{-1-a}(K_a)$ and $H^{-a}(Q_a)$ such that $d(Q_a) = \pm I_a$, we conclude

\begin{lem}\label{Ia}  For $0 \leq a \leq n-1$ we get
$$ H_D^\nu(Q_a) = \begin{cases}
  I_{a}    & \text{$\nu=-a$ }, \\
   0   & \text{otherwise}.
\end{cases}$$
\end{lem}

\medskip
{\bf Remark}. This result shows that for the $H_D^{\nu}$-cohomology there are do not exist long exact sequences attached to short exacts sequences in $\calR_n$. If these would exist, then $Q_1/L_1 \cong {\bf 1}$ would imply $H_D^{-1}(L_1)\cong H_D^{-1}(Q_1)$, in contrast to $H_D^{-1}(L_1) \cong I_1 \oplus {\bf 1}$ and $H^{-1}_D(Q_1) \cong I_1$.

\begin{cor}\label{killing}
$H_D^0(V)=0$ for $V=Q_a$ and $(Q_a^*)^\vee$ for $1 \leq a\leq n-1$.
\end{cor}

\medskip
Now we analyse in the case $a=1$ the nontrivial extension
$$  0 \to [0,...,0,-1] \to Q_1 \to {\bf 1} \to 0 \ .$$
Since $L_1^\vee \cong [0,...,0,-1]^\vee \cong Ber \otimes S^{n-1}$, also $V=(Q_1^*)^\vee$ defines  a nontrivial extension
$$  0 \to Ber\otimes S^{n-1}  \to V \to {\bf 1} \to 0 \ .$$

\begin{lem} \label{aKac} 
$V= (Q_1^*)^\vee$ defines a nontrivial extension
between ${\bf 1}$ and $Ber\otimes S^{n-1}$ in $\calR_n$
such that in $\calR_{n-1}$ the following holds
$$ H_D^\nu(V) \cong H^\nu(V) = \begin{cases}
  Ber\otimes S^{n-1}       & \text{$\nu=1$}, \\
   0   & \text{otherwise}.
\end{cases} $$
\end{lem}

\medskip
{\it Proof}. The statement about $H_D^\nu(V)$ follows immediately from lemma \ref{Ia}. We now calculate $H^{\nu}(V)$. Since the cohomology of the anti-Kac module $(V({\bf 1})^*)^\vee$ vanishes,
 $0\to (K_1^*)^\vee \to (V({\bf 1})^*)^\vee \to V \to 0$
gives
$$   H^{\ell -1}(V) \cong H^{\ell}((K_1^*)^\vee) \cong  H^{-\ell}(K_1^*)^\vee  \quad , \quad \text{ for all }  \ \ell \ .$$
$K_1^*$ is filtered with graded components $L_2,...,L_n$ so that the 
cohomology of  $K_1^*$ vanishes if the cohomology of the $L_i$ vanishes.
Hence $H^{-\ell}(K_1^*)=0$ unless $-\ell \notin \{-2,-3,...,-n\}$ and
$   H^\nu(V) =  0$  for all $ \nu \leq 0$ and  all $\nu\geq n $. On the other  hand
$H^\nu(Ber\otimes S^{n-1})=0 $ for $\nu\neq 1$ and $H^1(Ber\otimes S^{n-1})= {\bf 1} \oplus (Ber\otimes S^{n-1})$. Since $H^\nu(V)=0$, if $H^\nu({\bf 1})=0$ and $H^\nu(Ber\otimes S^{n-1})=0$,
therefore $H^\nu(V)=0$ unless $\nu=1$. \qed


\medskip
Applying $(n-1)$ times the functor $DS$ to $DS(V)\in \calR_{n-1}$,
the last lemma gives  

\begin{lem} If we apply $n$ times the functor $DS$ 
to $V= (Q_1^*)^\vee$ in $\calR_n$, we obtain that
$$   DS\circ DS \circ \cdots \circ DS(V) \ = \  \bigoplus_{\nu=0}^{n-2} \ k[-1-2\nu] \ 
$$  in $\calR_0$ is concentrated in the degrees $1,3,\cdots, 2n-3$.
\end{lem}

\medskip
The Leray type spectral sequences therefore imply the following result

\begin{cor} \label{aKac2}
For the module $V= (Q_1^*)^\vee$ in $\calR_n$, defining a nontrivial
extension between ${\bf 1}$ and $Ber\otimes S^{n-1}$, we have
$$ \fbox{$ DS_{n,0}^\ell(V) \ = \ 0 \ \text{ and }\  \omega_{n,0}^\ell(V) \ = \ 0  \quad \text{ for } \ell \leq 0 $} \ \ .$$
\end{cor}

\medskip



\section{Strict morphisms}\label{strictmorphisms}

\medskip\noindent
Recall the
functor $\omega: T_n \to svec_k$ defined by $\omega=\omega_{n,0}$.
A morphisms $q: V \to W$  in $T_n$ will be called a {\it strict epimorphism},
if the following holds
\begin{enumerate}
\item {\it $q$ is surjective}.
\item {\it $\omega(q)$ is surjective}.
\end{enumerate} 

\medskip
For a module $Z$ in $T_n$ and semisimple $L$ and 
$$ q: Z \twoheadrightarrow L $$
we make
the following 

\medskip
{\bf Assumption (S)}. {\it The induced morphism 
$$   \omega(q): \omega(Z) \to  \omega(L) $$ 
is surjective, i.e. $q$ is a strict epimorphism.}

\bigskip\noindent
Of course (S) holds for irreducible $Z$.
In the special case $L={\bf 1}$ condition
(S) is equivalent to $\omega(q)\neq 0$. We denote the cosocle of $Z$ by $C$.

\medskip
For any submodule $U\subseteq Kern(q)$ the map $q: Z\to L$ factorizes over the quotient $p: Z\to V=Z/U$
and induces the analogous morphism $q_V:  V \to L \hookrightarrow cosocle(Z/U)$. Hence
$$     q = q_V \circ p \quad , \quad  \omega_{n,i}(q) = \omega_{n,i}(q_V) \circ \omega_{n,i}(p) \ .$$
implies: $\omega_{n,i}(q)$ is surjective $\Longrightarrow \omega_{n,i}(q_V)$ is surjective.
For $i=0$ thus
\begin{itemize}
\item {\it If $Z$ is indecomposable, then $V$ is indecomposable}.
\item {\it Condition (S) for $q$ implies condition (S) for $q_V$}.
\item {\it $\omega(q_V)=0$ implies $\omega(q)=0$}.
\end{itemize}

\medskip
{\it Indecomposable $Z$}. Now assume $Z$ is indecomposable
and has upper Loewy length $m\geq 2$.
If $m\geq 3$, there exists a submodule $U \subset Z$
such that $V=Z/U$ has Loewy length 2 and such that $V$ 
again is indecomposable and satisfies assumption (S).  
So $V$ has Loewy length two and is indecomposable
with cosocle $C$. Then $(V,q_V)$ is a 
nontrivial extension
$$   0 \to S \to V \to C \to 0 \ $$
with semisimple socle $S$ decomposing into irreducible summands $S_\nu$
and cosocle $C$. The map $q$ is obtained from a projection map $pr_L:C \to L$
by composition with the canonical map $V \to C$. 
Since $V$ is indecomposable with cosocle $C$, all extensions $(V_\nu,q_\nu)$
obtained as pushouts 
$$ \xymatrix{0\ar[r] &  \oplus_\nu S_\nu \ar@{->>}[d]\ar[r] & V \ar@{->>}[d]_{\pi_\nu}\ar[r]^p & C \ar@{=}[d]\ar[r] & 0 \cr
0\ar[r] &  S_\nu \ar[r] & V_\nu \ar[r]^{p_\nu} & C \ar[r] & 0 } $$ 
must be nontrivial extensions. All $V_\nu$ again satisfy condition (S): Indeed $Im(\omega(q)) \subseteq Im(\omega(q_\nu))$
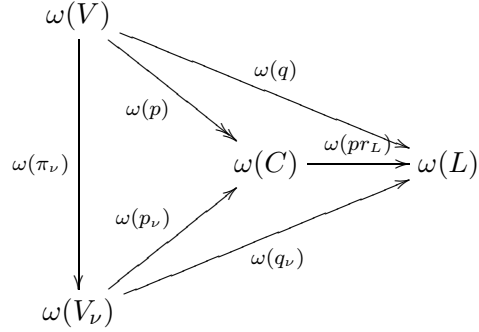
$$ \xymatrix@+0,5cm{\omega(V)\ar[rrd]^{\omega(q)}\ar@{->>}[rd]_{\omega(p)} 
\ar[dd]_{\omega(\pi_\nu)}& & \cr
& \omega(C)\ar[r]^-{\omega(pr_L)} & \omega(L)\cr 
\omega(V_\nu) \ar[rru]_{\omega(q_\nu)}\ar@{->}[ru]^{\omega(p_\nu)} & &} $$
The projection $pr_L: C \to L$ splits by an inclusion $i_L:L \to C$, since
$C$ is semisimple. Hence $C \cong L \oplus L'$ so that $pr_L$ and $i_L$
are considered as the canonical projection resp. inclusion for  the first summand.

\medskip
Since $V$ is indecomposable, $Ext^1(L,S_\nu)\neq 0$ holds for at least one $S_\nu$. Now divide by the submodule $U' \subset S$ 
generated by all $S_\nu$ with the property $Ext^1(L,S_\nu)=0$ and obtain
$V'=V/U'$. Then divide by the maximal submodule $U''$ of $L'$ that splits
in $V'$. Then $V'/U''$ is indecomposable and the map $q$ factorizes over 
this quotient and satisfies condition $S$. 

\medskip
{\it Resume}. Suppose $Z$ is indecomposable but not irreducible, $q:Z\to L$ satisfies condition (S),  the cosocle of $Z$ is $C=L \oplus L'$. Then there exists a quotient $V$ of $Z$ and a quotient
$\tilde L$ of $L'$ such that 
$$ \xymatrix{ 0 \ar[r] & S \ar[r] & V \ar[r]^-p & L \oplus \tilde L \ar[r] & 0 } $$
with 
\begin{itemize}
\item $V$ is indecomposable, 
\item $S$ is irreducible such that $Ext^1(L,S)\neq 0$
and $Ext^1(\tilde L,S)\neq 0$, 
\item the map $q=pr_L \circ p$ satisfies condition (S). 
\end{itemize}
The irreducible representations $X\not\cong {\bf 1}$
with the property $Ext^1(X,S)\neq 0$ will be called descendants of $S$.
 
\medskip
In the situation of the resume we get the extensions $E=E_S^L$ and $\tilde E= E_S^{\tilde L}$
defined by submodules of $V$. Hence $V/E_S^L \cong \tilde L$ and $V/E_S^{\tilde L} \cong L$
and we get  the following exact sequences

$$  \xymatrix{ \tilde L \ar@{=}[r] & \tilde L   &         \cr
\tilde E  \ar[u]\ar[r] &    V      \ar[u]\ar[r] &   L      \cr
 S   \ar[r]\ar[u] &    E      \ar[r]\ar[u] &    L    \ar@{=}[u] } $$

One of the potential candidates for
$S_\nu$ is the irreducible representation $L(\lambda -\mu)$ that
appears in the second upper Loewy level of the Kac module
$V(\lambda)$. 
Indeed this follows from
lemma \ref{KAC}, since $H_D(V(\lambda))=0$.
Since $Z_\nu$ is indecomposable, $Z_\nu$ is in this case 
a highest weight representation
of weight $\lambda$. This is clear, because all weights of $Z_\nu$ are 
in $\lambda - \sum_{\alpha\in\Delta_n} \mathbb Z \cdot \alpha$.
By corollary \ref{companion2}
a highest representation $V$ contains
a (nontrivial) highest weight subrepresentation $W$ of weight $\lambda-\mu$ only if 
$H_D(V)$ has trivial weight space $H_D(V)_{\overline \lambda}$.
For $V=Z_\nu$ as above this gives a contradiction, if $S_\nu = L(\lambda-\mu)$
occurs in the socle of $Z$. Indeed, notice that $H_D(L)$ contains $L(\overline\lambda)$
by lemma \ref{stable}. By condition (S) then also $H_D(Z)$ contains $L(\overline\lambda)$.
So by corollary \ref{companion2} $L(\lambda-\mu)$ is not contained in $Z_\nu$.  
This proves

\begin{lem}\label{minusmu}
Suppose $Z$ is an (indecomposable) module
with irreducible maximal atypical cosocle $L=L(\lambda)$. If $Z$ satisfies condition (S),
then the second layer of the upper Loewy filtration of $Z$ does not contain
the irreducible module $L(\lambda-\mu)$.
\end{lem}

A case of particular interest is $L={\bf 1}$. Fix some irreducible $S$
with the property $Ext_{\calR_n}(S,{\bf 1})\neq 0$. In section \ref{hooks} we will show for $L={\bf 1}$ that 
$\omega^0(q_E)=0$ (lemma \ref{omeganull}).

\medskip



\section{The module $R((n)^{n})$}\label{sec:n^n}

\medskip\noindent

We describe a certain maximal atypical mixed tensor for $n\geq 2$. 

\medskip

We recall some terminology from \cite{Brundan-Stroppel-1}. Given weights $\lambda, \mu \sim \alpha$ in the same block one can label the cup diagram $\lambda$ resp. the cap diagram $\mu$ with $\alpha$ to obtain $\underline{\lambda}\alpha$ resp. $\alpha\overline{\mu}$. These diagrams are by definition consistently oriented if and only if each cup resp cap has exactly one $\vee$ and one $\wedge$ and all the rays labelled $\wedge$ are to the left of all rays labelled $\vee$. Set $\lambda \subset \alpha$ iff $\lambda \sim \alpha$ and $\underline{\lambda} \alpha$ is consistently oriented. 

\medskip

A crossingless matching is a diagram obtained by drawing a cap diagram underneath a cup diagram and then joining rays according to some order-preserving bijection between the vertices. Given blocks $\Delta, \Gamma$ a $\Delta \Gamma$-matching is a crossingless matching $t$ such that the free vertices (not part of cups, caps or lines) at the bottom are exactly at the position as the vertices labelled $\circ$ or $\times$ in $\Delta$; and similarly for the top with $\Gamma$. Given a $\Delta \Gamma$-matching $t$ and $\alpha \in \Delta$ and $\beta \in \Gamma$, one can label the bottom line with $\alpha$ and the upper line with $\beta$ to obtain $\alpha t \beta$. $\alpha t \beta$ is consistently oriented if each cup resp cap has exactly one $\vee$ and one $\wedge$ and the endpoints of each line segment are labelled by the same symbol. Notation: $\alpha \rightarrow^t \beta$.

\medskip

For $t$ a crossingless $\Delta \Gamma$ and $\lambda \in \Delta, \ \mu \in \Gamma$ label the bottom and the upper line as usual. The \textit{lower reduction} $red(\underline{\lambda}t)$ is the cup diagram obtained from $\underline{\lambda}t$ by removing the bottom number line and all connected components that do not extend up to the top number line. 

\begin{thm} \cite{Brundan-Stroppel-5}, Thm 3.4. and \cite{Brundan-Stroppel-2}, Thm 4.11: In $K_0(\calR_n)$ the mixed tensor $R(\lambda)$ attached to the bipartition $\lambda$ satisfies \[ [ R(\lambda)] = \sum_{ \mu \subset \alpha \rightarrow^t \one, \ red(\underline{\mu}t) = \underline{\one} } [ L(\mu) ]\] where $t$ is a fixed matching determined by $\lambda$ between the block $\Gamma$ of $\one$ and the block $\Delta$ of $\lambda^{\dagger}$ \cite{Brundan-Stroppel-5}, 8.18. If $L(\mu)$ is a composition factor of $R(\lambda)$, its graded composition multiplicities are given by \[ \sum_{\mu} (q + q^{-1})^{n_{\mu}} [L(\mu)]\]where $n_{\mu}$ is the number of lower circles in $\underline{\mu}t$.
\end{thm}

\begin{lem} The module $R= R((n)^{n})$ in ${\calR}_{n+r+1}$, $r \geq 0$, has Loewy length $2n+1$ with socle
and cosocle equal to $\one$. We have $DS(R(n^n))) = R(n^n)$. If $r=0$, $DS(R) = P(\one)$. $R$ contains $\one$ with
multiplicity $2^{2n}$. It contains the irreducible module $L(h) = [n,1,\ldots,1,0,\ldots,0]$ (with $1$ occurring $n-1$-times) 
in the second Loewy layer. The multiplicity of $L(h)$ in $R$ is $2^{2(n-1)}$. It contains the module $[n,n,\ldots,n,0,\ldots,0]$ as the constituent of highest weight in the middle Loewy layer with multiplicity 1. It does not contain the modules $BS^{n-1} = [n,1,\ldots,1]$, $BS^n = [n+1,1,\ldots,1]$, $[n,1,\ldots,1,-1]$ and $[n,1,\ldots,1,-1,\ldots,-1]$ (with $1$ occurring $n-1$-times) as composition factors.
\end{lem}

\medskip{\it Proof}. The Loewy length of a mixed tensor is $2d(\lambda)+1$ (where $d(\lambda)$ is the number of caps)
and $d((n-1)^{n-1}) = n-1$ \cite{Heidersdorf-mixed-tensors}. The composition factors of $R$ are given as a sum $\sum_{\mu} (q + q^{-1})^{n_{\mu}}
[L(\mu)]$. For our choice of $\lambda = (n-1)^{n-1}$
the matching is given by \cite{Heidersdorf-mixed-tensors} (picture for $n=4$)

\bigskip
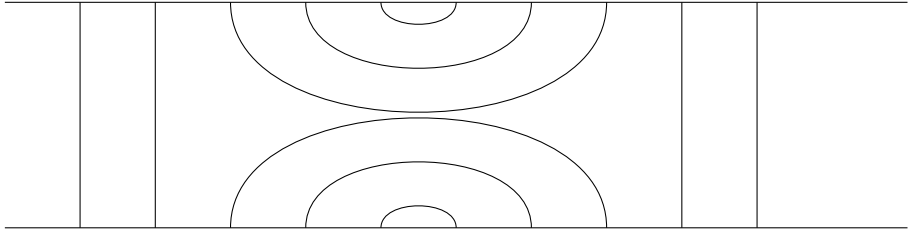
\begin{center}

  \begin{tikzpicture}
 \draw (-6,0) -- (6,0);
\foreach \x in {} 
     \draw (\x-.1, .2) -- (\x,0) -- (\x +.1, .2);
\foreach \x in {} 
     \draw (\x-.1, -.2) -- (\x,0) -- (\x +.1, -.2);
\foreach \x in {} 
     \draw (\x-.1, .1) -- (\x +.1, -.1) (\x-.1, -.1) -- (\x +.1, .1);

\begin{scope} [yshift = -3 cm]

 \draw (-6,0) -- (6,0);
\foreach \x in {} 
     \draw (\x-.1, .2) -- (\x,0) -- (\x +.1, .2);
\foreach \x in {} 
     \draw (\x-.1, -.2) -- (\x,0) -- (\x +.1, -.2);
\foreach \x in {} 
     \draw (\x-.1, .1) -- (\x +.1, -.1) (\x-.1, -.1) -- (\x +.1, .1);
\end{scope}

\draw [-,black,out=270,in=270](-1,0) to (0,0);
\draw [-,black,out=270,in=270](-2,0) to (1,0);
\draw [-,black,out=270,in=270](-3,0) to (2,0);

\draw [-,black,out=90,in=90](-1,-3) to (0,-3);
\draw [-,black,out=90,in=90](-2,-3) to (1,-3);
\draw [-,black,out=90,in=90](-3,-3) to (2,-3);

\draw [-,black,out=270,in=90](-4,0) to (-4,-3);
\draw [-,black,out=270,in=90](-5,0) to (-5,-3);
\draw [-,black,out=270,in=90](3,0) to (3,-3);
\draw [-,black,out=270,in=90](4,0) to (4,-3);

\end{tikzpicture}
\bigskip

\end{center}

with $n$ caps and where the rightmost vertex in a cap is at position $n$.
The irreducible module in the socle and cosocle is easily computed from the rules of
the section \ref{stable0}. The weight \[h =(n,1,\ldots,1,0|0,-1,\ldots,-1,n-n)\] easily seen to satisfy $h \rightarrow^t \one, \ red(\underline{h}t) = \underline{\one}$, hence occurs as a composition factor.
The number of lower circles in the lower reduction $\underline{h}t$ is $n-1$, hence $L(h)$ occurs with
multiplicity $2^{2(n-1)}$. If we number the Loewy layers starting with the socle by
$1,\ldots,2n+1$, $L(h)$ occurs in the $2k$-th Loewy layer ($k=1,\ldots,n$)
with multiplicity $\binom{n-1}{k-1}$. Likewise for $\one$ with $n_{\one} = n-1$. We note: A weight $\mu$ can only satisfy $red(\underline{\lambda}t) = \underline{\one}$ if the vertices $-n, -n-1,\ldots, -n-r$ (the first vertices left of the caps) are labelled by $\vee$. Hence:

\begin{itemize}
\item $BS^{n-1}$ does not occur as a composition factor. The vertex $-n$ is labelled by $\wedge$.
\item $[n,1,\ldots,1,-1]$ does not occur as a composition factor. The vertex $-n$ is labelled by $\wedge$
\item $[n+1,1,\ldots,1]$ does not occur as a composition factor since all composition factors $[\mu_1,\ldots,\mu_n]$ satisfy $\mu_1 \leq n$ since $[n,\ldots,n,0,\ldots,0]$ is the constituent of highest weight.
\end{itemize} \qed

\medskip

\textbf{Remark}. In particular the constituent $\one$ occurs with the same
multiplicity as in $P(\one) \in \calR_n$. 

\medskip

\textbf{Remark}. The module $R(n^n)$ can be obtained as follows. Let $\{ n^n \}$ be the covariant module to the partition $(n^n)$ and $\{ n^n\}^{\vee}$ its dual. Then $R(n^n)$ is the projection on the maximal atypical block of $\{ n^n \} \otimes \{ n^n\}^{\vee}$.

\medskip

\textbf{Example}. For $Gl(3|3)$ the Loewy structure of the module $R(2^2)$ is  \[ \begin{pmatrix}   [0,0,0] \\  [1,0,0] \oplus [2,1,0] \\ [2,0,0] \oplus
[2,-1,-1] \oplus [0,0,0] \oplus [0,0,0] \oplus [1,1,0] \oplus [2,2,0] \\ [1,0,0]
\oplus [2,1,0] \\ [0,0,0] \end{pmatrix} .\]

\medskip



\section{The basic hook representations $S$}\label{hooks}

\medskip
{\it The case $L={\bf 1}$}.  Suppose $Z$ has cosocle ${\bf 1}$ and the
projection $q:Z\to {\bf 1}$ satisfies condition (S). If $Z$ is not simple, we constructed
objects $V_\nu$ with cosocle ${\bf 1}$ and simple socle $S_\nu=Ker(q_\nu)$.
In this situation $Ext^1_{\calR_n}({\bf 1}, S_\nu)\neq 0$.

\begin{prop}\label{trivialextension}
For any nontrivial extension 
$$ \xymatrix{ 0 \ar[r] & S_\nu \ar[r] & V \ar[r] &  {\bf 1} \ar[r] &  0} $$
the vectorspace $  \omega_{n,0}^0(V)$ is zero (for simple $S_\nu$). Hence $\omega(q): \omega(V) \to \omega({\bf 1})$
is the zero map.
\end{prop}

\medskip
For the proof we use several lemmas. Finally lemma \ref{omeganull}
proves the proposition.

\begin{lem}\label{ext}
Up to isomorphism there are $n+1$ irreducible modules $L$ in $\calR_n$
such that $Ext^1({\bf 1},L)\neq 0$. They are 
\begin{enumerate}
\item $L_n(n)=Ber_n\otimes S^{n-1}$ and 
\item its dual $L_n(n)^\vee \cong [0,..,0,-1]$, and for 
\item $i=1,..,n-1$ the basic selfdual representations $$L_n(i)=[i,1,\cdots,1,0,\cdots,0]$$ (with $n-i$ entries $0$).
\end{enumerate}
In all cases $\dim(Ext^1(L,{\bf 1}))=1$.
Furthermore $$DS_{n,j}(L_n(i)) = L_{j}(i)$$ holds for $i < j \leq n$ and  
$$ \fbox{$ DS_{n,i}(L_n(i)) = L_{i}(i) \oplus L_{i}(i)^\vee \oplus Y $}     \ $$
where $Y\not\cong {\bf 1}$ is an irreducible module with $Ext^1({\bf 1},Y)=0$ and
sector structure
$$    [\vee_{-n},\wedge_{-n+1}]\wedge_{-n+2} [-n+3,...,n-2]\wedge_{n-1} [\vee_n,\wedge_{n+1}]                      \ .$$
\end{lem}

\medskip
{\bf Example}. $L_n(1)=S^1$.

\medskip{\it Proof}. $L^*\cong L$
for irreducible objects $L$ implies $Ext^1({\bf 1},L) \cong Ext^1(L,{\bf 1})$.
Furthermore $Ext^1(L,{\bf 1}) \cong Ext^1((L^*)^\vee,{\bf 1})$ and $L=L^*$, hence
$$Ext^1(L,{\bf 1}) \cong Ext^1(L^\vee,{\bf 1})\ .$$
By \cite{Brundan-Stroppel-2}, cor. 5.15 for $L=L(\lambda)$
$$\dim Ext^1(L(\lambda),{\bf 1}) = \dim Ext^1(V(\lambda),{\bf 1}) 
+ \dim Ext^1(V(0),L(\lambda)) \ $$
holds.  Since ${\bf 1}$ is a Kostant weight,  there exists a unique weight
$\lambda$ characterized by $\lambda \leq 0$ (Bruhat ordering) and 
$l(\lambda,0)=1$ in the notations of loc. cit. lemma 7.2, such that $\dim Ext^1(V(\lambda),{\bf 1})\neq 0$. One easily shows
$L(\lambda) \cong [0,..,0,-1]$.

\medskip 
On the other hand 
$ \dim Ext^1(V(0),L(\lambda)) \neq 0$ implies $0 < \lambda$
(see the explanations preceding loc. cit. (5.3) and loc. cit. lemma 5.2.(i)).
Then for any pair of adjacent labels $i,i+1$ of $\rho$ of type $i=\vee,i+1=\wedge$
we write $\rho\in \Lambda^{\vee,\wedge}$, if the labels of $\rho$ at $i,i+1$
are the same $i=\vee,i+1=\wedge$. Then     lemma 5.2(ii) of loc. cit. gives
$$ \dim(Ext^1_{\calR_n}(V({\rho}),L(\lambda)) = \begin{cases}
dim(Ext^1_{\calR_{n-1}}(V(\rho'),L(\lambda')) & \text{ if } \lambda \in \Lambda^{\vee,\wedge} \cr
\dim(Hom_{\calR_n}(V(\rho''),L(\lambda)) & \text{ otherwise } \cr
\end{cases} $$
Here $\lambda',\mu'$ are obtained from $\lambda,\mu$ by deleting $i,i+1$,
and $\rho''$ is obtained by transposing the labels at $i,i+1$.

\medskip
This shows our assertion, since for $$L(\rho)={\bf 1}$$ there is a unique
pair of such neighbouring indices for
$$ [\vee_{-n+1},...,\vee_0,\wedge_1,...,\wedge_n] \ ,$$
namely at the position $(i,i+1)=(0,1)$. We now assume $n\geq 2$. Then switching this pair gives $L_n(1)$
below. Freezing then also $(-1,.,.,2)$ gives $L_n(2)$ and so on. 
Hence applying this lemma of loc. cit. several times
will prove our first claim. 
Indeed, as long
as we freeze less than $n-2$ pairs, we end up
for every $j$ from $1,...,n-1$ with a representation $L_n(j)$. It has only one sector
$$ [\vee_{1-n},...,\vee_{-j-1}[\vee_{-j}\wedge_{-j+1}][\vee_{-j+2},...,\vee_0,\wedge_1,...,\wedge_{j-1}][\vee_j,\wedge_{j+1}]\wedge_{j+2},...,\wedge_n] \ .$$
In addition, if we freeze $n-1$ pairs we end up with $L_n(n)$
with the sector structure 
$$ [\vee_{2-n},\vee_{3-n}, ...., \wedge_{n-2},\wedge_{n-1}][\vee_n,\wedge_{n+1}] \ .$$
Indeed $L_n(n) \cong Ber_n \otimes S^{n-1}$. 

\medskip
The remaining assertions now follow from theorem \ref{mainthm},
since $L_{n+1}(n)$ has sectors 
$$S_1S_2S_3 =[\vee_{-n},\wedge_{-n+1}][-n+2,...,n-1][\vee_n,\wedge_{n+1}]\ .$$
Hence $DS(L_{n+1}(n)) = (Ber\otimes S^{n-1}) \oplus (Ber\otimes S^{n-1})^\vee \oplus 
Y$ for $Y$ with sector structure $[\vee_{-n},\wedge_{-n+1}][-n+3,...,n-2]\wedge_{n-1}[\vee_n,\wedge_{n+1}]$. 
\qed

\bigskip
{\it Basic cases}.
For a nontrivial extension 
$$ \xymatrix{ 0 \ar[r] & L_n(i) \ar[r] &  V \ar[r] &  {\bf 1} \ar[r] & 0 } $$
first suppose $S=L_n(i)$ is {\it basic}, so $i\in \{1,...,n-1\}$. 
Since $(L_n(i)^*)^\vee \cong L_n(i)^\vee \cong L_n(i)$
for $i<n$, $(V^*)^\vee$ again defines a nontrivial extension
$$ \xymatrix{ 0 \ar[r] & L_n(i) \ar[r] &  (V^*)^\vee \ar[r] &  {\bf 1} \ar[r] & 0 } \ .$$

\medskip
We now use $DS(L_n(i))=L_{n-1}(i)$ for $1\leq i<n-1$. 
In $T_{n-1}$ the induced long exact sequence  
$$ \xymatrix{ H^{-1}({\bf 1}) \ar[r] & H^0(S) \ar[r] &  H^0(V) \ar[r] &  H^0({\bf 1}) \ar[r] & H^1(S)  } $$
remains exact, since $H^\ell(S)=L_{n-1}(i)$ for $\ell=0$ and is zero otherwise and similarly 
$H^\ell({\bf 1})={\bf 1}$ for $\ell=0$ and is zero otherwise. 
In other words for {\it basic} $S$ we obtain from the given extension
in $\calR_n$ an exact sequence in $\calR_{n-1}$
$$ \xymatrix{ 0 \ar[r] & L_{n-1}(i) \ar[r] &  DS(V) \ar[r] &  {\bf 1} \ar[r] & 0 } \ .$$
Repeating this $n-i$ times we obtain an exact sequence 
$$ \xymatrix{ 0 \ar[r] & L_{i}(i) \oplus L_i(i)^\vee \oplus Y \ar[r] &  DS_{n,i}(V) \ar[r] &  {\bf 1} \ar[r] & 0 } \ .$$
Since $Ext^1({\bf 1},Y)=0$ this implies
$$   DS_{n,i}(V) = E \oplus Y $$
for some selfdual module $E$ defining an extension between ${\bf 1}$ and
$L_{i}(i) \oplus L_i(i)^\vee$.
We claim that this exact sequence does not split in $\calT_i$. 

\begin{prop} Suppose $r$ is an integer $\geq 0$. For an indecomposable module $V$ defining a nontrivial extension between ${\bf 1}$ and $L_{n+1+r}(n)$ in $\calR_{n+r+1}$, the object $(DS)^{\circ r+1} (V)$ decomposes into the direct sum of the irreducible module $Y$ from above and an 
indecomposable extension module $E$ in $\calR_n$.
\end{prop}

{\it Proof}. Note that any two such indecomposable extensions define isomorphic modules $V$, since the relevant $Ext$-groups are one-dimensional. We assume $r=0$ for simplicity.
Since the constituents $L_{n+1}(n)$ and ${\bf 1}$ of $V$ are basic, this implies
$DS(V) = H^0(V)=Y\oplus E$. If the module $E$ is not indecomposable, 
it is semisimple (for this use Tannaka duality). We proceed as follows:

\medskip
For the mixed tensor $R=R_{n^{n}}$ in $\calR_{n+1}$ we know that 
its image $DS(R_{n^n})$ is the projective hull
$P({\bf 1})$ of ${\bf 1}$ in $\calR_n$ and $P({\bf 1})$ is an indecomposable 
module with top ${\bf 1}$.
The module $R_{n^{n}}$ admits
as quotient an indecomposable module $V$ defining a nontrivial extension between ${\bf 1}$ (the top of $R$)
and the module $L_{n+1}(n)$ (which sits in the second layer of the Loewy filtration of $R$). 
Hence $R/K \cong V$ for some submodule $K$ of $R$. We claim that
$$    \xymatrix{  0  \ar[r] & H^0(K) \ar[r]^i & H^0(R) \ar[r]^p & E \oplus Y \ar[r] & H^1(K) \ar[r] & 0 } $$
is exact and $H^\nu(K)=0$ for $\nu\neq 0,1$. 
For this use $H^\bullet(V)=H^0(V)$ and $H^\nu(R)=P({\bf 1})$ for a unique $\nu$.
If $\nu\neq 0$, then $H^\nu(K) \to H^\nu(R)=P({\bf 1})$ would be surjective
and therefore $H^\nu(K) = P({\bf 1}) \oplus ?$. We exclude this later.
So suppose for the moment $\nu=0$. 

\medskip
The image of $p$ can not contain the irreducible module
$Y\not\cong 1$, since the top of $P({\bf 1})$ is ${\bf 1}$.
If $E$ splits, it is semisimple. Then the image of $p$ can not contain $E$ either,
since again this would contradict that $P({\bf 1})$ has top ${\bf 1}$. 
Therefore the image of $p$ is ${\bf 1}$ or zero, if $E$ splits.
This leads to a contradiction:

\medskip
Look at all constituents $X$ of $R$ with $Ber \otimes S^{n-1}$
in $H^\nu(X)$ for $\nu=-1,0,1$. These $X$ are isomorphic
to the following irreducible modules $X_{-1},X_0,X_1$ with $Ber \otimes S^{n-1}$ 
occuring  in $H^i(X_i)$ respectively: the basic module $X_0=L_{n+1}(n)$ with sector structure $[-n,-n+1][-n+2,...,n-1][n,n+1]$
and $X_1$ with
sector structure $[1-n[2-n,...,n-1][n,n+1]n+2]$ and
$X_{-1}$ with sector structure $[-n-1,-n]\boxminus [2-n,...,n-1][n,n+1]$.
Then $Ber\otimes S^{n-1}$ occurs in $H^i(X_i)$ for $i=0,\pm 1$.

\medskip
Let $F^i(.)$ denote the descending Loewy filtration. 
For a module $Z$ let $m(Z)$ denote the number of Jordan-H\"older
constituents of $Z$ that are isomorphic to $Ber^\otimes S^{n-1}$.
Next we use that for all $i$
$$   \fbox{$ X_1,X_{-1}\ \text{ does not occur in the }\  gr^i_F(R) $}  \ .$$
Indeed according to section \ref{sec:n^n} all irreducible constituents $[\lambda]$ satisfy
the property $\lambda_{n+1}=0$ except for one given by $[n,-n+1,...,-n+1]$.
Therefore  $m(H^{\pm 1}(gr^i_F(R)))=0$ and hence $m(H^1(F^i(X)))=0$.
Since also $m(H^{-1}(gr^i_F(R)))=0$, then
$$ H^{-1}(gr^i_F(R)) \to H^0(F^{i}(R)) \to H^0(F^{i-1}(R)) \to H^0(gr^i_F(R)) \to H^1(F^{i}(R)) $$
implies $  m(H^0(F^{i}(R))) = m(H^0(F^{i-1}(R))) + m(H^0(gr^i_F(R)))$. For small $i$
we have $F^i(R)=R$ and therefore
$$  m(H^0(R)) = \sum_i m(H^0(gr^i_F(R))) \ .$$
They same argument then applies for the submodule $K$ of $R$.
Hence $$ m(H^0(K)) = m(H^0(R)) - 1$$ by counting the multiplicities of $X_0$
in $K$ resp. $R$. Hence the image of $p$ must contain $Ber\otimes S^{n-1}$
and hence $E$ is an indecomposable quotient of $P({\bf 1})$.

\medskip
Now let us adress the assertion $\nu=0$ from above. If $\nu\neq 0$, then
$H^0(K) \cong P({\bf 1}) \oplus ?$ gives a contradiction using the same counting
argument. 

\medskip
In the case $r>0$ one uses the same kind of argument. Again the extension defined by $V$ in $\calR_{n+r+1}$ can be realized as a quotient of  
$R=R_{n^{n}}$ in $\calR_{n+r+1}$. The argument is modificatis modificandis the same.

\medskip
{\it The non-basic cases}.
For a nontrivial extension in $T_n$ of the form
$$ \xymatrix{ 0 \ar[r] & L_n(n) \ar[r] &  V \ar[r] &  {\bf 1} \ar[r] & 0 } \ .$$
we get a dual nontrivial extension
$$ \xymatrix{ 0 \ar[r] & [0,...,0,-1] \ar[r] &  (V^*)^\vee \ar[r] &  {\bf 1} \ar[r] & 0 } \ .$$
In lemma \ref{Qa} and lemma \ref{Ia} we defined $Q_a$, which for $a=1$ defines 
a nontrivial extension between ${\bf 1}$ and $L_n(n)^\vee$. Since $\dim(Ext^1({\bf 1}, L_n(n)^\vee))=1$, we get
$$  (V^*)^{\vee} \cong Q_1 \ .$$

By corollary \ref{aKac2} we get
$  DS_{n,0}^\ell(V) =  0$ and  $\omega_{n,0}^\ell(V) =  0$
for all $\ell \leq 0$. Similarly by duality 
$DS_{n,0}^{\ell}((V^*)^\vee)=0$ for $\ell \geq 0$.
This implies $\omega_{n,0}^{\ell}((V^*)^\vee)=0$ for $\ell \geq 0$.

Finally consider the nontrivial extension $V_i$ between $\bf 1$ and $L_n(i)$
in $\calR_n$ and the nontrivial extension
$DS_{n,i}(V_i)$ in $\calR_i$ from above.
It has the form $DS_{n,i}(V_i) = E \oplus Y$ for 
$$    0 \to L_i(i) \oplus L_i(i)^\vee  \to E \to {\bf 1} \to 0 \ .$$
The module $E = DS_{n,i}((V_i)/Y$
is the pullback of a nontrivial
extension of ${\bf 1}$ by $L_i(i)$ 
$$  0 \to L_i(i) \to   E_1 \to  {\bf 1} \to 0 $$
and of a nontrivial
extension of ${\bf 1}$ by $(L_i(i))^{\vee}$
$$  0 \to (L_i(i))^{\vee} \to   E_2 \to  {\bf 1} \to 0 \ .$$
 Hence there exists an exact sequence
$$    0 \to DS_{n,i}(V_i)/Y \to  E_1 \oplus E_2 \to {\bf 1} \to  0  \ $$ 
so that 
$$  \to  DS_{i,0}^{-1}({\bf 1}) \to DS_{i,0}^0(DS_{n,i}(V_i)/Y) \to  DS_{i,0}^0(E_1) \oplus DS_{i,0}^0(E_2) \to $$
 is exact.  Since $DS_{i,0}^0(E_1)=0$ and $DS_{i,0}^0(E_2)=0$ by 
  corollary \ref{aKac2}
 and since
  $DS_{i,0}^{-1}({\bf 1})=0$, therefore
$DS_{i,0}^0(DS_{n,i}(V_i)/Y)=0$. Hence
$\omega_{i,0}^0(DS_{n,i}(V_i))=\omega_{i,0}^0(Y)$.  
The Leray type spectral sequence 
$$  DS_{i,0}^p(DS_{n,i}^q(V_i)) \Longrightarrow DS_{n,0}^{p+q}(V_i) $$
degenerates, since $DS_{n,i}^q(V_i)=0$ for $q\neq 0$.
Therefore also $\omega^{\bullet}_{n,i}(V_i) \cong DS^{\bullet}_{n,i}(V_i)$.
One can now argue as in the proof of lemma \ref{insteps} 
to show
$$  \omega_{n,0}^0(V_i) = \omega_{i,0}^0(DS_{n,i}(V_i))=\omega_{i,0}(Y)\ .$$
Since the map
$$   \omega_{i,0}(q) : E \oplus Y \longrightarrow {\bf 1} $$
is trivial on the simple summand $Y\not\cong {\bf 1}$,  the next lemma follows.

\begin{lem}\label{omeganull}
For every nontrival extension $$\xymatrix{ 0\ar[r] & S \ar[r] & V \ar[r]^{q_V} & {\bf 1} \ar[r] & 0 }$$ 
of ${\bf 1}$ by a simple object $S$ in $\calR_n$
the map $\omega^0(q_V)$ vanishes.
\end{lem}

\medskip
The last lemma completes the proof of proposition \ref{trivialextension}.
This implies the following main result.

\medskip
\begin{cor}\label{splitting1}
Suppose $Z$ is indecomposable and $cosocle(Z)\cong {\bf 1}$.
If the quotient map $q:Z \to {\bf 1}$ is strict, then $q: Z\cong {\bf 1}$.
\end{cor}

{\it Remark.} A symmetric abelian tensor category in the sense of Deligne is semisimple if and only if $q:Z \to \one$, with cosocle of $Z$ isomorphic to $\one$, is an isomorphism.

\begin{cor}
For a nontrivial extension $V$ between ${\bf 1}$ and $L_n(n)$ or its dual $L_n(n)^\vee$
$$H_D^\nu(V)=0 \quad , \quad \nu\neq 1\ $$ holds, and hence the induced map $H_D(q): H_D(V) \to H_D({\bf 1})$
is trivial. 
\end{cor} 
 
\medskip







\end{document}